\documentclass[11pt,reqno]{amsart}
\usepackage{amsmath, amsthm, amssymb, amsfonts}
\usepackage{bbm}
\usepackage{bm}
\usepackage[utf8]{inputenc}
\usepackage{graphicx}
\usepackage{xcolor}
\usepackage{enumitem}
\usepackage{tikz}
\usepackage{tikz-cd}
\usepackage{pgfplots}
\usepackage{rotating}
\usepackage{float}
\usepackage{bookmark}
\usepackage{thmtools}
\usepackage{hyperref}
\usepackage{cleveref}

\pgfplotsset{compat=1.18}

\theoremstyle{plain}

\usepackage[backend=biber,sorting=nty,maxbibnames=99]{biblatex}

\numberwithin{equation}{section}

\makeatletter
\ifdefined\newcounteralias
  \let\orig@newcounteralias\newcounteralias
  \renewcommand*\newcounteralias[2]{%
    \def\reserved@a{#1}\def\reserved@b{#2}%
    \ifx\reserved@a\reserved@b\else\orig@newcounteralias{#1}{#2}\fi}
\fi
\makeatother

\newtheorem{theorem}{Theorem}[section]
\newtheorem{lemma}[theorem]{Lemma}
\newtheorem{proposition}[theorem]{Proposition}
\newtheorem{corollary}[theorem]{Corollary}
\newtheorem{claim}[theorem]{Claim}

\theoremstyle{definition}
\newtheorem{definition}[theorem]{Definition}
\newtheorem{notation}[theorem]{Notation}
\newtheorem{example}[theorem]{Example}
\newtheorem{qn}[theorem]{Question}

\newtheorem*{axiom}{Axiom}

\newcounter{subappendix}[section] % 0 = numbered subsections, 1 = lettered
\newcommand{\appsubsection}{%
  \ifnum\value{subappendix}=0 \setcounter{subsection}{0}\setcounter{subappendix}{1}\fi
  \subsection}
\renewcommand{\thesubsection}{\thesection.\ifnum\value{subappendix}>0 \Alph{subsection}\else\arabic{subsection}\fi}
\renewcommand{\theHsubsection}{\theHsection.\ifnum\value{subappendix}>0 \Alph{subsection}\else\arabic{subsection}\fi} % keeps hyperref anchors distinct from 2.1, 2.2, ...
\makeatletter
\newlength{\toclabelwidth}
\renewcommand{\tocsubsection}[3]{%
  \settowidth{\toclabelwidth}{0.0.\quad}%
  \indentlabel{\@ifnotempty{#2}{\ignorespaces#1 \makebox[\toclabelwidth][l]{#2.}}}#3}
\makeatother

\newcommand{\func}{\operatorname}

\renewcommand{\c}[1]{\langle #1 \rangle}

\newcommand{\bb}[1]{\left( #1 \right)}

\newcommand{\Z}{\mathbb{Z}}
\newcommand{\Q}{\mathbb{Q}}
\newcommand{\Po}{\cal{P}}
\newcommand{\N}{\mathbb{N}}
\newcommand{\F}{\cal{F}}

\renewcommand{\bar}[1]{\overline{#1}}

\newcommand{\eqbreak}{\phantom{{} = {}}}

\newcommand{\iffbreak}{\phantom{{} \iff {}}}

\newcommand{\cal}{\mathcal}
\renewcommand{\frak}{\mathfrak}

\renewcommand{\bf}{\mathbf}

\newcommand{\supp}{\func{supp}}

\newcommand{\Ob}{\func{Ob}}

\newcommand{\ZFC}{\mathsf{ZFC}}

\newcommand{\CH}{\mathsf{CH}}

\newcommand{\MA}{\mathsf{MA}}

\newcommand{\dom}{\func{dom}}
\newcommand{\ran}{\func{ran}}

\newcommand{\A}{\cal{A}}
\newcommand{\B}{\cal{B}}

\newcommand{\forces}{\Vdash}
\newcommand{\one}{\mathbbm{1}}
\newcommand{\restrictedto}{\mathord{\upharpoonright}}

\newcommand{\Coll}{\func{Coll}}

\newcommand{\FIN}{\bf{FIN}}
\newcommand{\U}{\cal{U}}
\newcommand{\V}{\bf{V}}

\newcommand{\X}{\cal{X}}
\newcommand{\Y}{\cal{Y}}

\newcommand{\lh}{\func{lh}}
\newcommand{\G}{\cal{G}}

\newcommand{\depth}{\func{depth}}
\newcommand{\K}{\cal{K}}
\newcommand{\fin}{\mathrm{fin}}
\newcommand{\D}{\cal{D}}

\newcommand{\E}{\cal{E}}

\newcommand{\height}{\func{height}}

\newcommand{\minmax}{\func{minmax}}

\NewCommandCopy\strokeL\L
\DeclareRobustCommand{\L}{\ifmmode\mathbf{L}\else\strokeL\fi}
\NewCommandCopy\pilcrowP\P
\DeclareRobustCommand{\P}{\ifmmode\mathbb{P}\else\pilcrowP\fi}
\NewCommandCopy\accH\H
\DeclareRobustCommand{\H}{\ifmmode\cal{H}\else\expandafter\accH\fi}
\NewCommandCopy\slashO\O
\DeclareRobustCommand{\O}{\ifmmode\cal{O}\else\slashO\fi}
\NewCommandCopy\Sec\S
\DeclareRobustCommand{\S}{\ifmmode\mathcal{S}\else\Sec\fi}

\newenvironment{midproof}[1][Proof]
  {\begin{proof}[#1]}
  {\end{proof}}

\makeatletter
\DeclareRobustCommand\bigop[1]{%
  \mathop{\vphantom{\sum}\mathpalette\bigop@{#1}}\slimits@
}
\newcommand{\bigop@}[2]{%
  \vcenter{%
    \sbox\z@{$#1\sum$}%
    \hbox{\resizebox{\ifx#1\displaystyle.9\fi\dimexpr\ht\z@+\dp\z@}{!}{$\m@th#2$}}%
  }%
}
\makeatother

\makeatletter
\DeclareRobustCommand{\sqsubsetneq}{\mathrel{\mathpalette\sqsubsetneq@\relax}}
\newcommand{\sqsubsetneq@}[2]{%
  \sbox\z@{$\m@th#1\sqsubseteq$}%
  \copy\z@
  \kern-.5\wd\z@
  \raise\dimexpr.03\wd\z@-\dp\z@\relax\hbox to\z@{\hss\rotatebox[origin=c]{40}{\rule[-.0275\wd\z@]{.42\wd\z@}{.055\wd\z@}}\hss}%
  \kern.5\wd\z@
}
\makeatother

\DeclareFontShape{OT1}{cmr}{bx}{sc}{<-> cmbcsc10}{}

\title[Projections between Ramsey quasi-orders]{Projections between Ramsey quasi-orders and topological Ramsey spaces}
\author{Spencer Unger, Clement Yung}

\begin{document}

\begin{abstract}
    A \emph{Ramsey quasi-order} is a quasi-order $(\cal{R},\leq)$ such that $\cal{R}$ is a Polish space, $\leq$ is analytic and every analytic subset of $\cal{R}$ is Ramsey. In this paper, we study the transfer of the analytic Ramsey property via monotone quasi-open maps, which are closely related to a family of forcing projections between quasi-orders (ranging from \emph{weak projections} to \emph{strong projections}). We apply this framework to the class of topological Ramsey spaces by providing a detailed list of examples of strong projections between topological Ramsey spaces, and show that, for many of these examples, there is no weak projection in the reverse direction. We also exhibit pairs of topological Ramsey spaces that do not weakly project into each other. Finally, partially confirming a conjecture of Dobrinen, we show that every space in a large class of topological Ramsey spaces either strongly projects into the Ellentuck space or contains a dense set of atoms.
\end{abstract}

\maketitle
\tableofcontents

\section{Introduction}
\label{sec:intro}
Consider two Polish spaces $\cal{R}$ and $\S$, each equipped with a meaningful notion of a ``Ramsey subset''. Let $\rho(\cal{R})$ denote the statement ``every analytic set is Ramsey'', a variant of Silver's infinite-dimensional Ramsey theorem for analytic sets \cite{S70}. This paper focuses on the following question:

\begin{qn}
    Assuming that both $\rho(\cal{R})$ and $\rho(\S)$ are true, can we prove or disprove that $\rho(\S)$ is a ``consequence'' of $\rho(\cal{R})$?
\end{qn}

\subsection{(Infinite) Ramsey theory} Let $\N = \{0,1,2,\dots\}$ denote the set of non-negative integers. Given $A \subseteq \N$ and some $n \in \N$, we let $[A]^n$ denote the set of all $n$-element subsets of $A$. Motivated by a problem about decision procedures for logical systems, Ramsey proved the original Ramsey theorem about countable sets:

\begin{theorem}[{\cite[Infinite Ramsey theorem]{R30}}]
\label{thm:infinite.Ramsey.theorem}
    For every positive integer $n$ and every finite colouring of $[\N]^n$, there exists some infinite \emph{monochromatic} subset $A \subseteq \N$ --- that is, every element of $[A]^n$ has the same colour.
\end{theorem}

\begin{corollary}[{\cite[Finite Ramsey theorem]{R30}}]
\label{cor:finite.Ramsey.theorem}
    For all positive integers $k,l,m$, there is a positive integer $n$ such that for every $a \in [\N]^n$ and every $k$-colouring of $[a]^l$, there exists some \emph{monochromatic} $b \in [a]^m$ --- that is, every element of $[b]^l$ has the same colour.
\end{corollary}

Results similar in flavour to the original Ramsey theorem existed prior to Ramsey's paper \cite{R30}, such as Schur's theorem \cite{S12} and van der Waerden's theorem \cite{V27}, both of which provided number theorists with useful insights. However, it was only in the 1970s that applications of Ramsey's ideas to other structures were extensively studied. One direction extended the finite Ramsey theorem to other structures, such as finite linearly ordered graphs (independently by Abramson and Harrington \cite{AH78}, and Ne\v{s}et\v{r}il and R\"{o}dl \cite{NR77}) and finite Boolean algebras (by Graham and Rothschild in \cite{GR71}). In fact, they proved what is now known as the Graham--Rothschild $n$-parameter-set theorem, a partition theorem about $n$-parameter words generalising the Hales--Jewett theorem in \cite{HJ63}. 

It was also in the 1970s that infinite-dimensional Ramsey theory, i.e. Ramsey theory surrounding the infinite Ramsey theorem for $n = \infty$, began to develop. While the infinite Ramsey theorem for $n = \infty$ was shown to be false by Erd\H{o}s and Rado in \cite{ER52} assuming the Axiom of Choice (in fact, the existence of a non-principal ultrafilter suffices --- see \cite{M75}), one may still ask whether the theorem holds for ``nice'' colourings, i.e. colourings of low descriptive complexity. Indeed, it was shown through the works of Galvin--Prikry \cite{GP73} and Silver \cite{S70} that the infinite Ramsey theorem for $n = \infty$ holds for colourings in the $\sigma$-algebra generated by analytic sets. Ellentuck later provided a purely combinatorial proof of Silver's theorem in \cite{E74} by considering the Ellentuck topology, which contrasts with Silver's original, deeply metamathematical proof. Ellentuck's ideas proved to be potent beyond the Ellentuck space --- Milliken extended Silver's theorem in \cite{Mil75} to the Milliken space (also called the Hindman space, see Example \ref{ex:Hindman.FIN}), showing that every analytic set of infinite block sequences of non-empty finite subsets of $\N$ is Ramsey. A similarly inspired theory was developed by Gowers in \cite{G92} and \cite{G02} to resolve several long-standing problems in Banach space theory. An abstract version of Ellentuck's theorem first appeared in \cite{CS84}, within a framework proposed by Carlson and Simpson that was also discussed in \cite{CS90}. This framework later formed the basis for the axioms \textbf{A1}--\textbf{A4} introduced by Todor\v{c}evi\'{c} in \cite{T10}. Carlson and Simpson also proved an infinite-dimensional extension of the Graham--Rothschild $n$-parameter-set theorem in \cite{C87}. 

\subsection{Comparing Ramsey theorems}
While some Ramsey theorems are obvious generalisations of others, such as the direct derivation of Schur's theorem from Hindman's theorem on finite sums in \cite{H74}, others are less straightforward. For instance, Pr\"{o}mel and Voigt showed in \cite{PV81} that one may use the Graham--Rothschild $n$-parameter-set theorem to prove the Ramsey property of finite linearly ordered graphs. Inspired by this proof, Ma\v{s}ulovi\'{c} proposed in \cite{M18} a category-theoretic framework for the transfer of the Ramsey property from one category to another. He showed that a \emph{pre-adjunction} between two categories transfers the Ramsey property (which we call the \emph{weak} Ramsey property, see Definition \ref{def:Ramsey.category}), and the proof of Pr\"{o}mel and Voigt in \cite{PV81} shows that there is a pre-adjunction between the category of finite linearly ordered graphs and the Graham--Rothschild category. Solecki took a related approach, deriving Ramsey theorems from a general pigeonhole principle that is preserved under \emph{interpretability} \cite{S13} and, more generally, \emph{modeling} \cite{S25}. Barto\v{s}ov\'{a} and Scow later provided a similar framework in the language of structural Ramsey theory in \cite{S21,BS24}, introducing the notion of a \emph{semi-retraction} between structures, and showed that semi-retractions transfer the Ramsey property between locally finite ordered structures. 

In \cite{T10}, Todor\v{c}evi\'{c} proved the pigeonhole property of several topological Ramsey spaces (Carlson--Simpson space \cite[Lemma 5.69]{T10}, infinite-dimensional vector subspaces \cite[Lemma 5.86]{T10}) by establishing a form of projection from the set of infinite sequences of left variable-words to the topological Ramsey space, then appealing to the infinite Hales--Jewett theorem for left variable-words. Similarly, in \cite[Theorem 4.50]{T10} he defined an explicit map:
\begin{align*}
    \pi_\infty : (\FIN_L^{[\infty]},\FIN_{Lv}^{[\infty]},\leq,\leq^o,r,s) \to (W_L^{[\infty]},W_{Lv}^{[\infty]},\leq,\leq^o,r,s),
\end{align*}
and asserted, without detailed explanation, that $\pi_\infty$ ``transfers the Ramsey-theoretic properties''. In this paper, we formalise this notion of transfer of Ramsey-theoretic properties by establishing a framework consisting of various types of maps between topological Ramsey spaces that transfer the Ramsey-theoretic properties of one to another. In much of this framework, we work with a far simpler structure that we call an \emph{analytic quasi-order} (Definition \ref{def:analytic.quasi-order}), i.e. a quasi-order $(\cal{R},\leq)$ such that $\cal{R}$ is Polish and ${\leq} \subseteq \cal{R}^2$ is analytic. We wish to transfer the analytic Ramsey property between analytic quasi-orders, i.e. the property that every analytic subset of $\cal{R}$ is Ramsey (Definition \ref{def:Ramsey.quasi-order}). A \emph{Ramsey quasi-order} is an analytic quasi-order with the analytic Ramsey property, and by the abstract Ellentuck theorem, every topological Ramsey space with a countable set of finite approximations is a Ramsey quasi-order (Lemma \ref{lem:A4.iff.local.Ramsey.properties}). Working with monotone maps between analytic quasi-orders, we identify two properties of interest: (continuous and cofinal) \emph{quasi-open} maps (Definition \ref{def:quasi-open}) transfer the analytic Ramsey property (Proposition \ref{prop:cofinal.qo.transfers.Ramsey}), while Alaoglu--Birkhoff convergent maps (Definition \ref{def:ab.convergent}) transfer generic filters (Lemma \ref{lem:ab.equivalence}). When the target analytic quasi-order satisfies the Ramsey property, under suitable structural assumptions on the analytic quasi-orders, these two properties coincide for a monotone map (Lemma \ref{lem:ab.iff.qo.TRS}). We remark that it is possible to transfer the (weak) analytic Ramsey property via non-monotone maps --- either via a family of maps, one into each cone of the target (Lemma \ref{lem:cts.qo.family.transfers.analytic.Ramsey.property}), or via a single map at the cost of a homogeneity assumption on the target quasi-order (Proposition \ref{prop:cts.qo.to.H.transfers.analytic.Ramsey.property}). A general framework that does not assume monotonicity appears to be considerably more involved, and we leave it as a direction for future work.

The main results of this paper consist of the (non-)existence of \emph{weak projections} (or \emph{WP}, i.e. monotone Alaoglu--Birkhoff convergent maps, see Definition \ref{def:weak.projection}) between topological Ramsey spaces. Natural examples of \emph{strong projections} (or \emph{SP}, see Definition \ref{def:strong.projection}) are plentiful --- we give a detailed collection of strong projections in \S\ref{sec:projection.examples}, most of which were studied in \cite{T10}. Dobrinen and Todor\v{c}evi\'{c} also considered a natural projection from $\cal{R}_1$ to $[\N]^\infty$ in \cite{DT14} (Example \ref{ex:Dobrinen.spaces.projections}).

\begin{theorem}
\label{thm:main}
    We have the following diagram of (non)-existence of projections:
    
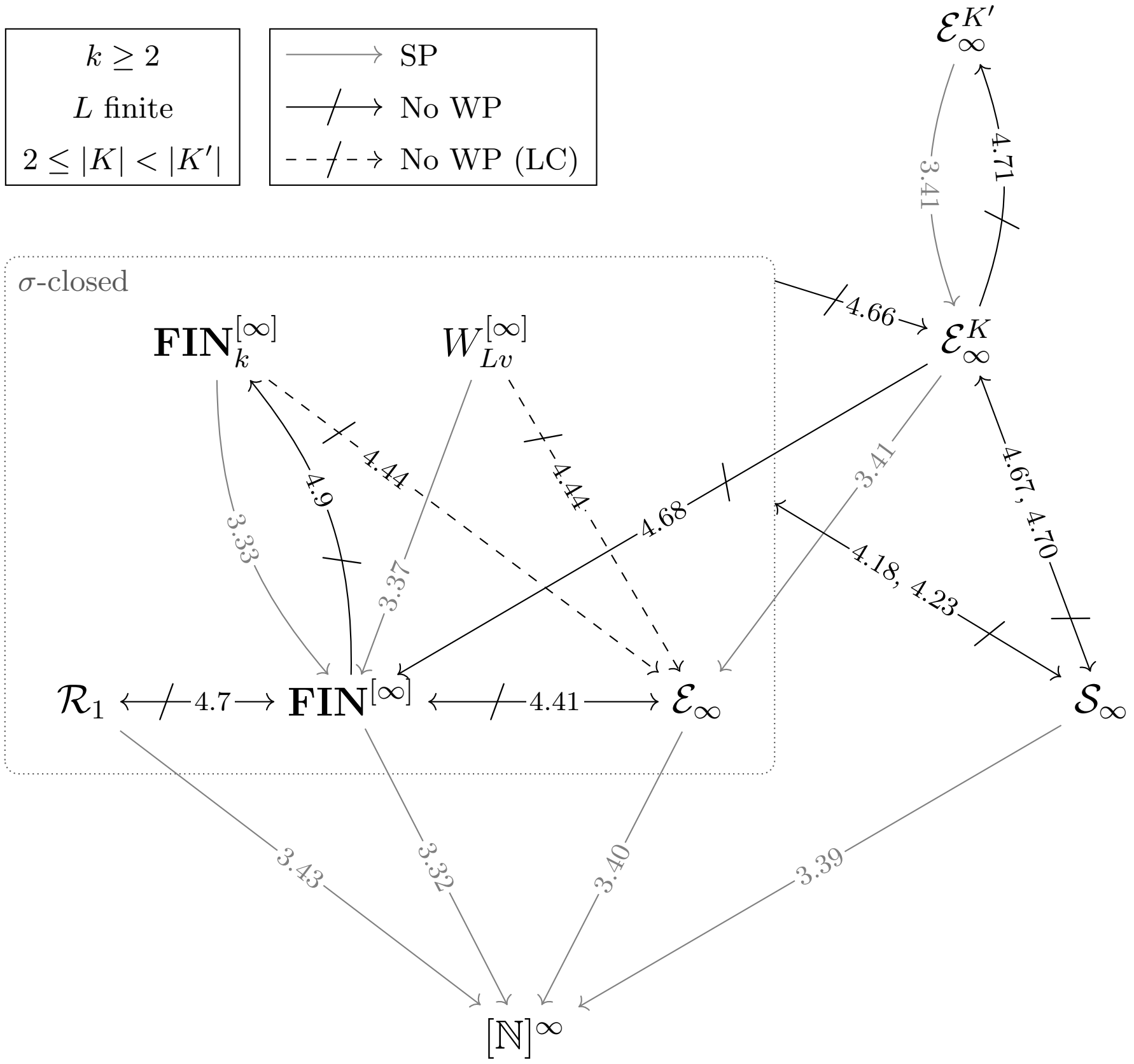
\begin{figure}[H]
        \centering
        \normalfont
        \begin{tikzpicture}[
            y=1.25cm,
            space/.style={font=\large},
            proj/.style={->, gray},
            slash/.style={sloped, inner sep=0pt},
            lbl/.style={font=\scriptsize, fill=white, inner sep=1pt, sloped},
            plbl/.style={lbl, text=gray}
        ]
            \node[space] (Ell) at (0,0) {$[\N]^\infty$};
            \node[space] (R1) at (-4.6,2.8) {$\cal{R}_1$};
            \node[space] (FIN) at (-1.8,2.8) {$\FIN^{\smash[t]{[\infty]}}$};
            \node[space] (FINk) at (-3.2,5.8) {$\FIN_k^{[\infty]}$};
            \node[space] (W) at (-0.4,5.8) {$W_{Lv}^{[\infty]}$};
            \node[space] (E) at (1.8,2.8) {$\E_\infty$};
            \node[space] (EK) at (4.6,5.8) {$\E_\infty^K$};
            \node[space] (EKp) at (4.6,8.4) {$\E_\infty^{K'}$};
            \node[space] (S) at (6.0,2.8) {$\S_\infty$};

            \draw[black!65, densely dotted, rounded corners] (-5.4,2.2) rectangle (2.6,6.5);
            \node[black!65, font=\footnotesize, anchor=north west] at (-5.4,6.5) {$\sigma$-closed};
            \node[draw, font=\footnotesize, inner sep=0pt, anchor=north west] (assumptions) at (-5.4,8.4) {\renewcommand{\arraystretch}{1.4}\begin{tabular}{@{\hspace{5pt}}c@{\hspace{5pt}}} $k \geq 2$ \\ $L$ finite \\ $2 \leq |K| < |K'|$ \end{tabular}};
            \node[draw, font=\footnotesize, inner sep=0pt, anchor=south west] (legend) at ([xshift=0.3cm]assumptions.south east) {\renewcommand{\arraystretch}{1.4}\begin{tabular}{@{\hspace{5pt}}p{1cm}@{\hspace{5pt}}l@{\hspace{5pt}}} & SP \\ & No WP \\ & No WP (LC) \end{tabular}};
            \path (legend.north west) -- (legend.south west) coordinate[pos=1/6] (l1) coordinate[pos=1/2] (l2) coordinate[pos=5/6] (l3);
            \draw[proj] ([xshift=5pt]l1) -- ++(1,0);
            \draw[->] ([xshift=5pt]l2) -- node[slash] {$/$} ++(1,0);
            \draw[->, dashed] ([xshift=5pt]l3) -- node[slash] {$/$} ++(1,0);

            \draw[proj] (R1) -- node[plbl] {\ref{ex:Dobrinen.spaces.projections}} (Ell);
            \draw[proj] (FIN) -- node[plbl] {\ref{ex:min.max.FIN}} (Ell);
            \draw[proj] (E) -- node[plbl] {\ref{ex:minimal.representatives}} (Ell);
            \draw[proj] (S) -- node[plbl] {\ref{ex:levels.strong.subtrees}} (Ell);
            \draw[proj] (W) -- node[plbl, pos=0.7] {\ref{ex:support.variable.words}} (FIN);
            \draw[proj] (FINk) to[bend right=20] node[plbl] {\ref{ex:tetris.operation}} (FIN);
            \draw[proj] (EK) -- node[plbl, pos=0.3] {\ref{ex:CS-variables.merging.alphabets}} (E);
            \draw[proj] (EKp) to[bend right=20] node[plbl] {\ref{ex:CS-variables.merging.alphabets}} (EK);

            \draw[<->] (R1) -- node[slash, pos=0.3] {$/$} node[lbl, pos=0.6] {\ref{thm:no.wp.R_alpha.and.FIN_k}} (FIN);
            \draw[<->] (FIN) -- node[slash, pos=0.3] {$/$} node[lbl, pos=0.55] {\ref{thm:no.wp.Carlson-Simpson.and.FIN}} (E);
            \draw[->] (FIN) to[bend right=20] node[slash, pos=0.35] {$/$} node[lbl, pos=0.6] {\ref{thm:no.wp.FIN.FIN_2}} (FINk);
            \draw[->, dashed] (FINk) -- node[slash, pos=0.18] {$/$} node[lbl, pos=0.3] {\ref{thm:no.wp.unconditional.Carlson-Simpson}} (E);
            \draw[->, dashed] (W) -- node[slash, pos=0.2] {$/$} node[lbl, pos=0.36] {\ref{thm:no.wp.unconditional.Carlson-Simpson}} (E);
            \draw[->] (2.6,6.3) -- node[slash, pos=0.4] {$/$} node[lbl, pos=0.63] {\ref{thm:no.wp.sigma-closed.CS-variables}} (EK);
            \draw[->] (EK) -- node[slash, pos=0.38] {$/$} node[lbl, pos=0.5] {\ref{thm:no.wp.CS-variables.FIN}} (FIN);
            \draw[<->] (S) -- node[slash, pos=0.25] {$/$} node[lbl, pos=0.55] {\ref{thm:no.wp.sigma-closed.trees}, \ref{cor:no.wp.strong.subtrees.sigma-closed.examples}} (2.6,4.45);
            \draw[<->] (S) -- node[slash, pos=0.18] {$/$} node[lbl, pos=0.57] {\ref{cor:no.wp.strong.subtrees.CS-variables}, \ref{thm:no.wp.CS-variables.trees}} (EK);
            \draw[->] (EK) to[bend right=20] node[slash, pos=0.35] {$/$} node[lbl, pos=0.62] {\ref{thm:CS-variables.alphabet.hierarchy}} (EKp);
        \end{tikzpicture}
        \label{fig:main}
    \end{figure}
    Here, LC denotes the existence of large cardinals sufficient for $\mathbf{\Sigma}_4^1$ generic absoluteness. Furthermore, $[\N]^\infty$ does not weakly project into any other topological Ramsey space in the diagram. 
\end{theorem}

Readers may refer to \S\ref{sec:non-existence.of.projections} for more non-existence results not mentioned above. It is worth noting that the proofs of the results vary significantly depending on the spaces of interest. Theorems \ref{thm:no.wp.Ellentuck.FIN}, \ref{thm:wp.Dobrinen.spaces.limit}, \ref{thm:no.wp.R_alpha.and.FIN_k} and \ref{thm:no.wp.FIN.FIN_2} are proved by an analysis of the initial Tukey structure of the generic Ramsey ultrafilters (as in Definition \ref{def:Ramsey.filter}) obtained by forcing with various spaces. The proofs of Theorems \ref{thm:no.wp.strong.subtrees.Sacks} and \ref{thm:no.wp.strong.subtrees.sigma-closed.projection} are forcing-theoretic and require us to prove the real degree minimality of forcing with strong subtrees (Lemma \ref{lem:strong.subtrees.forcing.minimal.real}). The proof of the $\E_\infty \nrightarrow \FIN^{[\infty]}$ direction of Theorem \ref{thm:no.wp.Carlson-Simpson.and.FIN} is purely combinatorial, while the converse direction (and hence Corollary \ref{cor:no.wp.Ellentuck.Carlson-Simpson}) relies on \"{O}zalp's classification theorem on the initial Tukey structure of stable ordered union ultrafilters up to Rudin--Keisler equivalence in \cite{O24} (Theorem \ref{thm:Tukey.below.souu.RK.classification}). Theorem \ref{thm:no.wp.unconditional.Carlson-Simpson} follows from a consistency result, due to Brendle, about a cardinal invariant inequality (Theorem \ref{thm:Carlson-Simpson.distributivity.omega_1}, \cite{B00}), followed by an appeal to the $\mathbf{\Pi}_4^1$ complexity of weak projections (hence the large cardinal hypothesis). Finally, Theorems \ref{thm:no.wp.CS-variables.trees} and \ref{thm:CS-variables.alphabet.hierarchy} are also forcing-theoretic, and require an analysis of the intermediate models of the generic extension by $(\E_\infty^L,\leq)$, as well as a canonisation theorem on analytic equivalence relations on $\E_\infty^L$ by Doucha \cite{D13}.

As hinted by the several examples of weak projections above, most topological Ramsey spaces are equipped with a natural projection to the Ellentuck space. Dobrinen conjectured that every topological Ramsey space contains a copy of the Ellentuck space \cite[\S2.1]{NF20}. In terms of projections, this asserts that every topological Ramsey space projects into the Ellentuck space, which fails in full generality, as $(\mathbb{S},\subseteq,r)$ does not weakly project into the Ellentuck space (Theorem \ref{thm:no.wp.Sacks.sigma-closed}). Such projections are usually first defined on finite approximations and respect initial segments, giving rise to a natural map on elements of $\cal{R}$. In \S\ref{sec:strong.projections}, we consider the notion of \emph{strong projections} (Definition \ref{def:strong.projection}). Strong projections are functions that are naturally extended from some map on finite approximations, and are projections on every basic clopen set and Ellentuck neighbourhood. Most natural projections are strong projections, including every example provided in \S\ref{sec:projection.examples}. We prove the following dichotomy result for \textbf{A2}-spaces (i.e. a closed triple $(\cal{R},\leq,r)$ satisfying Todor\v{c}evi\'{c}'s axioms \textbf{A1}--\textbf{A3}), confirming Dobrinen's conjecture for the atomless spaces in this class:

\begin{restatable*}{theorem}{thmdichotomy}
\label{thm:dichotomy}
    If $(\cal{R},\leq,r)$ is a depth-preserving depth-uniform \textbf{A2}-space with a $\leq$-maximal element $\mathbb{A}$, then exactly one of the following two cases holds:
    \begin{enumerate}
        \item $(\cal{R},\leq,r)$ strongly projects into the Ellentuck space.

        \item $(\cal{R},\leq,r)$ has a dense set of atoms. In particular, $(\cal{R},\leq,r)$ does not weakly project into the Ellentuck space.
    \end{enumerate}
\end{restatable*}

Wang showed in \cite{W25} that the set of codes for Ramsey positive $G_\delta$ subsets of a well-indexed topological Ramsey space is $\mathbf{\Sigma}_2^1$-complete, generalising Sabok's result in \cite{Sab12} for the Ellentuck space. As an application, we show that Wang's result applies to a class of depth-preserving depth-uniform topological Ramsey spaces (Corollary \ref{cor:dpdu.DR-nontrivial.implies.complete.codes}).

\subsection{Structure of this paper}
Our paper is organised as follows:
\begin{itemize}
    \item[(\S\ref{sec:Ramsey.quasi-orders.and.spaces})] We recap the axioms and notation of topological Ramsey spaces. We also consider spaces satisfying a fragment of these axioms, called \emph{(weak) \textbf{A2}-spaces}. In \S\ref{sec:examples}, we provide examples of (weak) \textbf{A2}-spaces (and topological Ramsey spaces). We introduce the notion of a \emph{Ramsey quasi-order} (Definition \ref{def:Ramsey.quasi-order}), which is an analytic quasi-order satisfying Silver's infinite-dimensional Ramsey theorem for analytic sets. We show that a topological Ramsey space satisfying Todor\v{c}evi\'{c}'s axioms \textbf{A1}--\textbf{A4} is a Ramsey quasi-order if $\cal{AR}$ is countable (Lemma \ref{lem:A4.iff.local.Ramsey.properties}). In \S\ref{sec:Ramsey.ultrafilters}, we consider a generalised notion of Ramsey ultrafilters (Definition \ref{def:Ramsey.filter}), which makes sense for any Ramsey quasi-order. We prove that for ``$\sigma$-closed'' Ramsey quasi-orders (Definition \ref{def:sigma-closed.quasi-order}), Ramsey ultrafilters exist under $\CH$ (Theorem \ref{thm:CH.implies.Ramsey.ultrafilter}), and forcing with the Ramsey quasi-order adds a Ramsey ultrafilter.
    
    \item[(\S\ref{sec:projections})] We consider the \emph{Alaoglu--Birkhoff convergent} and \emph{quasi-open} properties of maps between Ramsey quasi-orders, and show that continuous cofinal monotone quasi-open maps transfer the analytic Ramsey property (Proposition \ref{prop:cofinal.qo.transfers.Ramsey}). We show that monotone quasi-open maps are Alaoglu--Birkhoff convergent, hence transfer generic filters (Lemma \ref{lem:ab.equivalence}). We also show that under suitable assumptions, these two properties are equivalent (Lemma \ref{lem:ab.iff.qo.TRS}). In \S\ref{sec:projection.examples}, we include an extensive list of examples of such maps between topological Ramsey spaces.
    
    \item[(\S\ref{sec:non-existence.of.projections})] We prove the non-existence of weak projections between some topological Ramsey spaces introduced in \S\ref{sec:examples} (cf. Theorem \ref{thm:main}).
    
    \item[(\S\ref{sec:strong.projections})] We introduce and study the notion of \emph{strong projections} (Definition \ref{def:strong.projection}), and prove Theorem \ref{thm:dichotomy}. 
    
    \item[(\S\ref{sec:Ramsey.categories})] We discuss the relationship between Axiom \textbf{A4} of topological Ramsey spaces and the Ramsey property of categories.
\end{itemize}

\section{Ramsey quasi-orders and topological Ramsey spaces}
\label{sec:Ramsey.quasi-orders.and.spaces}
We recap the axioms of topological Ramsey spaces, and highlight several examples of topological Ramsey spaces that will be repeatedly discussed in this paper. These axioms also motivate the notion of \emph{Ramsey quasi-orders}, which, together with topological Ramsey spaces, form the two main frameworks of interest in this paper.  

\subsection{Axioms of topological Ramsey spaces}
\label{sec:TRS.axioms}
The axiomatisation of topological Ramsey spaces was initially studied by Carlson and Simpson in \cite{CS90}. In \cite{T10}, Todor\v{c}evi\'{c} presented the four axioms, labelled \textbf{A1}--\textbf{A4}, which are sufficient for a closed triple $(\cal{R},\leq,r)$ to be a topological Ramsey space. Since the introduction of these axioms, there has been extensive interest in generalising known results in Ramsey theory to topological Ramsey spaces --- see, for instance, \cite{DMN17,DN22,Y25,CD25}.

\subsubsection{Axioms \textbf{A1}--\textbf{A4}} Here, $\cal{R}$ is a non-empty set, $\leq$ is a quasi-order on $\cal{R}$, and $r : \cal{R} \times \omega \to \cal{AR}$ is a surjective function. We also define a sequence of maps $r_n : \cal{R} \to \cal{AR}$ by $r_n(A) := r(A,n)$ for all $A \in \cal{R}$. Let $\cal{AR}_n \subseteq \cal{AR}$ be the image of $r_n$ (i.e. $a \in \cal{AR}_n$ iff $a = r_n(A)$ for some $A \in \cal{R}$).

\begin{axiom}[\textbf{A1}, Sequencing]
\hfill
    \begin{enumerate}[label=(\arabic*)]
        \item $r_0(A) = \emptyset$ for all $A \in \cal{R}$.

        \item $A \neq B$ implies $r_n(A) \neq r_n(B)$ for some $n$.

        \item $r_n(A) = r_m(B)$ implies $n = m$ and $r_k(A) = r_k(B)$ for all $k < n$.
    \end{enumerate}
\end{axiom}

For each $a \in \cal{AR}$, let $\lh(a)$ denote the unique $n$ such that $a \in \cal{AR}_n$. By Axiom \textbf{A1}(3), this $n$ is well-defined. Axiom \textbf{A1}(3) also implies that if $a \in \cal{AR}_n$ and $k \leq n$, then $r_k(A) = r_k(B)$ for any $A,B \in \cal{R}$ with $r_n(A) = r_n(B) = a$, so we may unambiguously extend $r_k$ to $\{a \in \cal{AR} : \lh(a) \geq k\}$ by letting $r_k(a)$ denote this common value. For $a,b \in \cal{AR}$, we write $a \sqsubseteq b$ iff $a = r_{\lh(a)}(b)$, and for $B \in \cal{R}$, we write $a \sqsubseteq B$ iff $a = r_{\lh(a)}(B)$.

\begin{axiom}[\textbf{A2}, Finitisation]
    There is a quasi-ordering $\leq_\fin$ on $\cal{AR}$ such that:
    \begin{enumerate}[label=(\arabic*)]
        \item $\{a \in \cal{AR} : a \leq_\fin b\}$ is finite for all $b \in \cal{AR}$.

        \item $B \leq A$ iff $\forall n \exists m [r_n(B) \leq_\fin r_m(A)]$.

        \item $\forall a,b,c \in \cal{AR}[a \sqsubseteq b \wedge b \leq_\fin c \to \exists d \sqsubseteq c[a \leq_\fin d]]$.
    \end{enumerate}
    While the original presentation of the axioms in \cite{T10} does not assume this, we shall restrict our attention to triples satisfying the following slightly stronger variant of \textbf{A2}(3):
    \begin{enumerate}
        \item[(3')] $\forall a,b,c \in \cal{AR}[a \sqsubsetneq b \wedge b \leq_\fin c \to \exists d \sqsubsetneq c[a \leq_\fin d]]$.
    \end{enumerate}
\end{axiom}

We may define the \emph{Ellentuck neighbourhoods} as follows: For any $A \in \cal{R}$, $a \in \cal{AR}$ and $n \in \N$, we let:
\begin{align*}
    [a,A] &:= \{B \in \cal{R} : a \sqsubseteq B \wedge B \leq A\}, \\
    [n,A] &:= [r_n(A),A], \\
    [A] &:= [\emptyset,A], \\
    [a] &:= \{B \in \cal{R} : a \sqsubseteq B\}.
\end{align*}
We also define the depth function $\depth_A$, for $A \in \cal{R}$ and $a \in \cal{AR}$, by:
\begin{align*}
    \depth_A(a) := 
    \begin{cases}
        \min\{n < \omega : a \leq_\fin r_n(A)\}, &\text{if such $n$ exists}, \\
        \infty, &\text{otherwise}.
    \end{cases}
\end{align*}

We also define, for each $A \in \cal{R}$ and $n < \omega$:
\begin{align*}
    \cal{AR}\restrictedto A &:= \{a \in \cal{AR} : \exists m[a \leq_\fin r_m(A)]\}, \\
    \cal{AR}_n\restrictedto A &:= \cal{AR}_n \cap \cal{AR}\restrictedto A.
\end{align*}
If $a \in \cal{AR}\restrictedto A$, then:
\begin{align*}
    \cal{AR}\restrictedto[a,A] &:= \{b \in \cal{AR}\restrictedto A : a \sqsubseteq b\}, \\
    \cal{AR}\restrictedto[a] &:= \{b \in \cal{AR} : a \sqsubseteq b\}, \\
    \cal{AR}\restrictedto a &:= \{b \in \cal{AR} : b \leq_\fin a\},  \\
    \cal{AR}_n\restrictedto[a,A] &:= \cal{AR}_n \cap \cal{AR}\restrictedto[a,A], \\
    \cal{AR}_n\restrictedto[a] &:= \cal{AR}_n \cap \cal{AR}\restrictedto[a], \\
    \cal{AR}_n\restrictedto a &:= \cal{AR}_n \cap \cal{AR}\restrictedto a, \\
    r_n[a,A] &:= \{b \in \cal{AR}\restrictedto[a,A] : \lh(b) = n\}, \\
    r_n[A] &:= r_n[\emptyset,A].
\end{align*}
Note that for all $a \in \cal{AR}$, $\depth_A(a) < \infty$ iff $a \in \cal{AR}\restrictedto A$. We remark that \textbf{A2}(3') is equivalent to the statement that $\depth_A$ is \emph{strictly increasing} for all $A \in \cal{R}$ --- that is, for all $a,b \in \cal{AR}\restrictedto A$, if $a \sqsubsetneq b$ then $\depth_A(a) < \depth_A(b)$.

\begin{axiom}[\textbf{A3}, Amalgamation]
    Let $A \in \cal{R}$.
    \begin{enumerate}[label=(\arabic*)]
        \item For all $a \in \cal{AR}\restrictedto A$ and $B \in [\depth_A(a),A]$, $[a,B] \neq \emptyset$.

        \item If $B \leq A$ and $a \in \cal{AR}\restrictedto B$, then there exists $B' \in [\depth_A(a),A]$ such that $\emptyset \neq [a,B'] \subseteq [a,B]$.
    \end{enumerate}
\end{axiom}

\begin{axiom}[\textbf{A4}, Pigeonhole]
    For all $A \in \cal{R}$, $a \in \cal{AR}\restrictedto A$ and $\O \subseteq \cal{AR}_{\lh(a)+1}$, there exists some $B \in [\depth_A(a),A]$ such that $r_{\lh(a)+1}[a,B] \subseteq \O$ or $r_{\lh(a)+1}[a,B] \subseteq \O^c$.
\end{axiom}

Note that, by \textbf{A3}, the $B$ in \textbf{A4} may instead be taken to be an element of $[a,A]$. By Axiom \textbf{A1}, we may identify each element $A \in \cal{R}$ with a sequence of elements of $\cal{AR}$, via the map $A \mapsto (r_n(A))_{n<\omega}$. Therefore, we may identify $\cal{R}$ with a subset of $\cal{AR}^\N$.

We shall focus on the following topologies on $\cal{R}$.
\begin{itemize}
    \item The \emph{metric topology} induced by the first difference metric, where for $A,B \in \cal{R}$, $d(A,A) = 0$ and $d(A,B) = \frac{1}{2^n}$ for $A \neq B$, where $n$ is the least integer such that $r_n(A) \neq r_n(B)$. The sets of the form $[a]$ for $a \in \cal{AR}$ form a neighbourhood base of this topology. If $\cal{AR}$ is countable, then $\cal{R}$ is a Polish space under this metric topology.

    \item The \emph{Ellentuck topology} generated by open sets of the form $[a,A]$ for $A \in \cal{R}$ and $a \in \cal{AR}\restrictedto A$.

    \item The \emph{order topology} generated by open sets of the form $[A]$ for $A \in \cal{R}$.
\end{itemize}
Unless stated otherwise, all topological properties of subsets of $\cal{R}$ mentioned in this paper are with respect to the metric topology, except that we reserve the term \emph{dense open} for the order topology. That is, a subset $\D \subseteq \cal{R}$ is dense open (under the order topology) if for all $A \in \cal{R}$, there exists some $B \in \D$ such that $B \leq A$, and for all $A \in \D$ and $B \leq A$, $B \in \D$.

\begin{definition}
    We say that $(\cal{R},\leq,r)$ is a \emph{closed triple} if $\cal{R}$ is a metrically closed subset of $\cal{AR}^\N$.
\end{definition}

\begin{definition}
    Let $(\cal{R},\leq,r)$ be a closed triple satisfying \textbf{A1} and \textbf{A2}. A subset $\X \subseteq \cal{R}$ is \emph{Ramsey} if for all $A \in \cal{R}$ and $a \in \cal{AR}\restrictedto A$, there exists some $B \in [a,A]$ such that $[a,B] \subseteq \X$ or $[a,B] \subseteq \X^c$. If $[a,B] \subseteq \X^c$ always holds, then $\X$ is \emph{Ramsey null}.
\end{definition}

Todor\v{c}evi\'{c} subsequently showed in \cite{T10} the following theorem:

\begin{theorem}[{\cite[Abstract Ellentuck theorem]{T10}}]
\label{thm:abstract.Ellentuck.thm}
    Let $(\cal{R},\leq,r)$ be a closed triple satisfying \textbf{A1}--\textbf{A4}. Then a subset of $\cal{R}$ is Ramsey iff it has the property of Baire under the Ellentuck topology, and it is Ramsey null iff it is meagre under the Ellentuck topology. 
\end{theorem}

In particular, if $\cal{AR}$ is countable (so the metric topology is Polish), then every analytic set is Ramsey. A closed triple $(\cal{R},\leq,r)$ satisfying \textbf{A1} and \textbf{A2} is a \emph{topological Ramsey space} if it satisfies the conclusion of Theorem \ref{thm:abstract.Ellentuck.thm}. Theorem \ref{thm:abstract.Ellentuck.thm} thus asserts that every closed triple satisfying \textbf{A1}--\textbf{A4} is a topological Ramsey space.

On the other hand, the second author noted in \cite{Y25} that many examples satisfy most, but not all, of the axioms \textbf{A1}--\textbf{A4}, yet admit many interesting Ramsey-theoretic results. He thus introduced the notion of \emph{weak \textbf{A2}-spaces}, which generalise topological Ramsey spaces to a larger class of spaces, most notably countable vector spaces (see Example \ref{ex:vector.spaces}). He showed in the same paper that Kastanas Ramsey sets of weak \textbf{A2}-spaces serve as a common generalisation of Ramsey subsets of topological Ramsey spaces, and strategically Ramsey subsets of vector spaces introduced by Rosendal in \cite{R10} and \cite{R09}.

\begin{axiom}[\textbf{wA2}]
    There is a quasi-ordering $\leq_\fin$ on $\cal{AR}$ such that:
    \begin{enumerate}
        \item[(w1)] $\{a \in \cal{AR} : a \leq_\fin b\}$ is countable for all $b \in \cal{AR}$.

        \item[(2)] $B \leq A$ iff $\forall n \exists m [r_n(B) \leq_\fin r_m(A)]$.

        \item[(3)] $\forall a,b,c \in \cal{AR}[a \sqsubseteq b \wedge b \leq_\fin c \to \exists d \sqsubseteq c[a \leq_\fin d]]$.
    \end{enumerate}
\end{axiom}

\begin{definition}
    A triple $(\cal{R},\leq,r)$ is said to be:
    \begin{itemize}
        \item a \emph{weak \textbf{A2}-space} (or just \emph{\textbf{wA2}-space}) if $(\cal{R},\leq,r)$ is a closed triple satisfying Axioms \textbf{A1}, \textbf{wA2} and \textbf{A3};

        \item an \emph{\textbf{A2}-space} if it is a \textbf{wA2}-space satisfying \textbf{A2}.
    \end{itemize}
\end{definition}

\subsubsection{Useful lemmas}
We prove several consequences of the axioms \textbf{A1}--\textbf{A4} that we will need in later sections of the paper. The first lemma asserts that the axioms \textbf{A1}--\textbf{A3} allow us to perform diagonalisation over a family of dense open sets. The following definition is a variant of \cite[Definition 3.6]{DMN17}, and was introduced earlier by the second author in \cite{Y25} and \cite{Y26}.

\begin{definition}
\label{def:dense.open.diagonalisation}
    Let $A \in \cal{R}$ and $a \in \cal{AR}\restrictedto A$.
    \begin{itemize}
        \item A family of subsets $\vec{\D} = \{\D_b\}_{b \in \cal{AR}\restrictedto[a,A]}$ is \emph{dense open below $[a,A]$} if for all $b \in \cal{AR}\restrictedto[a,A]$, $\D_b$ is a $\leq$-downward closed subset of $[b,A]$, and for all $B \in [b,A]$, there exists some $C \in [b,B]$ such that $C \in \D_b$.
        
        \item Let $\vec{\D} = \{\D_b\}_{b \in \cal{AR}\restrictedto[a,A]}$ be dense open below $[a,A]$. We say that $B \in [a,A]$ \emph{diagonalises} $\vec{\D}$ if for all $b \in \cal{AR}\restrictedto[a,B]$, there exists some $A_b \in \D_b$ such that $[b,B] \subseteq [b,A_b]$.
    \end{itemize}
\end{definition}

\begin{lemma}
\label{lem:R.is.semiselective.improved}
    Let $(\cal{R},\leq,r)$ be a closed triple satisfying \textbf{A1}--\textbf{A3}, and let $A \in \cal{R}$. Let $\vec{\D} = \{\D_b\}_{b\in\cal{AR}\restrictedto A}$ be a family of subsets that is dense open below $[\emptyset,A]$. Then for all $a \in \cal{AR}\restrictedto A$, there exists some $B \in [a,A]$ such that for all $b \in \cal{AR}\restrictedto B$ with $\depth_B(b) \geq \lh(a)$, there exists some $A_b \in \D_b$ such that $[b,B] \subseteq [b,A_b]$.
\end{lemma}

\begin{proof}
    Fix some $A \in \cal{R}$ and $a \in \cal{AR}\restrictedto A$, and suppose that $\vec{\D} = \{\D_b\}_{b\in\cal{AR}\restrictedto A}$ is dense open below $[\emptyset,A]$. We shall define a fusion sequence $(A_n)_{n<\omega}$ in $[a,A]$, with $a_n = r_{\lh(a)+n}(A_n)$, such that $A_{n+1} \in [a_n,A_n]$: By \textbf{A3}(1), let $A_0 \in [a,A]$, and suppose that $A_n$ has been defined. Let $\{b_n^i : i < N\}$ enumerate the set of all $b \in \cal{AR}\restrictedto A_n$ such that $\depth_{A_n}(b) = \lh(a) + n$, which is finite by \textbf{A2}(1). We define another sequence $(A_n^i)_{i \leq N}$ in $[a_n,A_n]$ as follows: Let $A_n^0 := A_n$. Suppose that $A_n^i \in [a_n,A_n]$ has been defined. Then $r_j(A_n^i) = r_j(A_n)$ for all $j \leq \lh(a)+n$, so $\depth_{A_n^i}(b_n^i) = \depth_{A_n}(b_n^i) = \lh(a)+n$, and by \textbf{A3}(1) we may let $B_n^i \in [b_n^i,A_n^i]$. Since $\D_{b_n^i}$ is dense open in $[b_n^i,A]$, there exists some $A_{b_n^i} \in \D_{b_n^i}$ and $C_n^i \in [b_n^i,B_n^i]$ such that $[b_n^i,C_n^i] \subseteq [b_n^i,A_{b_n^i}]$. By \textbf{A3}(2), we may let $A_n^{i+1} \in [\depth_{A_n^i}(b_n^i),A_n^i] = [a_n,A_n^i]$ be such that $[b_n^i,A_n^{i+1}] \subseteq [b_n^i,C_n^i]$. We complete the induction by letting $A_{n+1} := A_n^N$. Since $\cal{R}$ is closed, $B := \lim_{n\to\infty} a_n$ is an element of $\cal{R}$, and $B \in [a,A]$ and $B \leq A_n^i$ for all $n$ and $i \leq N$ by \textbf{A2}(2). We shall show that $B$ is the required element of $[a,A]$.

    Let $b \in \cal{AR}\restrictedto B$ be such that $\depth_B(b) = \lh(a) + n$ for some $n \geq 0$. Since $r_j(B) = r_j(A_n)$ for all $j \leq \lh(a)+n$, we have that $\depth_{A_n}(b) = \lh(a)+n$, so $b = b_n^i$ for some $i$ in the induction step of the proof. Then $[b,B] \subseteq [b,A_n^{i+1}] \subseteq [b,C_n^i] \subseteq [b,A_b]$ and $A_b \in \D_b$, as desired.
\end{proof}

In particular (as $\depth_B(b) \geq \depth_B(a) = \lh(a)$ for all $b \in \cal{AR}\restrictedto[a,B]$), we have the following:

\begin{lemma}[{\cite[Lemma 3.16]{Y25}}]
\label{lem:R.is.semiselective}
    Let $(\cal{R},\leq,r)$ be an \textbf{A2}-space, and let $A \in \cal{R}$ and $a \in \cal{AR}\restrictedto A$. Then every family of subsets $\vec{\D} = \{\D_b\}_{b \in \cal{AR}\restrictedto[a,A]}$ that is dense open below $[a,A]$ has a diagonalisation.
\end{lemma}

The following lemma will be relevant to our motivation for defining Ramsey quasi-orders.

\begin{lemma}
\label{lem:reduction.is.G_delta}
    For any \textbf{wA2}-space $(\cal{R},\leq,r)$, the set $\leq$ is a $G_\delta$ subset of $\cal{R}^2$.
\end{lemma}

\begin{proof}
    This follows from the observation that, by \textbf{A2}:
    \begin{align*}
        {\leq} = \bigcap_{n<\omega} \bigcup_{m<\omega} \{(A,B) : r_n(A) \leq_\fin r_m(B)\},
    \end{align*}
    and the set $\{(A,B) : r_n(A) \leq_\fin r_m(B)\}$ is open in the metric topology of $\cal{R}^2$.
\end{proof}

Finally, we prove an abstract form of the combinatorial forcing used in \cite[\S4.3]{T10}, in which the notions of acceptance and rejection are left unspecified. Note that a localised variant of abstract combinatorial forcing was introduced for semiselective coideals of topological Ramsey spaces in \cite{DMN17} (i.e. coideals that satisfy a slightly weaker form of Lemma \ref{lem:R.is.semiselective}). Readers may also refer to \cite[Lemmas 6.6--6.9]{Y26} for a similar combinatorial forcing lemma, but localised to a semiselective coideal.

\begin{theorem}[Abstract combinatorial forcing]
\label{thm:abstract.combinatorial.forcing}
    Let $(\cal{R},\leq,r)$ be a closed triple satisfying \textbf{A1}--\textbf{A4}. Suppose that, for any $A \in \cal{R}$ and $a \in \cal{AR}$, we have defined the statements ``$A$ \emph{accepts} $a$'' and ``$A$ \emph{rejects} $a$'' (such that no $A$ both accepts and rejects $a$, and $A$ \emph{decides} $a$ iff $A$ accepts or rejects $a$), and the following statements hold:
    \begin{enumerate}
        \item If $a \notin \cal{AR}\restrictedto A$, then $A$ accepts $a$.

        \item If $A$ accepts $a$ and $B \leq A$, then $B$ accepts $a$.

        \item If $A$ rejects $a$, $B \leq A$ and $a \in \cal{AR}\restrictedto B$, then $B$ rejects $a$.

        \item If $A$ decides $a$ and $B \leq A$, then $B$ decides $a$.

        \item If $a \in \cal{AR}\restrictedto A$, then $B$ decides $a$ for some $B \in [a,A]$.

        \item If $A$ decides $a$, then for all $B \in [a,A]$, $B$ decides $a$ in the same way $A$ does.

        \item $A$ accepts $a$ if and only if $A$ accepts every $b \in r_{\lh(a)+1}[a,A]$.
    \end{enumerate}
    Then for any $A \in \cal{R}$ and $a \in \cal{AR}\restrictedto A$, if $A$ rejects $a$, then there exists some $B \in [a,A]$ that rejects every $b \in \cal{AR}\restrictedto[a,B]$.
\end{theorem}

\begin{proof}
    We fix some $A \in \cal{R}$ and $a \in \cal{AR}\restrictedto A$.

    \begin{claim}
    \label{claim:cf.no.reduction.accepts.all.ope}
        If $A$ rejects $a$, then there is no $B \in [a,A]$ that accepts every $b \in r_{\lh(a)+1}[a,B]$.
    \end{claim}

    \begin{midproof}
        Suppose, towards a contradiction, that some $B \in [a,A]$ accepts every $b \in r_{\lh(a)+1}[a,B]$. Since $a \in \cal{AR}\restrictedto B$, Property (7) implies that $B$ accepts $a$. On the other hand, $A$ rejects $a$, $B \leq A$ and $a \in \cal{AR}\restrictedto B$, so $B$ rejects $a$ by Property (3). This contradicts that no element both accepts and rejects $a$.
    \end{midproof}

    \begin{claim}
    \label{claim:cf.reduction.decides}
        For all $A \in \cal{R}$ and $a \in \cal{AR}\restrictedto A$, there exists some $B \in [a,A]$ that decides every $b \in \cal{AR}\restrictedto[a,B]$.
    \end{claim}

    \begin{midproof}
        For each $b \in \cal{AR}\restrictedto[a,A]$, let:
        \begin{align*}
            \D_b := \{B \in [b,A] : \text{$B$ decides $b$}\}.
        \end{align*}
        By Property (5), each $\D_b$ is dense in $[b,A]$, and by Property (4) it is $\leq$-downward closed in $[b,A]$, so $\{\D_b\}_{b\in\cal{AR}\restrictedto[a,A]}$ is dense open below $[a,A]$. Thus, by Lemma \ref{lem:R.is.semiselective}, there exists some $B \in [a,A]$ which diagonalises $\{\D_b\}_{b\in\cal{AR}\restrictedto[a,A]}$ --- that is, for all $b \in \cal{AR}\restrictedto[a,B]$, $[b,B] \subseteq [b,A_b]$ for some $A_b \in \D_b$. Since $A_b$ decides $b$ and $B \leq A_b$, $B$ also decides $b$ by Property (4) (in the same way as $A_b$, by Property (6)).
    \end{midproof}

    \begin{claim}
    \label{claim:cf.reduction.rejects.all.ope}
        If $A$ rejects $a$, then there exists some $B \in [a,A]$ that rejects every $b \in r_{\lh(a)+1}[a,B]$.
    \end{claim}

    \begin{midproof}
        By Claim \ref{claim:cf.reduction.decides}, let $A' \in [a,A]$ decide every $b \in r_{\lh(a)+1}[a,A']$. Define:
        \begin{align*}
            \O := \{b \in r_{\lh(a)+1}[a,A'] : \text{$A'$ rejects $b$}\}.
        \end{align*}
        By \textbf{A4}, there exists some $B \in [a,A']$ such that $r_{\lh(a)+1}[a,B] \subseteq \O$ or $r_{\lh(a)+1}[a,B] \subseteq \O^c$. It remains to show that the latter case is impossible. If $r_{\lh(a)+1}[a,B] \subseteq \O^c$, then for every $b \in r_{\lh(a)+1}[a,B]$, $A'$ does not reject $b$. Since $A'$ decides $b$, $A'$ accepts $b$, so $B$ also accepts $b$ by Property (2). This means that $B \in [a,A]$ accepts every $b \in r_{\lh(a)+1}[a,B]$, contradicting Claim \ref{claim:cf.no.reduction.accepts.all.ope} as $A$ rejects $a$.
    \end{midproof}

    We may now finish the proof. Again, by Claim \ref{claim:cf.reduction.decides}, we may let $A' \in [a,A]$ decide every $b \in \cal{AR}\restrictedto[a,A']$. For each $b \in \cal{AR}\restrictedto[a,A']$, we let:
    \begin{align*}
        \D_b := \{B \in [b,A'] : \text{$B$ accepts $b$ or $B$ rejects every $b' \in r_{\lh(b)+1}[b,B]$}\}.
    \end{align*}
    By Claim \ref{claim:cf.reduction.rejects.all.ope} (using that $A'$ decides every $b$, so that each $B \in [b,A']$ decides $b$ by Property (4)), each $\D_b$ is dense in $[b,A']$, and by Properties (2) and (3) it is $\leq$-downward closed in $[b,A']$. Therefore, $\{\D_b\}_{b\in\cal{AR}\restrictedto[a,A']}$ is dense open below $[a,A']$. By Lemma \ref{lem:R.is.semiselective}, we may let $B \in [a,A']$ diagonalise $\{\D_b\}_{b\in\cal{AR}\restrictedto[a,A']}$. By Properties (2) and (3), this implies that for all $b \in \cal{AR}\restrictedto[a,B]$, $B$ accepts $b$ or $B$ rejects every $b' \in r_{\lh(b)+1}[b,B]$. We shall show, by inducting on $\lh(b)$, that $B$ rejects $b$ for every $b \in \cal{AR}\restrictedto[a,B]$.

    For the base case, the only $b \in \cal{AR}\restrictedto[a,B]$ with $\lh(b) = \lh(a)$ is $b = a$. Since $A$ rejects $a$ and $B \in [a,A]$, $B$ rejects $a$ by Property (3). For the inductive step, suppose that $B$ rejects $b \in \cal{AR}\restrictedto[a,B]$. Since no element both accepts and rejects $b$, $B$ does not accept $b$, so $B$ rejects every $b' \in r_{\lh(b)+1}[b,B]$, as $B$ decides these $b'$ and $B$ diagonalises $\{\D_b\}_{b\in\cal{AR}\restrictedto[a,A']}$. This completes the induction, and $B \in [a,A]$ rejects every $b \in \cal{AR}\restrictedto[a,B]$, as required.
\end{proof}

\subsection{Examples}
\label{sec:examples}
We discuss several examples of \textbf{wA2}-spaces. All examples below, except for countable vector spaces (Example \ref{ex:vector.spaces}), Sacks forcing (Example \ref{ex:Sacks.forcing}), $\cal{R}_\alpha$ (Example \ref{ex:Dobrinen.space}) and the singleton space (Example \ref{ex:singleton.space}), may be found in \cite{T10}. Discussions of countable vector spaces and strategically Ramsey sets may be found in \cite{R10}, \cite{R09}, \cite{S18} and \cite{R20}, and of $\{(\cal{R}_\alpha,\leq_\alpha,r^\alpha) : \alpha < \omega_1\}$ may be found in \cite{DT14} and \cite{DT15}. In several of the examples below, we define only the finitisation $\leq_\fin$, and the quasi-order $\leq$ on $\cal{R}$ is the one determined by $\leq_\fin$ via \textbf{A2}(2), i.e. $B \leq A$ iff $\forall n \exists m[r_n(B) \leq_\fin r_m(A)]$.

\begin{example}[Natural numbers/Ellentuck space $\lbrack \N\rbrack^\infty$]
\label{ex:Ellentuck.space}
    Let $\cal{R} = [\N]^\infty$. For each $n$ and $A \in [\N]^\infty$, let $r_n(A)$ be the finite set containing the $n$ least elements of $A$ (so $\cal{AR} = [\N]^{<\infty}$). Let $\leq_\fin$ denote the subset relation $\subseteq$. It is easy to check that $([\N]^\infty,\subseteq,r)$ is a closed triple satisfying \textbf{A1}--\textbf{A4}.
\end{example}

\begin{example}[Hindman's space $\FIN^{[\infty]}$]
\label{ex:Hindman.FIN}
    Let $\FIN := [\N]^{<\infty} \setminus \{\emptyset\}$ be the set of non-empty finite subsets of $\N$. For $x,y \in \FIN$, we define:
    \begin{align*}
        x < y \iff \max(x) < \min(y).
    \end{align*}
    Note that $<$ is a partial order on $\FIN$. We let $\cal{R} := \FIN^{[\infty]}$ be the set of all infinite $<$-increasing sequences of elements of $\FIN$, and for each $A = (x_n)_{n<\omega} \in \FIN^{[\infty]}$, $r_N(A) := (x_n)_{n<N}$. We have that $\cal{AR} = \FIN^{[<\infty]}$, the set of all finite $<$-increasing sequences.

    Given $a = (x_n)_{n<N} \in \FIN^{[<\infty]}$, we let:
    \begin{align*}
        \c{a} := \{x_{n_0} \cup \cdots \cup x_{n_k} : 0 \leq n_0 < \cdots < n_k < N\}.
    \end{align*}
    For two $a,b \in \FIN^{[<\infty]}$, we write $a \leq_\fin b$ iff $\c{a} \subseteq \c{b}$. Then $(\FIN^{[\infty]},\leq,r)$ is a closed triple satisfying \textbf{A1}--\textbf{A4}. 

    $(\FIN^{[\infty]},\leq,r)$ was first introduced by Milliken in \cite{Mil75}. The fact that $(\FIN^{[\infty]},\leq,r)$ satisfies \textbf{A4} follows from Hindman's theorem (\cite[Corollary 3.3]{H74}, or \cite[Theorem 2.41]{T10}).
\end{example}

\begin{example}[Infinite block sequences/Gowers' space $\FIN_k^{[\infty]}$]
\label{ex:Gowers.FIN_k}
    Gowers' space is a generalisation of Hindman's space, introduced in Example \ref{ex:Hindman.FIN}. For each $k \geq 1$, let $\FIN_k$ be the set of all functions $x : \N \to \N$ with finite support (i.e. the set $\supp(x) := \{n \in \N : x(n) > 0\}$ is finite), such that $\ran(x) \subseteq \{0,\dots,k\}$ and $k \in \ran(x)$. For $x,y \in \FIN_k$, we also define the following:
    \begin{itemize}
        \item (Tetris operation) $T(x)$, defined by $T(x)(n) := \max\{x(n) - 1,0\}$. Note that $T(x) \in \FIN_{k-1}$ whenever $k \geq 2$.

        \item $x < y \iff \max(\supp(x)) < \min(\supp(y))$.

        \item If $x < y$, define $x + y \in \FIN_k$ by $(x + y)(n) := \max\{x(n),y(n)\}$. We adopt the convention that the zero function $\mathbf{0}$ (with empty support) is the additive identity, so that $x + \mathbf{0} = \mathbf{0} + x = x$.
    \end{itemize}
    It is easy to see that $<$ is a partial order on $\FIN_k$. We let $\cal{R} = \FIN_k^{[\infty]}$ be the set of all infinite $<$-increasing sequences, and for each $A = (x_n)_{n<\omega} \in \FIN_k^{[\infty]}$, $r_N(A) := (x_n)_{n<N}$. We have that $\cal{AR} = \FIN_k^{[<\infty]}$, the set of all finite $<$-increasing sequences.
    
    Given $a = (x_n)_{n<N} \in \FIN_k^{[<\infty]}$, we let:
    \begin{align*}
        \c{a} := \{T^{\lambda_0}(x_0) &+ \cdots + T^{\lambda_{N-1}}(x_{N-1}) : \\
        &\lambda_0,\dots,\lambda_{N-1} \in \{0,\dots,k\} \text{ and } \lambda_i = 0 \text{ for some } i < N\}.
    \end{align*}
    For two $a,b \in \FIN_k^{[<\infty]}$, we write $a \leq_\fin b$ iff $\c{a} \subseteq \c{b}$. Then $(\FIN_k^{[\infty]},\leq,r)$ is a closed triple satisfying \textbf{A1}--\textbf{A4}. Observe that $(\FIN_1^{[\infty]},\leq,r)$ is precisely Hindman's space $(\FIN^{[\infty]},\leq,r)$, as we can identify each $x \in \FIN_1$ with its support.
    
    $\FIN_k^{[\infty]}$ was first introduced by Gowers in \cite{G92} and \cite{G02} to resolve several long-standing problems in Banach space theory. The current formulation of $\FIN_k^{[\infty]}$ may be found in \cite{T10}. The fact that $(\FIN_k^{[\infty]},\leq,r)$ satisfies \textbf{A4} follows from Gowers' $\FIN_k$ theorem (\cite[Theorem 5]{G92}, or \cite[Theorem 2.22]{T10}).
\end{example}

\begin{example}[Infinite block sequences $\FIN_{\pm k}^{[\infty]}$]
\label{ex:Gowers.FIN_pmk}
    Consider a setup similar to Example \ref{ex:Gowers.FIN_k}. Instead, we let $\FIN_{\pm k}$ be the set of all functions $x : \N \to \Z$ with finite support, such that $\ran(x) \subseteq \{-k,\dots,k\}$, and ($k \in \ran(x)$ or $-k \in \ran(x)$). The tetris operation is now modified to:
    \begin{align*}
        T(x)(n) := 
        \begin{cases}
            x(n) - 1, &\text{if $x(n) > 0$}, \\
            x(n) + 1, &\text{if $x(n) < 0$}, \\
            0, &\text{otherwise}.
        \end{cases}
    \end{align*}
    The rest of the setup is similar. We let $\cal{R} = \FIN_{\pm k}^{[\infty]}$ be the set of all infinite $<$-increasing sequences, so $\cal{AR} = \FIN_{\pm k}^{[<\infty]}$ is the set of all finite $<$-increasing sequences.
    
    Given $a = (x_n)_{n<N} \in \FIN_{\pm k}^{[<\infty]}$, we let:
    \begin{align*}
        \c{a} := \{\varepsilon_0 T^{\lambda_0}(x_0) &+ \cdots + \varepsilon_{N-1} T^{\lambda_{N-1}}(x_{N-1}) : \lambda_0,\dots,\lambda_{N-1} \in \{0,\dots,k\}, \\
        &\varepsilon_0,\dots,\varepsilon_{N-1} \in \{\pm 1\} \text{ and } \lambda_i = 0 \text{ for some } i < N\}.
    \end{align*}
    We also let $\c{A} := \bigcup_{n<\omega} \c{r_n(A)}$. Then $\FIN_{\pm k}^{[\infty]}$ is an \textbf{A2}-space. However, $\FIN_{\pm k}^{[\infty]}$ does not satisfy \textbf{A4} --- for each $x \in \FIN_{\pm k}$, let:
    \begin{align*}
        n_x := \min\{n < \omega : x(n) \in \{\pm k\}\}.
    \end{align*}
    Now let:
    \begin{align*}
        Y := \{x \in \FIN_{\pm k} : x(n_x) = k\}.
    \end{align*}
    Then, for all $A \in \FIN_{\pm k}^{[\infty]}$, $\c{A} \cap Y \neq \emptyset$ and $\c{A} \cap Y^c \neq \emptyset$.

    $\FIN_{\pm k}^{[\infty]}$ was first introduced by Gowers in \cite{G92}, and the current formulation of $\FIN_{\pm k}^{[\infty]}$ may be found in \cite{T10}.
\end{example}

\begin{example}[Hales--Jewett space $W_{Lv}^{[\infty]}$]
\label{ex:Hales-Jewett}
    Let $L = \bigcup_{n<\omega} L_n$ be a countable increasing union of finite alphabets, and let $v \notin L$ be a variable. Let $W_{Lv}$ denote the set of all \emph{variable-words} over $L \cup \{v\}$, i.e. all finite strings of elements of $L \cup \{v\}$ in which $v$ occurs at least once. We also let $W_L$ denote the set of words over $L$ without any occurences of $v$. For each $x \in W_{Lv}$ and $\lambda \in L \cup \{v\}$, we let $x[\lambda] \in W_L \cup W_{Lv}$ denote the word in which all $v$'s occurring in $x$ are replaced by $\lambda$.
    
    Given $x_0,\dots,x_n \in W_{Lv}$, we write $(x_i)_{i<n} < x_n$ iff $\sum_{i<n} |x_i| < |x_n|$. A sequence $(x_n)_{n<N}$ is \emph{rapidly increasing} if $(x_i)_{i<n} < x_n$ for all $n < N$. Let $\cal{R} = W_{Lv}^{[\infty]}$ be the set of all infinite rapidly increasing sequences, and for each $A = (x_n)_{n<\omega} \in \cal{R}$, $r_N(A) := (x_n)_{n<N}$. We have that $\cal{AR} = W_{Lv}^{[<\infty]}$, the set of all finite rapidly increasing sequences.
    
    Given $a = (x_n)_{n<N} \in W_{Lv}^{[<\infty]}$, we let:
    \begin{align*}
        \c{a} := \{x_{n_0}[\lambda_0]^\frown\cdots^\frown x_{n_k}[\lambda_k] : \, &n_0 < \cdots < n_k < N, \, \lambda_i \in L_{n_i} \cup \{v\}, \\
        &\text{ and } \lambda_i = v \text{ for some } i\}.
    \end{align*}
    For two $a,b \in W_{Lv}^{[<\infty]}$, we write $a \leq_\fin b$ iff $\c{a} \subseteq \c{b}$, and $\c{a} \not\subseteq \c{c}$ for all $c \sqsubsetneq b$. Then $(W_{Lv}^{[\infty]},\leq,r)$ is a closed triple satisfying \textbf{A1}--\textbf{A4}.
    
    $W_{Lv}^{[\infty]}$ was based on the work of Hales and Jewett in \cite{HJ63}, and the current formulation of $W_{Lv}^{[\infty]}$ may be found in \cite{T10}. The fact that $(W_{Lv}^{[\infty]},\leq,r)$ satisfies \textbf{A4} follows from the infinite Hales--Jewett theorem (\cite[Theorem 4.21]{T10}).
\end{example}

\begin{example}[Block sequences of located variable-words $\FIN_{Lv}^{[\infty]}$]
\label{ex:block.sequences.located.variable-words}
    Once again, we let $L = \bigcup_{n<\omega} L_n$ be a countable increasing union of finite alphabets, and let $v \notin L$ be a variable. Let $\FIN_{Lv}$ be the set of \emph{located variable-words}, i.e. functions $x : F \to L \cup \{v\}$, where $F \subseteq \N$ is finite, such that $v \in \ran(x)$. We also let $\FIN_L$ be the set of located words $x : F \to L$ with $F \subseteq \N$ finite. For each $x \in \FIN_{Lv}$ with $\dom(x) = F$ and $\lambda \in L \cup \{v\}$, we let $x[\lambda] \in \FIN_L \cup \FIN_{Lv}$ be the function such that for all $n \in F$:
    \begin{align*}
        x[\lambda](n) := 
        \begin{cases}
            \lambda, &\text{if $x(n) = v$}, \\
            x(n), &\text{otherwise}.
        \end{cases}
    \end{align*}
    For $x,y \in \FIN_L \cup \FIN_{Lv}$, we define $x < y$ iff $\max(\dom(x)) < \min(\dom(y))$, and let $\cal{R} := \FIN_{Lv}^{[\infty]}$ be the set of all infinite $<$-increasing sequences of elements of $\FIN_{Lv}$. We have that $\cal{AR} = \FIN_{Lv}^{[<\infty]}$, the set of all finite $<$-increasing sequences. 

    Given $a = (x_n)_{n<N} \in \FIN_{Lv}^{[<\infty]}$, we let:
    \begin{align*}
        \c{a} := \{x_{n_0}[\lambda_0] \cup \cdots \cup x_{n_k}[\lambda_k] : \, &n_0 < \cdots < n_k < N, \, \lambda_i \in L_{n_i} \cup \{v\}, \\
        &\text{ and } \lambda_i = v \text{ for some } i\}.
    \end{align*}
    For two $a,b \in \FIN_{Lv}^{[<\infty]}$, we write $a \leq_\fin b$ iff $\c{a} \subseteq \c{b}$, and $\c{a} \not\subseteq \c{c}$ for all $c \sqsubsetneq b$. Then $(\FIN_{Lv}^{[\infty]},\leq,r)$ is a closed triple satisfying \textbf{A1}--\textbf{A4}.

    The fact that $(\FIN_{Lv}^{[\infty]},\leq,r)$ satisfies \textbf{A4} follows from \cite[Theorem 4.47]{T10}.
\end{example}

\begin{example}[Strong subtrees $\S_\infty(U)$]
\label{ex:strong.subtrees}
    Let $A \subseteq \omega^{<\omega}$ be a tree (i.e. a subset that is closed under initial segments). We introduce some terminology.
    \begin{itemize}
        \item A node $s \in A$ \emph{splits} in $A$ if there exist $m \neq n$ such that $s^\frown m \in A$ and $s^\frown n \in A$. For $n \geq 1$, we let $\func{split}_n(A)$ be the set of all $s \in A$ such that $s$ splits in $A$, and there are exactly $n-1$ values $k < \dom(s)$ such that $s\restrictedto k$ splits in $A$. We then let $\func{split}(A) := \bigcup_{n<\omega}\func{split}_n(A)$.

        \item A node $s \in A$ is \emph{terminal} in $A$ if $s^\frown n \notin A$ for all $n$.

        \item The \emph{height} of a non-empty tree $A$ is one more than the number of distinct levels of $A$ at which some node splits. In other words:
        \begin{align*}
            \height(A) := |\{\dom(s) : s \in \func{split}(A)\}| + 1.
        \end{align*}
        We set $\height(\emptyset) := 0$. 

        \item $A$ is \emph{perfect} if for all $s \in A$, either there exists some $t \sqsupseteq s$ such that $t$ splits in $A$, or there exists some $t \sqsupseteq s$ which is terminal in $A$.
    \end{itemize}
    We say that a perfect subtree $A$ is a \emph{strong subtree} if it satisfies the following conditions:
    \begin{itemize}
        \item For all $s \in A$, there are finitely many successors of $s$. That is, the set:
        \begin{align*}
            \func{succ}_A(s) := \{t \in A : t = s^\frown n \text{ for some } n < \omega\}
        \end{align*}
        is finite for all $s \in A$. 

        \item For all $s,t \in A$ such that $\dom(s) = \dom(t)$, $s$ splits in $A$ iff $t$ splits in $A$.

        \item If $s \in A$ and there exists some $t \in A$ such that $\dom(s) < \dom(t)$, then there exists some $u \in A$ such that $s \sqsubseteq u$ and $\dom(u) = \dom(t)$.
        
        \item If $A$ is of infinite height, then:
        \begin{align*}
            \sup\{|\func{succ}_A(s)| : s \in A\} < \infty.
        \end{align*}
    \end{itemize}
    Given a strong subtree $U$, we say that $A \subseteq U$ is a \emph{strong subtree of $U$} if $A$ is a strong subtree and, for every $s \in \func{split}(A)$ and every successor $s^\frown m$ of $s$ in $U$, there exists some $u \in A$ such that $s^\frown m \sqsubseteq u$ (i.e. every splitting node of $A$ splits fully with respect to $U$). Let $U$ be a strong subtree of infinite height, and let $\cal{R} = \S_\infty(U)$ be the set of strong subtrees of $U$ of infinite height. For each $A \in \S_\infty(U)$ and $n < \omega$, we define $r_0(A) := \emptyset$, and:
    \begin{align*}
        r_{n+1}(A) := \{s \in A : |\{t \in A : t \sqsubsetneq s \wedge t \text{ splits}\}| \leq n\}.
    \end{align*}
    We observe that if $A \in \S_\infty(U)$, then $r_n(A)$ is a strong subtree of height $n$, so $\cal{AR} = \S_{<\infty}(U)$ is the set of all strong subtrees of $U$ of finite height. We also write $\S_n(U) := \cal{AR}_n$, the set of strong subtrees of $U$ of height $n$. For $a,b \in \S_{<\infty}(U)$, we write $a \leq_\fin b$ iff $a \subseteq b$, and for all nodes $s \in a$ which are terminal in $a$, $s$ splits in $b$ or is terminal in $b$. Then $(\S_\infty(U),\subseteq,r)$ is a closed triple satisfying \textbf{A1}--\textbf{A4}.
    
    The space $\S_\infty(U)$ is also known as the Milliken space of strong subtrees (of $U$), and was first introduced by Milliken in \cite{M81} to generalise Silver's partition theorem for infinite trees. The fact that $(\S_\infty(U),\subseteq,r)$ satisfies \textbf{A4} (\cite[Lemma 6.1]{T10}) follows from the strong subtree variant of the Halpern--L\"{a}uchli theorem (\cite[Theorem 2.2]{M81}, or \cite[Theorem 3.2]{T10}). The current formulation may be found in \cite{T10}.
\end{example}

\begin{example}[Sacks forcing/Skewed perfect binary subtrees $\mathbb{S}$]
\label{ex:Sacks.forcing}
    We continue with the terminology introduced in Example \ref{ex:strong.subtrees}. We may define the shortlex well-ordering $\prec$ on $\omega^{<\omega}$, where for all $s,t \in \omega^{<\omega}$:
    \begin{align*}
        s \prec t \iff &\dom(s) < \dom(t) \\
        &\vee (\dom(s) = \dom(t) \wedge \exists i < \dom(s)(s\restrictedto i = t\restrictedto i \wedge s(i) < t(i))).
    \end{align*}
    Given a binary subtree $A \subseteq 2^{<\omega}$ and the $\prec$-increasing enumeration of its splitting nodes $\func{split}(A) = \{s_0,s_1,\dots\}$, we say that $A$ is \emph{skewed} if the following properties hold:
    \begin{enumerate}
        \item If $m < n$, then $\dom(s_m) < \dom(s_n)$. 
        
        \item If there are $N \leq \omega$ splitting nodes in $A$, then the map sending each $s_n$ (for $n < N$) to the $n^\text{th}$ $\prec$-least node of $2^{<\omega}$ is $\sqsubseteq$-preserving in both directions --- that is, $s_m \sqsubseteq s_n$ iff the $m^\text{th}$ $\prec$-least node of $2^{<\omega}$ is an initial segment of the $n^\text{th}$ $\prec$-least node. Furthermore, if $t^\frown i$ is an initial segment of the $n^\text{th}$ $\prec$-least node, where $t$ is the $m^\text{th}$ $\prec$-least node, then $s_m^\frown i \sqsubseteq s_n$.
    \end{enumerate}

    Let $\cal{R} = \mathbb{S}$ be the set of infinite skewed perfect binary subtrees of $2^{<\omega}$. For each $A \in \mathbb{S}$, if $\func{split}(A) = \{s_n : n < \omega\}$ is the $\prec$-increasing enumeration of splitting nodes, then we define $r_0(A) := \emptyset$, and:
    \begin{align*}
        r_{n+1}(A) := \{s \in A : \dom(s) \leq \dom(s_n)\}.
    \end{align*}
    Observe that $a := r_n(A)$ is a (finite) skewed perfect binary subtree whose splitting nodes are $s_0,\dots,s_{n-2}$, and $s_{n-1}$ is the unique terminal node in $a$ such that for any $B \in \mathbb{S}$ where $a \sqsubseteq B$, $s_{n-1}$ splits in $B$. We call $s_{n-1}$ the \emph{splitting terminal node} in $a$. We let $\cal{AR}_n = \mathbb{S}_n$ be the set of all finite skewed perfect binary subtrees of height $n$ whose terminal nodes all have the same length, and let $\cal{AR} = \mathbb{S}_{<\infty} = \bigcup_{n<\omega} \mathbb{S}_n$. For $a,b \in \mathbb{S}_{<\infty}$, we write $a \leq_\fin b$ iff $a \subseteq b$ and, if $a \neq \emptyset$, $a$ and $b$ have the same splitting terminal node. Then $(\mathbb{S},\subseteq,r)$ is a closed triple satisfying \textbf{A1}--\textbf{A4}.

    $(\mathbb{S},\subseteq)$ is a dense subposet of Sacks forcing, hence the given name. The proof that $(\mathbb{S},\subseteq,r)$ satisfies \textbf{A4} may be found in \cite[Polarized Theorem]{B81}. 
\end{example}

\begin{example}[Dual Ramsey theory/Carlson--Simpson Space $\E_\infty$]
\label{ex:Carlson-Simpson}
    Let $\cal{R} = \E_\infty$ denote the Carlson--Simpson space, i.e. the collection of equivalence relations $A$ on $\mathbb{N}$ such that $\mathbb{N}/A$ is infinite. For each $n$, let $A(n)$ be the equivalence class of $A$ with the $n^\text{th}$-least representative. We then let $p_n(A) := \min(A(n))$, and let $p(A) := \{p_n(A) : n < \omega\}$ (note that $p_0(A) = 0$ always). For each $A \in \E_\infty$, we let $r_n(A) := A\restrictedto p_n(A)$, i.e. the restriction of the equivalence relation $A$ to the \emph{domain} $\dom(A\restrictedto p_n(A)) := \{0,1,\dots,p_n(A) - 1\}$. We denote $\E_{<\infty} := \cal{AR}$ and $\E_n := \cal{AR}_n$ (note that this notation is non-standard --- see also the definition of $\E_n$ in the proof of \cite[Corollary 5.72]{T10}). For $a,b \in \E_{<\infty}$, we write $a \leq_\fin b$ iff $\dom(a) = \dom(b)$ and $a$ is coarser than $b$. We remark that for all $a \in \E_{<\infty}$, $\lh(a)$ is the number of equivalence classes in $a$. Then $(\E_\infty,\leq,r)$ is a closed triple satisfying \textbf{A1}--\textbf{A4}. 

    $\E_\infty$ was first introduced by Carlson and Simpson in \cite{CS84}. The current formulation may be found in \cite{T10}. The fact that $(\E_\infty,\leq,r)$ satisfies \textbf{A4} follows from \cite[Lemma 2.4]{CS84} (see also \cite[Lemma 5.69]{T10}).
\end{example}

\begin{example}[Carlson--Simpson space with variables $\E_\infty^L$ (or $W_\infty^L(\omega)$)]
\label{ex:CS-variables}
    We continue with the terminology introduced in Example \ref{ex:Carlson-Simpson}. Let $L$ be a finite non-empty alphabet, and let $\E_\infty^L$ be the set of equivalence relations $A$ on $\N \cup L$ such that $(\N \cup L)/A$ is infinite, and $A\restrictedto L$ is the equality relation. For each $n$, we let $A(n)$ be the equivalence class of $A$, among those disjoint from $L$, with the $n^\text{th}$-least representative. We then let $p_n(A) := \min(A(n))$, and let $p(A) := \{p_n(A) : n < \omega\}$. For each $A \in \E_\infty^L$, we let $r_0(A) := \emptyset$, and $r_{n+1}(A) := A\restrictedto(L \cup p_n(A))$ for $n < \omega$. We then denote $\E_{<\infty}^L := \cal{AR}$ and $\E_n^L := \cal{AR}_n$ (again, this notation is non-standard). The finitisation $\leq_\fin$ is defined similarly --- for $a,b \in \E_{<\infty}^L$, we write $a \leq_\fin b$ iff $\dom(a) = \dom(b)$ and $a$ is coarser than $b$. Note that if $L = \{w\}$, then by relabelling $w$ as $0$ and each $n \in \N$ as $n + 1$, $(\E_\infty^L,\leq,r)$ is equivalent to $(\E_\infty,\leq,r)$.

    The Carlson--Simpson space with variables has two possible reformulations. The first reformulation is similar to the Hales--Jewett space (Example \ref{ex:Hales-Jewett}) --- let $W_\infty^L(\omega)$ denote the set of all infinite variable words (which we call \emph{$\omega$-variable words}) over $L \cup \{v_i : i <\omega\}$ such that every $v_i$ occurs at least once, and if $i < j$, then the first occurrence of $v_i$ is before the first occurrence of $v_j$. Given $A \in W_\infty^L(\omega)$ and $n < \omega$, we let $A(n)$ denote the $n^\text{th}$ letter of the variable-word $A$. We set $r_0(A) := \emptyset$ and, for $n < \omega$:
    \begin{align*}
        r_{n+1}(A) := A\restrictedto m, \text{ where $m$ is the least integer such that $A(m) = v_n$}.
    \end{align*}
    Finally, for $A,B \in W_\infty^L(\omega)$, we write $A \leq B$ iff $A$ is obtained by substituting variables of $B$ with letters or variables.

    The second reformulation is a space of trees (cf. Examples \ref{ex:strong.subtrees} and \ref{ex:Sacks.forcing}). Let $L = \{w_0,\dots,w_{l-1}\}$ be a finite alphabet of size $l \geq 2$. For each $A \in W_\infty^L(\omega)$, recall that $v_k$ first occurs at position $p_k(A)$, and let $\Theta(A) \subseteq L^{<\omega}$ be the unique subtree such that $\emptyset \in \Theta(A)$ and, for each $s \in \Theta(A)$, the following properties hold:
    \begin{itemize}
        \item If $A(\dom(s)) = w_i$, then $s^\frown j \in \Theta(A)$ iff $j = i$.

        \item If $\dom(s) = p_k(A)$ for some $k$, then $s^\frown j \in \Theta(A)$ for all $j < l$.

        \item If $\dom(s) > p_k(A)$ and $A(\dom(s)) = v_k$, then $s^\frown j \in \Theta(A)$ iff $j = s(p_k(A))$.
    \end{itemize}
    Let $\cal{SE}_\infty^L := \{\Theta(A) : A \in W_\infty^L(\omega)\}$. Each $X = \Theta(A) \in \cal{SE}_\infty^L$ is a strong subtree of $L^{<\omega}$ whose splitting levels are exactly $\{p_k(A) : k < \omega\}$. To keep the tree formulation consistent with the equivalence-relation and variable-word formulations, we set $r_0(X) := \emptyset$ and, for $n < \omega$:
    \begin{align*}
        r_{n+1}(X) := \{s \in X : |\{t \in X : t \sqsubsetneq s \wedge t \text{ splits}\}| \leq n\}.
    \end{align*}
    Note that this is precisely the finite approximation used in Example \ref{ex:strong.subtrees}. We let $X \leq Y$ iff $X \subseteq Y$.

    The first reformulation of the Carlson--Simpson space was also introduced by Carlson and Simpson in \cite{CS84}. The current formulation may be found in \cite[\S5.6]{T10}, where $r_{n+1}(A)$ is instead defined to be $A\restrictedto(L \cup p_{n+1}(A))$. The fact that $(\E_\infty^L,\leq,r)$ satisfies \textbf{A4} follows from \cite[Lemma 2.4]{CS84} (see also \cite[Lemma 5.74]{T10}). That the three formulations are equivalent may be formalised by the assertion that there exist bijections $\E_\infty^L \to W_\infty^L(\omega)$ and $\E_\infty^L \to \cal{SE}_\infty^L$, both of which are \emph{strong projections} in either direction --- see Examples \ref{ex:CS-variables.reformulation.variable-words} and \ref{ex:CS-variables.reformulation.trees} respectively.
\end{example}

\begin{example}[Countable vector spaces $E^{[\infty]}$]
\label{ex:vector.spaces}
    Let $\mathbb{F}$ be a countable field, and let $E$ be an $\mathbb{F}$-vector space of dimension $\aleph_0$, with a distinguished Hamel basis $(e_n)_{n<\omega}$. Given $x \in E$, if $x = \sum_{n<\omega} \lambda_ne_n$, we write $\supp(x) := \{n < \omega : \lambda_n \neq 0\}$. We then define a partial order $<$ on $E \setminus \{0\}$ by:
    \begin{align*}
        x < y \iff \max(\supp(x)) < \min(\supp(y)).
    \end{align*}
    We let $\cal{R} := E^{[\infty]}$ denote the set of all infinite $<$-increasing sequences of non-zero vectors in $E$, and for each $A = (x_n)_{n<\omega} \in E^{[\infty]}$, define $r_N(A) := (x_n)_{n<N}$. $\cal{AR} = E^{[<\infty]}$ is the set of finite $<$-increasing sequences of non-zero vectors. Identifying each element of $\cal{AR}_1$ with a non-zero vector of $E$, for $a = (x_n)_{n<N} \in E^{[<\infty]}$ and $A = (x_n)_{n<\omega} \in E^{[\infty]}$ we let:
    \begin{align*}
        \c{a} &:= \func{span}\{x_n : n < N\} \setminus \{0\}, \\
        \c{A} &:= \func{span}\{x_n : n < \omega\} \setminus \{0\}.
    \end{align*}
    For two $a,b \in E^{[<\infty]}$, we write $a \leq_\fin b$ iff $\c{a} \subseteq \c{b}$. Then $(E^{[\infty]},\leq,r)$ is a \textbf{wA2}-space. Furthermore:
    \begin{itemize}
        \item $(E^{[\infty]},\leq,r)$ is an \textbf{A2}-space iff $\mathbb{F}$ is finite. If $\mathbb{F}$ is finite, then for all $a = (x_n)_{n<N},b = (y_m)_{m<M} \in E^{[<\infty]}$, $a \leq_\fin b$ iff each $x_n$ is an $\mathbb{F}$-linear combination of $y_0,\dots,y_{M-1}$, of which there are only finitely many. If $\mathbb{F}$ is infinite, then $(e_0 + \lambda e_1) \leq_\fin (e_0,e_1)$ for all $\lambda \in \mathbb{F}$, so \textbf{A2} fails.
        
        \item $(E^{[\infty]},\leq,r)$ satisfies \textbf{A4} iff $|\mathbb{F}| = 2$. If $|\mathbb{F}| = 2$, then \textbf{A4} follows from Hindman's theorem (\cite[Theorem 2.41]{T10}). If $|\mathbb{F}| > 2$, then define the set:
        \begin{align*}
            Y := \{x \in E : \text{$x = e_n + y$ for some $n$ and $e_n < y$}\}.
        \end{align*}
        It is not difficult to show that $\c{A} \cap Y \neq \emptyset$ and $\c{A} \cap Y^c \neq \emptyset$ for all $A \in E^{[\infty]}$.
    \end{itemize} 
    This Ramsey-theoretic framework of countable vector spaces was first introduced by Rosendal in \cite{R10} and \cite{R09}. He studied the set-theoretic properties of strategically Ramsey sets in this framework, a notion motivated by the Ramsey-theoretic methods employed by Gowers in \cite{G02}. Smythe studied the local Ramsey theory of this framework in \cite{S18}, extending some results by Rosendal to $\H$-strategically Ramsey sets, where $\H$ is a \emph{family} satisfying some combinatorial properties.
\end{example}

\begin{example}[$\cal{R}_\alpha$ for $\alpha < \omega_1$]
\label{ex:Dobrinen.space}
    $\{(\cal{R}_\alpha,\leq_\alpha,r^\alpha) : \alpha < \omega_1\}$ is a hierarchy of topological Ramsey spaces introduced in \cite{DT14} and \cite{DT15}. We give a description of $(\cal{R}_\alpha,\leq_\alpha,r^\alpha)$ for $\alpha < \omega$, and readers may refer to \cite{DT15} for the constructions of $(\cal{R}_\alpha,\leq_\alpha,r^\alpha)$ for $\alpha \geq \omega$. For $\alpha < \omega$ and $n < \omega$, we define a finite tree $\mathbb{A}_\alpha(n)$ by recursion on $\alpha$ as follows:
    \begin{align*}
        \mathbb{A}_0(n) &:= \{\emptyset,(n)\}, \\
        \mathbb{A}_{\alpha+1}(n) &:= \{\emptyset\} \cup \{(n)^\frown t : t \in \mathbb{A}_\alpha(m) \text{ for some } l_n^\alpha \leq m < l_{n+1}^\alpha\},
    \end{align*}
    where $l_0^\alpha := 0$ and $l_{n+1}^\alpha := l_n^\alpha + \max\{1,n+1-\alpha\}$ for all $n$. Let $\mathbb{A}_\alpha := \bigcup_{n<\omega} \mathbb{A}_\alpha(n)$. An element $A \in \cal{R}_\alpha$ is a tree of the form $\bigcup_{n<\omega} A(n)$ such that for some strictly increasing sequence $(k_n)_{n<\omega}$ of natural numbers, $A(n) \subseteq \mathbb{A}_\alpha(k_n)$ is a tree isomorphic to $\mathbb{A}_\alpha(n)$ for each $n$ (where two finite trees $a,b \subseteq \omega^{<\omega}$ are \emph{isomorphic} if there exists a bijection $\gamma : a \to b$ which is both $\sqsubseteq$-preserving and $\prec$-preserving in both directions, and $\prec$ is the shortlex ordering defined in Example \ref{ex:Sacks.forcing}). We call $A(n) = A \cap \mathbb{A}_\alpha(k_n)$ the \emph{$n^\text{th}$ subtree} of $A$, and denote $r_n^\alpha(A) := \bigcup_{i<n} A(i)$ (so $\cal{AR}$ is the set of all finite unions $\bigcup_{i<n} A(i)$ of such trees). Given $a,b \in \cal{AR}$, we define $b \leq_\fin a$ iff $b \subseteq a$ (so for all $A,B \in \cal{R}_\alpha$, $B \leq_\alpha A$ iff $B \subseteq A$). Note that $(\cal{R}_0,\leq_0,r^0)$ is equivalent to the Ellentuck space (Example \ref{ex:Ellentuck.space}), upon identifying each $A = \{\emptyset,(k_n) : n < \omega\} \in \cal{R}_0$ with the set $\{k_n : n < \omega\}$. Note that the tree $\mathbb{A}_\alpha(n)$ is, up to a relabelling of nodes, the structure $\mathbb{S}_\alpha(n)$ in \cite{DT15}, and while the elements of $\cal{R}_\alpha$ are originally defined as subtrees of a tree $\mathbb{T}_\alpha$, it is remarked before \cite[Definition 2.1]{DT15} that $\mathbb{T}_\alpha$ may be taken to be $\mathbb{S}_\alpha$ for $\alpha < \omega$, so no generality is lost.

    It was proven in \cite[Theorem 3.9]{DT14} and \cite[Theorem 3.14]{DT15} that for all $\alpha < \omega_1$, $(\cal{R}_\alpha,\leq_\alpha,r^\alpha)$ is a closed triple satisfying \textbf{A1}--\textbf{A4}, so each of them is a topological Ramsey space. The construction of these spaces was motivated by works in \cite{L89} and \cite{RT12}. For instance, the space $(\cal{R}_1,\leq_1,r^1)$ induces a natural variant of a Ramsey ultrafilter $\U_1$, which is isomorphic to the ultrafilter forced by Laflamme in \cite{L89}. By proving a canonisation theorem, Dobrinen and Todor\v{c}evi\'{c} classified all (Rudin--Keisler) isomorphism classes of ultrafilters Tukey reducible to $\U_1$.
\end{example}

\begin{example}[Singleton space]
\label{ex:singleton.space}
    Let $\cal{R} = \{\vec{0}\} := \{(0,0,\dots)\}$, the singleton containing the zero sequence. We define $r_n(A)$ to be the zero sequence $(0,\dots,0)$ of length $n$, and let $\leq_\fin$ be the equality relation. Then $(\{\vec{0}\},=,r)$ is a closed triple satisfying \textbf{A1}--\textbf{A4}. First introduced by the second author in \cite[Example 2.9]{Y25}, the singleton space serves as a pathological example of a topological Ramsey space. 
\end{example}

One may observe that the following examples are similar, in that they consist of all $<$-increasing sequences of elements of the countable set $\cal{AR}_1$, for some relation $<$ between $\cal{AR}$ and $\cal{AR}_1$:
\begin{itemize}
    \item Ellentuck space $[\N]^\infty$ (a subset $A \subseteq \N$ may be viewed as a strictly increasing sequence of integers).
    
    \item Hindman's space $\FIN^{[\infty]}$ and the spaces of infinite block sequences $\FIN_k^{[\infty]}$, $\FIN_{\pm k}^{[\infty]}$.

    \item Hales--Jewett space $W_{Lv}^{[\infty]}$ (where $L$ is finite).
    
    \item Infinite block sequences of located variable-words $\FIN_{Lv}^{[\infty]}$ (where $L$ is finite).
    
    \item Countable vector spaces $E^{[\infty]}$.
\end{itemize}

\begin{notation}
\label{not:first.approximation}
    Let $(\cal{R},\leq,r)$ be a \textbf{wA2}-space. For $A \in \cal{R}$ and $a \in \cal{AR}$, we write:
    \begin{align*}
        \c{a} &:= \{b \in \cal{AR}_1 : b \leq_\fin c \text{ for some } c \sqsubseteq a\}, \\
        \c{A} &:= \bigcup_{n<\omega} \c{r_n(A)}.
    \end{align*}
\end{notation}

These spaces satisfy the common property that for all $a,b \in \cal{AR}$, $a \leq_\fin b$ iff $\c{a} \subseteq \c{b}$ (and possibly with some additional depth-related restrictions). When we say that $(\cal{R},\leq,r)$ is \emph{unconditional}, we mean that it is one of the examples listed above. In this case, if $A = (x_n)_{n<\omega} \in \cal{R}$ is a sequence, then we write:
\begin{align*}
    A/n := (x_k)_{k \geq n}.
\end{align*}
In all of the above examples, $A/n$ is a well-defined element of $\cal{R}$. 

\subsection{Ramsey quasi-orders}
\label{sec:Ramsey.quasi-orders}
The motivation for the introduction of Ramsey quasi-orders is twofold. 
\begin{enumerate}
    \item While the axioms of topological Ramsey spaces serve as a standard abstract framework for Ramsey theory, we do not need the full power of these axioms for most of the results in this paper --- for instance, the finite approximations of elements of $\cal{R}$ are often not relevant. The notion of a Ramsey quasi-order serves as a much simplified version of the axioms \textbf{A1}--\textbf{A4}, whose properties are enough for most of the results in our study.
    
    \item In \cite{MS17}, Ma{\v{s}}ulovi{\'c} and Scow introduced the Ramsey property of a category, and showed that many known results in finite Ramsey theory can be realised as Ramsey properties of suitable categories. While the immediate generalisation to results in infinite-dimensional Ramsey theory fails, we show that Ramsey quasi-order is a natural attempt to generalise the categorical approach to infinite-dimensional Ramsey theory.
\end{enumerate}

\begin{definition}
\label{def:analytic.quasi-order}
    An \emph{analytic} (\emph{Borel}, $G_\delta$ etc.) \emph{quasi-order} is a quasi-order $(\cal{R},\leq)$ such that $\cal{R}$ is a Polish space, and $\leq$ is an analytic (Borel, $G_\delta$ etc.) subset of $\cal{R}^2$.
\end{definition}

We remark that Borel quasi-orders were also extensively studied in other contexts, particularly the study of countable Borel equivalence relations and Borel reducibility --- see, for instance,  \cite{F79, HS82, S84, HMS88, Lou89, RX24}.

Similar to the notation introduced in \S\ref{sec:TRS.axioms}, if $(\cal{R},\leq)$ is an analytic (Borel, $G_\delta$ etc.) quasi-order and $A \in \cal{R}$, we write:
\begin{align*}
    [A] := \{B \in \cal{R} : B \leq A\}.
\end{align*}
Note that $[A]$ is analytic (Borel, $G_\delta$ etc.), as it is the preimage of $\leq$ under the continuous map $B \mapsto (B,A)$.

\begin{definition}
\label{def:Ramsey.quasi-order}
    Let $(\cal{R},\leq)$ be an analytic quasi-order.
    \begin{enumerate}
        \item $(\cal{R},\leq)$ has the \emph{weak} (\emph{analytic}, \emph{Borel} etc.) {Ramsey property} if for every (analytic, Borel etc.) $\X \subseteq \cal{R}$, there exists some $A \in \cal{R}$ such that $[A] \subseteq \X$ or $[A] \subseteq \X^c$.

        \item $(\cal{R},\leq)$ has the (\emph{analytic}, \emph{Borel} etc.) {Ramsey property} if for every (analytic, Borel etc.) $\X \subseteq \cal{R}$ and $A \in \cal{R}$, there exists some $B \leq A$ such that $[B] \subseteq \X$ or $[B] \subseteq \X^c$.
        
        \item $(\cal{R},\leq)$ is a \emph{Ramsey quasi-order} if it has the analytic Ramsey property.
    \end{enumerate}
\end{definition}

By Lemma \ref{lem:reduction.is.G_delta} and the abstract Ellentuck theorem (Theorem \ref{thm:abstract.Ellentuck.thm}), we have that if $(\cal{R},\leq,r)$ is a closed triple satisfying \textbf{A1}--\textbf{A4} and $\cal{AR}$ is countable, then $(\cal{R},\leq)$ is a $G_\delta$ Ramsey quasi-order. More precisely, we have the following:

\begin{lemma}
\label{lem:A4.iff.local.Ramsey.properties}
    Let $(\cal{R},\leq,r)$ be an \textbf{A2}-space. The following are equivalent:
    \begin{enumerate}
        \item $(\cal{R},\leq,r)$ satisfies \textbf{A4}.
        
        \item For all $A \in \cal{R}$ and $a \in \cal{AR}\restrictedto A$, $([a,A],\leq)$ is a Ramsey quasi-order.

        \item For all $A \in \cal{R}$ and $a \in \cal{AR}\restrictedto A$, there exists some $B \in [a,A]$ such that $([a,B],\leq)$ is an analytic quasi-order that has the weak clopen Ramsey property.
    \end{enumerate}
\end{lemma}

\begin{proof}
    Note that $[a,A]$ is a closed subset of $\cal{R}$, $\leq$ is $G_\delta$ by Lemma \ref{lem:reduction.is.G_delta} and $\cal{AR}\restrictedto[a,A]$ is countable by \textbf{A2}, so $([a,A],\leq)$ is an analytic quasi-order.
    
    \underline{(1)$\implies$(2):} Fix some $A \in \cal{R}$ and $a \in \cal{AR}\restrictedto A$. By the abstract Ellentuck theorem, for all analytic $\X \subseteq [a,A]$ and $B \in [a,A]$, there exists some $C \in [a,B]$, such that $[a,C] \subseteq \X$ or $[a,C] \subseteq \X^c$, as desired. 

    \underline{(2)$\implies$(3):} Immediate.

    \underline{(3)$\implies$(1):} Fix some $A \in \cal{R}$, $a \in \cal{AR}\restrictedto A$ and $\O \subseteq \cal{AR}_{\lh(a)+1}$. Let $B \in [a,A]$ such that $([a,B],\leq)$ has the weak clopen Ramsey property. Define:
    \begin{align*}
        \X := \{C \in [a,B] : r_{\lh(a)+1}(C) \in \O\}.
    \end{align*}
    Observe that $\X$ is a clopen subset of $[a,B]$. Since $([a,B],\leq)$ has the weak clopen Ramsey property, there exists some $C \in [a,B]$ such that $[a,C] \subseteq \X$ or $[a,C] \subseteq \X^c$. Then $C \in [a,A]$ and $r_{\lh(a)+1}[a,C] \subseteq \O$ or $r_{\lh(a)+1}[a,C] \subseteq \O^c$ respectively.
\end{proof}

At the end of the subsection, we shall present Example \ref{ex:Rado.graph.diary}, a Ramsey quasi-order that is not a topological Ramsey space.

We now address the relationship between Ramsey quasi-orders and the Ramsey property of categories. Recall that a \emph{category} $\mathbf{C}$ consists of a class of \emph{objects} $\Ob(\mathbf{C})$, a set of \emph{morphisms} $\hom_{\mathbf{C}}(\A,\B)$ for each $\A \in \Ob(\mathbf{C})$, identity morphisms $\func{id}_\A \in \hom_{\mathbf{C}}(\A,\A)$ and a binary operation $\circ$ satisfying the following:
\begin{itemize}
    \item (Associativity) If $f \in \hom_{\mathbf{C}}(\A,\B)$, $g \in \hom_{\mathbf{C}}(\B,\cal{C})$ and $h \in \hom_{\mathbf{C}}(\cal{C},\D)$, then:
    \begin{align*}
        h \circ (g \circ f) = (h \circ g) \circ f.
    \end{align*}

    \item (Identity) For all $f \in \hom_{\mathbf{C}}(\A,\B)$:
    \begin{align*}
        \func{id}_\B \circ f = f = f \circ \func{id}_\A.
    \end{align*}
\end{itemize}
Recall also that given a set $X$ and positive integer $l$, an \emph{$l$-colouring} of $X$ is a function $\chi : X \to l$. Given $Y \subseteq X$, we say that $\chi$ is \emph{monochromatic} on $Y$ if $|\chi``[Y]| = 1$.

\begin{definition}[\cite{MS17}]
\label{def:Ramsey.category}
    Let $\mathbf{C}$ be a category. 
    \begin{enumerate}
        \item Given $\A,\B,\cal{C} \in \Ob(\mathbf{C})$ and a positive integer $k$, we write $\cal{C} \to (\B)_k^\A$ to denote that for all $k$-colouring $\chi$ of $\hom_{\mathbf{C}}(\A,\cal{C})$, there exists some $w \in \hom_{\mathbf{C}}(\B,\cal{C})$ such that $\chi$ is monochromatic on $w \circ \hom_{\mathbf{C}}(\A,\B) = \{w \circ v : v \in \hom_{\mathbf{C}}(\A,\B)\}$. 
        
        \item We say that $\mathbf{C}$ has the \emph{weak Ramsey property} if for all positive integer $k$ and $\A,\B \in \Ob(\mathbf{C})$ such that $\hom_{\mathbf{C}}(\A,\B) \neq \emptyset$, there exists some $\cal{C} \in \Ob(\mathbf{C})$ such that $\cal{C} \to (\B)_k^\A$ holds. We also say that $\mathbf{C}$ is a \emph{Ramsey category}.
    \end{enumerate}
\end{definition}

Note that in \cite{MS17} and \cite{M18}, Definition \ref{def:Ramsey.category} is simply called the \emph{Ramsey property}, but we shall call it the weak Ramsey property to be consistent with Definition \ref{def:Ramsey.quasi-order} (see also Lemma \ref{lem:Ramsey.monoid.equivalence}). Examples of Ramsey categories can be found in \cite{MS17} and \cite{M18}, but more can be derived from those in \S\ref{sec:examples}. 

\begin{example}[Finite subsets of $\N$, cf. Example \ref{ex:Ellentuck.space}]
\label{ex:Ellentuck.space.category}
    We define the category $\mathbf{E}$ by stipulating that $\Ob(\mathbf{E}) = \N$, and for $m,n \in \N$ such that $m \leq n$, we let $\hom_{\mathbf{E}}(m,n)$ be the set of strictly increasing functions $f : m \to n$. We let $\circ$ be the composition function, which is clearly associative. The identity morphism of $\hom_{\mathbf{E}}(n,n)$ is the identity function. We leave the readers to verify that the weak Ramsey property of $\mathbf{E}$ is precisely the finite Ramsey theorem (Corollary \ref{cor:finite.Ramsey.theorem}).
\end{example}

\begin{example}[Finite block sequences of $\FIN$, cf. Example \ref{ex:Hindman.FIN}]
\label{ex:Hindman.FIN.category}
    We define the category $\mathbf{BS}$ by stipulating that $\Ob(\mathbf{BS}) = \N$, and for $m,n \in \N$ such that $m \leq n$, we let $\hom_{\mathbf{BS}}(m,n)$ be the set of all block sequences $(x_i)_{i<m} \in \FIN^{[m]}$ such that $x_i \subseteq n$ for all $i$. If $l \leq m \leq n$, $b := (y_j)_{j<l} \in \hom_{\mathbf{BS}}(l,m)$ and $a := (x_i)_{i<m} \in \hom_{\mathbf{BS}}(m,n)$, then we define:
    \begin{align*}
        z_j := \bigcup_{i \in y_j} x_i,
    \end{align*}
    and let $a \circ b := (z_j)_{j<l} \in \hom_{\mathbf{BS}}(l,n)$. We leave readers to verify that $\circ$ is associative, and for each $n$, $(\{i\})_{i<n} \in \hom_{\mathbf{BS}}(n,n)$ is the identity morphism. 
    
    The weak Ramsey property of $\mathbf{BS}$ is precisely the following statement: For all positive integers $k,l,m$, there is a positive integer $n$ such that for every $a \in \FIN^{[n]}$ and every $k$-colouring $\chi$ of $\FIN^{[l]}\restrictedto a$, there exists some $b \in \FIN^{[m]}\restrictedto a$ such that $\chi$ is monochromatic on $\FIN^{[l]}\restrictedto b$. This is precisely a finite version of Hindman's theorem --- see, for instance, \cite[Theorem 2.6]{M75}.
\end{example}

We shall provide a more detailed account of constructing Ramsey categories from topological Ramsey spaces in \S\ref{sec:Ramsey.categories}.

Observe that in both examples above, we have that $\Ob(\mathbf{C}) = \N$, and each $n \in \N$ represents a sequence of length $n$. Thus, if we want to consider only infinite sequences, then $\Ob(\mathbf{C})= \{\infty\}$ is a singleton, and $\hom_{\mathbf{C}}(\infty,\infty)$ represents the set of infinite sequences of interest. Equivalently, we consider a \emph{monoid} $(M,\circ)$, i.e. a set $M$ and an associative binary operation $\circ$ with an identity element in $M$.

\begin{example}[Infinite subsets of $\N$, cf. Example \ref{ex:Ellentuck.space.category}]
\label{ex:Ellentuck.space.monoid}
    Let $M_{\mathbf{E}}$ be the set of all strictly increasing functions $f : \N \to \N$. We define the binary operation $\circ$ on $M_{\mathbf{E}}$ by function composition, which is clearly associative. The identity element is the identity function on $\N$. 

    Each $f_A \in M_{\mathbf{E}}$ may be identified with its range $A := \ran(f_A)$. Observe that if $f_B = f_A \circ f_C$, then $B = \ran(f_B) \subseteq \ran(f_A) = A$. Conversely, if $B \subseteq A$ and $A = \{m_n : n < \omega\}$ is its increasing enumeration, then $B = \{m_{n_k} : k < \omega\}$ for some increasing sequence $(n_k)_{k<\omega}$. Letting $f_C \in M_{\mathbf{E}}$ be the function for which $f_C(k) := n_k$ for all $k$, we have that $f_B = f_A \circ f_C$. Therefore, for all infinite subsets $A,B \subseteq \N$:
    \begin{align*}
        B \subseteq A \iff \exists C \in [\N]^\infty \, f_B = f_A \circ f_C.
    \end{align*}
    We also define $A \circ C := \ran(f_A \circ f_C)$, so $B = A \circ C$ iff $f_B = f_A \circ f_C$.
\end{example}

\begin{example}[Infinite block sequences of $\FIN$, cf. Example \ref{ex:Hindman.FIN.category}]
\label{ex:Hindman.FIN.monoid}
    Let $M_{\mathbf{BS}} = \FIN^{[\infty]}$ be the set of all infinite block sequences of elements of $\FIN$. The binary operation $\circ$ can be defined similarly to Example \ref{ex:Hindman.FIN.category} --- if $B := (y_j)_{j<\omega} \in \FIN^{[\infty]}$ and $A := (x_i)_{i<\omega} \in \FIN^{[\infty]}$, we let $z_j := \bigcup_{i \in y_j} x_i$, and we define $A \circ B := (z_j)_{j<\omega} \in \FIN^{[\infty]}$. 
    
    An observation on the relationship between the reduction relation $\leq$ and the operation $\circ$, similar to that of Example \ref{ex:Ellentuck.space.monoid}, can be made as follows: Let $A,B \in \FIN^{[\infty]}$. If $B = A \circ C$ for some $C \in \FIN^{[\infty]}$, then every element in the sequence $B$ is a union of elements of the sequence $A$, so $B \leq A$. Conversely, if $B = (y_j)_{j<\omega}$, $(x_i)_{i<\omega}$ and $B \leq A$, then for every $j$ there exists a finite set $w_j \subseteq \N$ such that $y_j = \bigcup_{i \in w_j} x_i$. Since $B$ is an infinite block sequence, $\max(w_j) < \min(w_{j+1})$ for all $j$. Therefore, we may let $C := (w_k)_{k<\omega}$, and we have that $B = A \circ C$. This gives the equivalence:
    \begin{align*}
        B \leq A \iff \exists C \in \FIN^{[\infty]} \, B = A \circ C.
    \end{align*}
\end{example}

In general, given an arbitrary monoid $(M,\circ)$, we may define a relation $\leq$ on $M$ where for all $A,B \in M$:
\begin{align*}
    B \leq A \iff \exists C \in M \, B = A \circ C.
\end{align*}
Observe that $\leq$ is a quasi-order on $M$:
\begin{itemize}
    \item (Transitivity) Suppose that $A_2 \leq A_1$ and $A_1 \leq A_0$. Let $B_0,B_1 \in M$ be such that $A_2 = A_1 \circ B_1$ and $A_1 = A_0 \circ B_0$. By the associativity of $\circ$:
    \begin{align*}
        A_2 = A_1 \circ B_1 = (A_0 \circ B_0) \circ B_1 = A_0 \circ (B_0 \circ B_1),
    \end{align*}
    so $B_0 \circ B_1$ witnesses that $A_2 \leq A_0$.

    \item (Reflexivity) For any $A \in M$, since $A = A \circ 1$ where $1$ is the identity element in $M$, $1$ witnesses that $A \leq A$.
\end{itemize}
Furthermore, if $(M,\circ)$ is a Polish monoid, i.e. a topological monoid such that the underlying topology is Polish, then $\leq$ is an analytic subset of $M^2$. Indeed, letting $\pi : M^3 \to M^2$ be the projection $\pi(B,A,C) := (B,A)$, we have:
\begin{align*}
    {\leq} &= \{(B,A) : \exists C \in M \, B = A \circ C\} \\
    &= \pi``\{(B,A,C) \in M^3 : B = A \circ C\}.
\end{align*}
The set $\{(B,A,C) \in M^3 : B = A \circ C\}$ is the preimage of the closed set $\{(x,y) \in M^2 : x = y\}$ under the continuous map $(B,A,C) \mapsto (B,A \circ C)$, so it is a closed subset of the Polish space $M^3$. Therefore, $\leq$ is a continuous image of a closed subset of the Polish space $M^3$, so it is analytic. Consequently, every Polish monoid $(M,\circ)$ gives rise to an analytic quasi-order $(M,\leq)$.

\begin{lemma}
\label{lem:Ramsey.monoid.equivalence}
    Let $(M,\circ)$ be a monoid, realised as a category $\mathbf{C}$ with $\Ob(\mathbf{C}) = \{\infty\}$ and $M = \hom_{\mathbf{C}}(\infty,\infty)$. The following are equivalent:
    \begin{enumerate}
        \item $\mathbf{C}$ has the weak Ramsey property.
        
        \item $\infty \to (\infty)_k^\infty$ for all $k \geq 2$.
        
        \item $\infty \to (\infty)_2^\infty$.
        
        \item $(M,\circ)$ has the weak Ramsey property.
    \end{enumerate}
\end{lemma}

\begin{proof}
    (1)$\iff$(2) follows from the definition of the weak Ramsey property. (2)$\implies$(3) is immediate, while (3)$\implies$(2) follows from applying $\infty \to (\infty)_2^\infty$ $k - 1$ times to a $k$-colouring. (4) is equivalent to (3) by identifying each $2$-colouring with a $2$-partition, $\X$ and $\X^c$.
\end{proof}

Assuming the Axiom of Choice, $(M_{\mathbf{E}},\circ)$ does not satisfy the weak Ramsey property, as there is a non-Ramsey subset of $[\N]^\infty$. A similar construction shows that $(M_{\mathbf{BS}},\circ)$ also does not satisfy the weak Ramsey property. However, both monoids satisfy the weak Ramsey property \emph{for analytic sets}, i.e. every analytic $\X \subseteq M$ satisfies the statement in Lemma \ref{lem:Ramsey.monoid.equivalence}(4) --- since the quasi-orders $(M_{\mathbf{E}},\leq)$ and $(M_{\mathbf{BS}},\leq)$ coincide with $([\N]^\infty,\subseteq)$ and $(\FIN^{[\infty]},\leq)$ respectively, this follows from the abstract Ellentuck theorem (Theorem \ref{thm:abstract.Ellentuck.thm}). Therefore, $(M_{\mathbf{E}},\leq)$ and $(M_{\mathbf{BS}},\leq)$ are analytic quasi-orders that satisfy the weak Ramsey property (of course, the abstract Ellentuck theorem implies that they are Ramsey quasi-orders). 

We conclude this subsection with an example of a Polish monoid whose induced analytic quasi-order is a Ramsey quasi-order that is not a topological Ramsey space.

\begin{example}[A diary that codes the Rado graph]
\label{ex:Rado.graph.diary}
    In \cite[\S8.3]{DZ23}, Dobrinen and Zucker constructed a diary $T$ that codes the Rado graph. Let $\func{DEmb}(T,T)$ be the set of diary embeddings of $T$ into itself, equipped with the composition $\circ$ and the topology of pointwise convergence. Then $(\func{DEmb}(T,T),\circ)$ is a Polish monoid, and the induced analytic quasi-order $(\func{DEmb}(T,T),\leq)$ satisfies the weak analytic Ramsey property by \cite[Theorem 7.4]{DZ23}. Since $(\func{DEmb}(T,T),\leq)$ satisfies property (H) by Proposition \ref{prop:Polish.monoid.property.H}, it is a Ramsey quasi-order by Lemma \ref{lem:cts.qo.family.transfers.analytic.Ramsey.property}. However, it is not a topological Ramsey space, as it was noted in \cite[\S8.3]{DZ23} that the Ellentuck topology of $\func{DEmb}(T,T)$ is discrete.
\end{example}

\subsection{Ramsey ultrafilters}
\label{sec:Ramsey.ultrafilters}
Notions of Ramsey ultrafilters for topological Ramsey spaces have been proposed in various works --- see the examples following Definition \ref{def:Ramsey.filter}. The relationships between Ramsey ultrafilters, both within a single topological Ramsey space and across different ones, have been extensively studied, and will prove very useful in understanding when projections between topological Ramsey spaces exist. We propose a notion of Ramsey ultrafilters for Ramsey quasi-orders and prove some of their basic properties. In this section, we fix a Ramsey quasi-order $(\cal{R},\leq)$.

\begin{definition}
\label{def:filter}
    A \emph{filter} $\F \subseteq \cal{R}$ is a subset that is a filter in the quasi-order $(\cal{R},\leq)$. In other words, it satisfies the following two properties:
    \begin{enumerate}
        \item If $A \in \F$ and $A \leq B$, then $B \in \F$.
        
        \item If $A,B \in \F$, then there exists some $C \in \F$ such that $C \leq A$ and $C \leq B$.
    \end{enumerate}
    If $\F$ is maximal (i.e. there is no filter $\G \subseteq \cal{R}$ such that $\F \subsetneq \G$), then $\F$ is an \emph{ultrafilter}.
\end{definition}

We remark that this notion of ultrafilter differs from several other existing notions of ultrafilters of topological Ramsey spaces.
\begin{itemize}
    \item If $\U \subseteq \cal{R}$ is an ultrafilter (as in Definition \ref{def:filter}), then $\{\c{A} : A \in \U\}$ (see Notation \ref{not:first.approximation}) does not necessarily generate a first-approximation ultrafilter on $\cal{AR}_1$ (as introduced in \cite{M07}). For instance, by Zorn's lemma, there is an ultrafilter in $(\FIN^{[\infty]},\leq)$ (as in Definition \ref{def:filter}). However, a first-approximation ultrafilter on $\FIN$ is an ordered-union ultrafilter, whose existence is not provable in $\ZFC$.
    
    \item On the other hand, it is unclear whether a first-approximation ultrafilter on $\cal{AR}_1$ generates an ultrafilter in $\cal{R}$. If $U$ is an ultrafilter on $\cal{AR}_1$ generated by sets of the form $\c{A}$ for $A \in \cal{R}$, we may consider the set:
    \begin{align*}
        \U := \{A \in \cal{R} : \c{A} \in U\}.
    \end{align*}
    However, it is not true in general that for any $A,B \in \U$, if $\c{A} \subseteq \c{B}$ then $A \leq B$. The Carlson--Simpson space $\E_\infty$ is an example. In fact, there exist two incompatible $A,B \in \E_\infty$ such that $\c{A} = \c{B}$ --- let $2\N + 1$ be the set of all odd positive integers, and let $\Gamma : \N^2 \to 2\N + 1$ be a bijection such that $\Gamma(i,j) \geq 2\max\{i,j\}$ for all $i,j \in \N$. We then let:
    \begin{align*}
        (2m,n) \in A &\iff n = \Gamma(m,j) \text{ for some $j \in \N$}, \\
        (2m,n) \in B &\iff n = \Gamma(i,m) \text{ for some $i \in \N$}.
    \end{align*}
    Then the equivalence relation on $\N$ generated by $A \cup B$ is the trivial relation on $\N$ (i.e. $m \sim n$ for all $m,n \in \N$). On the other hand, for each $m \in \N$, let $e_m$ be the single-class equivalence relation on $m = \{0,\dots,m-1\}$ (i.e. $(i,j) \in e_m$ for all $i,j < m$). Then:
    \begin{align*}
        \c{A} = \c{B} = \{e_{2m} : m \in \N\}.
    \end{align*}

    \item It is clear that if $\U \subseteq \cal{R}$ is an ultrafilter as in \cite[Definition 6.1 (and Definition 3.1)]{DMN17}, then it is an ultrafilter as in Definition \ref{def:filter}. It seems unlikely that the two notions are equivalent. An abstract version of Ramsey filters was also considered in \cite{DMT17}. 
\end{itemize}
From here onward, all uses of the term ``ultrafilter'' refer to Definition \ref{def:filter}. 

\begin{definition}
\label{def:F-Ramsey}
    Let $\F \subseteq \cal{R}$ be a filter. A subset $\X \subseteq \cal{R}$ is \emph{$\F$-Ramsey} if there exists some $A \in \F$ such that $[A] \subseteq \X$ or $[A] \subseteq \X^c$.
\end{definition}

We remark that this definition does \emph{not} coincide with the definition of $\H$-Ramsey for topological Ramsey spaces --- see, for instance, \cite[Definition 7.2]{T10} and \cite{DMN17}.

\begin{definition}
\label{def:Ramsey.filter}
    A filter $\F \subseteq \cal{R}$ is \emph{Ramsey} if every analytic subset of $\cal{R}$ is $\F$-Ramsey.
\end{definition}

\begin{lemma}
\label{lem:Ramsey.filter.is.ultrafilter}
    Every Ramsey filter is an ultrafilter.
\end{lemma}

\begin{proof}
    Let $\F \subseteq \cal{R}$ be a Ramsey filter, and let $\G \subseteq \cal{R}$ be a filter such that $\F \subseteq \G$. Let $A \in \G$. Then $[A]$ is analytic, so there exists some $B \in \F$ such that $[B] \subseteq [A]$ or $[B] \cap [A] = \emptyset$. Since $\F \subseteq \G$, we have that $B \in \G$, so $[B] \cap [A] = \emptyset$ is impossible as $A \in \G$ and $\G$ is a filter. This implies that $B \leq A$, so $A \in \F$ as $\F$ is upward-closed. Therefore $\G \subseteq \F$, which gives $\F = \G$, i.e. $\F$ is maximal.
\end{proof}

Filters for which every clopen set is $\F$-Ramsey satisfy a finite-dimensional version of the pigeonhole principle:

\begin{lemma}
\label{lem:Ramsey.ultrafilter.finite.pigeonhole}
    Let $(\cal{R},\leq,r)$ be a closed triple satisfying \textbf{A1}--\textbf{A4}, and let $\F \subseteq \cal{R}$ be a filter such that every clopen subset of $\cal{R}$ is $\F$-Ramsey. For any $n$, if $\O \subseteq \cal{AR}_n$, then there exists some $A \in \F$ such that $r_n[A] \subseteq \O$ or $r_n[A] \subseteq \O^c$.
\end{lemma}

\begin{proof}
    Let $n < \omega$, and let $\O \subseteq \cal{AR}_n$. Let:
    \begin{align*}
        \X := \{A \in \cal{R} : r_n(A) \in \O\}.
    \end{align*}
    Observe that $\X = \bigcup_{a \in \O}[a]$ and $\X^c = \bigcup_{a \in \cal{AR}_n \setminus \O}[a]$, so $\X$ is clopen. Therefore, there exists some $A \in \F$ such that $[A] \subseteq \X$ or $[A] \subseteq \X^c$. By the definition of $\X$, we have that $r_n[A] \subseteq \O$ or $r_n[A] \subseteq \O^c$, as desired.
\end{proof}

Ramsey ultrafilters of various topological Ramsey spaces have been extensively studied.

\begin{example}[Ellentuck space $\lbrack\N\rbrack^\infty$]
\label{ex:selective.ultrafilter}
    A Ramsey ultrafilter $\U \subseteq [\N]^\infty$ is an ultrafilter on $\N$ which is Ramsey (or selective) in the usual sense. If $\U$ is a selective ultrafilter on $\N$ (see, e.g., \cite[Definition 7.3]{T10}), then by the local Ellentuck theorem (\cite[Theorem 7.19]{T10}), $\U$ is Ramsey as in Definition \ref{def:Ramsey.filter}. Conversely, every Ramsey ultrafilter is selective --- see, for instance, \cite[Theorem 0.10]{M77}. 
\end{example}

\begin{example}[Hindman space $\FIN^{[\infty]}$]
\label{ex:stable.ordered-union.ultrafilter}
    A Ramsey ultrafilter $\U \subseteq \FIN^{[\infty]}$ is precisely a stable ordered-union ultrafilter. This was proved by Blass in \cite[Theorem 4.2]{B87}, where Definition \ref{def:Ramsey.filter} coincides with the \emph{infinitary partition property}.
\end{example}

\begin{example}[Carlson--Simpson space $\E_\infty$]
\label{ex:Ramseyan.ultrafilter}
    Filters of $\E_\infty$ were first studied by Matet in \cite{M86}. He showed that a Ramsey ultrafilter $\U \subseteq \E_\infty$ exists under $\CH$ (see \cite[Theorem 8.1]{M86}). A similar version of ultrafilters of $\E_\infty$ was considered by Halbeisen in \cite{H01,H17}, which he called \emph{Ramseyan ultrafilters}. A Ramseyan ultrafilter, as defined in \cite{H17} is an ultrafilter $\U \subseteq \E_\infty$ (as in Definition \ref{def:filter}) such that for all $n$ and $\O \subseteq \E_n$, there exists some $A \in \U$ such that $r_n[A] \subseteq \O$ or $r_n[A] \subseteq \O^c$ (see also the definition provided in \cite{H01}, which is formally stronger). By Lemma \ref{lem:Ramsey.ultrafilter.finite.pigeonhole}, every Ramsey ultrafilter is Ramseyan.
\end{example}

\begin{example}[$\cal{R}_\alpha$]
\label{ex:R_alpha.Ramsey.ultrafilter}
    Filters of $\cal{R}_\alpha$ were first introduced by Dobrinen and Todor\v{c}evi\'{c} in \cite{DT14} and \cite{DT15} to study the Tukey theory of ultrafilters. They considered a variant of Ramsey ultrafilter, where an ultrafilter $\U \subseteq \cal{R}_\alpha$ is \emph{Ramsey for $\cal{R}_\alpha$} if \emph{it satisfies the abstract Nash-Williams Theorem}, i.e. for any front $\O \subseteq \cal{AR}$ and partition $\O = \O_0 \sqcup \O_1$, there is some $A \in \U$ such that $\cal{AR}\restrictedto A \cap \O_i = \emptyset$ for some $i < 2$.

    One can show that if $\U \subseteq \cal{R}_\alpha$ is a filter such that every clopen subset of $\cal{R}_\alpha$ is $\U$-Ramsey, then it satisfies the abstract Nash-Williams Theorem. In fact, this holds for any closed triple satisfying \textbf{A1}--\textbf{A4}.
\end{example}

\begin{proof}
    Let $\U \subseteq \cal{R}_\alpha$ be a filter such that every clopen subset of $\cal{R}_\alpha$ is $\U$-Ramsey. Let $\O \subseteq \cal{AR}$ be a front, and let $\O_0 \sqcup \O_1 = \O$ be a partition of this front. Since $\O$ is a front, for all $A \in \cal{R}_\alpha$, there exists a unique $a(A) \in \O$ such that $a(A) \sqsubseteq A$. We let:
    \begin{align*}
        \X := \{A \in \cal{R}_\alpha : a(A) \in \O_0\}.
    \end{align*}

    $\X$ is a clopen set, so there exists some $A \in \U$ such that $[A] \subseteq \X$ or $[A] \subseteq \X^c$. Then $\cal{AR}\restrictedto A \cap \O_1 = \emptyset$ or $\cal{AR}\restrictedto A \cap \O_0 = \emptyset$ respectively.
\end{proof}

Di Prisco, Mijares and Nieto also proposed abstract notions of \emph{selective} ultrafilters and semiselective coideals for topological Ramsey spaces in \cite{DMN17}.

\subsubsection{Existence of Ramsey ultrafilters}
In all the previous examples, it was shown that Ramsey ultrafilters exist under $\CH$, and there is a natural forcing that adds a Ramsey ultrafilter. These examples satisfy a common property, in which the partial order $(\cal{R},\leq)$ shares an isomorphic separative quotient with a $\sigma$-closed quasi-order.

\begin{definition}[{\cite[Definition 52]{DN22}}]
\label{def:sigma-closed.quasi-order}
    A triple $(\cal{R},\leq,\leq^*)$ is a \emph{$\sigma$-closed} (Borel, Ramsey etc.) quasi-order if the following hold:
    \begin{enumerate}
        \item $(\cal{R},\leq^*)$ is a $\sigma$-closed (Borel, Ramsey etc.) quasi-order coarsening $(\cal{R},\leq)$.
        
        \item $(\cal{R},\leq^*)$ and $(\cal{R},\leq)$ share the same separative quotient.
    \end{enumerate}
\end{definition}

Dobrinen and Navarro-Flores introduced the notion of a \emph{$\sigma$-closed topological Ramsey space} in \cite{DN22}, which is a closed triple $(\cal{R},\leq,r)$ satisfying \textbf{A1}--\textbf{A4} for which $(\cal{R},\leq,\leq^*)$ is a $\sigma$-closed Ramsey quasi-order. If one defines a relation $\leq^*$ on $\cal{R}$ by stipulating that $A \leq^* B$ iff $[a,A] \subseteq [a,B]$ for some $a \in \cal{AR}\restrictedto A$, then $(\E_\infty,\leq^*)$ and $(\cal{R}_1,\leq^*)$ are $\sigma$-closed --- see \cite[Proposition 4.2]{M86} and \cite[Lemma 1.2]{L89}. Additionally, all unconditional spaces are $\sigma$-closed, witnessed by the relation $A \leq^* B$ iff $A/n \leq B$ for some $n$. On the other hand, forcing with strong subtrees $(\S_\infty(U),\subseteq)$, Sacks forcing $(\mathbb{S},\subseteq)$ or Carlson--Simpson space with variable words $(\E_\infty^L,\leq,r)$ (where $|L| \geq 2$) adds new reals, so they do not share the same separative quotient with any $\sigma$-closed quasi-order.

\begin{theorem}
\label{thm:CH.implies.Ramsey.ultrafilter}
    Assume $\CH$. If $(\cal{R},\leq,\leq^*)$ is a $\sigma$-closed Ramsey quasi-order, then for all $A \in \cal{R}$, there is a Ramsey ultrafilter $\U \subseteq \cal{R}$ such that $A \in \U$.
\end{theorem}

\begin{proof}
    Let $\{\X_\alpha : \alpha < \frak{c}\}$ be an enumeration of all analytic subsets of $\cal{R}$. We fix a bijection $\Gamma : \{(\beta,\gamma) : \gamma \leq \beta < \frak{c}\} \to \frak{c}$, and note that we may choose $\Gamma$ so that if $\Gamma(\beta,\gamma) = \alpha$, then $\beta \leq \alpha$. We define a $\leq^*$-decreasing sequence $(A_\alpha)_{\alpha<\frak{c}}$ in $\cal{R}$ as follows: Let $A_0 \in \cal{R}$ be arbitrary. If $A_\alpha$ has been defined, then by the analytic Ramsey property of $(\cal{R},\leq)$, there exists some $A_{\alpha+1} \leq A_\alpha$ such that $[A_{\alpha+1}] \subseteq \X_\alpha$ or $[A_{\alpha+1}] \subseteq \X_\alpha^c$. Furthermore, if $\alpha = \Gamma(\beta,\gamma)$, then we shall choose $A_{\alpha+1}$ such that $A_{\alpha+1} \leq A_\beta$ and $A_{\alpha+1} \leq A_\gamma$ (this is possible as $A_\alpha \leq^* A_\beta \leq^* A_\gamma$, so they are compatible in $(\cal{R},\leq^*)$, hence they are also compatible in $(\cal{R},\leq)$ as they have isomorphic separative quotients). Finally, if $\alpha$ is a limit ordinal and $A_\beta$ has been defined for $\beta < \alpha$, then we let $A_\alpha \in \cal{R}$ be any element such that $A_\alpha \leq^* A_\beta$ for all $\beta < \alpha$ (which is possible as $\alpha$ is countable and $(\cal{R},\leq^*)$ is $\sigma$-closed). 
    
    Now let:
    \begin{align*}
        \U := \{A \in \cal{R} : A_\alpha \leq A \text{ for some } \alpha < \frak{c}\}.
    \end{align*}
    We first show that $\U$ is a filter. It is clear that $\U$ is upward-closed. If $A,B \in \U$, then $A_\beta \leq A$ and $A_\gamma \leq B$ for some $\beta,\gamma < \frak{c}$. Assuming that $\gamma \leq \beta$, we may let $\alpha := \Gamma(\beta,\gamma)$, and we have that $A_{\alpha+1} \leq A$ and $A_{\alpha+1} \leq B$. Finally, to see that $\U$ is Ramsey, if $\X_\alpha$ is an analytic subset of $\cal{R}$ (since every analytic subset of $\cal{R}$ appears in the enumeration $\{\X_\alpha : \alpha < \frak{c}\}$), we have that $A_{\alpha+1} \in \U$, and $[A_{\alpha+1}] \subseteq \X_\alpha$ or $[A_{\alpha+1}] \subseteq \X_\alpha^c$. Finally, since $A_0$ may be arbitrarily chosen and $A_0 \in \U$, the theorem is proved.
\end{proof}

\begin{lemma}
\label{lem:generic.filter.is.Ramsey}
    If $(\cal{R},\leq,\leq^*)$ is a $\sigma$-closed Ramsey quasi-order, then forcing with $(\cal{R},\leq^*)$ adds a Ramsey ultrafilter.
\end{lemma}

\begin{proof}
    Let $\V$ be a model of set theory, and let $\U$ be $(\cal{R},\leq^*)$-generic over $\V$. Since $(\cal{R},\leq)$ and $(\cal{R},\leq^*)$ have isomorphic separative quotients and $(\cal{R},\leq^*)$ is $\sigma$-closed, no new elements of $\cal{R}$ are added in the generic extension. Therefore, every analytic subset of $\cal{R}$ is in the ground model (as analytic subsets of $\cal{R}$ are encoded by some real). Thus, it suffices to show that for every analytic $\X \subseteq \cal{R}$ in the ground model, there exists some $A \in \U$ such that $[A] \subseteq \X$ or $[A] \subseteq \X^c$.
    
    Let $\X \subseteq \cal{R}$ be analytic. We define (in the ground model):
    \begin{align*}
        \D := \{A \in \cal{R} : [A] \subseteq \X \text{ or } [A] \subseteq \X^c\}.
    \end{align*}
    Since every analytic set is Ramsey and $(\cal{R},\leq)$ and $(\cal{R},\leq^*)$ have isomorphic separative quotients, $\D$ is a dense subset of $(\cal{R},\leq^*)$ in the ground model. Therefore, $\U \cap \D \neq \emptyset$, i.e. there exists some $A \in \U$ such that $[A] \subseteq \X$ or $[A] \subseteq \X^c$.
\end{proof}

See also the discussion after \cite[Definition 47]{DMT17}, which gives analogues of Theorem \ref{thm:CH.implies.Ramsey.ultrafilter} and Lemma \ref{lem:generic.filter.is.Ramsey} for a variant of \emph{selective ultrafilters}, the latter due to Mijares \cite{M07}.

\appsubsection{Abstract orders and abstract Ramsey spaces}
\label{sec:ARS}
In \cite{T10}, topological Ramsey spaces were presented as a subclass of \emph{abstract Ramsey spaces}. An abstract Ramsey space is a structure of the form $(\cal{R},\S,\leq,\leq^o,r,s)$, where $\S$ is closed in $\cal{AS}^\N$, satisfying four axioms \textbf{A1}--\textbf{A4} that generalise those presented in \S\ref{sec:TRS.axioms}. We refer readers to \cite[\S4.2]{T10} for a description of the axioms.

\begin{definition}
    We say that $(\cal{R},\S,\leq,\leq^o,r,s)$ is a \emph{closed sextuple} if $\cal{R}$ is a metrically closed subset of $\cal{AR}^\N$ and $\S$ is a metrically closed subset of $\cal{AS}^\N$.
\end{definition}

We remark that Todor\v{c}evi\'{c} does not assume that $\cal{R}$ is closed in $\cal{AR}^\N$ in \cite{T10}, and the abstract Ramsey theorem below does not require it.

\begin{theorem}[{\cite[Abstract Ramsey theorem]{T10}}]
\label{thm:abstract.Ramsey.thm}
    Let $(\cal{R},\S,\leq,\leq^o,r,s)$ be a closed sextuple satisfying \textbf{A1}--\textbf{A4}. Then a subset of $\cal{R}$ is $\S$-Ramsey iff it is $\S$-Baire, and it is $\S$-Ramsey null iff it is $\S$-meagre.
\end{theorem}

In particular, if $\cal{AR}$ is countable (so the metric topology is Polish), then every analytic subset of $\cal{R}$ is $\S$-Ramsey. A closed sextuple $(\cal{R},\S,\leq,\leq^o,r,s)$ satisfying \textbf{A1} and \textbf{A2} is an \emph{abstract Ramsey space} if it satisfies the conclusion of Theorem \ref{thm:abstract.Ramsey.thm}. Theorem \ref{thm:abstract.Ramsey.thm} thus asserts that every closed sextuple satisfying \textbf{A1}--\textbf{A4} is an abstract Ramsey space. Two important examples of abstract Ramsey spaces were illustrated in \cite{T10}, which we shall describe below.

\begin{example}[Hales--Jewett space]
    We continue from the setup in Example \ref{ex:Hales-Jewett}. We may extend the definition of a rapidly increasing sequence to sequences of elements in $W_L$, and let $W_L^{[\infty]}$ and $W_L^{[<\infty]}$ denote the set of all infinite and finite rapidly increasing sequences of elements of $W_L$ respectively. Given $x = (w_n)_{n<N} \in W_{Lv}^{[<\infty]}$, we let:
    \begin{align*}
        \c{x}_{Lv} := \{w_{n_0}[\lambda_0]^\frown\cdots^\frown w_{n_k}[\lambda_k] : \, &n_0 < \cdots < n_k < N, \, \lambda_i \in L_{n_i} \cup \{v\}, \\
        &\text{ and } \lambda_i = v \text{ for some } i\}, \\
        \c{x}_L := \{w_{n_0}[\lambda_0]^\frown\cdots^\frown w_{n_k}[\lambda_k] : \, &n_0 < \cdots < n_k < N, \, \lambda_i \in L_{n_i}\}.
    \end{align*}
    Given $a = (w_n)_{n<N} \in W_L^{[<\infty]}$, we also let:
    \begin{align*}
        \c{a}_L := \{w_{n_0}^\frown\cdots^\frown w_{n_k} : n_0 < \cdots < n_k < N\}.
    \end{align*}
    We also define $\c{X}_L := \bigcup_{n<\omega} \c{s_n(X)}_L$ and $\c{X}_{Lv} := \bigcup_{n<\omega} \c{s_n(X)}_{Lv}$ for $X \in W_{Lv}^{[\infty]}$. Under this notation, we have that for all $x,y \in W_{Lv}^{[<\infty]}$, $x \leq_\fin y$ iff $\c{x}_{Lv} \subseteq \c{y}_{Lv}$, and for all $y' \sqsubsetneq y$, $\c{x}_{Lv} \nsubseteq \c{y'}_{Lv}$. Given $a \in W_L^{[<\infty]}$ and $x \in W_{Lv}^{[<\infty]}$, we write $a \leq_\fin^o x$ iff $\c{a}_L \subseteq \c{x}_L$, and for all $x' \sqsubsetneq x$, $\c{a}_L \nsubseteq \c{x'}_L$. Then $(W_L^{[\infty]},W_{Lv}^{[\infty]},\leq,\leq^o,r,s)$ is a closed sextuple satisfying \textbf{A1}--\textbf{A4} --- see \cite[\S4.3]{T10}.
\end{example}

\begin{example}[Block sequences of located variable-words]
    We continue from the setup in Example \ref{ex:block.sequences.located.variable-words}. Let $\FIN_L^{[\infty]}$ and $\FIN_L^{[<\infty]}$ denote the set of infinite and finite $<$-increasing sequences of elements of $\FIN_L$. Given $x = (w_n)_{n<N} \in \FIN_{Lv}^{[<\infty]}$, we let:
    \begin{align*}
        \c{x}_{Lv} := \{w_{n_0}[\lambda_0]\cup\cdots\cup w_{n_k}[\lambda_k] : \, &n_0 < \cdots < n_k < N, \, \lambda_i \in L_{n_i} \cup \{v\}, \\
        &\text{ and } \lambda_i = v \text{ for some } i\}, \\
        \c{x}_L := \{w_{n_0}[\lambda_0]\cup\cdots\cup w_{n_k}[\lambda_k] : \, &n_0 < \cdots < n_k < N, \, \lambda_i \in L_{n_i}\}.
    \end{align*}
    Given $a = (w_n)_{n<N} \in \FIN_L^{[<\infty]}$, we also let:
    \begin{align*}
        \c{a}_L := \{w_{n_0} \cup \cdots \cup w_{n_k} : n_0 < \cdots < n_k < N\}.
    \end{align*}
    We also define $\c{X}_L := \bigcup_{n<\omega} \c{s_n(X)}_L$ and $\c{X}_{Lv} := \bigcup_{n<\omega} \c{s_n(X)}_{Lv}$ for $X \in \FIN_{Lv}^{[\infty]}$. Under this notation, we have that for all $x,y \in \FIN_{Lv}^{[<\infty]}$, $x \leq_\fin y$ iff $\c{x}_{Lv} \subseteq \c{y}_{Lv}$, and for all $y' \sqsubsetneq y$, $\c{x}_{Lv} \nsubseteq \c{y'}_{Lv}$. Given $a \in \FIN_L^{[<\infty]}$ and $x \in \FIN_{Lv}^{[<\infty]}$, we write $a \leq_\fin^o x$ iff $\c{a}_L \subseteq \c{x}_L$, and for all $x' \sqsubsetneq x$, $\c{a}_L \nsubseteq \c{x'}_L$. Then $(\FIN_L^{[\infty]},\FIN_{Lv}^{[\infty]},\leq,\leq^o,r,s)$ is a closed sextuple satisfying \textbf{A1}--\textbf{A4} --- see \cite[Theorem 4.48]{T10}.
\end{example}

We shall now introduce natural analogues of analytic/Ramsey quasi-orders for abstract Ramsey spaces.

\begin{definition}
    An \emph{abstract order} is a quadruple $(\P,\Q,\leq,\leq^o)$ that satisfies the following:
    \begin{enumerate}
        \item $(\Q,\leq)$ is a quasi-order.

        \item ${\leq^o} \subseteq \P \times \Q$.

        \item If $p \leq^o q$ and $q \leq q'$, then $p \leq^o q'$.
    \end{enumerate}
\end{definition}

\begin{definition}
    An \emph{analytic} (\emph{Borel}, $G_\delta$ etc.) \emph{abstract order} is an abstract order $(\cal{R},\S,\leq,\leq^o)$ such that $\cal{R}$ and $\S$ are Polish spaces, $\leq$ is an analytic (Borel, $G_\delta$ etc.) subset of $\S^2$, and $\leq^o$ is an analytic (Borel, $G_\delta$ etc.) subset of $\cal{R} \times \S$.
\end{definition}

If $(\cal{R},\S,\leq,\leq^o)$ is an analytic (Borel, $G_\delta$ etc.) abstract order and $X \in \S$, then we write:
\begin{align*}
    [X] := \{A \in \cal{R} : A \leq^o X\},
\end{align*}
which is an analytic (Borel, $G_\delta$ etc.) subset of $\cal{R}$.

\begin{definition}
    Let $(\cal{R},\S,\leq,\leq^o)$ be an analytic abstract order.
    \begin{enumerate}
        \item $(\cal{R},\S,\leq,\leq^o)$ has the \emph{weak} (\emph{analytic}, \emph{Borel} etc.) \emph{Ramsey property} if for every (analytic, Borel etc.) $\X \subseteq \cal{R}$, there exists some $X \in \S$ such that $[X] \neq \emptyset$, and $[X] \subseteq \X$ or $[X] \subseteq \X^c$.

        \item $(\cal{R},\S,\leq,\leq^o)$ has the (\emph{analytic}, \emph{Borel} etc.) \emph{Ramsey property} if for every (analytic, Borel etc.) $\X \subseteq \cal{R}$ and $X \in \S$, there exists some $Y \leq X$ such that $[Y] \neq \emptyset$, and $[Y] \subseteq \X$ or $[Y] \subseteq \X^c$.

        \item $(\cal{R},\S,\leq,\leq^o)$ is a \emph{Ramsey abstract order} if it has the analytic Ramsey property.
    \end{enumerate}
\end{definition}

The following lemma is the abstract Ramsey space analogue of Lemma \ref{lem:A4.iff.local.Ramsey.properties}. The proof is very similar to that of Lemma \ref{lem:A4.iff.local.Ramsey.properties} --- we leave the details to the reader.

\begin{lemma}
\label{lem:A4.iff.local.Ramsey.properties.ARS}
    Let $(\cal{R},\S,\leq,\leq^o,r,s)$ be a closed sextuple satisfying \textbf{A1}--\textbf{A3} such that $\cal{AR}$ is countable. The following are equivalent:
    \begin{enumerate}
        \item $(\cal{R},\S,\leq,\leq^o,r,s)$ satisfies \textbf{A4}.

        \item For all $X \in \S$ and $a \in \cal{AR}\restrictedto X$, $([a,X],[\depth_X(a),X],\leq,\leq^o)$ is a Ramsey abstract order.
        
        \item For all $X \in \S$ and $a \in \cal{AR}\restrictedto X$, there exists some $Y \in [\depth_X(a),X]$ such that $([a,Y],[\depth_Y(a),Y],\leq,\leq^o)$ is an analytic abstract order satisfying the weak clopen Ramsey property.
    \end{enumerate}
\end{lemma}

\section{Projections between quasi-orders}
\label{sec:projections}
This section focuses on functions between quasi-orders that allow us to ``compare'' one quasi-order with another. In the context of analytic Ramsey properties, we explain the relationship between maps that transfer the (weak) analytic Ramsey property and maps that transfer genericity. We also provide an extensive list of examples of such functions between topological Ramsey spaces and briefly discuss the relationship between these functions and Ramsey ultrafilters.

For this section, $(\cal{R},\leq_\cal{R})$ and $(\S,\leq_\S)$ always denote analytic quasi-orders, while $(\P,\leq_\P)$ and $(\Q,\leq_\Q)$ always denote quasi-orders with no additional topological properties. In contexts where the quasi-order of interest is obvious, we use $\leq$ to denote the order in place of the subscripted order.

\begin{notation}
    Let $(\P,\leq)$ be a quasi-order. If $\cal{K} \subseteq \P$, we write:
    \begin{align*}
        \cal{K}^\uparrow &:= \{p \in \P : p' \leq p \text{ for some } p' \in \E\}, \\
        \cal{K}^\downarrow &:= \{p \in \P : p \leq p' \text{ for some } p' \in \E\}.
    \end{align*}
    Given $p \in \P$, we also write $p^\uparrow := \{p\}^\uparrow$ and $p^\downarrow := \{p\}^\downarrow$.
\end{notation}

For $A \in \cal{R}$ in an analytic quasi-order, $[A] = A^\downarrow$ (see the paragraph before Definition \ref{def:Ramsey.quasi-order}), but we reserve the notation $[A]$ for analytic quasi-orders.

\subsection{Maps preserving the analytic Ramsey property}
\label{sec:quasi-open}
As mentioned in the introduction, Ma\v{s}ulovi\'{c} introduced the notion of \emph{pre-adjunctions} between categories, and showed that pre-adjunctions transfer the weak Ramsey property. We recap the definition. 

\begin{definition}[{\cite[Definition 3.1]{M18}}]
\label{def:pre-adjunction}
    Let $\mathbf{C},\mathbf{D}$ be two categories. A pair of maps:
    \begin{align*}
        F : \Ob(\mathbf{D}) \rightleftarrows \Ob(\mathbf{C}) : G
    \end{align*}
    is a \emph{pre-adjunction between $\mathbf{C}$ and $\mathbf{D}$} if there is a family of maps:
    \begin{align*}
        \Phi_{\D,\cal{C}} : \hom_{\mathbf{C}}(F(\D),\cal{C}) \to \hom_{\mathbf{D}}(\D,G(\cal{C}))
    \end{align*}
    indexed by the family $\{(\D,\cal{C}) \in \Ob(\mathbf{D}) \times \Ob(\mathbf{C}) : \hom_{\mathbf{C}}(F(\D),\cal{C}) \neq \emptyset\}$ satisfying the following:
    \begin{itemize}
        \item[(PA)] For all $\cal{C} \in \Ob(\mathbf{C})$, $\D,\E \in \Ob(\mathbf{D})$, $u \in \hom_{\mathbf{C}}(F(\D),\cal{C})$ and $f \in \hom_{\mathbf{D}}(\E,\D)$, there exists some $v \in \hom_{\mathbf{C}}(F(\E),F(\D))$ such that $\Phi_{\D,\cal{C}}(u) \circ f = \Phi_{\E,\cal{C}}(u \circ v)$.
    \end{itemize}
\end{definition}

\begin{theorem}[{\cite[Theorem 3.2]{M18}}]
\label{thm:category.pre-adjunction.transfer.weak.Ramsey.property}
    Let $\mathbf{C},\mathbf{D}$ be two categories. If $\mathbf{C}$ has the weak Ramsey property and there is a pre-adjunction between $\mathbf{C}$ and $\mathbf{D}$, then $\mathbf{D}$ also has the weak Ramsey property.
\end{theorem}

As explained in \S\ref{sec:Ramsey.quasi-orders}, categories capturing infinite-dimensional Ramsey theory of topological Ramsey spaces only have a single object $\Ob(\mathbf{C}) = \{\infty\}$, so such a category is equivalent to a monoid. In this case, given two monoids $(M,\circ_M)$ and $(N,\circ_N)$, Definition \ref{def:pre-adjunction} collapses to a single map $\Phi : M \to N$ satisfying the following:

\begin{definition}
\label{def:pre-adjunction.monoid}
    Let $(M,\circ_M)$ and $(N,\circ_N)$ be two monoids. A function $\Phi : M \to N$ is a \emph{pre-adjunction} if for all $A \in M$ and $X \in N$, there exists some $B \in M$ such that $\Phi(A) \circ_N X = \Phi(A \circ_M B)$.
\end{definition}

As discussed in \S\ref{sec:Ramsey.quasi-orders}, each monoid $(M,\circ)$ induces a natural quasi-order $\leq$ on $M$ by:
\begin{align*}
    B \leq A \iff \exists C \in M \, B = A \circ C,
\end{align*}
and if $(M,\circ)$ is a Polish monoid then $(M,\leq)$ is an analytic quasi-order. Definition \ref{def:pre-adjunction.monoid} may now be rephrased in the following manner:

\begin{lemma}
    Let $(M,\circ_M)$ and $(N,\circ_N)$ be two monoids. A function $\Phi : M \to N$ is a pre-adjunction iff for all $A \in M$, $[\Phi(A)] \subseteq \Phi``[A]$.
\end{lemma}

This motivates the following natural generalisation of pre-adjunctions to quasi-orders:

\begin{definition}
\label{def:pre-adjunction.quasi-order}
    A function $f : \P \to \Q$ is a \emph{pre-adjunction} if for all $p \in \P$, $f(p)^\downarrow \subseteq f``(p^\downarrow)$.
\end{definition}

This motivates a topological version of Theorem \ref{thm:category.pre-adjunction.transfer.weak.Ramsey.property} for analytic quasi-orders. In fact, the weaker \emph{quasi-open} property suffices.

\begin{definition}
\label{def:quasi-open}
    A function $f : \P \to \Q$ is \emph{quasi-open} if for all $p \in \P$, there exists some $q \in \Q$ such that $q^\downarrow \subseteq f``(p^\downarrow)$. In other words, $f$ is quasi-open under the order topologies of $(\P,\leq_\P)$ and $(\Q,\leq_\Q)$.
\end{definition}

\begin{lemma}
\label{lem:cts.qo.transfers.weak.analytic.Ramsey.property}
    If $(\cal{R},\leq_\cal{R})$ satisfies the weak analytic Ramsey property and there exists a continuous quasi-open $f : \cal{R} \to \S$, then $(\S,\leq_\S)$ satisfies the weak analytic Ramsey property.
\end{lemma}

\begin{proof}
    Let $\Y \subseteq \S$ be analytic. Since $f$ is continuous, $f^{-1}(\Y)$ is analytic, so by the weak analytic Ramsey property of $(\cal{R},\leq_\cal{R})$ there exists some $A \in \cal{R}$ such that $[A] \subseteq f^{-1}(\Y)$ or $[A] \subseteq (f^{-1}(\Y))^c = f^{-1}(\Y^c)$. This implies that $f``[A] \subseteq \Y$ or $f``[A] \subseteq \Y^c$. Since $f$ is quasi-open, there exists some $X \in \S$ such that $[X] \subseteq f``[A]$. Then $[X] \subseteq \Y$ or $[X] \subseteq \Y^c$, as desired.
\end{proof}

\subsubsection{Transferring the analytic Ramsey property}
To transfer the analytic Ramsey property, it suffices to have a collection of maps from $(\cal{R},\leq_\cal{R})$ to $([X],\leq_\S)$, one for each $X \in \S$.

\begin{lemma}
\label{lem:cts.qo.family.transfers.analytic.Ramsey.property}
    If $(\cal{R},\leq_\cal{R})$ satisfies the weak analytic Ramsey property, and for all $X \in \S$, there exists a continuous quasi-open $f_X : \cal{R} \to [X]$, then $(\S,\leq_\S)$ satisfies the analytic Ramsey property.
\end{lemma}

\begin{proof}
    Let $\Y \subseteq \S$ be analytic, and let $X \in \S$. Since $f_X$ is continuous as a map $\cal{R} \to \S$, $f_X^{-1}(\Y)$ is analytic, so by the weak analytic Ramsey property of $(\cal{R},\leq_\cal{R})$ there exists some $A \in \cal{R}$ such that $[A] \subseteq f_X^{-1}(\Y)$ or $[A] \subseteq (f_X^{-1}(\Y))^c = f_X^{-1}(\Y^c)$. This implies that $f_X``[A] \subseteq \Y$ or $f_X``[A] \subseteq \Y^c$. Since $f_X$ is quasi-open, there exists some $Y \in [X]$ such that $[Y] \subseteq f_X``[A]$. Then $Y \leq X$, and $[Y] \subseteq \Y$ or $[Y] \subseteq \Y^c$. Since $X$ is arbitrary, $(\S,\leq_\S)$ has the analytic Ramsey property.
\end{proof}

If, however, we wish to transfer the analytic Ramsey property via a single map $f : \cal{R} \to \S$, then there are two natural approaches:
\begin{itemize}
    \item Restrict $(\S,\leq_\S)$ to analytic quasi-orders that satisfy \emph{property (H)}, a weak form of cone homogeneity. 

    \item Require the continuous quasi-open map $f : \cal{R} \to \S$ to be \emph{monotone} and \emph{cofinal}.
\end{itemize}

\begin{definition}
\label{def:property.H}
    An analytic quasi-order $(\cal{R},\leq)$ satisfies \emph{property (H)} if for all $A \in \cal{R}$, there is a continuous quasi-open map $f : \cal{R} \to [A]$.
\end{definition}

\begin{proposition}
\label{prop:cts.qo.to.H.transfers.analytic.Ramsey.property}
    Suppose that $(\S,\leq_\S)$ satisfies property (H). If $(\cal{R},\leq_\cal{R})$ satisfies the weak analytic Ramsey property and there exists a continuous quasi-open $f : \cal{R} \to \S$, then $(\S,\leq_\S)$ satisfies the analytic Ramsey property.
\end{proposition}

\begin{proof}
    This follows from Lemma \ref{lem:cts.qo.family.transfers.analytic.Ramsey.property} --- for any $X \in \S$, by property (H) there exists some continuous quasi-open $g_X : \S \to [X]$. Then $g_X \circ f : \cal{R} \to [X]$ is continuous quasi-open.
\end{proof}

Any analytic quasi-order induced by a Polish monoid satisfies property (H):

\begin{proposition}
\label{prop:Polish.monoid.property.H}
    Let $(M,\circ)$ be a Polish monoid, and let $(M,\leq)$ be the analytic quasi-order defined by:
    \begin{align*}
        B \leq A \iff \exists C \in M \, B = A \circ C.
    \end{align*}
    Then for all $A \in M$, there is a continuous monotone pre-adjunction (in particular, quasi-open) $g_A : M \to [A]$.
\end{proposition}

\begin{proof}
    Define $g_A : M \to [A]$ by $g_A(B) := A \circ B$. By the definition of $\leq$, $g_A$ is well-defined. $g_A$ is continuous as $\circ$ is continuous. If $C \leq B$, then $C = B \circ D$ for some $D \in M$, so:
    \begin{align*}
        A \circ C = A \circ (B \circ D) = (A \circ B) \circ D.
    \end{align*}
    Therefore $A \circ C \leq A \circ B$, so $g_A$ is monotone. Finally, to see that $g_A$ is a pre-adjunction --- let $C \in M$, and let $B := g_A(C) = A \circ C \in [A]$. If $B' \in [B]$, then $B' = B \circ D$ for some $D \in M$, so $B' = A \circ (C \circ D)$. Then $C \circ D \leq C$, and $g_A(C \circ D) = B'$. Therefore, $[g_A(C)] \subseteq g_A``[C]$, as required.
\end{proof}

Since the operations $\circ$ in Examples \ref{ex:Ellentuck.space.monoid} and \ref{ex:Hindman.FIN.monoid} are continuous, $(M_{\mathbf{E}},\circ)$ and $(M_{\mathbf{BS}},\circ)$ are Polish monoids. In particular, the Ellentuck space $([\N]^\infty,\subseteq,r)$ and Hindman space $(\FIN^{[\infty]},\leq,r)$ satisfy property (H). 

\begin{definition}
    A function $f : \P \to \Q$ is \emph{monotone} if $f(p') \leq f(p)$ whenever $p' \leq p$.
\end{definition}

\begin{definition}
    \hfill
    \begin{enumerate}
        \item A subset $\D \subseteq \P$ is \emph{cofinal} (or \emph{dense}) if for all $p \in \P$, there exists some $p' \in \D$ such that $p' \leq p$.
        
        \item A function $f : \P \to \Q$ is \emph{cofinal} if for every cofinal subset $\D$ of $\P$, $f``(\D)$ is a cofinal subset of $\Q$.
    \end{enumerate}
\end{definition}

Note that a monotone map is cofinal iff the range is cofinal.

\begin{proposition}
\label{prop:cofinal.qo.transfers.Ramsey}
    Suppose that $(\cal{R},\leq_\cal{R})$ is a Ramsey quasi-order, and there exists a continuous cofinal monotone quasi-open map $f : \cal{R} \to \S$. Then $(\S,\leq_\S)$ is a Ramsey quasi-order.
\end{proposition}

\begin{proof}
    Let $\X \subseteq \S$ be analytic, and let $X \in \S$. Since $f$ is continuous, $f^{-1}(\X) \subseteq \cal{R}$ is analytic. By the cofinality of $f$, we may let $A \in \cal{R}$ be such that $f(A) \leq X$. Since $(\cal{R},\leq_\cal{R})$ has the analytic Ramsey property, there exists some $B \leq A$ such that $[B] \subseteq f^{-1}(\X)$ or $[B] \subseteq (f^{-1}(\X))^c = f^{-1}(\X^c)$. This implies that $f``[B] \subseteq \X$ or $f``[B] \subseteq \X^c$. By the quasi-open property of $f$, let $Y \in \S$ be such that $[Y] \subseteq f``[B]$. Then $[Y] \subseteq \X$ or $[Y] \subseteq \X^c$. Furthermore, there exists some $C \in [B]$ such that $f(C) = Y$, and:
    \begin{align*}
        Y = f(C) \leq f(B) \leq f(A) \leq X,
    \end{align*}
    as needed.
\end{proof}

While we believe that all examples in \S\ref{sec:examples} satisfy property (H), a general framework that does not assume monotonicity appears to be considerably more involved, so we shall focus on monotone maps. An important consequence of monotonicity is automatic continuity: we shall prove in \S\ref{sec:automatic.continuity} that, under suitable assumptions, every monotone function from a topological Ramsey space to an analytic quasi-order is continuous on a cone (Lemma \ref{lem:monotone.automatic.continuity}).

Note that in the hypothesis of Proposition \ref{prop:cofinal.qo.transfers.Ramsey}, we have, in particular, that $\S$ is \emph{Tukey reducible} to $\cal{R}$.

\begin{definition}
\label{def:Tukey.reducible}
    Given quasi-orders $(\P,\leq_\P)$ and $(\Q,\leq_\Q)$, $\Q$ is \emph{Tukey reducible} to $\P$, denoted $\Q \leq_T \P$, if there is a cofinal map $f : \P \to \Q$. We say that $\P$ and $\Q$ are \emph{Tukey equivalent}, denoted $\P \equiv_T \Q$, if $\Q \leq_T \P$ and $\P \leq_T \Q$. We write $\Q <_T \P$ if $\Q \leq_T \P$ and $\Q \not\equiv_T \P$.
\end{definition}

We shall see in the next subsection and \S\ref{sec:non-existence.of.projections} that the Tukey theory of ultrafilters plays an important role in the existence of certain maps between quasi-orders.

\subsection{Maps preserving generic filters}
The key observation that motivates this subsection is that quasi-open maps are defined for arbitrary quasi-orders, and that monotone quasi-open maps send generic filters to generic filters (Lemmas \ref{lem:ab.equivalence} and \ref{lem:qo.implies.ab}). For the remainder of this paper, we shall write ``ctm of $\ZFC$'' to denote a countable transitive model of a sufficiently large fragment of $\ZFC$.

\subsubsection{Alaoglu--Birkhoff convergence}
\begin{definition}[\cite{AB40,TZ99}]
\label{def:ab.convergent}
    A function $f : \P \to \Q$ is \emph{Alaoglu--Birkhoff convergent} if for every dense open $\D \subseteq \Q$, $f^{-1}(\D)$ contains a dense open subset of $\P$.
\end{definition}

\begin{lemma}
\label{lem:ab.equivalence}
    Let $M$ be a ctm of $\ZFC$ with $\P,\Q \in M$. For any monotone function $f : \P \to \Q$ in $M$, the following properties are equivalent:
    \begin{enumerate}
        \item $M \models \text{$f$ is Alaoglu--Birkhoff convergent}$.

        \item If $G$ is a $\P$-generic filter over $M$, then $H := f``(G)^\uparrow$ is $\Q$-generic over $M$ (and $H \leq_T G$).
    \end{enumerate}
\end{lemma}

\begin{proof}
    \underline{(1)$\implies$(2):} We work in $M[G]$. Since $f$ is monotone, it is clear that $H$ is a filter, so it remains to show that $H$ is $\Q$-generic over $M$. Let $\D$ be a dense open subset of $\Q$ in $M$. By our hypothesis, $f^{-1}(\D)$ contains a dense open subset of $\P$ in $M$, so there exists some $p \in G \cap f^{-1}(\D)$. This implies that $f(p) \in H \cap \D$, so $H \cap \D \neq \emptyset$. Since $\D$ is arbitrary, we have that $H$ is $\Q$-generic over $M$. Finally, $f\restrictedto G$ is a monotone map whose range is cofinal in $H$, so $H \leq_T G$.

    \underline{(2)$\implies$(1):} Let $\D$ be a dense open subset of $\Q$ in $M$, fix some $p \in \P$, and let $G$ be a $\P$-generic filter over $M$ with $p \in G$. We work in $M[G]$. Since $H$ is $\Q$-generic over $M$, $H \cap \D \neq \emptyset$. In other words, there exists some $p' \in G$ such that $f(p') \in \D$, and since $G$ is a filter, we may choose $p'$ such that $p' \leq p$. Then $p' \in G \cap f^{-1}(\D)$, so $f^{-1}(\D)$ is dense in $\P$ as $p$ is arbitrary. Furthermore, since $f$ is monotone and $\D$ is open, if $p \in f^{-1}(\D)$ and $p' \leq p$, then $f(p') \leq f(p) \in \D$ so $f(p') \in \D$. Thus, $f^{-1}(\D)$ contains a dense open subset of $\P$.
\end{proof}

Since we are focusing on monotone maps, we reserve a term for monotone Alaoglu--Birkhoff convergent maps.

\begin{definition}
\label{def:weak.projection}
    A map $f : \P \to \Q$ is a \emph{weak projection} if it is monotone and Alaoglu--Birkhoff convergent.
\end{definition}

Projections maps between posets under various definitions are well-studied. The basic definitions and facts about them appear to be well-known to the community from the early days of the development of forcing, although we don't have a concrete reference for this.

\begin{lemma}
\label{lem:monotone.preimage.dense.implies.ab}
    Let $f : \P \to \Q$ be monotone, and suppose that for all dense open $\D \subseteq \Q$, $f^{-1}(\D)$ is dense. Then $f$ is Alaoglu--Birkhoff convergent.
\end{lemma}

\begin{proof}
    Let $\D \subseteq \Q$ be dense open. For any $p \in \P$, let $r(p) \in f^{-1}(\D)$ be any element such that $r(p) \leq p$, which exists as $f^{-1}(\D)$ is dense. Define:
    \begin{align*}
        \K := \{s \in \P : \exists p \in \P \, s \leq r(p)\}.
    \end{align*}
    Since $r(p) \in \K$ for all $p \in \P$, $\K$ is dense. It is clear that $\K$ is open. Furthermore, if $s \leq r(p)$, then since $\D$ is open and $f$ is monotone, we have that $f(s) \leq f(r(p)) \in \D$, so $f(s) \in \D$. Therefore, $\K \subseteq f^{-1}(\D)$ and $\K$ is dense open.
\end{proof}

\begin{lemma}
\label{lem:qo.implies.ab}
    If $f : \P \to \Q$ is quasi-open, then for all dense subsets $\D$ of $\Q$, $f^{-1}(\D)$ is a dense subset of $\P$. Therefore, every monotone quasi-open map is Alaoglu--Birkhoff convergent.
\end{lemma}

\begin{proof}
    Let $\D \subseteq \Q$ be dense, and let $p \in \P$. Since $f$ is quasi-open, there exists some $q \in \Q$ such that $q^\downarrow \subseteq f``(p^\downarrow)$. Since $\D$ is dense in $\Q$, there exists some $q' \in \D$ such that $q' \leq q$. Then $q' \in q^\downarrow \subseteq f``(p^\downarrow)$, so there exists some $p' \leq p$ such that $f(p') = q' \in \D$. Then $p' \in f^{-1}(\D)$ and $p' \leq p$. 

    If $f$ is also monotone, then $f$ is Alaoglu--Birkhoff convergent by Lemma \ref{lem:monotone.preimage.dense.implies.ab}.
\end{proof}

The converse holds for continuous maps into Ramsey quasi-orders:

\begin{lemma}
\label{lem:cts.ab.implies.qo}
    Suppose that $(\S,\leq_\S)$ is a Ramsey quasi-order. If $f : \cal{R} \to \S$ is continuous and Alaoglu--Birkhoff convergent, then $f$ is quasi-open.
\end{lemma}

\begin{proof}
    Let $A \in \cal{R}$, and let:
    \begin{align*}
        \D := \{X \in \S : [X] \subseteq f``[A] \text{ or } [X] \cap f``[A] = \emptyset\}.
    \end{align*}
    Note that $f``[A]$ is a continuous image of an analytic set, so it is analytic. Therefore, by the analytic Ramsey property of $(\S,\leq_\S)$, $\D$ is a dense open subset of $\S$. Since $f$ is Alaoglu--Birkhoff convergent, $f^{-1}(\D)$ contains a dense open subset of $\cal{R}$, so there exists some $B \leq A$ such that $[f(B)] \subseteq f``[A]$ or $[f(B)] \cap f``[A] = \emptyset$. However, the latter case is not possible as $f(B) \in f``[A]$.
\end{proof}

Therefore, assuming that $(\S,\leq_\S)$ has the analytic Ramsey property, being Alaoglu--Birkhoff convergent and being quasi-open are equivalent properties of continuous monotone maps. This observation is crucial for our study of weak projections between topological Ramsey spaces in \S\ref{sec:non-existence.of.projections}, where many proofs rely on the forcing-theoretic properties of the quasi-orders underlying the topological Ramsey spaces.

\subsubsection{Projections}
We conclude this subsection with a mention of a related property: the notion of projection between quasi-orders in the sense of forcing. 

\begin{definition}
\label{def:projection}
    A function $f : \P \to \Q$ is a \emph{projection} if it satisfies the following properties:
    \begin{enumerate}
        \item $f$ is monotone.
        
        \item $f$ is cofinal.
        
        \item For all $p \in \P$ and $q \leq f(p)$, there exists some $p' \leq p$ such that $f(p') \leq q$.
    \end{enumerate}
\end{definition}

We remark that the standard definition of a projection replaces condition (2) with $f(\one_\P) = \one_\Q$, where $\one_\P$ and $\one_\Q$ are the $\leq$-maximal elements in $\P$ and $\Q$ respectively (see \cite{E20} and \cite{K23}). Since the quasi-orders of interest in this paper need not have a $\leq$-maximal element, we replace this condition with cofinality. Note, however, that these two definitions are essentially equivalent --- if $\P$ and $\Q$ have $\leq$-maximal elements and $f : \P \to \Q$ is a projection (as in Definition \ref{def:projection}), then the map $f^* : \P \to \Q$ defined by:
\begin{align*}
    f^*(p) = 
    \begin{cases}
        \one_\Q, &\text{if $p = \one_\P$}, \\
        f(p), &\text{otherwise}
    \end{cases}
\end{align*}
satisfies the standard definition of a projection (i.e. properties (1) and (3) of Definition \ref{def:projection}, and $f^*(\one_\P) = \one_\Q$). 

\begin{lemma}
\label{lem:projection.implies.ab}
    Every projection is Alaoglu--Birkhoff convergent.
\end{lemma}

\begin{proof}
    Let $f : \P \to \Q$ be a projection. Let $\D \subseteq \Q$ be a dense open subset. By Lemma \ref{lem:monotone.preimage.dense.implies.ab}, it suffices to show that $f^{-1}(\D)$ is dense. Fix any $p \in \P$. Since $\D$ is dense in $\Q$, there exists some $q \in \D$ such that $q \leq f(p)$. Since $f$ is a projection, there exists some $p' \leq p$ such that $f(p') \leq q \in \D$. Then $f(p') \in \D$, so $p' \in f^{-1}(\D)$ and $p' \leq p$.
\end{proof}

Therefore, every projection is a weak projection, and by Lemma \ref{lem:cts.ab.implies.qo}, every continuous projection into a Ramsey quasi-order is quasi-open.

\subsection{Automatic continuity}
\label{sec:automatic.continuity}
It was proven in \cite[Theorem 20]{DT11} that if $\U \subseteq [\omega]^\omega$ is a p-point, and $f : \U \to \cal{V}$ is a monotone map such that $\cal{V} = f``(\U)^\uparrow$ is an arbitrary ultrafilter, then there is a continuous map $f^* : \Po(\omega) \to \Po(\omega)$ such that $\cal{V} = (f^*)``(\U)^\uparrow$. Dobrinen and Todor\v{c}evi\'{c} proved this by first showing that there is some $A \in \U$ such that $f\restrictedto[A] : \U\restrictedto A \to \cal{V}$ is continuous. We prove an analogous result for monotone maps $f : \cal{R} \to \S$ under suitable assumptions.

\begin{definition}
\label{def:epsilon_0.property}
    Let $\{U_n\}_{n<\omega}$ be a countable base for $\cal{R}$. An analytic quasi-order $(\cal{R},\leq)$ has the \emph{$(\varepsilon_0)$-property with respect to $\{U_n\}_{n<\omega}$} if for all $n$, $A \in U_n$ and $B \in [A] \setminus U_n$, $[B] \subseteq U_n^c$. We also say that $(\cal{R},\leq)$ has the \emph{$(\varepsilon_0)$-property} if it has the $(\varepsilon_0)$-property with respect to some countable base.
\end{definition}

Let $(\cal{R},\leq,r)$ be a \textbf{wA2}-space such that $\cal{AR}$ is countable, so that $\{[a]\}_{a \in \cal{AR}}$ is a countable base for the metric topology on $\cal{R}$. Then $(\cal{R},\leq,r)$ has the $(\varepsilon_0)$-property with respect to $\{[a]\}_{a \in \cal{AR}}$ iff for all $a \in \cal{AR}$ and $A \in [a]$, if $B \leq A$ and $a \not\sqsubseteq B$ then $a \notin \cal{AR}\restrictedto B$. It is easy to verify that all topological Ramsey spaces in \S\ref{sec:examples} satisfy this property. We give a sketch of the proof for some examples below.
\begin{enumerate}
    \item $(\FIN^{[\infty]},\leq,r)$ satisfies the $(\varepsilon_0)$-property. Since $B \leq A$, if $\c{a} \subseteq \c{A}$ and $a \not\sqsubseteq B$ then $\c{a} \not\subseteq \c{B}$, so $a \notin \FIN^{[<\infty]}\restrictedto B$. A similar argument can be used to show that all unconditional spaces, except for countable vector spaces, satisfy the $(\varepsilon_0)$-property.
    
    \item $(\E_\infty,\leq,r)$ satisfies the $(\varepsilon_0)$-property. Since $B \leq A$, if $a \in \E_{<\infty}\restrictedto A$ and $a \not\sqsubseteq B$, then there exist some $m,n \in \dom(a)$ such that $(m,n) \notin a$ but $(m,n) \in B$. Then for any $C \leq B$, $(m,n) \in C$ so $a \not\sqsubseteq C$. Thus, $a \notin \E_{<\infty}\restrictedto B$.
    
    \item $(\S_\infty(U),\subseteq,r)$ and $(\mathbb{S},\subseteq,r)$ satisfy the $(\varepsilon_0)$-property. If $a \sqsubseteq A$, $B \subseteq A$ and $a \in \cal{AR}\restrictedto B$, then $a \subseteq r_{\depth_B(a)}(B)$, and $a$ and $r_{\depth_B(a)}(B)$ have the same top level, so $a = r_{\depth_B(a)}(B) \sqsubseteq B$. 
\end{enumerate}
On the other hand, if $E$ is a countable vector space over a field $\mathbb{F}$ and $|\mathbb{F}| \neq 2$, then $(E^{[\infty]},\leq,r)$ does not have the $(\varepsilon_0)$-property with respect to $\{[a]\}_{a \in E^{[<\infty]}}$. Indeed, if $A = (x_n)_{n<\omega} \in E^{[\infty]}$, $\lambda \in \mathbb{F} \setminus \{0,1\}$ and $B := (\lambda x_0)^\frown(x_n)_{n \geq 1}$, then $B \in [A] \setminus [(x_0)]$ but $(x_0) \in E^{[<\infty]}\restrictedto B$.

\begin{definition}
    Let $(\cal{R},\leq,r)$ be a \textbf{wA2}-space such that $\cal{AR}$ is countable. We say that $(\cal{R},\leq,r)$ satisfies the \emph{$(\varepsilon_0)$-property} if it has the $(\varepsilon_0)$-property with respect to $\{[a]\}_{a \in \cal{AR}}$.
\end{definition}

\begin{lemma}
\label{lem:monotone.automatic.continuity}
    Let $(\cal{R},\leq,r)$ be an \textbf{A2}-space, and let $(\S,\leq_\S)$ be an analytic quasi-order satisfying the $(\varepsilon_0)$-property. If $f : \cal{R} \to \S$ is monotone, then for all $A \in \cal{R}$, there exists some $B \leq A$ such that $f\restrictedto[B]$ is continuous.
\end{lemma}

\begin{proof}
    Let $\{U_n\}_{n<\omega}$ be a countable base for $\S$ such that $(\S,\leq_\S)$ has the $(\varepsilon_0)$-property with respect to $\{U_n\}_{n<\omega}$. Fix some $A \in \cal{R}$. For each $B \leq A$, $a \in \cal{AR}\restrictedto B$ and $n$, we write:
    \begin{itemize}
        \item $\varphi_0(a,B,n)$, if $f(C) \in U_n$ for all $C \in [a,B]$.
        
        \item $\varphi_1(a,B,n)$, if $f(C) \notin U_n$ for all $C \in [a,B]$.
        
        \item $\varphi(a,B,n)$, if $\varphi_0(a,B,n)$ or $\varphi_1(a,B,n)$.
    \end{itemize}
    For each $a \in \cal{AR}\restrictedto A$, let:
    \begin{align*}
        \K_{a,n} &:= \{B \in [a,A] : \varphi(a,B,n)\}, \\
        \K_a &:= \bigcap_{n \leq \lh(a)} \K_{a,n}.
    \end{align*}

    \begin{claim}
        $\K_{a,n}$ is dense open in $[a,A]$ for all $a \in \cal{AR}\restrictedto A$ and $n$ (so $\K_a$, being a finite intersection of such sets, is dense open in $[a,A]$).
    \end{claim}

    \begin{midproof}
        We first note that $\K_{a,n}$ is open in $[a,A]$, i.e. $\leq$-downward closed: if $\varphi_0(a,B,n)$ holds and $B' \in [a,B]$, then $[a,B'] \subseteq [a,B]$, so $f(C) \in U_n$ for all $C \in [a,B']$, i.e. $\varphi_0(a,B',n)$ holds. The same applies to $\varphi_1$.

        For density, let $B \in [a,A]$ and $n < \omega$, and let:
        \begin{align*}
            \X := \{D \in [a,B] : f(D) \in U_n\}.
        \end{align*}
        We show that there exists some $C \in [a,B]$ such that $[a,C] \subseteq \X$ or $[a,C] \subseteq \X^c$ --- then $\varphi_0(a,C,n)$ or $\varphi_1(a,C,n)$ holds respectively, so $C \in \K_{a,n}$ and $C \leq B$, as required.

        If $[a,D] \subseteq \X^c$ for some $D \in [a,B]$, then we are done by taking $C := D$. Otherwise, for all $D \in [a,B]$ there exists some $E \in [a,D]$ such that $f(E) \in U_n$ --- in particular, there exists some $D \in [a,B]$ such that $f(D) \in U_n$. We claim that $[a,D] \subseteq \X$. Indeed, if $f(E) \notin U_n$ for some $E \in [a,D]$, then $E \leq D$, so $f(E) \leq f(D)$ by monotonicity. Since $f(D) \in U_n$ and $f(E) \notin U_n$, by the $(\varepsilon_0)$-property, we have that $[f(E)] \subseteq U_n^c$, while $f(G) \notin U_n$ for all $G \in [a,E]$ and so $[a,E] \subseteq \X^c$. But $E \in [a,B]$, so this contradicts our assumption. Therefore, $[a,D] \subseteq \X$, and we take $C := D$.
    \end{midproof}

    By Lemma \ref{lem:R.is.semiselective}, there exists some $B \leq A$ that diagonalises $(\K_a)_{a\in \cal{AR}\restrictedto A}$. Thus, for all $a \in \cal{AR}\restrictedto B$ and $n \leq \lh(a)$, $\varphi(a,B,n)$ holds. We shall show that for all $n$:
    \begin{align*}
        (f\restrictedto[B])^{-1}(U_n) = \bigcup\{[a,B] : \varphi_0(a,B,n)\}.
    \end{align*}
    This is a metrically open subset of $[B]$ for all $n$, and since $\{U_n\}_{n<\omega}$ is a base for $\S$, it follows that $f\restrictedto[B]$ is continuous.

    If $C \in [a,B]$ and $\varphi_0(a,B,n)$ holds, then immediately we have that $f(C) \in U_n$. Conversely, let $C \in [B]$ such that $f(C) \in U_n$. Note that $r_n(C) \in \cal{AR}\restrictedto B$ (as $C \leq B$), so $\varphi(r_n(C),B,n)$ holds. Since $C \in [r_n(C),B]$ and $f(C) \in U_n$, $\varphi_0(r_n(C),B,n)$ must necessarily hold. Therefore, $C$ belongs to $[r_n(C),B]$, which is one of the sets in the union.
\end{proof}

Consequently, we can conclude that for monotone maps, the Alaoglu--Birkhoff convergent property and quasi-open property are equivalent under suitable assumptions. More precisely:

\begin{lemma}
\label{lem:ab.iff.qo.TRS}
    Let $(\cal{R},\leq,r)$ be an \textbf{A2}-space, and let $(\S,\leq_\S)$ be a Ramsey quasi-order satisfying the $(\varepsilon_0)$-property. For any monotone map $f : \cal{R} \to \S$, $f$ is Alaoglu--Birkhoff convergent iff $f$ is quasi-open.
\end{lemma}

\begin{proof}
    By Lemma \ref{lem:qo.implies.ab}, every quasi-open map is Alaoglu--Birkhoff convergent, so it suffices to show that the converse also holds. Let $A \in \cal{R}$. By Lemma \ref{lem:monotone.automatic.continuity}, there exists some $B \leq A$ such that $f\restrictedto[B]$ is continuous. Then $f\restrictedto[B] : [B] \to \S$ is continuous, monotone and Alaoglu--Birkhoff convergent, so by Lemma \ref{lem:cts.ab.implies.qo}, $f\restrictedto[B]$ is quasi-open. Thus, there exists some $C \leq B$ such that $[f(C)] \subseteq f``[B] \subseteq f``[A]$. Since $A$ is arbitrary, $f$ is quasi-open.
\end{proof}

The automatic continuity of monotone maps leads to improvements of earlier results. Most importantly, we can remove the continuity requirement in Proposition \ref{prop:cofinal.qo.transfers.Ramsey}:

\begin{proposition}
\label{prop:cofinal.qo.transfers.Ramsey.TRS}
    Let $(\cal{R},\leq,r)$ be a closed triple satisfying \textbf{A1}--\textbf{A4}, and let $(\S,\leq_\S)$ be an analytic quasi-order satisfying the $(\varepsilon_0)$-property. If there is a cofinal monotone quasi-open map $f : \cal{R} \to \S$, then $(\S,\leq_\S)$ has the analytic Ramsey property.
\end{proposition}

\begin{proof}
    Let $\X \subseteq \S$ be analytic, and let $X \in \S$. Since $f$ is cofinal, there exists some $A \in \cal{R}$ such that $f(A) \leq X$. By Lemma \ref{lem:monotone.automatic.continuity}, there exists some $B \leq A$ such that $f\restrictedto[B]$ is continuous. Then $(f\restrictedto[B])^{-1}(\X)$ is an analytic subset of $[B]$, so by \textbf{A4} there exists some $C \leq B$ such that $[C] \subseteq (f\restrictedto[B])^{-1}(\X)$ or $[C] \subseteq (f\restrictedto[B])^{-1}(\X^c)$. This implies that $(f\restrictedto[B])``[C] \subseteq \X$ or $(f\restrictedto[B])``[C] \subseteq \X^c$. By the quasi-open property of $f$, there exists some $Y \in \S$ such that $[Y] \subseteq (f\restrictedto[B])``[C] = f``[C]$. Then $Y \leq f(C) \leq f(A) \leq X$ and $[Y] \subseteq \X$ or $[Y] \subseteq \X^c$, as desired.
\end{proof}

\subsection{Examples of projections}
\label{sec:projection.examples}
We provide an extensive list of projections between topological Ramsey spaces. In fact, every example in this subsection satisfies the following stronger properties:

\begin{lemma}
\label{lem:monotone.ab.stronger.properties}
    Let $(\P,\leq_\P)$ and $(\Q,\leq_\Q)$ be two quasi-orders, with $\leq$-maximal elements $\one_\P$ and $\one_\Q$ respectively. Let $f : \P \to \Q$ be a monotone function, and suppose that:
    \begin{itemize}
        \item If $\one_\P$ is $\leq$-maximal in $\P$, then $f(\one_\P)$ is $\leq$-maximal in $\Q$.
        
        \item For all $p \in \P$, $f``(p^\downarrow) = f(p)^\downarrow$.
    \end{itemize}
    Then $f$ is a quasi-open projection (so it is Alaoglu--Birkhoff convergent).
\end{lemma}

\begin{proof}
    $f$ is quasi-open, as $f(p)^\downarrow \subseteq f``(p^\downarrow)$ for all $p \in \P$. If $q \leq f(p)$, then $q \in f(p)^\downarrow = f``(p^\downarrow)$, so there exists some $p' \leq p$ such that $f(p') = q$. Applying this argument to any $\leq$-maximal element $f(\one_\P)$ shows that $f$ is cofinal. Thus, $f$ is a projection. 
\end{proof}

Note that if $f$ is monotone, then the containment $f``(p^\downarrow) \subseteq f(p)^\downarrow$ follows immediately. Thus, it suffices to just check that $f(p)^\downarrow \subseteq f``(p^\downarrow)$, i.e. that $f$ is a pre-adjunction (Definition \ref{def:pre-adjunction.quasi-order}).

\begin{example}[min and max of $\FIN$]
\label{ex:min.max.FIN}
    For each $x \in \FIN$, let $\min(x)$ and $\max(x)$ denote the smallest and largest elements of $x$ respectively. We can extend this map to $\min,\max : \FIN^{[\infty]} \to [\N]^\infty$ by:
    \begin{align*}
        \min((x_n)_{n<\omega}) &:= \{\min(x_n) : n<\omega\}, \\
        \max((x_n)_{n<\omega}) &:= \{\max(x_n) : n<\omega\}.
    \end{align*}
    Then these maps satisfy the properties of Lemma \ref{lem:monotone.ab.stronger.properties}.
\end{example}

\begin{proof}
    Note that $(\{n\})_{n<\omega}$ is the $\leq$-maximal element in $\FIN^{[\infty]}$, and:
    \begin{align*}
        \min((\{n\})_{n<\omega}) = \max((\{n\})_{n<\omega}) = \N.
    \end{align*}
    
    Observe that both maps are monotone --- if $A = (x_n)_{n<\omega} \in \FIN^{[\infty]}$, $B = (y_n)_{n<\omega} \in \FIN^{[\infty]}$ and $B \leq A$, then for each $n$, $y_n = x_{n_0} \cup \cdots \cup x_{n_k}$ for some $n_0 < \cdots < n_k$. This implies that $\min(y_n) = \min(x_{n_0})$ and $\max(y_n) = \max(x_{n_k})$. Therefore, $\min(B) \subseteq \min(A)$ and $\max(B) \subseteq \max(A)$.

    Now suppose that $\{k_n : n < \omega\}$ is the increasing enumeration of $\min(A)$, and $X \in [\min(A)]$. Then $X = \{k_{n_i} : i < \omega\}$ where $(n_i)_{i<\omega}$ is increasing, so we have that $\min((x_{n_i})_{i<\omega}) = X$. Therefore, $\min``[A] = [\min(A)]$. A similar argument can be used to show that $\max``[A] = [\max(A)]$.
\end{proof}

\begin{example}[Tetris operation of $\FIN_k$]
\label{ex:tetris.operation}
    Recall the \emph{tetris operation} from Example \ref{ex:Gowers.FIN_k}, defined by $T(x)(n) := \max\{x(n)-1,0\}$ for all $x \in \FIN_k$, so that $T(x) \in \FIN_{k-1}$. We can extend this map to $T : \FIN_k^{[\infty]} \to \FIN_{k-1}^{[\infty]}$ by:
    \begin{align*}
        T((x_n)_{n<\omega}) := (T(x_n))_{n<\omega}.
    \end{align*}
    Then this map satisfies the properties of Lemma \ref{lem:monotone.ab.stronger.properties}.
\end{example}

\begin{proof}
    Define $e_n^k \in \FIN_k^{[\infty]}$ so that $\supp(e_n^k) = \{n\}$ and $e_n^k(n) = k$. Then $(e_n^k)_{n<\omega}$ is the $\leq$-maximal element in $\FIN_k^{[\infty]}$, and $T(e_n^k) = e_n^{k-1}$. Therefore, $T$ preserves the $\leq$-maximal element.

    We show that $T$ is monotone. Suppose that $A = (x_n)_{n<\omega} \in \FIN_k^{[\infty]}$, $B = (y_n)_{n<\omega} \in \FIN_k^{[\infty]}$ and $B \leq A$. Then for each $n$, $y_n \in \c{A}$, so $y_n = T^{\lambda_0}(x_{n_0}) + \cdots + T^{\lambda_j}(x_{n_j})$ for some $n_0 < \cdots < n_j$, and $\lambda_i = 0$ for some $i$. Then:
    \begin{align*}
        T(y_n) &= T(T^{\lambda_0}(x_{n_0}) + \cdots + T^{\lambda_j}(x_{n_j})) \\
        &= T^{\lambda_0}(T(x_{n_0})) + \cdots + T^{\lambda_j}(T(x_{n_j}))
    \end{align*}
    and $\lambda_i = 0$, so $T(y_n) \in \c{T(A)}$. Therefore, $T(B) \leq T(A)$.

    Now suppose that $A = (x_n)_{n<\omega}$ and $X = (z_n)_{n<\omega} \leq T(A)$ in $\FIN_{k-1}^{[\infty]}$. Then for each $n$, we have that:
    \begin{align*}
        z_n &= T^{\lambda_{n,0}}(T(x_{j_{n,0}})) + \cdots + T^{\lambda_{n,N_n}}(T(x_{j_{n,N_n}})) \\
        &= T\bb{T^{\lambda_{n,0}}(x_{j_{n,0}}) + \cdots + T^{\lambda_{n,N_n}}(x_{j_{n,N_n}})}
    \end{align*}
    for some $j_{n,0} < \cdots < j_{n,N_n}$. Furthermore, $j_{n,N_n} < j_{n+1,0}$ for all $n$. Thus, we may let:
    \begin{align*}
        B := (T^{\lambda_{n,0}}(x_{j_{n,0}}) + \cdots + T^{\lambda_{n,N_n}}(x_{j_{n,N_n}}))_{n<\omega}.
    \end{align*}
    Then $B \in \FIN_k^{[\infty]}$ (as for all $n$, $\lambda_{n,i} = 0$ for some $i \leq N_n$), $B \leq A$ and $T(B) = X$. Thus, $T``[A] = [T(A)]$. 
\end{proof}

\begin{example}[Tetris operation of $\FIN_{\pm k}$]
\label{ex:tetris.operation.FIN_pmk}
    A variant of the tetris operation for $\FIN_{\pm k}$ was also introduced in Example \ref{ex:Gowers.FIN_pmk}, where $T : \FIN_{\pm k} \to \FIN_{\pm(k - 1)}$ is defined by:
    \begin{align*}
        T(x)(n) := 
        \begin{cases}
            x(n) - 1, &\text{if $x(n) > 0$}, \\
            x(n) + 1, &\text{if $x(n) < 0$}, \\
            0, &\text{otherwise}.
        \end{cases}
    \end{align*}
    We can similarly extend the map to $T : \FIN_{\pm k}^{[\infty]} \to \FIN_{\pm(k-1)}^{[\infty]}$ by:
    \begin{align*}
        T((x_n)_{n<\omega}) := (T(x_n))_{n<\omega}.
    \end{align*}
    Then $T$ is monotone and $T``[A] = [T(A)]$ for all $A \in \FIN_{\pm k}^{[\infty]}$, which can be proved similarly to Example \ref{ex:tetris.operation}. 
\end{example}

\begin{example}[Absolute value map of $\FIN_{\pm k}$]
\label{ex:absolute.value.FIN_pmk}
    We define a map $v : \FIN_{\pm k} \to \FIN_k$ by:
    \begin{align*}
        v(x)(n) := |x(n)|.
    \end{align*}
    We can extend the map to $v : \FIN_{\pm k}^{[\infty]} \to \FIN_k^{[\infty]}$ by:
    \begin{align*}
        v((x_n)_{n<\omega}) := (v(x_n))_{n<\omega}.
    \end{align*}
    Then this map satisfies the properties of Lemma \ref{lem:monotone.ab.stronger.properties}.
\end{example}

\begin{proof}
    We defined, in the proof of Example \ref{ex:tetris.operation}, the function $e_n^k \in \FIN_{\pm k}$. Note that $e_n^k,-e_n^k \in \FIN_{\pm k}$, and $v(e_n^k) = v(-e_n^k) = e_n^k$. We make the following observations:
    \begin{enumerate}
        \item If $x,y \in \FIN_{\pm k}$ and $x < y$, then $v(x + y) = v(x) + v(y)$.
        
        \item For all $x \in \FIN_{\pm k}$, $v(T(x)) = T(v(x))$.
    \end{enumerate}
    Every $\leq$-maximal element in $\FIN_k^{[\infty]}$ is of the form $(\varepsilon_ne_n^k)_{n<\omega}$ for some sequence $(\varepsilon_n)_{n<\omega} \in \{-1,1\}^\omega$, and $v((\varepsilon_ne_n^k)_{n<\omega}) = (e_n^k)_{n<\omega}$, the $\leq$-maximal element in $\FIN_k^{[\infty]}$.

    To see that $v$ is monotone, let $A = (x_n)_{n<\omega},B = (y_n)_{n<\omega} \in \FIN_{\pm k}^{[\infty]}$, and suppose that $B \leq A$. Then for all $n$, $y_n = T^{\lambda_0}(x_{n_0}) + \cdots + T^{\lambda_j}(x_{n_j})$ for some $n_0 < \cdots < n_j$, and $\lambda_i = 0$ for some $i$. By the observations above, we have that:
    \begin{align*}
        v(y_n) &= v\bb{T^{\lambda_0}(x_{n_0}) + \cdots + T^{\lambda_j}(x_{n_j})} \\
        &= T^{\lambda_0}(v(x_{n_0})) + \cdots + T^{\lambda_j}(v(x_{n_j})),
    \end{align*}
    so $v(y_n) \in \c{v(A)}$ for all $n$. Therefore, $v(B) \leq v(A)$.

    Now suppose that $A = (x_n)_{n<\omega}$ and $X = (z_n)_{n<\omega} \leq v(A)$ in $\FIN_k^{[\infty]}$. Observe that $\FIN_k \subseteq \FIN_{\pm k}$ and $\FIN_k^{[\infty]} \subseteq \FIN_{\pm k}^{[\infty]}$. Since $v(A) \leq A$ (as $v(x_n) \in \c{A}$ for all $n$), we have that $X \in [v(A)] \subseteq [A]$ in $\FIN_{\pm k}^{[\infty]}$ and $v(X) = X$ (as $\ran(z_n) \subseteq \N$ for all $n$). Thus, $v``[A] = [v(A)]$.
\end{proof}

\begin{example}[Supports of vectors]
\label{ex:support.vectors}
    Recall from Example \ref{ex:vector.spaces} that if $E$ is a countable-dimensional vector space over a countable field $\mathbb{F}$, with a Hamel basis $(e_n)_{n<\omega}$, then we defined the support of a vector $x = \sum_{n<\omega} \lambda_ne_n \in E$ by $\supp(x) := \{n < \omega : \lambda_n \neq 0\}$. This gives us a map $\supp : E \setminus \{0\} \to \FIN$, which induces a map $\supp : E^{[\infty]} \to \FIN^{[\infty]}$ by:
    \begin{align*}
        \supp((x_n)_{n<\omega}) = (\supp(x_n))_{n<\omega}.
    \end{align*}
    Then this map satisfies the properties of Lemma \ref{lem:monotone.ab.stronger.properties}.
\end{example}

\begin{proof}
    Note that $(e_n)_{n<\omega}$ is the $\leq$-maximal element in $E^{[\infty]}$, and:
    \begin{align*}
        \supp((e_n)_{n<\omega}) = (\{n\})_{n<\omega},
    \end{align*}
    which is the $\leq$-maximal element in $\FIN^{[\infty]}$.

    To see that $\supp$ is monotone --- if $A = (x_n)_{n<\omega} \in E^{[\infty]}$, $B = (y_n)_{n<\omega} \in E^{[\infty]}$ and $B \leq A$, then for each $n$, $y_n = \lambda_{n_0}x_{n_0} + \cdots + \lambda_{n_k}x_{n_k}$ for some $n_0 < \cdots < n_k$ and $\lambda_{n_i} \neq 0$ for all $i$. Then:
    \begin{align*}
        \supp(y_n) &= \supp(\lambda_{n_0}x_{n_0} + \cdots + \lambda_{n_k}x_{n_k}) \\
        &= \supp(\lambda_{n_0}x_{n_0}) \cup \cdots \cup \supp(\lambda_{n_k}x_{n_k}) \\
        &= \supp(x_{n_0}) \cup \cdots \cup \supp(x_{n_k}),
    \end{align*}
    so $\supp(B) \leq \supp(A)$.

    Now suppose that $A = (x_n)_{n<\omega} \in E^{[\infty]}$ and $X = (y_n)_{n<\omega} \leq \supp(A)$. Then for each $n$, $y_n = \supp(x_{j_{n,0}}) \cup \cdots \cup \supp(x_{j_{n,N_n}})$ for some $j_{n,0} < \cdots < j_{n,N_n} < j_{n+1,0}$. Then let:
    \begin{align*}
        B := (x_{j_{n,0}} + \cdots + x_{j_{n,N_n}})_{n<\omega},
    \end{align*}
    and we have that $B \in E^{[\infty]}$, $B \leq A$ and $\supp(B) = X$. Thus, $\supp``[A] = [\supp(A)]$.
\end{proof}

\begin{example}[Supports of variable-words]
\label{ex:support.variable.words}
    Let $A = (x_n)_{n<\omega} \in W_{Lv}^{[\infty]}$ be a sequence of variable-words in the Hales--Jewett space (Example \ref{ex:Hales-Jewett}). Note that for all $B \in W_{Lv}^{[\infty]}$, $B \leq A$ iff $\c{B} \subseteq \c{A}$. For each $x \in \c{A}$, we may let $\supp_A(x) := \{n_0,\dots,n_k\}$, where $x = x_{n_0}[\lambda_0]^\frown\cdots^\frown x_{n_k}[\lambda_k]$ for some $n_0 < \cdots < n_k$. Since $A$ is rapidly increasing, $\supp_A$ is well-defined, and if $B = (y_n)_{n<\omega} \leq A$ then $\supp_A(y_n) < \supp_A(y_{n+1})$ for all $n$. Therefore, for every $A \in W_{Lv}^{[\infty]}$, we have a map $\supp_A : [A] \to \FIN^{[\infty]}$ defined by:
    \begin{align*}
        \supp_A((y_n)_{n<\omega}) := (\supp_A(y_n))_{n<\omega}.
    \end{align*}
    $\supp_A$ shares many properties with $\supp$ in Example \ref{ex:support.vectors}, so one may show similarly that $\supp_A$ satisfies the properties of Lemma \ref{lem:monotone.ab.stronger.properties}.
\end{example}

\begin{example}[Located variable-words to variable-words]
\label{ex:unlocate.variable-words}
    Recall that $\FIN_{Lv}^{[\infty]}$ is the set of infinite block sequences of located variable-words (Example \ref{ex:block.sequences.located.variable-words}). We say that $A = (x_n)_{n<\omega} \in \FIN_{Lv}^{[\infty]}$ is \emph{rapidly increasing} if for all $n$, $\sum_{i<n} |\dom(x_i)| < |\dom(x_n)|$. Readers may verify that the set of rapidly increasing sequences is dense open in $\FIN_{Lv}^{[\infty]}$.

    In \cite[Theorem 4.50]{T10}, a map $\zeta : \FIN_L \cup \FIN_{Lv} \to W_L \cup W_{Lv}$ (originally denoted as $\pi$) is defined as follows: Let $x \in \FIN_L \cup \FIN_{Lv}$ with $\dom(x) = \{n_i : i < l\} \subseteq \N$, where $n_0 < \cdots < n_{l-1}$. We define $\zeta(x)$ to be the word of length $l$ such that the $i^\text{th}$ letter of $\zeta(x)$ is $x(n_i)$. Note that $x \in \FIN_L$ iff $\zeta(x) \in W_L$, and $x \in \FIN_{Lv}$ iff $\zeta(x) \in W_{Lv}$. We fix some rapidly increasing $A \in \FIN_{Lv}^{[\infty]}$, and extend the map to $\zeta$ to the domain $[A]$ by:
    \begin{align*}
        \zeta((x_n)_{n<\omega}) := (\zeta(x_n))_{n<\omega}.
    \end{align*}
    Then $\zeta : [A] \to [\zeta(A)]$ satisfies the properties of Lemma \ref{lem:monotone.ab.stronger.properties}. We remark that, in general, we may define $\zeta$ on the dense open set of rapidly increasing sequences in $\FIN_{Lv}^{[\infty]}$, and $\zeta$ would be a cofinal map such that $\zeta``[B] = [\zeta(B)]$ for all rapidly increasing $B \in \FIN_{Lv}^{[\infty]}$.
\end{example}

\begin{proof}
    It is easy to verify from the definitions that for all $x,y \in \FIN_{Lv}$ with $x < y$, and $\lambda,\mu \in L \cup \{v\}$:
    \begin{align*}
        \zeta(x[\lambda] \cup y[\mu]) = \zeta(x)[\lambda]^\frown\zeta(y)[\mu]. \tag{$\star$}
    \end{align*}
    That $\zeta : [A] \to [\zeta(A)]$ preserves the $\leq$-maximal element follows from the choice of the codomain of $\zeta$. We show that $\zeta$ is monotone. Let $B = (y_n)_{n<\omega},C = (z_n)_{n<\omega} \in [A]$, and suppose that $C \leq B$. Then for all $n$, $z_n = y_{n_0}[\lambda_{n_0}] \cup \cdots \cup y_{n_j}[\lambda_{n_j}]$ for some $n_0 < \cdots < n_j$, and $\lambda_{n_i} \in L_{n_i} \cup \{v\}$ for all $i$, with $\lambda_{n_i} = v$ for some $i$. Applying the identity $(\star)$, we have that for all $n$:
    \begin{align*}
        \zeta(z_n) &= \zeta(y_{n_0}[\lambda_{n_0}] \cup \cdots \cup y_{n_j}[\lambda_{n_j}]) \\
        &= \zeta(y_{n_0})[\lambda_{n_0}]^\frown\cdots^\frown\zeta(y_{n_j})[\lambda_{n_j}],
    \end{align*}
    so $\zeta(C) \leq \zeta(B)$.

    Now suppose that $B = (y_n)_{n<\omega} \in [A]$, and $X = (w_n)_{n<\omega} \in [\zeta(B)]$. For each $n$, $w_n = \zeta(y_{j_{n,0}})[\lambda_{j_{n,0}}]^\frown\cdots^\frown\zeta(y_{j_{n,N_n}})[\lambda_{j_{n,N_n}}]$ for some $j_{n,0} < \cdots < j_{n,N_n} < j_{n+1,0}$. Then let:
    \begin{align*}
        C := (y_{j_{n,0}}[\lambda_{j_{n,0}}] \cup \cdots \cup y_{j_{n,N_n}}[\lambda_{j_{n,N_n}}])_{n<\omega},
    \end{align*}
    and we have that $C \in \FIN_{Lv}^{[\infty]}$, $C \leq B$, and $\zeta(C) = X$ by $(\star)$. Thus, $\zeta``[B] = [\zeta(B)]$. 
\end{proof}

\begin{example}[Levels of strong subtrees]
\label{ex:levels.strong.subtrees}
    Recall that $\S_\infty(U)$ is the set of strong subtrees of $U$ (Example \ref{ex:strong.subtrees}). Let $\{m_n : n < \omega\}$ be the increasing enumeration of the set of all $m$'s such that $\dom(s) = m$ for some splitting node $s \in U$. We then define a map $S : \S_\infty(U) \to [\N]^\infty$ by:
    \begin{align*}
        S(A) := \{n \in \N : \dom(s) = m_n \text{ for some splitting node } s \in A\}.
    \end{align*}
    In other words, $S(A)$ is the set of all $n$'s such that $A$ splits at level $m_n$. Then $S$ satisfies the properties of Lemma \ref{lem:monotone.ab.stronger.properties}.
\end{example}

\begin{proof}
    Note that $U$ is the $\subseteq$-maximal element in $\S_\infty(U)$. For notational convenience, we may assume that $m_n = n$ for all $n$, which implies that $S(U) = \N$. To see that $S$ is monotone, suppose that $B \leq A$ in $\S_\infty(U)$. Then every splitting node in $B$ is also a splitting node in $A$, so $S(B) \subseteq S(A)$. It remains to show that $S``[A] = [S(A)]$ for all $A \in \S_\infty(U)$ (which also implies that $S$ is surjective, as $S(A) = \N$ for some $A \in \S_\infty(U)$).

    Let $X = \{n_k : k <\omega\} \in [S(A)]$. We shall construct a binary subtree of $A$ that splits at exactly the levels in $X$. We define $s_t$ for $t \in 2^{<\omega}$ inductively as follows: Let $s_\emptyset \in A$ be any splitting node such that $\dom(s_\emptyset) = n_0$ (which exists, as $n_0 \in S(A)$). If $s_t$ has been defined and is a splitting node in $A$, then let $s_{t^\frown0},s_{t^\frown1} \in A$ be two distinct splitting nodes in $A$ extending $s_t$, such that $\dom(s_{t^\frown 0}) = \dom(s_{t^\frown 1}) = n_{|t|+1}$. Let:
    \begin{align*}
        B := \{u \in A : u \sqsubseteq s_t \text{ for some } t \in 2^{<\omega}\}.
    \end{align*}
    Observe that $B \in \S_\infty(U)$, as $s_t$ is a splitting node in $B$ for all $t \in 2^{<\omega}$, and if $u \in B$ and $\dom(u) = n_k$ then $u = s_t$ for some $t \in 2^{<\omega}$ such that $|t| = k$. We have that $B \leq A$ and $S(B) = X$, as desired. 
\end{proof}

\begin{example}[Minimal representatives of equivalence classes]
\label{ex:minimal.representatives}
    Recall that for each $A \in \E_\infty$ (the Carlson--Simpson space, see Example \ref{ex:Carlson-Simpson}), $p(A) = \{p_n(A) : n < \omega\}$ is the set of minimal representatives of equivalence classes of $A$. We define $\N^+ := \N \setminus \{0\}$, and $p^* : \E_\infty \to [\N^+]^\infty$ by $p^*(A) := p(A) \setminus \{0\}$. By identifying $[\N^+]^\infty$ with $[\N]^\infty$ via a bijection between $\N$ and $\N^+$, this map satisfies the properties of Lemma \ref{lem:monotone.ab.stronger.properties}. 
\end{example}

\begin{proof}
    The $\leq$-maximal element in $\E_\infty$ is the equality relation $E$ (i.e. $(m,n) \in E$ iff $m = n$). Then $p_n(E) = n$ for all $n$, so $p^*(E) = \N^+$. Thus, $p^*$ preserves $\leq$-maximal elements. 

    Let $A,B \in \E_\infty$, and suppose that $B \leq A$. This means that every equivalence class of $B$ is a union of equivalence classes of $A$. Thus, if $k = p_n(B) = \min(B(n))$ ($B(n)$ is the equivalence class of $B$ with the $n^\text{th}$-least minimal representative) and $B(n) = \bigcup_{i \in I} A(n_i)$ for some $I \subseteq \omega$, then $k$ is minimal in $\bigcup_{i \in I} A(n_i)$. Assuming that $n_0$ is minimal among the $n_i$'s, we have that $k = \min(A(n_0))$, so $k \in p(A)$. Therefore, $p(B) \subseteq p(A)$, which implies that $p^*(B) \subseteq p^*(A)$.

    Now suppose that $p^*(A) = \{k_n : 1 \leq n < \omega\}$ (so $k_n = p_n(A)$), and let $X = \{k_{n_i} : 1 \leq i < \omega\} \in [p^*(A)]$. Let $B$ be the unique equivalence relation on $\N$ such that for all $i \geq 1$, $\bigcup_{m<n_1} A(m)$ and $\bigcup_{n_i \leq m < n_{i+1}} A(m)$ are equivalence classes of $B$. Observe that $B(i) = \bigcup_{n_i \leq m < n_{i+1}} A(m)$, so for all $i \geq 1$:
    \begin{align*}
        \min(B(i)) = \min\bb{\bigcup_{n_i \leq m < n_{i+1}} A(m)} = \min(A(n_i)) = k_{n_i}.
    \end{align*}
    Therefore, $B \leq A$ and $p^*(B) = X$, as desired. 
\end{proof}

\begin{example}[Merging letters of Carlson--Simpson spaces]
\label{ex:CS-variables.merging.alphabets}
    Let $L$ be a finite non-empty alphabet, and consider the Carlson--Simpson space with this alphabet (Example \ref{ex:CS-variables}). Let $\cal{L}$ be a partition of $L$, which we view as a finite alphabet. We define $\sigma_{\cal{L}} : \E_\infty^L \to \E_\infty^{\cal{L}}$ by stipulating that for all $A \in \E_\infty^L$, $\sigma_{\cal{L}}(A)$ is the equivalence relation on $\N \cup \cal{L}$ such that for all $m,n \in \N$ and $u \in \cal{L}$:
    \begin{small}
        \begin{align*}
            (u,n) \in \sigma_{\cal{L}}(A) &\iff \exists w \in u \, ((w,n) \in A), \\
            (m,n) \in \sigma_{\cal{L}}(A) &\iff (m,n) \in A \vee \exists u \in \cal{L} \, ((u,m) \in \sigma_{\cal{L}}(A) \wedge (u,n) \in \sigma_{\cal{L}}(A)).
        \end{align*}
    \end{small}
    In other words, $\sigma_{\cal{L}}(A)$ is obtained from $A$ by merging, for each $u \in \cal{L}$, the equivalence classes of all $w \in u$ into one, which is then identified with $u$. Then $\sigma_{\cal{L}}$ satisfies the properties of Lemma \ref{lem:monotone.ab.stronger.properties}. From this map, we can also define $\sigma : \E_\infty^L \to \E_\infty$ by $\sigma(A) := \sigma_{\{L\}}(A)\restrictedto\N$, and it's easy to verify that $\sigma$ also satisfies the properties of Lemma \ref{lem:monotone.ab.stronger.properties}.
\end{example}

\begin{proof}
    It is clear that $\sigma_{\cal{L}}$ is monotone, and that it sends the equality relation on $\N \cup L$ to the equality relation on $\N \cup \cal{L}$. Let $A \in \E_\infty^L$, and fix a function $\tau : \cal{L} \to L$ such that $\tau(u) \in u$ for all $u \in \cal{L}$. If $X \leq \sigma_{\cal{L}}(A)$, then, writing $[x]_Y$ for the equivalence class of $x$ in an equivalence relation $Y$, we let $B$ be the equivalence relation on $\N \cup L$ whose equivalence classes are:
    \begin{itemize}
        \item $[\tau(u)]_A \cup ([u]_X \setminus [u]_{\sigma_{\cal{L}}(A)})$ for each $u \in \cal{L}$;

        \item $[w]_A$ for each $w \in L \setminus \ran(\tau)$;

        \item the equivalence classes of $X$ disjoint from $\cal{L}$.
    \end{itemize}
    Then $B \in \E_\infty^L$, $B \leq A$ and $\sigma_{\cal{L}}(B) = X$.
\end{proof}

\begin{example}[Minimal representatives of equivalence classes disjoint from $L$]
\label{ex:minimal.representatives.not.L}
    It follows from Examples \ref{ex:minimal.representatives} and \ref{ex:CS-variables.merging.alphabets} that there is a map $\E_\infty^L \to [\N]^\infty$ that satisfies the properties of Lemma \ref{lem:monotone.ab.stronger.properties}. There is, however, a simpler example: recall from Example \ref{ex:CS-variables} that given some $A \in \E_\infty^L$, 
    $p(A) = \{p_n(A) : n <\omega\}$ is the set of minimal representatives of classes disjoint from $L$. Unlike the case for Example \ref{ex:minimal.representatives}, $p_0(A)$ need not always be $0$, so we can prove, in a way similar to that of Example \ref{ex:minimal.representatives}, that $p$ satisfies the properties of Lemma \ref{lem:monotone.ab.stronger.properties}.
\end{example}

\begin{example}[Projections between $\cal{R}_\alpha$'s]
\label{ex:Dobrinen.spaces.projections}
    Let $\alpha < \omega$, and recall the definition of $(\cal{R}_\alpha,\leq_\alpha,r^\alpha)$ from Example \ref{ex:Dobrinen.space}. We define a map $\pi_\alpha : \cal{R}_{\alpha+1} \to \cal{R}_\alpha$ by stipulating that for all $A \in \cal{R}_{\alpha+1}$:
    \begin{align*}
        \pi_\alpha(A) := \{\emptyset\} \cup \{(s(i) - 1)_{i<\dom(s)} : s \in A \setminus A(0) \text{ and } \dom(s) \leq \alpha + 1\}.
    \end{align*}
    The projection $\pi_\alpha$ may be visualised as follows: given a tree $A \in \cal{R}_{\alpha+1}$ with its root at the bottom (i.e. positioned in the same way as \cite[Figure 5]{DT15}), $\pi_\alpha$ removes the topmost layer and the leftmost branch of $A$, followed by a relabelling of the remaining branches (by subtraction of $1$). Then $\pi_\alpha$ satisfies the properties of Lemma \ref{lem:monotone.ab.stronger.properties}. Consequently, by composing these projections, we have that $(\cal{R}_\alpha,\leq_\alpha,r^\alpha)$ projects onto $(\cal{R}_\beta,\leq_\beta,r^\beta)$ for $\beta < \alpha < \omega$.
\end{example}

\begin{proof}
    For each $s \in \omega^{<\omega}$, we write $s^- := (s(i) - 1)_{i<\dom(s)}$, and for each tree $a \subseteq \omega^{<\omega}$, we write $a^- := \{s^- : s \in a \text{ and } \dom(s) \leq \alpha + 1\}$. Under this notation, we have that for all $A \in \cal{R}_{\alpha+1}$:
    \begin{align*}
        \pi_\alpha(A) = \bigcup_{n<\omega} A(n+1)^-.
    \end{align*}
    Readers may verify from the definitions that $l_{n+1}^{\beta+1} = l_n^\beta + 1$ for all $\beta,n < \omega$, so it follows by induction on $\alpha$ that $\mathbb{A}_{\alpha+1}(n+1)^- = \mathbb{A}_\alpha(n)$ for all $n$. In particular, $\pi_\alpha$ sends the $\leq_{\alpha+1}$-maximal element $\mathbb{A}_{\alpha+1}$ to the $\leq_\alpha$-maximal element $\mathbb{A}_\alpha$. Since the map $s \mapsto s^-$ is $\sqsubseteq$-preserving and $\prec$-preserving in both directions, if $a \subseteq \mathbb{A}_{\alpha+1}(k+1)$ is a tree isomorphic to $\mathbb{A}_{\alpha+1}(n+1)$, then $a^- \subseteq \mathbb{A}_\alpha(k)$ is a tree isomorphic to $\mathbb{A}_\alpha(n)$, so $\pi_\alpha$ is well-defined. To see that $\pi_\alpha$ is monotone --- from the definition above, we have that if $b \subseteq a \subseteq \omega^{<\omega}$, then $b^- \subseteq a^-$. Therefore, if $A,B \in \cal{R}_{\alpha+1}$ and $B \leq_{\alpha+1} A$, then for each $m$, $B(m+1) \subseteq A(n+1)$ for some $n$, so $\pi_\alpha(B) \subseteq \pi_\alpha(A)$. It remains to show that $\pi_\alpha``[A] = [\pi_\alpha(A)]$ for all $A \in \cal{R}_{\alpha+1}$. 
    
    Let $A \in \cal{R}_{\alpha+1}$, and let $X \in [\pi_\alpha(A)]$. Let $b := \{s \in A : s^- \in X\}$, and let $\{s_j : j < \omega\}$ be the $\prec$-increasing enumeration of the set of all $s \in A(0) \cup b$ such that $\dom(s) = \alpha + 1$. Observe that each $s_j$ has at least $j+1$ successors in $A$ --- indeed, for all $k$, the $k^\text{th}$ $\prec$-least node of length $\alpha+1$ in $\mathbb{A}_{\alpha+1}$ has exactly $l_{k+1}^0 - l_k^0 = k+1$ successors, so the same holds for $A$, and $s_j$ is the $k^\text{th}$ $\prec$-least node of length $\alpha+1$ in $A$ for some $k \geq j$. Thus, we may let $c_j$ be the set of the $j+1$ $\prec$-least successors of $s_j$ in $A$. Let:
    \begin{align*}
        B := A(0) \cup b \cup \bigcup_{j<\omega} c_j.
    \end{align*}
    Note that $B \in \cal{R}_{\alpha+1}$ --- indeed, $b^- = X \in \cal{R}_\alpha$, and each $s_j$ is the $j^\text{th}$ $\prec$-least node of length $\alpha+1$ in $B$, with exactly $j+1$ successors in $B$. Then $B \leq_{\alpha+1} A$ (as $B \subseteq A$) and $\pi_\alpha(B) = X$, as desired.
\end{proof}

We shall revisit some of these examples in \S\ref{sec:strong.projections} as \emph{strong projections}.

\subsubsection{Biprojections}
We stated in Example \ref{ex:CS-variables} that the Carlson--Simpson space with variables has three equivalent formulations. The notion of \emph{biprojection} formalises this claim. 

\begin{definition}
\label{def:biprojection}
    Let $(\P,\leq_\P)$ and $(\Q,\leq_\Q)$ be two quasi-orders, with $\leq$-maximal elements $\one_\P$ and $\one_\Q$ respectively. A \emph{biprojection} is a bijective function $f : \P \to \Q$ such that:
    \begin{itemize}
        \item $f(\one_\P) = \one_\Q$.
        
        \item For all $p \in \P$, $f``(p^\downarrow) = f(p)^\downarrow$.
    \end{itemize}
\end{definition}

Note that if $f$ is a biprojection, then for all $q \in \Q$, $(f^{-1})``(q^\downarrow) = f^{-1}(q)^\downarrow$.

\begin{example}[Variable-word reformulation of Carlson--Simpson spaces with variables]
\label{ex:CS-variables.reformulation.variable-words}
    Let $L$ be a finite non-empty alphabet, and consider the first two formulations of the Carlson--Simpson space with variables $\E_\infty^L$ and $W_\infty^L(\omega)$. Define a map $\Phi : \E_\infty^L \to W_\infty^L(\omega)$ such that for all $A \in \E_\infty^L$, $\Phi(A)$ is the unique $\omega$-variable-word such that for all $n$:
    \begin{align*}
        \Phi(A)(n) = 
        \begin{cases}
            w, &\text{if $w \in L$ and $(w,n) \in A$}, \\
            v_i, &\text{if $(p_i(A),n) \in A$}.
        \end{cases}
    \end{align*}
    Then $\Phi$ is a biprojection.
\end{example}

\begin{proof}
    \underline{$\Phi$ is well-defined:} Fix some $A \in \E_\infty^L$. It is clear that $\Phi(A)$ is an $\omega$-variable-word, but we need to show that if $i < j$, then the first occurrence of $v_i$ comes before the first occurrence of $v_j$. Let $n$ be the least integer such that $\Phi(A)(n) = v_i$. Then $n$ is the least integer such that $(p_i(A),n) \in A$. But $p_i(A)$ is the least representative of its own equivalence class, so $n = p_i(A)$. Therefore, $v_i$ first occurs at position $p_i(A)$ for every $i$, and since $i < j$, $p_i(A) < p_j(A)$.

    \underline{$\Phi$ is monotone:} Let $A,B \in \E_\infty^L$ such that $B \leq A$. It suffices to show that the following hold:
    \begin{itemize}
        \item For all $n \in \N$, if $\Phi(A)(n) = w \in L$, then $\Phi(B)(n) = w$.

        \item For all $m,n \in \N$, if $\Phi(A)(m) = \Phi(A)(n)$ then $\Phi(B)(m) = \Phi(B)(n)$.
    \end{itemize}
    If $\Phi(A)(n) = w \in L$, then $(w,n) \in A$. Since $B$ coarsens $A$, we have that $(w,n) \in B$ as well, so $\Phi(B)(n) = w$. If $\Phi(A)(m) = \Phi(A)(n)$, then $m$ and $n$ lie in the same equivalence class of $A$, i.e. $(m,n) \in A$. Since $B$ coarsens $A$, we have that $(m,n) \in B$, so $\Phi(B)(m) = \Phi(B)(n)$.

    \underline{$\Phi$ is bijective and $\Phi^{-1}$ is monotone:} For each $X \in W_\infty^L(\omega)$, let $\Psi(X)$ be the unique equivalence relation consisting of equivalence classes of the form $\{n \in \N : X(n) = v_i\}$ for some $i$, and $\{w\} \cup \{n \in \N : X(n) = w\}$ for some $w \in L$. The map $\Psi : W_\infty^L(\omega) \to \E_\infty^L$ is well-defined, as every $v_i$ occurs at least once in $X$, so $\Psi(X)$ has infinitely many equivalence classes, and $(w,w') \notin \Psi(X)$ for all $w \neq w'$ in $L$. We shall show that $\Psi$ is monotone, and $\Psi = \Phi^{-1}$.
    
    If $Y \leq X$ in $W_\infty^L(\omega)$, then for all $m,n \in \N$, $X(m) = X(n)$ implies that $Y(m) = Y(n)$. If $(w,n) \in \Psi(X)$, then $X(n) = w$, which implies that $Y(n) = w$ as well. Therefore, $(w,n) \in \Psi(Y)$, so $\Psi(Y)$ coarsens $\Psi(X)$.
    
    It remains to show that $\Psi = \Phi^{-1}$. We first show that for all $A \in \E_\infty^L$, $(\Psi \circ \Phi)(A) = A$. For all $m,n \in \N$, we have that $(m,n) \in (\Psi \circ \Phi)(A)$ iff $\Phi(A)(m) = \Phi(A)(n)$ iff $(m,n) \in A$, and for all $w \in L$ and $n \in \N$, $(w,n) \in (\Psi \circ \Phi)(A)$ iff $\Phi(A)(n) = w$ iff $(w,n) \in A$. Therefore, $(\Psi \circ \Phi)(A) = A$. Conversely, let $X \in W_\infty^L(\omega)$. For all $n \in \N$, $(\Phi \circ \Psi)(X)(n) = w$ iff $(w,n) \in \Psi(X)$ iff $X(n) = w$, and for all $m,n \in \N$, $(\Phi \circ \Psi)(X)(m) = (\Phi \circ \Psi)(X)(n)$ iff $(m,n) \in \Psi(X)$ iff $X(m) = X(n)$. Thus, $(\Phi \circ \Psi)(X) = X$.

    \underline{$\Phi$ is a biprojection:} That $\Phi$ sends the $\leq$-maximal element in $\E_\infty^L$ (i.e. the equality relation on $L \cup \N$) to the $\leq$-maximal element in $W_\infty^L(\omega)$ (i.e. the variable word $v_0v_1v_2v_3\dots$) follows from the fact that $\Phi$ is monotone and surjective. Let $A \in \E_\infty^L$ and $X \leq \Phi(A)$. Since $\Phi$ is surjective, we have that $\Phi(B) = X$ for some $B \in \E_\infty^L$. Since $\Phi^{-1}$ is monotone, we have that $\Phi(B) \leq \Phi(A)$ implies that $B \leq A$, as needed. This also implies that $\Psi$ satisfies the properties of Lemma \ref{lem:monotone.ab.stronger.properties}, as for all $X \in W_\infty^L(\omega)$:
    \begin{align*}
        \Phi``[\Psi(X)] = [X] &\implies (\Phi^{-1})``[X] = [\Psi(X)] \implies \Psi``[X] = [\Psi(X)].
    \end{align*}
\end{proof}

\begin{example}[Tree representation of Carlson--Simpson space with variables]
\label{ex:CS-variables.reformulation.trees}
    Let $L = \{w_0,\dots,w_{d-1}\}$ be a finite alphabet of size $d \geq 2$. We consider the last two formulations of the Carlson--Simpson space with variables $W_\infty^L(\omega)$ and $\cal{SE}_\infty^L$. In Example \ref{ex:CS-variables}, we defined a map $\Theta : W_\infty^L(\omega) \to \cal{SE}_\infty^L$ by stipulating that for all $A \in W_\infty^L(\omega)$, $\Theta(A) \subseteq d^{<\omega}$:
    \begin{itemize}        
        \item If $A(\dom(s)) = w_i$, then $s^\frown j \in \Theta(A)$ iff $j = i$.
        
        \item Note that for each $k$, $v_k$ first occurs at position $p_k(A)$.
        \begin{itemize}
            \item If $\dom(s) = p_k(A)$, then $s^\frown j \in \Theta(A)$ for all $j < d$.
            
            \item If $\dom(s) > p_k(A)$ and $A(\dom(s)) = v_k$, then $s^\frown j \in \Theta(A)$ iff $j = s(p_k(A))$.
        \end{itemize}
    \end{itemize}
    Then $\Theta$ is a biprojection. 
\end{example}

\begin{proof}
    It is clear that $\Theta$ is well-defined.

    \underline{$\Theta$ is monotone:} Suppose that $B \leq A$. We shall show that for all $s \in \Theta(B)$, $s \in \Theta(A)$ by induction on $|s|$. The base case $|s| = 0$ is clear. For the induction step, suppose that for all $s \in \Theta(B)$ such that $|s| \leq n$, $s \in \Theta(A)$. Let $i < d$ be such that $s^\frown i \in \Theta(B)$.
    \begin{itemize}
        \item Suppose that $A(n) = v_k$ for some $k$. Let $m \leq n$ be the first position at which $v_k$ occurs in $A$. Since $B \leq A$ and $A(m) = A(n) = v_k$, we have that $B(m) = B(n)$. This implies that $i = s(m)$, and since $s \in \Theta(A)$, we have that $s^\frown i \in \Theta(A)$.
        
        \item Suppose that $A(n) = w_j$ for some $j$. Since $B \leq A$, $B(n) = w_j$, so $j = i$. Therefore, $s^\frown i \in \Theta(A)$ as $s \in \Theta(A)$.
    \end{itemize}

    \underline{$\Theta$ is bijective:} $\Theta$ is surjective as $\cal{SE}_\infty^L$ is defined to be the image of $\Theta$. To see that $\Theta$ is injective, let $A,B \in W_\infty^L(\omega)$, and let $n$ be the smallest integer for which $A(n) \neq B(n)$. 
    \begin{itemize}
        \item Suppose that $A(n) = w_i \neq w_j = B(n)$. Then for all $s \in \Theta(A)$ and $t \in \Theta(B)$ such that $|s| = |t| \geq n + 1$, we have that $s(n) = i$ and $t(n) = j$. Therefore, $\Theta(A) \neq \Theta(B)$.
        
        \item Suppose that $A(n) = w_i$ and $B(n) = v_l$. Let $p \leq n$ be the first position at which $v_l$ occurs in $B$. If $p = n$, then every node of length $n$ splits in $\Theta(B)$, while no node of length $n$ splits in $\Theta(A)$ (as $A(n) = w_i$), so $\Theta(A) \neq \Theta(B)$. If $p < n$, then $A$ and $B$ agree below $n$, so $\Theta(A)$ and $\Theta(B)$ have the same nodes of length $\leq n$, and every node of length $p$ splits (as $p$ is the first occurrence of $v_l$). Pick a node $s$ of length $n$ with $s(p) = i' \neq i$ (which exists as $d \geq 2$). Then $s^\frown j \in \Theta(A)$ iff $j = i$, while $s^\frown j \in \Theta(B)$ iff $j = s(p) = i'$. Therefore, $\Theta(A) \neq \Theta(B)$.
        
        \item Suppose for a contradiction that $A(n) = v_k \neq v_l = B(n)$ and $\Theta(A) = \Theta(B)$. We may assume that $k < l$. Let $m_k < n$ be the first position at which $v_k$ occurs in $B$. Then $A(m_k) = B(m_k) = v_k$, so every $s \in \Theta(A) = \Theta(B)$ such that $|s| = n$ is not a splitting node. This, in turn, implies that $v_l$ does not first occur at position $n$ in $B$, so we may let $m_l < n$ be the first position at which $v_l$ occurs in $B$. Every node in $\Theta(B)$ of length either $m_k$ or $m_l$ splits, so there exists some $s \in \Theta(B)$ such that $|s| = n$, $s(m_k) = 0$ and $s(m_l) = 1$. Then we have that $s^\frown i \in \Theta(A)$ iff $i = 0$, while $s^\frown i \in \Theta(B)$ iff $i = 1$, contradicting that $\Theta(A) = \Theta(B)$.
    \end{itemize}

    \underline{$\Theta^{-1}$ is monotone:} Let $A,B \in W_\infty^L(\omega)$, and suppose that $\Theta(B) \subseteq \Theta(A)$. Recall that $B \leq A$ iff the following hold:
    \begin{enumerate}[label=(\roman*)]
        \item For all $n \in \N$, $A(n) = w \in L$ implies $B(n) = w$.

        \item For all $m,n \in \N$, $A(m) = A(n)$ implies $B(m) = B(n)$.
    \end{enumerate}
    For (i), suppose that $A(n) = w_i$. Then every node of length $n + 1$ in $\Theta(A)$ has $s(n) = i$, so the same holds in $\Theta(B) \subseteq \Theta(A)$. This implies that $B(n) = w_i$, for otherwise if $B(n)$ is a variable, then there exists some node $s$ such that $s(n) \neq i$. It remains to show that (ii) holds. 

    \begin{claim}
        Let $A \in W_\infty^L(\omega)$ and $m,n \in \N$. Then $A(m) = A(n)$ iff for all $s \in \Theta(A)$ such that $|s| = n + 1$, $s(n) = s(m)$.
    \end{claim}

    \begin{midproof}
        Suppose that $A(m) = A(n)$. We may assume that $m < n$.
        \begin{itemize}
            \item If $A(m) = A(n) = w_i$, then for all $s,t \in \Theta(A)$ such that $|s| = n + 1$ and $|t| = m + 1$, $s(n) = t(m) = i$. In particular, $s(n) = (s\restrictedto(m+1))(m) = s(m)$. 
            
            \item If $A(m) = A(n) = v_k$, then let $p \leq m$ be the first position at which $v_k$ occurs in $A$. Then for all $s,t \in \Theta(A)$ such that $|s| = n + 1$ and $|t| = m + 1$, $s(n) = s(p)$ and $t(m) = t(p)$. In particular, $s(n) = s(m)$.
        \end{itemize}

        Conversely, suppose $A(m) \neq A(n)$. 
        \begin{itemize}
            \item If $A(m) = w_i$ and $A(n) = w_j$, where $i \neq j$, then for any $s \in \Theta(A)$ such that $|s| = n + 1$, $s(n) = j$ and $s(m) = i$. 
            
            \item If $A(m) = v_k$ and $A(n) = w_j$, then let $p \leq m$ be the first position at which $v_k$ occurs in $A$. Let $s \in \Theta(A)$ be any node such that $|s| = n + 1$ and $s(p) = i \neq j$, which exists as every node in $\Theta(A)$ of length $p$ splits. Then $s(m) = s(p) = i$ and $s(n) = j$.

            \item If $A(m) = v_k$ and $A(n) = v_l$ with $k \neq l$, then let $p \leq m$ and $q \leq n$ be the first positions at which $v_k$ and $v_l$ occur in $A$, respectively. Since $k \neq l$, we have $p \neq q$, and both $p$ and $q$ are splitting levels of $\Theta(A)$, so there exists some $s \in \Theta(A)$ such that $|s| = n + 1$ and $s(p) \neq s(q)$. Then $s(m) = s(p) \neq s(q) = s(n)$.

            \item A similar argument holds for the case where $A(m) = w_i$ and $A(n) = v_l$.
        \end{itemize}
    \end{midproof}

    \underline{$\Theta$ is a biprojection:} The argument presented in the proof of Example \ref{ex:CS-variables.reformulation.variable-words} works.
\end{proof}

\subsubsection{Ramsey ultrafilters}
All examples listed in this subsection also send Ramsey ultrafilters to Ramsey ultrafilters:

\begin{lemma}
\label{lem:cts.monotone.pre-adjunction.preserves.Ramsey.ultrafilters}
    Let $f : \cal{R} \to \S$ be a continuous monotone pre-adjunction. If $\U \subseteq \cal{R}$ is a Ramsey ultrafilter, then $\cal{V} := f``(\U)^\uparrow$ is also a Ramsey ultrafilter.
\end{lemma}

\begin{proof}
    Since $f$ is monotone, $\cal{V}$ is a filter, so by Lemma \ref{lem:Ramsey.filter.is.ultrafilter}, it suffices to show that every analytic subset of $\S$ is $\cal{V}$-Ramsey. Let $\X \subseteq \S$ be analytic, so $f^{-1}(\X)$ is analytic as $f$ is continuous. Since $\U$ is a Ramsey ultrafilter, there exists some $A \in \U$ such that $[A] \subseteq f^{-1}(\X)$ or $[A] \subseteq (f^{-1}(\X))^c = f^{-1}(\X^c)$. Then $f(A) \in \cal{V}$, and $[f(A)] = f``[A] \subseteq \X$ or $[f(A)] = f``[A] \subseteq \X^c$ respectively.
\end{proof}

It is unclear whether Lemma \ref{lem:cts.monotone.pre-adjunction.preserves.Ramsey.ultrafilters} holds for continuous monotone quasi-open maps. One can show that $\cal{V}$ is an ultrafilter and that every Borel subset of $\S$ is $\cal{V}$-Ramsey, but we do not know whether the same holds for every analytic subset of $\S$.

\subsection{Inducing projections from maps between finite approximations}
\label{sec:induce.projections}
Observe from \S\ref{sec:projection.examples} that all maps there are defined in two steps --- first, a map $\alpha : \cal{AR} \to \cal{AS}$ respecting $\sqsubseteq$ is defined, and then $f : \cal{R} \to \S$ is defined by stipulating that:
\begin{align*}
    f(A) := \lim_{n\to\infty} \alpha(r_n(A)).
\end{align*}
Conversely, given a weak projection $f$, it is often helpful to find some $\alpha : \cal{AR} \to \cal{AS}$ that respects $\sqsubseteq$ and defines $f$ in this way. This subsection shows that, under suitable assumptions, this can be done on a cone below any given element.

Given a \textbf{wA2}-space $(\cal{R},\leq,r)$ and $A \in \cal{R}$, we write:
\begin{align*}
    \cal{R}^{\leq\infty} &:= \cal{R} \cup \cal{AR}, \\
    [A]^{\leq\infty} &:= [A] \cup \cal{AR}\restrictedto A.
\end{align*}
There is a natural metric topology on $\cal{R}^{\leq\infty}$, defined as follows: Note that for $A \in \cal{AR}$, $r_n(A)$ is defined only for $n \leq \lh(A)$. For distinct $A,B \in \cal{R}^{\leq\infty}$, we let $d(A,B) := \frac{1}{2^n}$, where $n$ is the least integer such that $r_n(A) \neq r_n(B)$, or exactly one of $r_n(A)$ and $r_n(B)$ is undefined. In particular, the following hold:
\begin{itemize}
    \item If $A,B \in \cal{R}$ and $A \neq B$, then $n$ is the least integer such that $r_n(A) \neq r_n(B)$.

    \item If $A \in \cal{AR}$, $B \in \cal{R}$ and $A \sqsubseteq B$, then $n = \lh(A) + 1$.

    \item If $A,B \in \cal{AR}$ and $A \sqsubsetneq B$, then $n = \lh(A) + 1$.
\end{itemize}
It is clear that the metric topology on $\cal{R}$ is precisely the subspace topology on $\cal{R}$ inherited from $\cal{R}^{\leq\infty}$ under this metric.

\begin{definition}
\label{def:sqsubseteq-preserving}
    A function $\alpha : \cal{AR} \to \S^{\leq\infty}$ is \emph{$\sqsubseteq$-preserving} if $\alpha(a) \sqsubseteq \alpha(b)$ whenever $a \sqsubseteq b$.
\end{definition}

\begin{definition}
\label{def:induce}
    Let $\alpha : \cal{AR} \to \S^{\leq\infty}$ be a function. If $f : \cal{R} \to \S^{\leq\infty}$ is a function such that for all $A \in \cal{R}$:
    \begin{align*}
        f(A) = \lim_{n\to\infty} \alpha(r_n(A)),
    \end{align*}
    then $\alpha$ \emph{induces} $f$.
\end{definition}

\begin{definition}
\label{def:epsilon_1.property}
    We say that a \textbf{wA2}-space $(\cal{R},\leq,r)$ satisfies the \emph{$(\varepsilon_1)$-property} if for all $A \in \cal{R}$ and $a \in \cal{AR}\restrictedto A$, $|[a,A]| > 1$.
\end{definition}

It is easy to see that all examples in \S\ref{sec:examples} (except for the singleton space in Example \ref{ex:singleton.space}) satisfy the $(\varepsilon_1)$-property. We also say that $(\cal{R},\leq,r)$ satisfies the \emph{$(\varepsilon)$-properties} if it satisfies both the $(\varepsilon_0)$-property (see Definition \ref{def:epsilon_0.property}) and the $(\varepsilon_1)$-property.

\begin{proposition}
\label{prop:ab.induced.by.sqsubseteq-preserving}
    Let $(\cal{R},\leq_\cal{R},r)$ and $(\S,\leq_\S,s)$ be two closed triples satisfying \textbf{A1}--\textbf{A4}, and suppose that $(\S,\leq_\S,s)$ satisfies the $(\varepsilon)$-properties. If $f : \cal{R} \to \S$ is a weak projection, then for all $A \in \cal{R}$, there exist some $B \leq A$ and some $\sqsubseteq$-preserving $\alpha : \cal{AR}\restrictedto B \to \cal{AS}\restrictedto f(B)$ that induces $f\restrictedto[B]$.
\end{proposition}

We first consider arbitrary monotone maps:

\begin{proposition}
\label{prop:nm.induced.by.sqsubseteq-preserving}
    Let $(\cal{R},\leq_\cal{R},r)$ and $(\S,\leq_\S,s)$ be two closed triples satisfying \textbf{A1}--\textbf{A4}, and suppose that $(\S,\leq_\S,s)$ satisfies the $(\varepsilon_0)$-property.
    \begin{enumerate}
        \item If $f : \cal{R} \to \S$ is monotone, then for all $A \in \cal{R}$, there exist some $B \leq A$ and some $\sqsubseteq$-preserving $\alpha : \cal{AR}\restrictedto B \to [f(B)]^{\leq\infty}$ that induces $f\restrictedto[B]$. In particular, $f\restrictedto[B]$ is continuous.
        
        \item If, furthermore, $|f``[a,C]| > 1$ for all $a \in \cal{AR}\restrictedto A$ and $C \in [a,A]$, then $\ran(\alpha) \subseteq \cal{AS}\restrictedto f(B)$.
    \end{enumerate}
\end{proposition}

\begin{proof}
    Fix some $A \in \cal{R}$. By Lemma \ref{lem:monotone.automatic.continuity}, we may assume that $f\restrictedto[A] : [A] \to [f(A)]$ is continuous. We shall pick up where we left off in the proof of Lemma \ref{lem:monotone.automatic.continuity} with the countable base $\{[x]\}_{x \in \cal{AS}\restrictedto f(A)}$. We write:
    \begin{enumerate}
        \item $\varphi_0(a,B,x)$, iff $f(C) \in [x]$ for all $C \in [a,B]$, iff $x \sqsubseteq f(C)$ for all $C \in [a,B]$.

        \item $\varphi_1(a,B,x)$, iff $f(C) \notin [x]$ for all $C \in [a,B]$, iff $x \not\sqsubseteq f(C)$ for all $C \in [a,B]$.
        
        \item $\varphi(a,B,x)$, iff $\varphi_0(a,B,x)$ or $\varphi_1(a,B,x)$.
    \end{enumerate}
    At the end of the proof, we concluded that for all $x \in \cal{AS}\restrictedto f(A)$:
    \begin{align*}
        (f\restrictedto[A])^{-1}([x]) = \bigcup\{[a,A] : \varphi_0(a,A,x)\}.
    \end{align*}
    For each $a \in \cal{AR}\restrictedto A$, we let:
    \begin{align*}
        \G_a := \{x \in \cal{AS} : \varphi_0(a,A,x)\}.
    \end{align*}
    Note that $\G_a \neq \emptyset$, as $\emptyset \in \G_a$ for all $a \in \cal{AR}\restrictedto A$. Observe that for all $x,y \in \G_a$, $x \sqsubseteq y$ or $y \sqsubseteq x$ --- if we take any $B \in [a,A]$, then $\varphi_0(a,A,x)$ and $\varphi_0(a,A,y)$ imply that $x \sqsubseteq f(B)$ and $y \sqsubseteq f(B)$, so $x$ and $y$ are $\sqsubseteq$-comparable. Furthermore, from the definition of $\varphi_0(a,A,x)$, we have that if $a \sqsubseteq b$ then $\G_a \subseteq \G_b$.

    \begin{claim}
    \label{claim:G.tends.to.infinity}
        If $a_0 \sqsubseteq a_1 \sqsubseteq \cdots$ in $\cal{AR}\restrictedto A$ are such that $\lim_{n\to\infty} \lh(a_n) = \infty$, then $\lim_{n\to\infty} |\G_{a_n}| = \infty$.
    \end{claim}

    \begin{midproof}
        Since $a_n \sqsubseteq a_{n+1}$, we have $\G_{a_n} \subseteq \G_{a_{n+1}}$, so $(|\G_{a_n}|)_{n<\omega}$ is non-decreasing. If $\G_{a_n}$ is infinite for some $n$, then we are done. Otherwise, for each $n$ let $x_n$ be the unique $\sqsubseteq$-maximal element of $\G_{a_n}$. Suppose, for a contradiction, that $|\G_{a_n}| \not\to \infty$ as $n \to \infty$. Then there exist some $x \in \cal{AS}$ and $N$ such that $x_n = x$ for all $n \geq N$. Let $B := \lim_{n\to\infty} a_n$, and let $y := r_{\lh(x)+1}(f(B))$. For each $n$, since $y \notin \G_{a_n}$, there exists some $B_n \in [a_n,A]$ such that $y \not\sqsubseteq f(B_n)$. But $\lim_{n\to\infty} B_n = B$ in the metric topology, and since $f\restrictedto[A]$ is continuous, we have that $(f^{-1}[y,f(A)])^c$ is closed in $[A]$ and $B_n \in (f^{-1}[y,f(A)])^c$ for all $n$. Therefore, $B \notin f^{-1}[y,f(A)]$, contradicting that $y \sqsubseteq f(B)$. 
    \end{midproof}

    Thus, we may define $\alpha : \cal{AR}\restrictedto A \to [f(A)]^{\leq\infty}$ by stipulating that:
    \begin{align*}
        \alpha(a) :=
        \begin{cases}
            \text{the unique $\sqsubseteq$-maximal element of $\G_a$}, & \text{if $\G_a$ is finite}, \\
            X, & \text{if $\G_a$ is infinite},
        \end{cases}
    \end{align*}
    where, in the latter case, $X \in [f(A)]$ is the unique element such that $x \sqsubseteq X$ for all $x \in \G_a$. By Claim \ref{claim:G.tends.to.infinity}, for all $C \in [A]$ we have $\lim_{n\to\infty} \lh(\alpha(r_n(C))) = \infty$, so $\alpha$ induces $f$. Moreover, since $\G_a \subseteq \G_b$ whenever $a \sqsubseteq b$, we have $\alpha(a) \sqsubseteq \alpha(b)$, so $\alpha$ is $\sqsubseteq$-preserving.
    
    Suppose, furthermore, that $|f``[a,A]| > 1$ for all $a \in \cal{AR}\restrictedto A$. To show that $\ran(\alpha) \subseteq \cal{AS}$, it suffices to show that $\G_a$ is finite for all $a \in \cal{AR}\restrictedto A$. Otherwise, there exist $x_0 \sqsubseteq x_1 \sqsubseteq \cdots$ in $\G_a$ such that $X := \lim_{n\to\infty} x_n \in \S$. This implies that $x_n \sqsubseteq f(B)$ for all $B \in [a,A]$, so $f``[a,A] = \{X\}$, a contradiction.
\end{proof}

Note that if $f$ is a weak projection, then by Lemma \ref{lem:ab.iff.qo.TRS}, $f``[A]$ is not Ellentuck meagre for all $A \in \cal{R}$. In particular, $|f``[A]| > 1$ for all $A \in \cal{R}$ if every singleton subset of $\S$ is Ellentuck meagre. The following lemma shows that this observation may be extended to all Ellentuck neighbourhoods below some $B \in \cal{R}$.

\begin{lemma}
\label{lem:function.automatic.non-meagre}
    Let $(\cal{R},\leq,r)$ be a closed triple satisfying \textbf{A1}--\textbf{A4}, and let $(P,\tau)$ be a Baire space (i.e. every non-empty open subset of $P$ is non-meagre). Let $f : \cal{R} \to P$ be any function. Let $A \in \cal{R}$, and suppose that $f``[B]$ is not meagre for all $B \leq A$. Then there exists some $B \leq A$ such that for all $b \in \cal{AR}\restrictedto B$ and $C \in [b,B]$, $f``[b,C]$ is not meagre.
\end{lemma}

\begin{proof}
    We shall appeal to Theorem \ref{thm:abstract.combinatorial.forcing} and use a combinatorial forcing argument. Given $A \in \cal{R}$ and $a \in \cal{AR}$, we say that:
    \begin{enumerate}[label=(\roman*)]
        \item $A$ \emph{accepts} $a$ if $f``[a,A]$ is meagre.
        
        \item $A$ \emph{rejects} $a$ if $a \in \cal{AR}\restrictedto A$ and no $B \in [a,A]$ accepts $a$.
    \end{enumerate}
    It follows from the definitions that it is not possible for $A$ to both accept and reject $a$. We check that statements (1)--(7) of Theorem \ref{thm:abstract.combinatorial.forcing} hold.
    \begin{enumerate}
        \item If $a \notin \cal{AR}\restrictedto A$, then $[a,A] = \emptyset$ is meagre.

        \item If $B \leq A$, then $f``[a,B] \subseteq f``[a,A]$, so $f``[a,B]$ is meagre if $f``[a,A]$ is.

        \item If $B \leq A$, then $[a,B] \subseteq [a,A]$, so for all $C \in [a,B] \subseteq [a,A]$, $f``[a,C]$ is not meagre as $A$ rejects $a$.

        \item This follows from (1), (2) and (3).
        
        \item If $A$ rejects $a$, then any $B \in [a,A]$ rejects $a$ by (3), and so decides $a$. Otherwise, there exists some $B \in [a,A]$ that accepts $a$, so $B$ decides $a$.
        
        \item This follows from (2) and (3).
        
        \item If $b \in r_{\lh(a)+1}[a,A]$, then $[b,A] \subseteq [a,A]$, so if $f``[a,A]$ is meagre then so is $f``[b,A] \subseteq f``[a,A]$. Conversely, observe that for all $B \in [a,A]$:
        \begin{align*}
            f``[a,B] &= f``\bb{\bigcup_{b \in r_{\lh(a)+1}[a,B]} [b,B]} = \bigcup_{b \in r_{\lh(a)+1}[a,B]} f``[b,B].
        \end{align*}
        Thus, if $f``[b,B]$ is meagre for all $b \in r_{\lh(a)+1}[a,B]$, then $f``[a,B]$ is a meagre set as it is a countable union of meagre sets. 
    \end{enumerate}
    Since $f``[B]$ is not meagre for all $B \leq A$, $A$ rejects $\emptyset$, so by Theorem \ref{thm:abstract.combinatorial.forcing} there exists some $B \leq A$ that rejects every $b \in \cal{AR}\restrictedto B$. In other words, $f``[b,C]$ is not meagre for all $b \in \cal{AR}\restrictedto B$ and $C \in [b,B]$. 
\end{proof}

Therefore, Proposition \ref{prop:nm.induced.by.sqsubseteq-preserving} applies to all weak projections, which gives Proposition \ref{prop:ab.induced.by.sqsubseteq-preserving}:

\begin{proof}[Proof of Proposition \ref{prop:ab.induced.by.sqsubseteq-preserving}]
    By Lemma \ref{lem:ab.iff.qo.TRS}, $f$ is quasi-open, so for all $A \in \cal{R}$, $f``[A]$ is not Ellentuck meagre as it contains $[f(B)]$ for some $B \leq A$. By Lemma \ref{lem:function.automatic.non-meagre}, there exists some $B \leq A$ such that $f``[b,C]$ is not Ellentuck meagre for all $b \in \cal{AR}\restrictedto B$ and $C \in [b,B]$. Since $(\S,\leq_\S,s)$ satisfies the $(\varepsilon_1)$-property, every singleton subset of $\S$ is Ellentuck meagre, so $|f``[b,C]| > 1$ for all such $b,C$, and the conclusion follows from Proposition \ref{prop:nm.induced.by.sqsubseteq-preserving}.
\end{proof}

\appsubsection{Non-monotone pre-adjunctions}
If a pre-adjunction $f$ is also monotone, then $f``(p^\downarrow) = f(p)^\downarrow$ for all $p \in \P$. In particular, every function satisfying the properties of Lemma \ref{lem:monotone.ab.stronger.properties} is a pre-adjunction, so every example in \S\ref{sec:projection.examples} is a pre-adjunction. However, a pre-adjunction is not necessarily monotone --- in fact, it is possible that $f : \cal{R} \to \S$ is a continuous pre-adjunction between two topological Ramsey spaces, but for every $A \in \cal{R}$, $f\restrictedto[A]$ is not monotone.

\begin{example}
    Define $\min_{\geq 1} : \FIN_2 \to \N$ by:
    \begin{align*}
        \min\nolimits_{\geq 1}(x) := \min\{m \in \N : x(m) \geq 1\} = \min(\supp(x)).
    \end{align*}
    We can extend this map to $\min_{\geq 1} : \FIN_2^{[\infty]} \to [\N]^\infty$ by:
    \begin{align*}
        \min\nolimits_{\geq 1}((x_n)_{n<\omega}) := \{\min\nolimits_{\geq 1}(x_n) : n < \omega\}.
    \end{align*}
    Then $\min_{\geq 1}$ is a continuous pre-adjunction that preserves the $\leq$-maximal element, but for all $A \in \FIN_2^{[\infty]}$, there exist $C \leq B \leq A$ such that for all $D \leq C$, $\min_{\geq 1}(D) \cap \min_{\geq 1}(B) = \emptyset$ --- in particular, $\min_{\geq 1}\restrictedto[A]$ is not monotone.
\end{example}

\begin{proof}
    It is clear that $\min_{\geq 1} : \FIN_2^{[\infty]} \to [\N]^\infty$ is continuous. The $\leq$-maximal element of $\FIN_2^{[\infty]}$ is $(e_n^2)_{n<\omega}$, where $e_n^2 \in \FIN_2$ is the function such that $\supp(e_n^2) = \{n\}$ and $e_n^2(n) = 2$. Then $\min_{\geq 1}(e_n^2) = n$ for all $n$, so $\min_{\geq 1}$ preserves the $\leq$-maximal element as $\min_{\geq 1}((e_n^2)_{n<\omega}) = \N$.

    We show that $\min_{\geq 1}$ is a pre-adjunction. Let $A = (x_n)_{n<\omega} \in \FIN_2^{[\infty]}$, and let $X \in [\min_{\geq 1}(A)]$. Then there exists an increasing sequence $(n_k)_{k<\omega}$ such that $X = \{\min_{\geq 1}(x_{n_k}) : k < \omega\}$. Let $B := (x_{n_k})_{k<\omega}$, and we have that $B \leq A$ and $\min_{\geq 1}(B) = X$. Therefore, $[\min_{\geq 1}(A)] \subseteq \min_{\geq 1}``[A]$. 

    For the final property, fix some $A = (x_n)_{n<\omega} \in \FIN_2^{[\infty]}$. Let $y_n := T(x_{2n}) + x_{2n+1}$, and let $B := (y_n)_{n<\omega}$. Then $B \in \FIN_2^{[\infty]}$, as for all $n$:
    \begin{align*}
        \max(\supp(y_n)) &= \max(\supp(x_{2n+1})) \\
        &< \min(\supp(T(x_{2n+2}))) \\
        &= \min(\supp(y_{n+1})).
    \end{align*}
    It is clear that $B \leq A$, and:
    \begin{align*}
        \min\nolimits_{\geq1}(B) = \{\min\nolimits_{\geq1}(y_n) : n < \omega\} = \{\min(\supp(T(x_{2n}))) : n < \omega\}.
    \end{align*}
    In particular, every element of $\min_{\geq1}(B)$ belongs to $\supp(x_{2n})$ for some $n$. Now let $z_n := T(y_{2n}) + y_{2n+1}$. By the same argument as above, $C := (z_n)_{n<\omega} \in \FIN_2^{[\infty]}$ and $C \leq B \leq A$. Observe that $T(y_n) = T(x_{2n+1})$ for all $n$, so for all $n$:
    \begin{align*}
        \min\nolimits_{\geq1}(z_n) &= \min\nolimits_{\geq1}(T(y_{2n})) = \min\nolimits_{\geq1}(T(x_{4n+1})), \\
        \min\nolimits_{\geq1}(T(z_n)) &= \min\nolimits_{\geq1}(T(y_{2n+1})) = \min\nolimits_{\geq1}(T(x_{4n+3})). 
    \end{align*}
    If $D = (w_m)_{m<\omega} \leq C$, then for all $m$ we have that $\min_{\geq1}(w_m) = \min_{\geq1}(z_n)$ or $\min_{\geq1}(w_m) = \min_{\geq1}(T(z_n))$ for some $n$. This implies that $\min_{\geq1}(w_m) = \min_{\geq1}(T(x_l)) \in \supp(x_l)$ for some odd $l$. Since the supports of the $x_n$'s are pairwise disjoint, $\min_{\geq1}(w_m) \notin \min_{\geq1}(B)$. Therefore, $\min_{\geq1}(D) \cap \min_{\geq1}(B) = \emptyset$.
\end{proof}

A similar proof shows that for any $k \geq 2$ and $1 \leq i < k$, the function $\min_{\geq i} : \FIN_k \to \N$ defined by:
\begin{align*}
    \min\nolimits_{\geq i}(x) := \min\{m \in \N : x(m) \geq i\}
\end{align*}
satisfies the same properties. 

\appsubsection{Projections between abstract orders}
\label{sec:projection.abstract.orders}
Let $(\cal{R}_0,\S_0,\leq_0,\leq_0^o)$ and $(\cal{R}_1,\S_1,\leq_1,\leq_1^o)$ be two analytic abstract orders. As the analytic Ramsey property is a statement about analytic subsets of $\cal{R}_0$ and $\cal{R}_1$, the set of continuous functions $f : \cal{R}_0 \to \cal{R}_1$ is the natural class of functions to consider. Since $\cal{R}_0$ and $\cal{R}_1$ are not equipped with any quasi-orders on them, there is no meaningful definition of a ``dense open'' subset of $\cal{R}$, or a ``monotone'' $f : \cal{R}_0 \to \cal{R}_1$. However, we can extend the notion of the quasi-open property to abstract orders.

\begin{definition}
    Let $(\P_0,\Q_0,\leq_0,\leq_0^o)$ and $(\P_1,\Q_1,\leq_1,\leq_1^o)$ be two abstract orders, and let $f : \P_0 \to \P_1$ be a function. We say that $f$ is \emph{quasi-open} if for all $q \in \Q_0$, there exists some $q' \in \Q_1$ such that:
    \begin{align*}
        \emptyset \neq \{p' \in \P_1 : p' \leq_1^o q'\} \subseteq f``(\{p \in \P_0 : p \leq_0^o q\}).
    \end{align*}
\end{definition}

\begin{lemma}
\label{lem:cts.qo.transfers.weak.analytic.Ramsey.property.ARS}
    Let $(\cal{R}_0,\S_0,\leq_0,\leq_0^o)$ and $(\cal{R}_1,\S_1,\leq_1,\leq_1^o)$ be two analytic abstract orders. If $(\cal{R}_0,\S_0,\leq_0,\leq_0^o)$ has the weak analytic Ramsey property and there exists a continuous quasi-open $f : \cal{R}_0 \to \cal{R}_1$, then $(\cal{R}_1,\S_1,\leq_1,\leq_1^o)$ has the weak analytic Ramsey property.
\end{lemma}

\begin{proof}
    Let $\Y \subseteq \cal{R}_1$ be analytic. Since $f$ is continuous, $f^{-1}(\Y)$ is analytic, so by our hypothesis there exists some $X \in \S_0$ such that $[X] \subseteq f^{-1}(\Y)$ or $[X] \subseteq (f^{-1}(\Y))^c = f^{-1}(\Y^c)$. This implies that $f``[X] \subseteq \Y$ or $f``[X] \subseteq \Y^c$. Since $f$ is quasi-open, there exists some $X' \in \S_1$ such that $\emptyset \neq [X'] \subseteq f``[X]$. Then $[X'] \subseteq \Y$ or $[X'] \subseteq \Y^c$, as desired.
\end{proof}

The following example, which is an extension of the map defined in Example \ref{ex:unlocate.variable-words}, is the map used in \cite[Theorem 4.50]{T10} that transfers the Ramsey-theoretic properties of $(\FIN_L^{[\infty]},\FIN_{Lv}^{[\infty]},\leq,\leq^o,r,s)$ to $(W_L^{[\infty]},W_{Lv}^{[\infty]},\leq,\leq^o,r,s)$.

\begin{example}[Located variable-words to variable-words]
    Recall, from Example \ref{ex:unlocate.variable-words}, that we defined a map $\zeta : \FIN_L \cup \FIN_{Lv} \to W_L \cup W_{Lv}$ by letting, for all $x \in \FIN_L \cup \FIN_{Lv}$ such that $\dom(x) = \{n_i : i < l\}$ where $n_0 < \cdots < n_{l-1}$, $\zeta(x)$ be the word of length $l$ such that the $i^\text{th}$ letter of $\zeta(x)$ is $x(n_i)$. Let $X \in \FIN_{Lv}^{[\infty]}$ be rapidly increasing, so $\zeta(X) \in W_{Lv}^{[\infty]}$ (as in Example \ref{ex:unlocate.variable-words}). Then define $\zeta_L : [X] \to [\zeta(X)]$ by letting, for all $(w_n)_{n<\omega} \in [X]$:
    \begin{align*}
        \zeta_L((w_n)_{n<\omega}) := (\zeta(w_n))_{n<\omega}.
    \end{align*}
    Then it's clear that $\zeta_L : [X] \to [\zeta(X)]$ is continuous, and we can show, in almost the same way as the proof of Example \ref{ex:unlocate.variable-words}, that for all $Y \in [\emptyset_{\FIN_{Lv}^{[\infty]}},X]$, $\zeta_L``[Y] = [\zeta(Y)]$. Therefore, $\zeta_L$ is a continuous quasi-open map from $([X],[\emptyset_{\FIN_{Lv}^{[\infty]}},X],\leq,\leq^o)$ to $([\zeta(X)],[\emptyset_{W_{Lv}^{[\infty]}},\zeta(X)],\leq,\leq^o)$.
\end{example}

In fact, one may show that $\zeta_L$ is similar to a \emph{strong projection} (Definition \ref{def:strong.projection}), in the following sense: Extending $\zeta_L$ to $\FIN_L^{[<\infty]}$ by letting $\zeta_L((w_n)_{n<N}) := (\zeta(w_n))_{n<N}$, we can show that for all rapidly increasing $X \in \FIN_{Lv}^{[\infty]}$ and $a \in \FIN_L^{[<\infty]}\restrictedto X$:
\begin{align*}
    \zeta_L``[a,X] = [\zeta_L(a),\zeta(X)].
\end{align*}
Therefore, for all rapidly increasing $X \in \FIN_{Lv}^{[\infty]}$ and $a \in \FIN_L^{[<\infty]}\restrictedto X$, $\zeta_L\restrictedto[a,X] : [a,X] \to [\zeta_L(a),\zeta(X)]$ is a continuous quasi-open map from $([a,X],[\depth_X(a),X],\leq,\leq^o)$ to $([\zeta_L(a),\zeta(X)],[\depth_{\zeta(X)}(\zeta_L(a)),\zeta(X)],\leq,\leq^o)$. Furthermore, for every $X' \in W_{Lv}^{[\infty]}$ and $a' \in W_L^{[<\infty]}\restrictedto X'$, there exist a rapidly increasing $X \in \FIN_{Lv}^{[\infty]}$ and $a \in \FIN_L^{[<\infty]}\restrictedto X$ such that $\zeta(X) = X'$ and $\zeta_L(a) = a'$ --- indeed, if $X' = (p_n')_{n<\omega} \in W_{Lv}^{[\infty]}$, then let $p_n \in \FIN_{Lv}$ be the located variable-word such that $\dom(p_n) = \{M_n,\dots,M_{n+1}-1\}$, where $M_n = \sum_{i<n} |p_i'|$, and $p_n(M_n + i)$ is the letter/variable at the $i^\text{th}$ position of $p_n'$. Write $a' = (q_i')_{i<N}$, and for each $i < N$:
\begin{align*}
    q_i' = p_{j_{i,0}}'[\lambda_{j_{i,0}}]^\frown\cdots^\frown p_{j_{i,N_i}}'[\lambda_{j_{i,N_i}}]
\end{align*}
for some $j_{i,0} < \cdots < j_{i,N_i} < j_{i+1,0}$ and $\lambda_{j_{i,k}} \in L_{j_{i,k}}$ for all $k \leq N_i$. Let:
\begin{align*}
    q_i := p_{j_{i,0}}[\lambda_{j_{i,0}}] \cup \cdots \cup p_{j_{i,N_i}}[\lambda_{j_{i,N_i}}],
\end{align*}
and $a := (q_i)_{i<N}$. By property $(\star)$ in the proof of Example \ref{ex:unlocate.variable-words}, we have that $a \in \FIN_L^{[<\infty]}\restrictedto X$ and $\zeta_L(a) = a'$.

Todor\v{c}evi\'{c}'s assertion that $\zeta_L$ ``transfers the Ramsey-theoretic properties'' may now be formalised as follows: Assume that $(\FIN_L^{[\infty]},\FIN_{Lv}^{[\infty]},\leq,\leq^o,r,s)$ satisfies \textbf{A4}, and that $(W_L^{[\infty]},W_{Lv}^{[\infty]},\leq,\leq^o,r,s)$ is a closed sextuple satisfying \textbf{A1}--\textbf{A3} (note that $\FIN_L^{[<\infty]}$ and $W_L^{[<\infty]}$ are countable). By Lemma \ref{lem:A4.iff.local.Ramsey.properties.ARS}, it suffices to show that for all $X' \in W_{Lv}^{[\infty]}$ and $a' \in W_L^{[<\infty]}\restrictedto X'$, $([a',X'],[\depth_{X'}(a'),X'],\leq,\leq^o)$ has the weak analytic Ramsey property. Let $X \in \FIN_{Lv}^{[\infty]}$ be rapidly increasing and $a \in \FIN_L^{[<\infty]}\restrictedto X$ such that $\zeta_L(a) = a'$ and $\zeta(X) = X'$. Then $\zeta_L\restrictedto[a,X] : [a,X] \to [\zeta_L(a),\zeta(X)]$ is continuous quasi-open. Since $([a,X],[\depth_X(a),X],\leq,\leq^o)$ has the weak analytic Ramsey property by Lemma \ref{lem:A4.iff.local.Ramsey.properties.ARS}, by Lemma \ref{lem:cts.qo.transfers.weak.analytic.Ramsey.property.ARS}, so does $([a',X'],[\depth_{X'}(a'),X'],\leq,\leq^o)$, as desired.

\section{Non-existence of projections}
\label{sec:non-existence.of.projections}
We shall now prove various non-existence results for weak projections between topological Ramsey spaces, including those mentioned in Theorem \ref{thm:main}.

\subsection{Non-existence by Tukey type theory}
\label{sec:non-existence.by.Tukey.type.theory}
Let $(\cal{R},\leq_\cal{R},\leq_\cal{R}^*,r)$ and $(\S,\leq_\S,\leq_\S^*,s)$ be two $\sigma$-closed topological Ramsey spaces. By Lemma \ref{lem:generic.filter.is.Ramsey}, forcing with $(\cal{R},\leq_\cal{R}^*)$ and $(\S,\leq_\S^*)$ adds Ramsey ultrafilters. This allows us to apply existing results about Tukey types of Ramsey ultrafilters of various topological Ramsey spaces to prove the non-existence of weak projections.

Recall that $(\FIN^{[\infty]},\leq,r)$ is the Hindman space of infinite block sequences --- see Example \ref{ex:Gowers.FIN_k}. Given a Ramsey ultrafilter $\U \subseteq \FIN^{[\infty]}$, i.e. a stable ordered-union ultrafilter (see Example \ref{ex:stable.ordered-union.ultrafilter}), we let $\U_{\min}$, $\U_{\max}$, $\U_{\minmax}$ and $\U_\FIN$ be the ultrafilters on $\N$, $\N$, $\N \times \N$ and $\FIN$ generated by the sets of the form $\min(A)$, $\max(A)$, $\{(\min(x),\max(x)) : x \in \c{A}\}$ and $\c{A}$ respectively, where $A \in \U$. Note that $\U_\FIN \equiv_T \U$. It is clear that $\U_{\min} \leq_T \U_{\minmax} \leq_T \U_\FIN$ and $\U_{\max} \leq_T \U_{\minmax} \leq_T \U_\FIN$, and it is known that $\U_{\min}$ and $\U_{\max}$ are selective ultrafilters that are Tukey incomparable --- see the remark after \cite[Corollary 69]{DT11}, and \cite[Theorem 10]{RT12}. \"{O}zalp extended this observation to a complete classification of the initial Tukey structure of $\U$:

\begin{theorem}[{\cite[Theorem 5.12]{O24}}]
\label{thm:souu.initial.Tukey.structure}
    Let $\U \subseteq \FIN^{[\infty]}$ be a Ramsey ultrafilter, and let $\cal{V}$ be any ultrafilter on $\N$ such that $\cal{V} \leq_T \U$. Then either $\cal{V}$ is principal, or $\cal{V} \equiv_T \cal{W}$ for some $\cal{W} \in \{\U_{\min},\U_{\max},\U_{\minmax},\U_\FIN\}$. 
\end{theorem}

This gives us the following theorem:

\begin{theorem}
\label{thm:no.wp.Ellentuck.FIN}
    $([\N]^\infty,\subseteq,r)$ does not weakly project into $(\FIN^{[\infty]},\leq,r)$.
\end{theorem}

\begin{proof}
    Suppose that there is a weak projection $f : [\N]^\infty \to \FIN^{[\infty]}$. Let $M$ be the transitive collapse of a countable elementary submodel of $H_\theta$ for some $\theta$ large enough, where $M$ contains $f$ and everything else relevant. Let $\U$ be $([\N]^\infty,\subseteq)$-generic over $M$. Then $\U$ is a selective ultrafilter on $\N$, and by Lemmas \ref{lem:generic.filter.is.Ramsey} and \ref{lem:ab.equivalence}, $\cal{V} := f``(\U)^\uparrow$ is a stable ordered-union ultrafilter (see also Example \ref{ex:stable.ordered-union.ultrafilter}) such that $\cal{V} \leq_T \U$. As mentioned above, $\cal{V}_{\min}$ and $\cal{V}_{\max}$ are Tukey incomparable. However, it was also shown in \cite{RT12} that selective ultrafilters on $\omega$ are Tukey-minimal. Since $\cal{V}_{\min} \leq_T \U$ and $\cal{V}_{\max} \leq_T \U$ and $\U$ is selective, we have that:
    \begin{align*}
        \cal{V}_{\min} \equiv_T \U \equiv_T \cal{V}_{\max},
    \end{align*}
    a contradiction.
\end{proof}

We shall later provide an alternative, forcing-free proof of Theorem \ref{thm:no.wp.Ellentuck.FIN} --- see the paragraph after the proof of the $\E_\infty \nrightarrow \FIN^{[\infty]}$ direction of Theorem \ref{thm:no.wp.Carlson-Simpson.and.FIN}.

\begin{corollary}
\label{cor:no.wp.Ellentuck.FIN_k}
    For all $1 \leq k < \omega$, $([\N]^\infty,\subseteq,r)$ does not weakly project into $(\FIN_k^{[\infty]},\leq,r)$.
\end{corollary}

\begin{proof}
    If $f : [\N]^\infty \to \FIN_k^{[\infty]}$ is a weak projection, then the map $T^{k-1} \circ f : [\N]^\infty \to \FIN^{[\infty]}$ (where $T$ is the \emph{tetris operation}, see Example \ref{ex:tetris.operation}) is also a weak projection, contradicting Theorem \ref{thm:no.wp.Ellentuck.FIN}.
\end{proof}

\begin{corollary}
\label{cor:no.wp.Ellentuck.Hales-Jewett}
    For any countable increasing union of finite alphabets $L = \bigcup_{n<\omega} L_n$ with $v \notin L$, $([\N]^\infty,\subseteq,r)$ does not weakly project into $(W_{Lv}^{[\infty]},\leq,r)$.
\end{corollary}

\begin{proof}
    If $f : [\N]^\infty \to W_{Lv}^{[\infty]}$ is a weak projection, then the map $\supp_{f(\N)} \circ f : [\N]^\infty \to \FIN^{[\infty]}$ (see Example \ref{ex:support.variable.words} for the definition of $\supp_A$ for $A \in W_{Lv}^{[\infty]}$) is also a weak projection, contradicting Theorem \ref{thm:no.wp.Ellentuck.FIN}.
\end{proof}

In \cite{DT15}, Dobrinen and Todor\v{c}evi\'{c} considered an ultrafilter $\U_\alpha \subseteq \cal{R}_\alpha$ (as in Definition \ref{def:filter}) that is Ramsey for $\cal{R}_\alpha$, and such that for any front $\F$ on $\cal{R}_\alpha$ and equivalence relation $E$ on $\F$, there is some $A \in \U_\alpha$ such that $E$ is canonical on $\F\restrictedto A = \F \cap \cal{AR}\restrictedto A$ (for the definition of a canonical equivalence relation on a front on $\cal{R}_\alpha$, see \cite[Definition 4.10 and Theorem 4.12]{DT15}).

\begin{lemma}
\label{lem:generic.implies.U_alpha}
    Let $M$ be a model of $\ZFC$, and let $\U_\alpha$ be $(\cal{R}_\alpha,\leq_\alpha)$-generic over $M$. Then, in $M[\U_\alpha]$, $\U_\alpha \subseteq \cal{R}_\alpha$ is an ultrafilter satisfying the following properties:
    \begin{enumerate}
        \item $\U_\alpha$ is Ramsey for $\cal{R}_\alpha$.

        \item $\U_\alpha$ is canonical for fronts on $\cal{R}_\alpha$.
    \end{enumerate}
\end{lemma}

\begin{proof}
    (1) follows from Lemma \ref{lem:generic.filter.is.Ramsey} and Example \ref{ex:R_alpha.Ramsey.ultrafilter}, so it remains to show that (2) holds. Working in $M[\U_\alpha]$, let $\F$ be a front on $\cal{R}_\alpha$, and let $E$ be an equivalence relation on $\F$. Since forcing with $(\cal{R}_\alpha,\leq_\alpha^*)$ adds no countable sets and both $\F$ and $E$ are countable, we have $\F, E \in M$. By \cite[Theorem 4.12]{DT15}, the set:
    \begin{align*}
        \D := \{A \in \cal{R}_\alpha : \text{$E$ is canonical on $\F\restrictedto A$}\}
    \end{align*}
    is a dense subset of $\cal{R}_\alpha$ in $M$. By the genericity of $\U_\alpha$, we may let $A \in \U_\alpha \cap \D$. Then $E$ is canonical on $\F\restrictedto A$, as desired.
\end{proof}

\begin{theorem}
\label{thm:wp.Dobrinen.spaces.limit}
    Let $\beta,\alpha < \omega_1$. If $(\cal{R}_\alpha,\leq_\alpha,r^\alpha)$ weakly projects into $(\cal{R}_\beta,\leq_\beta,r^\beta)$, then $\alpha = \gamma + \beta$ for some ordinal $\gamma$.
\end{theorem}

In particular, if $\alpha < \beta$, then $(\cal{R}_\alpha,\leq_\alpha,r^\alpha)$ does not weakly project into $(\cal{R}_\beta,\leq_\beta,r^\beta)$. This theorem also shows that the assertion in Example \ref{ex:Dobrinen.spaces.projections} that $(\cal{R}_\alpha,\leq_\alpha,r^\alpha)$ weakly projects into $(\cal{R}_\beta,\leq_\beta,r^\beta)$ for all $\beta < \alpha < \omega$ does not extend to all countable ordinals --- for instance, $(\cal{R}_\omega,\leq_\omega,r^\omega)$ does not weakly project into $(\cal{R}_\beta,\leq_\beta,r^\beta)$ for any $1 \leq \beta < \omega$. Note that $(\cal{R}_0,\leq_0,r^0)$ is equivalent to the Ellentuck space (see Example \ref{ex:Dobrinen.space}), so this theorem implies that the Ellentuck space does not weakly project into $(\cal{R}_1,\leq_1,r^1)$.

\begin{proof}
    Suppose that $f : \cal{R}_\alpha \to \cal{R}_\beta$ is a weak projection. Once again, let $M$ be the transitive collapse of a countable elementary submodel of some $H_\theta$, where $M$ contains $f$ and everything else relevant. Let $\U_\alpha$ be $(\cal{R}_\alpha,\leq_\alpha)$-generic over $M$. By Lemmas \ref{lem:generic.filter.is.Ramsey}, \ref{lem:ab.equivalence} and \ref{lem:generic.implies.U_alpha}, $\U_\beta := f``(\U_\alpha)^\uparrow$ is an ultrafilter of $\cal{R}_\beta$ that is Ramsey for $\cal{R}_\beta$ and canonical for fronts on $\cal{R}_\beta$. By \cite[Theorem 5.13]{DT15}, the Tukey types of non-principal ultrafilters below $\U_\alpha$ form a decreasing chain of order type $(\alpha + 1)^*$, and the Tukey types of non-principal ultrafilters below $\U_\beta$ form a decreasing chain of order type $(\beta + 1)^*$. Since $\U_\beta$ is Tukey below $\U_\alpha$, the latter chain is a downward-closed subset of the former, so $\beta + 1$ is the order type of a final segment of $\alpha + 1$. Therefore, $\alpha + 1 = \gamma + (\beta + 1)$ for some ordinal $\gamma$, so $\alpha = \gamma + \beta$.
\end{proof}

These results may be improved by examining the Tukey classes below ultrafilters of $\FIN^{[\infty]}$ and $\cal{R}_\alpha$ more closely.

\begin{theorem}
\label{thm:no.wp.R_alpha.and.FIN_k}
\hfill
\begin{enumerate}
    \item For all $1 \leq \alpha < \omega_1$, $(\FIN^{[\infty]},\leq,r)$ does not weakly project into $(\cal{R}_\alpha,\leq_\alpha,r^\alpha)$.
    
    \item For all $\alpha < \omega_1$ and $k \geq 1$, $(\cal{R}_\alpha,\leq_\alpha,r^\alpha)$ does not weakly project into $(\FIN_k^{[\infty]},\leq,r)$.
\end{enumerate}
\end{theorem}

\begin{proof}
    For both items, let $M$ be the transitive collapse of a countable elementary submodel of some $H_\theta$ containing the relevant weak projection and everything else relevant, and force over $M$ as in the proofs of Theorems \ref{thm:no.wp.Ellentuck.FIN} and \ref{thm:wp.Dobrinen.spaces.limit}.

    \underline{$\FIN^{[\infty]} \nrightarrow \cal{R}_\alpha$:} Let $1 \leq \alpha < \omega_1$, and suppose for a contradiction that $(\FIN^{[\infty]},\leq,r)$ weakly projects into $(\cal{R}_\alpha,\leq_\alpha,r^\alpha)$, so in the $(\FIN^{[\infty]},\leq)$-generic extension of $M$, we have that $\U_\alpha \leq_T \cal{V}$, where $\cal{V}$ is a stable ordered-union ultrafilter and $\U_\alpha \subseteq \cal{R}_\alpha$ is an ultrafilter that is Ramsey for $\cal{R}_\alpha$ and canonical for fronts on $\cal{R}_\alpha$ (by Lemmas \ref{lem:ab.equivalence} and \ref{lem:generic.implies.U_alpha}, as in the proof of Theorem \ref{thm:wp.Dobrinen.spaces.limit}, $\U_\alpha$ arises as the image of the generic under the weak projection, hence is itself $\cal{R}_\alpha$-generic over $M$). Since $\U_\alpha$ is non-principal (as $\c{A}$ is infinite for all $A \in \U_\alpha$), we have that by Theorem \ref{thm:souu.initial.Tukey.structure}, either $\U_\alpha \equiv_T \cal{V}_{\min}$ or $\U_\alpha \equiv_T \cal{V}_{\max}$, or $\cal{V}_{\min} <_T \U_\alpha$ and $\cal{V}_{\max} <_T \U_\alpha$.
    \begin{itemize}
        \item Since $\cal{V}_{\min}$ and $\cal{V}_{\max}$ are selective ultrafilters on $\omega$, if $\U_\alpha \equiv_T \cal{V}_{\min}$ or $\U_\alpha \equiv_T \cal{V}_{\max}$, then we may repeat the proof of Theorem \ref{thm:wp.Dobrinen.spaces.limit} and obtain a contradiction (i.e. $\U_0 \equiv_T \U_\alpha$, where $\U_0$ is the unique selective ultrafilter Tukey below $\U_\alpha$).

        \item If $\cal{V}_{\min} <_T \U_\alpha$ and $\cal{V}_{\max} <_T \U_\alpha$, then since the Tukey types of non-principal ultrafilters below $\U_\alpha$ form a chain \cite[Theorem 5.13]{DT15}, $\cal{V}_{\min}$ and $\cal{V}_{\max}$ are Tukey comparable. This contradicts \cite[Theorem 10]{RT12}.
    \end{itemize}

    \underline{$\cal{R}_\alpha \nrightarrow \FIN_k^{[\infty]}$:} If $(\cal{R}_\alpha,\leq_\alpha,r^\alpha)$ weakly projects into $(\FIN_k^{[\infty]},\leq,r)$, then it also weakly projects into $(\FIN^{[\infty]},\leq,r)$, by composing with $T^{k-1}$ as in the proof of Corollary \ref{cor:no.wp.Ellentuck.FIN_k}. Thus, in the $(\cal{R}_\alpha,\leq_\alpha)$-generic extension of $M$, we have that $\cal{V} \leq_T \U_\alpha$, where $\cal{V}$ is a stable ordered-union ultrafilter and $\U_\alpha \subseteq \cal{R}_\alpha$ is the generic ultrafilter, which is Ramsey for $\cal{R}_\alpha$ and canonical for fronts on $\cal{R}_\alpha$ (by Lemma \ref{lem:generic.implies.U_alpha}). Then $\cal{V}_{\min}$ and $\cal{V}_{\max}$ are non-principal ultrafilters Tukey below $\U_\alpha$ that are Tukey incomparable \cite[Theorem 10]{RT12}. But by \cite[Theorem 5.13]{DT15}, the Tukey types of non-principal ultrafilters below $\U_\alpha$ form a chain, a contradiction.
\end{proof}

As the initial Tukey structure of a Ramsey ultrafilter of $\FIN_k^{[\infty]}$ is unknown for $k \geq 2$, the analogue of Theorem \ref{thm:no.wp.R_alpha.and.FIN_k}(1) for $\FIN_k^{[\infty]}$ remains open:

\begin{qn}
    Let $2 \leq k < \omega$ and $1 \leq \alpha < \omega_1$. Does $(\FIN_k^{[\infty]},\leq,r)$ weakly project into $(\cal{R}_\alpha,\leq_\alpha,r^\alpha)$?
\end{qn}

\subsection{\texorpdfstring{Projections between $\FIN_k^{[\infty]}$'s}{Projections between Gowers spaces}}
\label{sec:projections.Gowers.spaces}
Recall from Example \ref{ex:tetris.operation} the tetris operation $T : \FIN_k \to \FIN_{k-1}$, defined by:
\begin{align*}
    T(x)(m) = \max\{x(m)-1,0\} \text{ for all $m < \omega$}.
\end{align*}
For $A = (x_n)_{n<\omega} \in \FIN_k^{[\infty]}$, we let $T(A) := (T(x_n))_{n<\omega}$. Then $T : \FIN_k^{[\infty]} \to \FIN_{k-1}^{[\infty]}$ is a projection (in fact, it is a \emph{strong projection} --- see Definition \ref{def:strong.projection}). The main theorem of this subsection asserts that this claim cannot be reversed in general:

\begin{theorem}
\label{thm:no.wp.FIN.FIN_2}
    $(\FIN^{[\infty]},\leq,r)$ does not weakly project into $(\FIN_2^{[\infty]},\leq,r)$.
\end{theorem}

We shall prove this theorem by analysing the Tukey structure below Ramsey ultrafilters of $\FIN_k^{[\infty]}$. For $k = 1$, Ramsey ultrafilters are precisely the stable ordered-union ultrafilters, and their Tukey structure has been resolved by \"{O}zalp in \cite{O24} using the tools of Tukey theory of ultrafilters in \cite{DT11} and \cite{RT12}, and the canonisation theory by Klein and Spinas in \cite{KS05}. L\'{o}pez-Abad described the canonical equivalence relations on $\FIN_k^{[\infty]}$ in \cite{L07}, but the Tukey theory of Ramsey ultrafilters of $\FIN_k^{[\infty]}$ remains largely unexplored.

For an ultrafilter $\U \subseteq \FIN_k^{[\infty]}$, we write $T(\U) := T``(\U)^\uparrow \subseteq \FIN_{k-1}^{[\infty]}$. Note that if $\U$ is a Ramsey ultrafilter, then $T(\U)$ is also a Ramsey ultrafilter of $\FIN_{k-1}^{[\infty]}$ by Lemma \ref{lem:cts.monotone.pre-adjunction.preserves.Ramsey.ultrafilters}. While we will not give a full description of the Tukey structure below a Ramsey ultrafilter of $\FIN_k^{[\infty]}$, we shall prove the following.

\begin{theorem}
\label{thm:TU.Tukey.above.U}
    Let $k \geq 2$, and let $\U \subseteq \FIN_k^{[\infty]}$ be a Ramsey ultrafilter. Then $T(\U) <_T \U$.
\end{theorem}

To prove this theorem, we need a localised version of Lemma \ref{lem:monotone.automatic.continuity} for Ramsey ultrafilters of $\FIN_k^{[\infty]}$.

\begin{lemma}
\label{lem:monotone.automatic.continuity.localised.FIN_k}
    Let $k \geq 1$, let $\U \subseteq \FIN_k^{[\infty]}$ be a Ramsey ultrafilter, and let $(\S,\leq_\S)$ be an analytic quasi-order satisfying the $(\varepsilon_0)$-property. If $f : \U \to \S$ is monotone, then for all $A \in \U$, there exists some $B \in \U\restrictedto A$ such that $f\restrictedto(\U\restrictedto B)$ is continuous.
\end{lemma}

\begin{proof}[Sketch of proof]
    One can repeat the proof of Lemma \ref{lem:monotone.automatic.continuity} inside $\U$, after checking that $[a,A] \cap \U \neq \emptyset$ for all $A \in \U$ and $a \in \FIN_k^{[<\infty]}\restrictedto A$, and that the conclusion of Lemma \ref{lem:R.is.semiselective} holds relative to $\U$. We omit the details.
\end{proof}

\begin{lemma}
\label{lem:tetris.non-trivial.quotient}
    Let $k \geq 2$, and let $\U \subseteq \FIN_k^{[\infty]}$ be a Ramsey ultrafilter. For all $A \in \U$, there exist incompatible $B_0,B_1 \leq A$ such that $T(B_0) = T(B_1)$ and $B_0 \in \U$.
\end{lemma}

\begin{proof}
    Fix any $A = (x_n)_{n<\omega} \in \U$. For each $n$, let:
    \begin{align*}
        m_n := \min\{m \in \supp(x_n) : x_n(m) = k\}.
    \end{align*}
    Define:
    \begin{align*}
        P := \{x \in \FIN_k : \forall n < \omega \, (m_n \in \supp(x) \to x(m_n) \geq 2)\}.
    \end{align*}
    Since $\U$ is a Ramsey ultrafilter, there exists some $B_0 \in \U\restrictedto A$ such that $\c{B_0} \subseteq P$ or $\c{B_0} \subseteq P^c$. We write $B_0 = (y_i)_{i<\omega}$.
    
    \begin{claim}
        $\c{B_0} \cap P^c \neq \emptyset$. Therefore, $\c{B_0} \subseteq P^c$.
    \end{claim}

    \begin{midproof}
        Let $l := \min\{m \in \supp(y_0) : y_0(m) = k\}$. Since $B_0 \leq A$, $l = m_n$ for some $n$. Then $y := T^{k-1}(y_0) + y_1 \in \c{B_0} \cap P^c$, as $m_n \in \supp(y)$ but $y(m_n) = 1$.
    \end{midproof}
    
    This implies that for all $i$, $y_i \in P^c$ --- that is, for all $i$ there exists some $n_i'$ such that $y_i(m_{n_i'}) = 1$. Since $B_0 \leq A$, for each $i$ we may write:
    \begin{align*}
        y_i = T^{\lambda_{n_i}}(x_{n_i}) + \cdots + T^{\lambda_{n_{i+1}-1}}(x_{n_{i+1}-1}).
    \end{align*}
    Note that $n_i \leq n_i' < n_{i+1}$, and $\lambda_{n_i'} = k - 1$. We define:
    \begin{align*}
        z_i &:= y_i - T^{k-1}(x_{n_i'}) \\
        &:= T^{\lambda_{n_i}}(x_{n_i}) + \cdots + T^{\lambda_{n_i'-1}}(x_{n_i'-1}) \\
        &\eqbreak + T^{\lambda_{n_i'+1}}(x_{n_i'+1}) + \cdots + T^{\lambda_{n_{i+1}-1}}(x_{n_{i+1}-1}).
    \end{align*}
    Note that $z_i \in \FIN_k$, as $\lambda_{n_i'} = k-1 \neq 0$. Then, for all $i$:
    \begin{align*}
        T(z_i) = T\bb{y_i - T^{k-1}(x_{n_i'})} = T(y_i) - T^k(x_{n_i'}) = T(y_i).
    \end{align*}
    Therefore, if we let $B_1 := (z_i)_{i<\omega}$, then $T(B_0) = T(B_1)$. It is also clear that $B_1 \leq A$. We shall show that $\c{B_0} \cap \c{B_1} = \emptyset$, which implies that they are incompatible, completing the proof.

    Suppose otherwise, so there exists some $w \in \c{B_0} \cap \c{B_1}$. Since $w \in \c{B_0}$, $w\restrictedto\supp(y_i) = y_i$ for some $i$. In particular, $w(m_{n_i'}) = y_i(m_{n_i'}) = 1$, so $m_{n_i'} \in \supp(w)$. On the other hand, $n_i \leq n_i' < n_{i+1}$ and each $z_j$ omits exactly the summand $T^{k-1}(x_{n_j'})$, so $m_{n_i'} \notin \supp(z_j)$ for all $j$. But $w \in \c{B_1}$, so this gives $m_{n_i'} \notin \supp(w)$, a contradiction.
\end{proof}

We may now return to the proof of Theorem \ref{thm:TU.Tukey.above.U}.

\begin{proof}[Proof of Theorem \ref{thm:TU.Tukey.above.U}]
    We write $(\cal{R},\leq,r)$ for $(\FIN_k^{[\infty]},\leq,r)$ and $(\S,\leq,s)$ for $(\FIN_{k-1}^{[\infty]},\leq,r)$, and let $\cal{V} := T(\U)$, which is a Ramsey ultrafilter of $\S$. Since $T$ is monotone, we have that $\cal{V} \leq_T \U$. Suppose for a contradiction that $\U \leq_T \cal{V}$, so there exists some cofinal function $h : \cal{V} \to \U$. Define $f : \cal{V} \to \Po(\cal{AR})$ by:
    \begin{align*}
        f(Y) := \bigcup\{\cal{AR}\restrictedto h(Y') : Y' \in \cal{V}\restrictedto Y\}.
    \end{align*}
    Then $f(Y') \subseteq f(Y)$ whenever $Y' \leq Y$, and $\cal{AR}\restrictedto h(Y) \subseteq f(Y)$ for all $Y \in \cal{V}$. Note that, by the cofinality of $h$, for all $B \in \U$, there exists some $Y \in \cal{V}$ such that $f(Y) \subseteq \cal{AR}\restrictedto B$.

    We equip $\Po(\cal{AR})$ with the product topology. Then $(\Po(\cal{AR}),\subseteq)$ is an analytic quasi-order satisfying the $(\varepsilon_0)$-property with respect to the base consisting of the sets $\{W \subseteq \cal{AR} : F \subseteq W \wedge W \cap G = \emptyset\}$ for finite $F,G \subseteq \cal{AR}$. By Lemma \ref{lem:monotone.automatic.continuity.localised.FIN_k} and the definition of $\cal{V}$, there exists some $A \in \U$ such that $f\restrictedto(\cal{V}\restrictedto T(A))$ is continuous. For each $x \in \cal{AS}\restrictedto T(A)$, let:
    \begin{align*}
        \alpha(x) := \bigcap\{f(Y) : Y \in \cal{V}\restrictedto[x,T(A)]\},
    \end{align*}
    and define $\bar{f} : [T(A)] \to \Po(\cal{AR})$ by:
    \begin{align*}
        \bar{f}(X) := \bigcup_{n<\omega} \alpha(s_n(X)).
    \end{align*}
    Then:
    \begin{itemize}
        \item $\bar{f}$ and $f$ coincide on $\cal{V}\restrictedto T(A)$. Indeed, if $X \in \cal{V}\restrictedto T(A)$, then $\alpha(s_n(X)) \subseteq f(X)$ for all $n$, and if $a \in f(X)$, then since $f\restrictedto(\cal{V}\restrictedto T(A))$ is continuous, there exists some $n$ such that $a \in f(Y)$ for all $Y \in \cal{V}\restrictedto[s_n(X),T(A)]$, so $a \in \alpha(s_n(X))$.

        \item For all $B \in \cal{R}$, $\bar{f}^{-1}(\Po(\cal{AR}\restrictedto B))$ is a closed subset of $[T(A)]$, as $\bar{f}(X) \subseteq \cal{AR}\restrictedto B$ iff $\alpha(s_n(X)) \subseteq \cal{AR}\restrictedto B$ for all $n$.
    \end{itemize}

    We consider the following set:
    \begin{align*}
        \X := \{B \leq A : \exists C \leq B \, \bar{f}``[T(C)] \subseteq \Po(\cal{AR}\restrictedto B)\}.
    \end{align*}

    \begin{claim}
    \label{claim:bar.f.image.analytic}
        $\X$ is analytic.
    \end{claim}

    \begin{midproof}
        We first note that for all $X \leq T(A)$ and $B \leq A$:
        \begin{align*}
            \bar{f}``[X] \subseteq \Po(\cal{AR}\restrictedto B) \iff \alpha``(\cal{AS}\restrictedto X) \subseteq \Po(\cal{AR}\restrictedto B).
        \end{align*}
        \begin{itemize}
            \item Suppose that $\bar{f}``[X] \subseteq \Po(\cal{AR}\restrictedto B)$. Let $x \in \cal{AS}\restrictedto X$, and let $Y \in [x,X]$. Then $Y \leq X$, so $\bar{f}(Y) \in \Po(\cal{AR}\restrictedto B)$. Since $x \sqsubseteq Y$, we have that $\alpha(x) \subseteq \bar{f}(Y)$, so $\alpha(x) \in \Po(\cal{AR}\restrictedto B)$. 

            \item Conversely, suppose that $\alpha``(\cal{AS}\restrictedto X) \subseteq \Po(\cal{AR}\restrictedto B)$. Let $Y \leq X$, so $\alpha(s_m(Y)) \in \Po(\cal{AR}\restrictedto B)$ for all $m$. Then $\bar{f}(Y) = \bigcup_{m<\omega} \alpha(s_m(Y)) \in \Po(\cal{AR}\restrictedto B)$.
        \end{itemize}
        Therefore:
        \begin{align*}
            \X = \{B \leq A : \exists C \leq B \, \alpha``(\cal{AS}\restrictedto T(C)) \subseteq \Po(\cal{AR}\restrictedto B)\}.
        \end{align*}
        Observe that the statement $\alpha``(\cal{AS}\restrictedto X) \subseteq \Po(\cal{AR}\restrictedto B)$ is arithmetical w.r.t. some real parameters, as:
        \begin{gather*}
            \alpha``(\cal{AS}\restrictedto X) \subseteq \Po(\cal{AR}\restrictedto B) \\
            \Updownarrow \\
            \forall x \in \cal{AS}\restrictedto T(A) \, \forall a \in \alpha(x) \, (\exists m[x \leq_\fin s_m(X)] \to \exists n[a \leq_\fin r_n(B)]).
        \end{gather*}
        Therefore, $\X$ is analytic.
    \end{midproof}
    
    \begin{claim}
    \label{claim:X.is.avoidable}
        There exists some $A_0 \in \U$ such that $[A_0] \subseteq \X^c$.
    \end{claim}

    \begin{midproof}
        Suppose for a contradiction that $[A_0] \not\subseteq \X^c$ for all $A_0 \in \U$. Since $\X$ is analytic and $\U$ is a Ramsey ultrafilter, there exists some $A_0 \in \U$ such that $[A_0] \subseteq \X$. By Lemma \ref{lem:tetris.non-trivial.quotient}, let $B_0,B_1 \leq A_0$ be incompatible reductions such that $T(B_0) = T(B_1) =: Y$ and $B_0 \in \U$. Note that $Y \in \cal{V}$.

        For each $i < 2$, let $\D_i := \bar{f}^{-1}(\Po(\cal{AR}\restrictedto B_i)) \cap [Y]$, which is Borel. Furthermore, $\D_i$ is dense below $Y$ --- if $Z \leq Y$, then since $T``[B_i] = [T(B_i)] = [Y]$ by Example \ref{ex:tetris.operation}, we have that $Z = T(B)$ for some $B \leq B_i$. Since $B \leq B_i \leq A_0$ and $[A_0] \subseteq \X$, there exists some $C \leq B$ such that $\bar{f}``[T(C)] \subseteq \Po(\cal{AR}\restrictedto B) \subseteq \Po(\cal{AR}\restrictedto B_i)$. Then $T(C) \leq Z$ and $T(C) \in \D_i$, as needed. Since $\cal{V}$ is a Ramsey ultrafilter and $Y \in \cal{V}$, there is some $Z_i \in \cal{V}\restrictedto Y$ such that $[Z_i] \subseteq \D_i$ or $[Z_i] \subseteq \D_i^c$. The latter case is not possible as $\D_i$ is dense, so $\bar{f}``[Z_i] \subseteq \Po(\cal{AR}\restrictedto B_i)$.

        Since $\cal{V}$ is a filter, we may let $Z \in \cal{V}$ such that $Z \leq Z_0$ and $Z \leq Z_1$. Then:
        \begin{align*}
            \bar{f}``[Z] \subseteq \bar{f}``[Z_0] \cap \bar{f}``[Z_1] \subseteq \Po(\cal{AR}\restrictedto B_0) \cap \Po(\cal{AR}\restrictedto B_1).
        \end{align*}
        But $Z \in \cal{V}$, so $\cal{AR}\restrictedto h(Z) \subseteq f(Z) = \bar{f}(Z)$. This implies that $h(Z) \in [B_0] \cap [B_1]$, contradicting that $B_0$ and $B_1$ are incompatible.
    \end{midproof}

    Therefore, we may fix some $A' \in \U\restrictedto A$ such that $[A'] \subseteq \X^c$. Recall that $\bar{f}^{-1}(\Po(\cal{AR}\restrictedto A'))$ is Borel. Since $\cal{V}$ is a Ramsey ultrafilter and $T(A') \in \cal{V}$, there exists some $Y \in \cal{V}\restrictedto T(A')$ such that $\bar{f}``[Y] \subseteq \Po(\cal{AR}\restrictedto A')$ or $\bar{f}``[Y] \subseteq \Po(\cal{AR}\restrictedto A')^c$. We shall show that either case leads to a contradiction.
    \begin{itemize}
        \item Suppose that $\bar{f}``[Y] \subseteq \Po(\cal{AR}\restrictedto A')$. Since $Y \in \cal{V}$, there is some $Z \in \U$ with $T(Z) \leq Y$. Let $B \in \U$ be a common reduction of $Z$ and $A'$, so $B \leq A'$ and $T(B) \leq Y$. Then $\bar{f}``[T(B)] \subseteq \bar{f}``[Y] \subseteq \Po(\cal{AR}\restrictedto A')$, so in particular, $A' \in \X$. This contradicts $[A'] \subseteq \X^c$.

        \item Suppose that $\bar{f}``[Y] \subseteq \Po(\cal{AR}\restrictedto A')^c$. Since $A' \in \U$, there is some $Z \in \cal{V}$ with $f(Z) \subseteq \cal{AR}\restrictedto A'$. As $\cal{V}$ is a filter, fix $Z' \in \cal{V}$ with $Z' \leq Z$ and $Z' \leq Y$. On the one hand, $Z' \leq Z$ gives $f(Z') \subseteq f(Z) \subseteq \cal{AR}\restrictedto A'$. On the other hand, $Z' \leq Y$ gives $\bar{f}(Z') \in \bar{f}``[Y] \subseteq \Po(\cal{AR}\restrictedto A')^c$, so $\bar{f}(Z') \nsubseteq \cal{AR}\restrictedto A'$. But $Z' \leq Y \leq T(A')$ and $Z' \in \cal{V}$, so $\bar{f}(Z') = f(Z')$, a contradiction.
    \end{itemize}
    Therefore, $\U \not\leq_T \cal{V}$, and hence $\cal{V} <_T \U$.
\end{proof}

We may now prove Theorem \ref{thm:no.wp.FIN.FIN_2} using the same argument as in the proof of Theorem \ref{thm:no.wp.Ellentuck.FIN}.

\begin{proof}
    Suppose that there is a weak projection $f : \FIN^{[\infty]} \to \FIN_2^{[\infty]}$. Let $M$ be the transitive collapse of a countable elementary submodel of $H_\theta$ for some $\theta$ large enough, where $M$ contains $f$ and everything else relevant. Let $\U$ be $(\FIN^{[\infty]},\leq)$-generic over $M$. Then $\U$ is a stable ordered-union ultrafilter, i.e. a Ramsey ultrafilter of $\FIN^{[\infty]}$. By Lemmas \ref{lem:generic.filter.is.Ramsey} and \ref{lem:ab.equivalence}, $\cal{V} := f``(\U)^\uparrow$ is a Ramsey ultrafilter of $\FIN_2^{[\infty]}$ such that $\cal{V} \leq_T \U$. By the work of \"{O}zalp in \cite{O24}, $\U$ has exactly four Tukey types of non-principal ultrafilters below it. However, $T(\cal{V}) <_T \cal{V}$ by Theorem \ref{thm:TU.Tukey.above.U}, so $\cal{V}$ has at least five Tukey types of non-principal ultrafilters below it --- the four Tukey types below $T(\cal{V})$ (as it is a stable ordered-union ultrafilter), and $\cal{V}$ itself. Therefore, $\cal{V} \leq_T \U$ is not possible.
\end{proof}

Since the initial Tukey structure of a Ramsey ultrafilter of $\FIN_k^{[\infty]}$ is unknown for $k \geq 2$, the following question remains open.

\begin{qn}
\label{qn:no.wp.FIN_k.FIN_kp1}
    Let $k \geq 2$. Does $(\FIN_k^{[\infty]},\leq,r)$ weakly project into $(\FIN_{k+1}^{[\infty]},\leq,r)$?
\end{qn}

We note that we do not need the full classification of Tukey types below a stable ordered-union ultrafilter, but merely that there are finitely many Tukey types of ultrafilters below a stable ordered-union ultrafilter. Consequently, for a fixed $k \geq 2$, if there is some finite $n$ such that every Ramsey ultrafilter of $\FIN_k^{[\infty]}$ has $n$ Tukey types of ultrafilters below it, then the same argument proves that there is no weak projection from $\FIN_k^{[\infty]}$ to $\FIN_{k+1}^{[\infty]}$. This motivates the following question:

\begin{qn}
\label{qn:FIN_k.Ramsey.ultrafilter.finite.initial.Tukey.structure}
    Let $k \geq 2$, and let $\U$ be a Ramsey ultrafilter of $\FIN_k^{[\infty]}$. Are there only finitely many Tukey types of ultrafilters Tukey below $\U$?
\end{qn}

Answers to the $k = 2$ case of both Questions \ref{qn:no.wp.FIN_k.FIN_kp1} and \ref{qn:FIN_k.Ramsey.ultrafilter.finite.initial.Tukey.structure} seem within reach. \"{O}zalp's determination of the initial Tukey structure of Ramsey ultrafilters of $\FIN^{[\infty]}$ in \cite{O24} relies on the canonisation theorem for Borel functions on $\FIN^{[\infty]}$ by Klein and Spinas \cite{KS05}, and \"{O}zalp recently extended this canonisation theorem to $\FIN_2^{[\infty]}$ in \cite{O26}.

\subsection{Strong subtrees and Sacks forcing}
\label{sec:projections.trees}
Forcing with $(\S_\infty(U),\subseteq)$ (Example \ref{ex:strong.subtrees}) or $(\mathbb{S},\subseteq)$ (Example \ref{ex:Sacks.forcing}) adds a real. On the other hand, if $(\cal{R},\leq,\leq^*,r)$ is a $\sigma$-closed closed triple satisfying \textbf{A1}--\textbf{A4}, then forcing with $(\cal{R},\leq)$ does not add any reals. Therefore, by Lemma \ref{lem:ab.equivalence}, $(\cal{R},\leq,r)$ does not weakly project into $(\S_\infty(U),\subseteq,r)$ or $(\mathbb{S},\subseteq,r)$. However, it is possible to provide a purely combinatorial proof.

\begin{theorem}
\label{thm:no.wp.sigma-closed.trees}
    Let $(\cal{R},\leq,\leq^*,r)$ be a $\sigma$-closed closed triple satisfying \textbf{A1}--\textbf{A4}. Then it does not weakly project into $(\S_\infty(U),\subseteq,r)$ or $(\mathbb{S},\subseteq,r)$.
\end{theorem}

\begin{proof}
    Let $(\S,\subseteq,r)$ be either $(\S_\infty(U),\subseteq,r)$ or $(\mathbb{S},\subseteq,r)$, and suppose that $f : \cal{R} \to \S$ is a weak projection. Without loss of generality, we may assume that $s^\frown 0,s^\frown 1 \in U$ for every splitting node $s$ of $U$. For each $n$, let:
    \begin{align*}
        \X_n := \{X \in \S : \exists s \in X(|s| = n \text{ and } s^\frown 0 \in X)\}.
    \end{align*}
    We define a $\leq$-decreasing sequence $(A_n)_{n<\omega}$ in $\cal{R}$ as follows: Since both $(\S_\infty(U),\subseteq,r)$ and $(\mathbb{S},\subseteq,r)$ satisfy the $(\varepsilon_0)$-property (see the paragraph after Definition \ref{def:epsilon_0.property}), by Lemma \ref{lem:monotone.automatic.continuity}, we may let $A_0 \in \cal{R}$ be any element such that $f\restrictedto[A_0]$ is continuous. Suppose that $A_n$ has been defined. Since $\X_n$ is a clopen subset of $\S$, $f^{-1}(\X_n) \cap [A_0]$ is clopen in $[A_0]$, hence Ramsey in $[A_0]$, so there exists some $A_{n+1} \leq A_n$ such that $[A_{n+1}] \subseteq f^{-1}(\X_n)$ or $[A_{n+1}] \subseteq (f^{-1}(\X_n))^c$. This implies that $f``[A_{n+1}] \subseteq \X_n$ or $f``[A_{n+1}] \subseteq \X_n^c$ respectively. This completes the construction of the sequence.

    Since $(\cal{R},\leq^*)$ is $\sigma$-closed, we may let $A^* \in \cal{R}$ such that $A^* \leq^* A_n$ for all $n$. Since $f$ is quasi-open by Lemma \ref{lem:ab.iff.qo.TRS}, we may let $Y \in \S$ such that $[Y] \subseteq f``[A^*]$, and let $n$ be the least integer such that $Y$ has a splitting node of length $n$.

    Let $t \in Y$ be the (unique) splitting node of length $n$. By construction, either $[A_{n+1}] \subseteq f^{-1}(\X_n)$ or $[A_{n+1}] \subseteq (f^{-1}(\X_n))^c$.
    \begin{itemize}
        \item First suppose that $[A_{n+1}] \subseteq f^{-1}(\X_n)$, so that $f``[A_{n+1}] \subseteq \X_n$. Let $Z \leq Y$ be such that $t^\frown 1 \in Z$ but $t^\frown 0 \notin Z$. Then the only node of length $n$ in $Z$ is $t$. Since $[Y] \subseteq f``[A^*]$, we have that $Z = f(B)$ for some $B \leq A^* \leq^* A_{n+1}$. Since $B$ and $A_{n+1}$ are compatible in $(\cal{R},\leq^*)$, and $(\cal{R},\leq^*)$ and $(\cal{R},\leq)$ share the same separative quotient, they are compatible in $(\cal{R},\leq)$. Thus, there exists some $B' \leq B$ such that $B' \leq A_{n+1}$. Let $Z' := f(B')$. Then $Z' \in f``[A_{n+1}] \subseteq \X_n$, so there exists some $s \in Z'$ with $|s| = n$ and $s^\frown 0 \in Z'$. Since $B' \leq B$, we have that $Z' \subseteq Z$, so $s$ is a node of length $n$ in $Z$. Hence $s = t$, and $t^\frown 0 \in Z' \subseteq Z$, contradicting $t^\frown 0 \notin Z$.
        
        \item Now suppose instead that $[A_{n+1}] \subseteq (f^{-1}(\X_n))^c$, so that $f``[A_{n+1}] \subseteq \X_n^c$. Let $Z \leq Y$ be such that $t^\frown 0$ is the only successor of $t$ in $Z$, so that again the only node of length $n$ in $Z$ is $t$. As before, $Z = f(B)$ for some $B \leq A^*$, and there exists some $B' \leq B$ with $B' \leq A_{n+1}$. Let $Z' := f(B')$. Then $Z' \subseteq Z$ and $Z' \in f``[A_{n+1}] \subseteq \X_n^c$. Since $Z' \in \S$, its unique node $t$ of length $n$ must have a successor in $Z'$. As $Z' \subseteq Z$, this successor is $t^\frown 0$. Hence $t^\frown 0 \in Z'$, so $Z' \in \X_n$, contradicting $Z' \in \X_n^c$.
    \end{itemize}
\end{proof}

The minimality of Sacks forcing may be used to show that it does not weakly project into a large class of topological Ramsey spaces.

\begin{theorem}
\label{thm:no.wp.Sacks.sigma-closed}
    Let $(\cal{R},\leq,\leq^*,r)$ be a $\sigma$-closed closed triple satisfying \textbf{A1}--\textbf{A4} such that for all $X \in \cal{R}$, there exist incompatible $X_0,X_1 \leq X$. Then $(\mathbb{S},\subseteq,r)$ does not weakly project into $(\cal{R},\leq,r)$.
\end{theorem}

\begin{proof}
    Suppose that $f : \mathbb{S} \to \cal{R}$ is a weak projection. Let $M$ be the transitive collapse of a countable elementary submodel of some $H_\theta$ containing $f$ and everything else relevant, and let $G$ be an $\mathbb{S}$-generic filter over $M$. By Lemma \ref{lem:ab.equivalence}, $\U := f``(G)^\uparrow \in M[G]$ is $(\cal{R},\leq)$-generic over $M$. Since $(\cal{R},\leq)$ and $(\cal{R},\leq^*)$ share the same separative quotient, $\U$ generates the same extension as a $(\cal{R},\leq^*)$-generic filter, which adds no reals as $(\cal{R},\leq^*)$ is $\sigma$-closed. Then $M \subsetneq M[\U] \subsetneq M[G]$ (as $M[\U]$ and $M$ share the same reals), contradicting the minimality of Sacks forcing.
\end{proof}

\begin{corollary}
\label{cor:no.wp.strong.subtrees.Sacks}
    $(\mathbb{S},\subseteq,r)$ does not weakly project into $(\S_\infty(U),\subseteq,r)$.
\end{corollary}

\begin{proof}
    Let $S : \S_\infty(U) \to [\N]^\infty$ be the weak projection provided in Example \ref{ex:levels.strong.subtrees}. If $f : \mathbb{S} \to \S_\infty(U)$ is a weak projection, then so is $S \circ f : \mathbb{S} \to [\N]^\infty$, contradicting Theorem \ref{thm:no.wp.Sacks.sigma-closed}.
\end{proof}

Since $(\S_\infty(U),\subseteq)$ is not a minimal forcing (due to Example \ref{ex:levels.strong.subtrees}), we cannot use the same argument to show that $(\S_\infty(U),\subseteq,r)$ does not weakly project into $(\mathbb{S},\subseteq,r)$. However, this fact turns out to be true for $\S_\infty := \S_\infty(2^{<\omega})$:

\begin{theorem}
\label{thm:no.wp.strong.subtrees.Sacks}
    $(\S_\infty,\subseteq,r)$ does not weakly project into $(\mathbb{S},\subseteq,r)$.
\end{theorem}

While Theorem \ref{thm:no.wp.Sacks.sigma-closed} does not hold for strong subtrees (see Example \ref{ex:levels.strong.subtrees}), a similar result is true:

\begin{theorem}
\label{thm:no.wp.strong.subtrees.sigma-closed.projection}
    Let $(\cal{R},\leq,\leq^*,r)$ be a $\sigma$-closed \textbf{A2}-space, and suppose that there is a projection $\pi : \cal{R} \to [\N]^\infty$ such that for all $X \in \cal{R}$, there exist incompatible $X_0,X_1 \leq X$ such that $\pi(X_0) = \pi(X_1)$. Then $(\S_\infty,\subseteq,r)$ does not weakly project into $(\cal{R},\leq,r)$.
\end{theorem}

Consequently:

\begin{corollary}
\label{cor:no.wp.strong.subtrees.sigma-closed.examples}
    $(\S_\infty,\subseteq,r)$ does not weakly project into the following:
    \begin{enumerate}
        \item $(\FIN^{[\infty]},\leq,r)$.
        
        \item $(\FIN_k^{[\infty]},\leq,r)$ for $k \geq 2$.
        
        \item $(\FIN_{\pm k}^{[\infty]},\leq,r)$ for $k \geq 2$.
        
        \item $(W_{Lv}^{[\infty]},\leq,r)$, where $L$ is finite.
        
        \item $(\cal{R}_1,\leq_1,r^1)$.
        
        \item $(\E_\infty,\leq,r)$.
    \end{enumerate}
\end{corollary}

\subsubsection{Forcing properties of $(\S_\infty,\subseteq)$}
To prove these theorems, we need to first study the forcing-theoretic properties of $(\S_\infty,\subseteq)$. The forcing $(\S_\infty,\subseteq)$ is also known as \emph{uniform Sacks forcing}, and ideals related to $(\S_\infty,\subseteq)$ were previously studied in \cite{R93} and \cite{B95}. The key lemma is the real-degree minimality of $(\S_\infty,\subseteq)$:

\begin{lemma}
\label{lem:strong.subtrees.forcing.minimal.real}
    Let $M$ be a ctm of $\ZFC$, and let $G$ be $\S_\infty$-generic over $M$. If $x \in 2^\omega \cap (M[G] \setminus M)$, then $M[x] = M[G]$.
\end{lemma}

Theorem \ref{thm:no.wp.strong.subtrees.Sacks} then quickly follows from the minimality of Sacks forcing. However, Theorem \ref{thm:no.wp.strong.subtrees.sigma-closed.projection} requires deeper analysis of the forcing $(\S_\infty,\subseteq)$, which we shall discuss after the proof of Theorem \ref{thm:no.wp.strong.subtrees.Sacks}.

\begin{notation}
\label{not:tree.restriction}
    Let $A \subseteq \omega^{<\omega}$ be a tree.
    \begin{enumerate}
        \item We let $\func{br}(A)$ be the set of all branches of $A$.

        \item For each $s \in A$, we write:
        \begin{align*}
            A\restrictedto s := \{t \in A : t \sqsubseteq s \vee s \sqsubseteq t\}.
        \end{align*}
    \end{enumerate}
\end{notation}

\begin{lemma}
\label{lem:tree.generic.real}
    Let $M$ be a ctm of $\ZFC$, and let $(\P,\subseteq) \in M$ be a partial order consisting of non-empty trees $A \subseteq \omega^{<\omega}$ without terminal nodes such that for all $A \in \P$ and $s \in A$, $A\restrictedto s \in \P$. Let $G$ be $\P$-generic over $M$, and let:
    \begin{align*}
        x_G := \bigcup\{s \in \omega^{<\omega} : \forall A \in G \, s \in A\}.
    \end{align*}
    Then $x_G \in \omega^\omega \cap M[G]$ and for all $A \in \P$, $A \in G$ iff $x_G \in \func{br}(A)$, so $M[G] = M[x_G]$.
\end{lemma}

Note that for all $A \in \S_\infty$ and $s \in A$, $A\restrictedto s \in \S_\infty$, so we may apply Lemma \ref{lem:tree.generic.real} to $(\S_\infty,\subseteq)$.

\begin{proof}
    For each $n$, define:
    \begin{align*}
        \D_n := \{A \in \P : \text{$A$ has exactly one $s \in A$ with $\dom(s) = n$}\}.
    \end{align*}
    It's clear that $\D_n$ is open, and it is dense as it contains $A\restrictedto t$ for all $A \in \P$ and $t \in A$ with $\dom(t) \geq n$. Since $G$ is $\P$-generic, we may let $A_n \in G \cap \D_n$, and let $t_n \in A_n$ be the unique node of length $n$. Since $G$ is a filter and $\D_{n+1}$ is open, we may assume that $A_{n+1} \subseteq A_n$ for all $n$, so $t_n \sqsubseteq t_{n+1}$ and hence $\bigcup_{n<\omega} t_n \in \omega^\omega$. We shall show that:
    \begin{align*}
        x_G = \bigcup_{n<\omega} t_n.
    \end{align*}
    If $A \in G$, then let $B \in G$ such that $B \subseteq A \cap A_n$. Then $B \in \D_n$, so $t_n \in B \subseteq A$ as $B \subseteq A_n$ and $B$ contains a unique node of length $n$. Conversely, if $s \in \omega^{<\omega}$ is of length $n$ and $s \in A$ for all $A \in G$, then $s \in A_n$, so $s = t_n$ is the unique node in $A_n$ of length $n$.

    Now let $A \in \P$. If $A \in G$, then $x_G \in \func{br}(A)$ by the definition of $x_G$. Conversely, suppose that $x_G \in \func{br}(A)$. Let:
    \begin{align*}
        \cal{K} := \{B \in \P : \text{$B \subseteq A$ or $B \cap A$ is finite}\}.
    \end{align*}
    It's clear that $\cal{K}$ is open. Observe that $\cal{K}$ is dense: if $B \in \P$ and $B \nsubseteq A$, then let $s \in B \setminus A$, and we have that $B\restrictedto s \subseteq B$ and every element of $(B\restrictedto s) \cap A$ is an initial segment of $s$, so $(B\restrictedto s) \cap A$ is finite. Since $G$ is $\P$-generic over $M$, we may let $B \in G \cap \cal{K}$. Since $x_G \in \func{br}(B) \cap \func{br}(A)$, every initial segment of $x_G$ lies in $B \cap A$, so $B \cap A$ is infinite. Therefore, $B \subseteq A$, so $A \in G$.
\end{proof}

Let $\U$ be an ultrafilter on $\N$. We let $\S_\U := \S_\U(2^{<\omega})$ be the set of all strong binary subtrees $A \subseteq 2^{<\omega}$ such that $S(A) \in \U$ (where $S(A)$ is the set of splitting levels of $A$, see Example \ref{ex:levels.strong.subtrees}). Letting $\dot{\U}$ be the canonical name for a $[\N]^\infty$-generic filter, we consider the two-step iterated forcing $[\N]^\infty * \S_{\dot{\U}}$ under the usual iteration order.

\begin{lemma}
\label{lem:strong.subtrees.forcing.product}
    $(\S_\infty,\subseteq)$ is forcing equivalent to $([\N]^\infty * \S_{\dot{\U}},\leq)$.
\end{lemma}

\begin{proof}
    Recall the function $S : \S_\infty \to [\N]^\infty$ sending a strong subtree to the set of splitting levels (see Example \ref{ex:levels.strong.subtrees}). We define a function $h : \S_\infty \to [\N]^\infty * \S_{\dot{\U}}$ by:
    \begin{align*}
        h(A) := (S(A),\check{A}),
    \end{align*}
    where $\check{A}$ is the canonical $[\N]^\infty$-name for the strong subtree $A$. Observe that $h$ is well-defined, as if $\dot{\U}$ is the canonical name for the $([\N]^\infty,\subseteq)$-generic ultrafilter, then $S(A) \forces S(A) \in \dot{\U}$, so $S(A) \forces \check{A} \in \S_{\dot{\U}}$. We shall show that $h$ is a dense embedding.
    \begin{itemize}
        \item Let $A,B \in \S_\infty$. Since inclusion of subtrees is an absolute statement, we have that:
        \begin{align*}
            B \subseteq A &\iff S(B) \subseteq S(A) \wedge B \subseteq A \\
            &\iff S(B) \subseteq S(A) \wedge S(B) \forces \check{B} \subseteq \check{A} \\
            &\iff (S(B),\check{B}) \leq (S(A),\check{A}).
        \end{align*}

        \item Let $(X,\dot{A}) \in [\N]^\infty * \S_{\dot{\U}}$. Since forcing with $([\N]^\infty,\subseteq)$ does not add reals, no new binary subtrees are added, so there exist some $X' \leq X$ and $B \in \S_\infty$ such that $X' \forces \check{B} = \dot{A}$. Since $X' \forces \check{B} \in \S_{\dot{\U}}$, we have that $X' \forces S(\check{B}) \in \dot{\U}$, so $X' \subseteq^* S(B)$. We then let $C$ be any strong subtree of $B$ such that $S(C) = X' \cap S(B)$ (which is infinite as $X' \subseteq^* S(B)$). Then $(S(C),\check{C}) \in [\N]^\infty * \S_{\dot{\U}}$ as $S(C) \forces S(C) \in \dot{\U}$ (so $S(C) \forces \check{C} \in \S_{\dot{\U}}$), and $(S(C),\check{C}) \leq (X',\dot{A}) \leq (X,\dot{A})$. Therefore, the image under $h$ is dense.
    \end{itemize}
\end{proof}

\begin{lemma}
\label{lem:strong.subtrees.forcing.minimal}
    Let $M$ be a ctm of $\ZFC$. Let $G$ be a $([\N]^\infty * \S_{\dot{\U}})$-generic filter over $M$, and let:
    \begin{align*}
        \U := \{Y \in [\N]^\infty : (Y,\dot{A}) \in G \text{ for some } \dot{A}\}
    \end{align*}
    be the projection of $G$ to the first coordinate, so $\U$ is a $[\N]^\infty$-generic filter over $M$ (cf. Lemma \ref{lem:strong.subtrees.forcing.product}). If $X \in M[G] \setminus M[\U]$ is a set of ordinals, then $M[G] = M[X]$.
\end{lemma}

\begin{proof}
    The proof of this lemma is similar to the proof of the minimality of Sacks forcing, but some modifications are required. Let $\P := [\N]^\infty * \S_{\dot{\U}}$, and let $\dot{X}$ be a $\P$-name for $X$. For notational convenience, we shall suppress checkmarks in names (so we write $(Y,A) \in \P$ instead of $(Y,\check{A}) \in \P$ for a strong subtree $A$). Since $X \notin M[\U]$, there exists some $(\bar{Y},\bar{A}) \in G$, with $\bar{A} \in \S_\infty$, such that $(\bar{Y},\bar{A}) \forces \dot{X} \notin M[\dot{\U}]$, where $\dot{\U}$ is the canonical name for the $[\N]^\infty$-generic filter.
    
    \begin{claim}
    \label{claim:ordinal-determining.A*}
        In $M$, for any $(Y_0,A_0) \in \P$ with $(Y_0,A_0) \leq (\bar{Y},\bar{A})$, there exist some $(Y^*,A^*) \leq (Y_0,A_0)$ in $\P$ and a set of ordinals $\{\gamma_s : s \in \func{split}(A^*)\}$ such that for all $s \in \func{split}(A^*)$, $(Y^*,A^*\restrictedto(s^\frown 0))$ and $(Y^*,A^*\restrictedto(s^\frown 1))$ force different truth values of the statement $\gamma_s \in \dot{X}$.
    \end{claim}
    
    \begin{midproof}
        We shall define $A^*$ as the limit of a fusion sequence $(A_n)_{n<\omega}$ of strong subtrees, starting with the given $A_0$. We also define in parallel a $\subseteq$-decreasing sequence $(Y_n)_{n<\omega}$ of infinite subsets of $\N$ such that $Y_n \subseteq^* S(A_n)$ for all $n$, starting with the given $Y_0$, and let $Y^* := S(A^*) \cap Y_0$ so that $Y^*$ is infinite and $Y^* \subseteq^* Y_n$ for all $n$. Assume that $A_n$ and $Y_n$ have been defined and $Y_n \subseteq^* S(A_n)$, and let $\{s_i : i < N\}$ enumerate the set $\func{split}_n(A_n)$, i.e. the set of splitting nodes of $A_n$ in the $n^\text{th}$ level. We shall define a $\subseteq$-decreasing sequence $(A_n^i)_{i \leq N}$ in $[n,A_n]$, and a $\subseteq$-decreasing sequence $(Y_n^i)_{i \leq N}$ of subsets of $Y_n$, then complete the induction by letting $A_{n+1} := A_n^N$ and $Y_{n+1} := Y_n^N$.

        Let $A_n^0 := A_n$ and $Y_n^0 := Y_n$, and suppose that $A_n^i$ and $Y_n^i$ have been defined. We first show that the following statement holds:
        \begin{enumerate}
            \item[$(\star)$] There are some ordinal $\gamma_{s_i}$, some infinite $Y_n^{i+1} \subseteq Y_n^i$ and strong subtrees $B_n^i \subseteq A_n^i\restrictedto({s_i}^\frown 0)$ and $C_n^i \subseteq A_n^i\restrictedto({s_i}^\frown 1)$ such that $Y_n^{i+1} \subseteq S(B_n^i) \cap S(C_n^i)$, and $(Y_n^{i+1},B_n^i),(Y_n^{i+1},C_n^i)$ decide $\gamma_{s_i} \in \dot{X}$ oppositely.
        \end{enumerate}
        If $(\star)$ holds, then we can complete the induction by letting $A_n^{i+1} \in [n,A_n^i]$ be such that $A_n^{i+1}\restrictedto({s_i}^\frown 0) \subseteq B_n^i$, $A_n^{i+1}\restrictedto({s_i}^\frown 1) \subseteq C_n^i$, and the splitting levels of $A_n^{i+1}$ above its $n^\text{th}$ splitting level are precisely the elements of $Y_n^{i+1}$, which is possible as $Y_n^{i+1} \subseteq S(B_n^i) \cap S(C_n^i)$ is infinite. Note that ${s_i}^\frown 0,{s_i}^\frown 1 \in A_n^{i+1}$, so every $s_i$ remains a splitting node, and $\func{split}_n(A_{n+1}) = \func{split}_n(A_n)$. We see that $(Y^*,A^*)$ works --- if $s \in \func{split}(A^*)$, then $s \in \func{split}_n(A^*) = \func{split}_n(A_n)$ for some $n$, so $s = s_i$ for some $i$ in the $n^\text{th}$ induction step. Since $(Y_n^{i+1},B_n^i)$ and $(Y_n^{i+1},C_n^i)$ decide $\gamma_{s_i} \in \dot{X}$ differently, $Y^* \subseteq^* Y_n^{i+1}$, $A^*\restrictedto({s_i}^\frown 0) \subseteq A_n^{i+1}\restrictedto({s_i}^\frown 0) \subseteq B_n^i$ and $A^*\restrictedto({s_i}^\frown 1) \subseteq A_n^{i+1}\restrictedto({s_i}^\frown 1) \subseteq C_n^i$, we have that for all infinite $Z \subseteq Y^*$ such that $Z \subseteq Y_n^{i+1}$ (of which there are densely many below $Y^*$), $(Z,A^*\restrictedto({s_i}^\frown 0))$ and $(Z,A^*\restrictedto({s_i}^\frown 1))$ decide $\gamma_{s_i} \in \dot{X}$ differently. Therefore, no extension of $(Y^*,A^*\restrictedto({s_i}^\frown 0))$ or of $(Y^*,A^*\restrictedto({s_i}^\frown 1))$ forces the opposite truth value, so $(Y^*,A^*\restrictedto({s_i}^\frown 0))$ and $(Y^*,A^*\restrictedto({s_i}^\frown 1))$ also decide $\gamma_{s_i} \in \dot{X}$ differently, as needed.

        It remains to show that $(\star)$ holds. Suppose for a contradiction that $(\star)$ fails. For each ordinal $\gamma$ we let:
        \begin{align*}
            \D_\gamma := \{Z \subseteq Y_n^i : (Z,A_n^i\restrictedto s_i) \text{ decides } \gamma \in \dot{X}\}.
        \end{align*}
        We show that $\D_\gamma$ is dense below $Y_n^i$. Fix any $Y \subseteq Y_n^i$, and let $(Z,B) \leq (Y,A_n^i\restrictedto({s_i}^\frown0))$ be any condition that decides $\gamma \in \dot{X}$. We shall show that for all $(W,C) \leq (Z,A_n^i\restrictedto s_i)$ that decides $\gamma \in \dot{X}$, $(W,C)$ decides $\gamma \in \dot{X}$ in the same way as $(W,B)$ (as $(W,B) \leq (Z,B)$, so $(W,B)$ also decides $\gamma \in \dot{X}$ in the same way as $(Z,B)$). This implies that $(Z,A_n^i\restrictedto s_i)$ decides $\gamma \in \dot{X}$, so $Z \subseteq Y$ and $Z \in \D_\gamma$, completing the proof. We fix any $(W,C) \leq (Z,A_n^i\restrictedto s_i)$ that decides $\gamma \in \dot{X}$, and consider two cases.
        \begin{itemize}
            \item Suppose that ${s_i}^\frown1 \in C$. Let $C' := C\restrictedto({s_i}^\frown1)$. Then $(W,C') \leq (W,C)$, and $(W,C')$ decides $\gamma \in \dot{X}$ in the same way as $(W,C)$. Let $W' := W \cap S(B) \cap S(C')$, which is cofinite in $W$, as $W \subseteq^* S(B)$ and $W \subseteq^* S(C')$. Then $B \subseteq A_n^i\restrictedto({s_i}^\frown0)$, $C' \subseteq A_n^i\restrictedto({s_i}^\frown1)$ and $W' \subseteq S(B) \cap S(C')$, so by the negation of $(\star)$, $(W',B)$ and $(W',C')$ decide $\gamma \in \dot{X}$ in the same way. Since $(W',B) \leq (W,B)$ and $(W',C') \leq (W,C')$, it follows that $(W,B)$ and $(W,C)$ decide $\gamma \in \dot{X}$ in the same way.

            \item Suppose that $C \subseteq A_n^i\restrictedto({s_i}^\frown0)$. Let $(W',D) \leq (W,A_n^i\restrictedto({s_i}^\frown1))$ be a condition that decides $\gamma \in \dot{X}$. Applying the argument in the previous case twice, we have that $(W',C)$ decides $\gamma \in \dot{X}$ in the same way as $(W',D)$, which decides $\gamma \in \dot{X}$ in the same way as $(W',B)$. Since $(W',C) \leq (W,C)$ and $(W',B) \leq (W,B)$, it follows that $(W,C)$ and $(W,B)$ decide $\gamma \in \dot{X}$ in the same way.
        \end{itemize}
        Therefore, for every ordinal $\gamma$, $\D_\gamma$ is a dense open subset of $[\N]^\infty$ below $Y_n^i$, and $\D_\gamma \in M$. This implies that:
        \begin{align*}
            (Y_n^i,A_n^i\restrictedto s_i) \forces \dot{X} = \{\gamma \in \mathrm{Ord} : \exists Y \in \dot{\U} \, (Y \subseteq Y_n^i \wedge (Y,A_n^i\restrictedto s_i) \forces \gamma \in \dot{X})\},
        \end{align*}
        which implies that $(Y_n^i,A_n^i\restrictedto s_i) \forces \dot{X} \in M[\dot{\U}]$. This is a contradiction, as $(Y_n^i,A_n^i\restrictedto s_i) \leq (\bar{Y},\bar{A})$, and $(\bar{Y},\bar{A}) \forces \dot{X} \notin M[\dot{\U}]$.
    \end{midproof}

    By Claim \ref{claim:ordinal-determining.A*}, the set of conditions $(Y^*,A^*) \in \P$ satisfying the properties of the claim is dense below $(\bar{Y},\bar{A})$, and this set lies in $M$. Since $(\bar{Y},\bar{A}) \in G$ and $G$ is $\P$-generic over $M$, we may fix some such $(Y^*,A^*) \in G$. We define $s_0 \sqsubseteq s_1 \sqsubseteq \cdots$ in $\func{split}(A^*)$ as follows: Let $s_0$ be the $\sqsubseteq$-minimal splitting node in $A^*$. If $s_n \in \func{split}(A^*)$ has been defined, then let $s_{n+1}$ be the $\sqsubseteq$-minimal splitting node in $A^*$ extending ${s_n}^\frown{i_n}$, where $i_n < 2$ is such that $\gamma_{s_n} \in X$ iff $(Y^*,A^*\restrictedto({s_n}^\frown i_n)) \forces \gamma_{s_n} \in \dot{X}$. Now let $x_G^* := \bigcup_{n<\omega} s_n$. Note that $A^*$ and the function $s \mapsto \gamma_s$ lie in $M$, as Claim \ref{claim:ordinal-determining.A*} is proved in $M$. Therefore, $x_G^* \in M[X]$. Finally, recall from Lemmas \ref{lem:strong.subtrees.forcing.product} and \ref{lem:tree.generic.real} the generic real $x_G \in M[G]$, which satisfies that for all $A \in \S_\infty \cap M$, $(S(A),A) \in G$ iff $x_G$ is a branch of $A$.
    
    \begin{claim}
        In $M[G]$, $x_G = x_G^*$.
    \end{claim}

    \begin{midproof}
        Suppose otherwise. Let $n$ be the least integer such that $s_{n+1} \not\sqsubseteq x_G$ (note that $s_0 \sqsubseteq x_G$ --- since $(Y^*,A^*) \leq (S(A^*),A^*)$, we have $(S(A^*),A^*) \in G$, so $x_G$ is a branch of $A^*$). We may as well assume that $i_n = 0$ and $(Y^*,A^*\restrictedto({s_n}^\frown 0)) \forces \gamma_{s_n} \in \dot{X}$, so that $\gamma_{s_n} \in X$ in $M[G]$. Since there are no other splitting nodes between $s_n$ and $s_{n+1}$, and $s_{n+1} \not\sqsubseteq x_G$, this implies that ${s_n}^\frown1 \sqsubseteq x_G$. Therefore, $x_G$ is a branch of $A^*\restrictedto({s_n}^\frown1)$, so $(S(A^*\restrictedto({s_n}^\frown1)),A^*\restrictedto({s_n}^\frown1)) \in G$. Since $G$ is a filter also containing $(Y^*,A^*)$, we have that $(Y^*,A^*\restrictedto({s_n}^\frown1)) \in G$. But $(Y^*,A^*\restrictedto({s_n}^\frown1)) \forces \gamma_{s_n} \notin \dot{X}$, contradicting that $\gamma_{s_n} \in X$.
    \end{midproof}

    Therefore, we may define $G$ in $M[X]$ as the upward closure of the set of all conditions $(S(A),A)$, where $A \in \S_\infty \cap M$ and $x_G^*$ is a branch of $A$, completing the proof.
\end{proof}

In particular, applying Lemma \ref{lem:strong.subtrees.forcing.minimal} to a real, we can prove Lemma \ref{lem:strong.subtrees.forcing.minimal.real}:

\begin{proof}[Proof of Lemma \ref{lem:strong.subtrees.forcing.minimal.real}]
    Let $\U := S``(G)$, and note that $\U$ is $[\N]^\infty$-generic over $M$. Since no reals are added in $M[\U]$, we have that $x \notin M[\U]$. By Lemma \ref{lem:strong.subtrees.forcing.minimal} (and that $(\S_\infty,\subseteq)$ and $([\N]^\infty * \S_{\dot{\U}},\leq)$ are forcing equivalent by Lemma \ref{lem:strong.subtrees.forcing.product}), we have that $M[x] = M[G]$. 
\end{proof}

We are now ready to prove Theorem \ref{thm:no.wp.strong.subtrees.Sacks}.

\begin{proof}[Proof of Theorem \ref{thm:no.wp.strong.subtrees.Sacks}]
    Suppose that there is a weak projection $f : \S_\infty \to \mathbb{S}$. Let $M$ be the transitive collapse of a countable elementary submodel of some $H_\theta$ containing $f$ and everything else relevant. Let $G$ be $\S_\infty$-generic over $M$. By Lemma \ref{lem:ab.equivalence}, there exists a Sacks real $x \in M[G]$ over $M$, and by Lemma \ref{lem:strong.subtrees.forcing.minimal.real}, we have that $M[x] = M[G]$. But since $\U := S``(G)$ is $[\N]^\infty$-generic over $M$, we have that $M \subsetneq M[\U] \subsetneq M[x]$, contradicting the minimality of Sacks forcing.
\end{proof}

\subsubsection{Proof of Theorem \ref{thm:no.wp.strong.subtrees.sigma-closed.projection}}
We dedicate the rest of the subsection to proving Theorem \ref{thm:no.wp.strong.subtrees.sigma-closed.projection}. We first prove a combinatorial forcing lemma regarding continuous functions $\Phi : [a,A] \to \N$ for some closed triple $(\cal{R},\leq,r)$ satisfying \textbf{A1}--\textbf{A4}. Given $a,b \in \cal{AR}$, we define:
\begin{align*}
    b \preceq_\fin a \iff b \leq_\fin a' \text{ for some } a' \sqsubseteq a.
\end{align*}
Note that $\preceq_\fin$ is a quasi-order on $\cal{AR}$, and every $a \in \cal{AR}$ has only finitely many $b \preceq_\fin a$.

\begin{definition}
\label{def:amalgamable.cover}
    Let $A \in \cal{R}$, $a \sqsubseteq A$ and $b \preceq_\fin a$. We say that a family $\{H_i\}_{i \in I} \subseteq \cal{R}$ is an \emph{amalgamable cover for $b$ in $[a,A]$} if the following statements hold:
    \begin{enumerate}
        \item For all $B \in [b,A]$ such that $\depth_A(r_{\lh(b)+1}(B)) > \lh(a)$, we have that $B \leq H_i$ for some $i \in I$.

        \item If $A' \in [a,A]$, $B \in [b,A']$ and $B \leq H_i$, then there exists some $B' \in [a,A']$ such that $[b,B'] \cap [H_i] \subseteq [b,B]$.
    \end{enumerate}
\end{definition}

\begin{definition}
\label{def:finite.amalgamable.covers}
    We say that $(\cal{R},\leq,r)$ has \emph{finite amalgamable covers} if for all $A \in \cal{R}$, $a \sqsubseteq A$ and $b \preceq_\fin a$, there is a finite amalgamable cover for $b$ in $[a,A]$.
\end{definition}

Note that if $(\cal{R},\leq,r)$ has finite amalgamable covers, then $\lh(a) \leq \depth_A(a)$ for all $A \in \cal{R}$ and $a \in \cal{AR}\restrictedto A$, with equality if $a \sqsubseteq A$. Several examples in \Sec\ref{sec:examples} have finite amalgamable covers --- we discuss a few examples. 

\begin{example}
\label{ex:Ellentuck.has.FACs}
    The Ellentuck space $([\N]^\infty,\subseteq,r)$ has finite amalgamable covers.
\end{example}

\begin{proof}
    Let $A \in [\N]^\infty$, $a \sqsubseteq A$ and $b \subseteq a$ (note that $b \preceq_\fin a$ iff $b \subseteq a$). Define:
    \begin{align*}
        H := b \cup (\N \setminus (\max(a) + 1)),
    \end{align*}
    where $\max(\emptyset) := -1$. We shall show that $\{H\}$ is an amalgamable cover for $b$ in $[a,A]$.
    \begin{enumerate}
        \item Let $B \in [b,A]$ such that $\depth_A(r_{\lh(b)+1}(B)) > \lh(a)$. Then $\min(B \setminus b) > \max(a)$, so $B \setminus b \subseteq \N \setminus (\max(a) + 1)$. Thus, $B = b \cup (B \setminus b) \subseteq H$.
        
        \item Let $A' \in [a,A]$ and $B \in [b,A']$ such that $B \subseteq H$. Again, $\min(B \setminus b) > \max(a)$, so we may let $B' := a \cup (B \setminus b)$, and we have that $B' \in [a,A']$ and $[b,B'] \cap [H] = [b,B]$.
    \end{enumerate}
\end{proof}

\begin{example}
\label{ex:strong.subtrees.has.FACs}
    The space of strong binary subtrees $(\S_\infty,\subseteq,r)$ has finite amalgamable covers.
\end{example}

\begin{proof}
    Let $A \in \S_\infty$, $a \sqsubseteq A$ and $b \preceq_\fin a$. Let:
    \begin{align*}
        \F &:= \{c \in r_{\lh(b)+1}[b,A] : \depth_A(c) = \lh(a) + 1\}.
    \end{align*}
    Note that $\F$ is finite. For each $c \in \F$, let:
    \begin{align*}
        A\restrictedto c := \bigcup\{A\restrictedto t : t \text{ is a terminal node of } c\}.
    \end{align*}
    Note that $A\restrictedto c \in \S_\infty$ for all $c \in \F$. We shall show that $\{A\restrictedto c\}_{c \in \F}$ is an amalgamable cover for $b$ in $[a,A]$. We fix $\ell := \max\{\dom(s) : s \in r_{\lh(a)+1}(A)\}$.
    \begin{enumerate}
        \item Let $B \in [b,A]$ be such that $\depth_A(r_{\lh(b)+1}(B)) > \lh(a)$. Let $c := \{s \in B : \dom(s) \leq \ell\}$. Since $\depth_A(r_{\lh(b)+1}(B)) > \lh(a)$, the terminal nodes of $r_{\lh(b)+1}(B)$ have length at least $\ell$, so $c \in \F$ and $B \subseteq A\restrictedto c$.
        
        \item Let $A' \in [a,A]$ and $B \in [b,A']$ such that $B \subseteq A\restrictedto c$ for some $c \in \F$. Since $b \sqsubseteq B$, $B \subseteq A\restrictedto c$ and $\lh(c) = \lh(b) + 1$, we have that $c \subseteq B$. Let $B' \in [a,A']$ be such that every $s \in B'$ extending a terminal node of $c$ lies in $B$. Then $B' \cap A\restrictedto c \subseteq B$, so $[b,B'] \cap [A\restrictedto c] \subseteq [b,B]$.
    \end{enumerate}
\end{proof}

We remark that the Hindman space $(\FIN^{[\infty]},\leq,r)$ does not have finite amalgamable covers --- see the end of this subsection for a discussion. On the other hand, we shall show, in the next subsection, that the Carlson--Simpson space $(\E_\infty,\leq,r)$ has finite amalgamable covers (Example \ref{ex:Carlson-Simpson.has.FACs}), which is needed in the proofs of Theorems \ref{thm:no.wp.Carlson-Simpson.and.FIN} and \ref{thm:no.wp.CS-variables.FIN}. 

Let $(\cal{R},\leq,r)$ be a closed triple satisfying \textbf{A1}--\textbf{A4}, let $A \in \cal{R}$ and $a_0 \in \cal{AR}\restrictedto A$, and let $\Phi : [a_0,A] \to \N$ be a function. For $B \in [a_0,A]$ and $a \in \cal{AR}\restrictedto[a_0,B]$, we say that $B$ \emph{accepts} $a$ if $\Phi$ is constant on $[a,B]$, that $B$ \emph{rejects} $a$ if no $B' \in [a,B]$ accepts $a$, and that $B$ \emph{decides} $a$ if $B$ accepts or rejects $a$. Note that property (7) of Theorem \ref{thm:abstract.combinatorial.forcing} fails for this notion of acceptance, as $B$ may accept every $b \in r_{\lh(a)+1}[a,B]$ without accepting $a$, so Theorem \ref{thm:abstract.combinatorial.forcing} does not apply.

\begin{lemma}
\label{lem:decision}
    Let $(\cal{R},\leq,r)$ be a closed triple satisfying \textbf{A1}--\textbf{A4} that has finite amalgamable covers. Let $A \in \cal{R}$ and $a_0 \in \cal{AR}\restrictedto A$, and let $\Phi : [a_0,A] \to \N$ be continuous (where $\N$ is equipped with the discrete topology). Then there exists some $B \in [a_0,A]$ such that the following statements hold:
    \begin{enumerate}
        \item $B$ decides every $a \in \cal{AR}\restrictedto[a_0,B]$.

        \item If $a$ and $a'$ are $\sqsubseteq$-minimal among the elements of $\cal{AR}\restrictedto[a_0,B]$ accepted by $B$, and $\Phi$ takes the same value on $[a,B]$ and $[a',B]$, then $\depth_B(a) = \depth_B(a')$.
    \end{enumerate}
\end{lemma}

For the case where $(\cal{R},\leq,r)$ is the Ellentuck space, Lemma \ref{lem:decision} follows from the Pudl\'{a}k--R\"{o}dl theorem on irreducible mappings (see \cite[Theorem 2]{PR82} and \cite[Theorem 1.57]{T10}). Readers may refer to \cite{PSV86} and \cite{V14} for canonisation results of a similar kind for functions on the Carlson--Simpson space and on strong subtrees respectively.

\begin{proof}
    For each $a \in \cal{AR}\restrictedto[a_0,A]$, the set of all $B \in [a,A]$ that decide $a$ is dense open in $[a,A]$, so by Lemma \ref{lem:R.is.semiselective}, there exists some $A_1 \in [a_0,A]$ that decides every $a \in \cal{AR}\restrictedto[a_0,A_1]$. We extend $\Phi$ to all $a \in \cal{AR}\restrictedto[a_0,A_1]$ accepted by $A_1$, by letting $\Phi(a)$ be the constant value of $\Phi$ on $[a,A_1]$. For each $n \in \N$ and $a,b \in \cal{AR}\restrictedto[a_0,A_1]$ such that $b \preceq_\fin a$ and $A_1$ rejects $b$, we define:
    \begin{small}
        \begin{align*}
            \D_{a,b,n} := \{C \in [a,A_1] : \forall C' \in [b,C](\depth_C(r_{\lh(b)+1}(C')) > \lh(a) \to \Phi(C') \neq n)\}.
        \end{align*}
    \end{small}

    \begin{claim}
    \label{claim:deep.reduction.avoids.constant}
        $\D_{a,b,n}$ is dense open in $[a,A_1]$. 
    \end{claim}    

    \begin{midproof}
        If $C \in [a,A_1]$, $C'' \in [a,C]$ and $e \in \cal{AR}\restrictedto C''$, then $\depth_{C''}(e) > \lh(a)$ iff $\depth_C(e) > \lh(a)$, so $\D_{a,b,n}$ is open. Let $B \in [a,A_1]$, and let $\{H_i\}_{i<k}$ be a finite amalgamable cover for $b$ in $[a,B]$. We define a $\leq$-decreasing sequence $C_0 \geq C_1 \geq \cdots \geq C_k$ in $[a,B]$ such that $\Phi(C') \neq n$ for all $i < k$ and $C' \in [b,C_{i+1}] \cap [H_i]$ as follows: Let $C_0 := B$, and suppose that $C_i$ has been defined. If $[b,C_i] \cap [H_i] = \emptyset$, then we let $C_{i+1} := C_i$. Otherwise, let $C_i' \in [b,C_i]$ be such that $C_i' \leq H_i$. Since $\Phi$ is continuous, the set $\X := \{E \in [b,C_i'] : \Phi(E) = n\}$ is closed, hence Ramsey, so there exists some $D_i \in [b,C_i']$ such that $[b,D_i] \subseteq \X$ or $[b,D_i] \subseteq \X^c$. The former is impossible as $A_1$ rejects $b$, so $\Phi(E) \neq n$ for all $E \in [b,D_i]$. Since $D_i \leq H_i$, there exists some $C_{i+1} \in [a,C_i]$ such that $[b,C_{i+1}] \cap [H_i] \subseteq [b,D_i]$. Then $C_k \in \D_{a,b,n}$: if $C' \in [b,C_k]$ and $\depth_{C_k}(r_{\lh(b)+1}(C')) > \lh(a)$, then $\depth_B(r_{\lh(b)+1}(C')) > \lh(a)$, so $C' \leq H_i$ for some $i < k$, and hence $C' \in [b,C_k] \cap [H_i] \subseteq [b,C_{i+1}] \cap [H_i]$.
    \end{midproof}
    
    For each $a \in \cal{AR}\restrictedto[a_0,A_1]$, let $\D_a$ be the intersection of all $\D_{a,b,\Phi(c)}$ such that $b,c \in \cal{AR}\restrictedto[a_0,A_1]$, $b,c \preceq_\fin a$, $A_1$ rejects $b$ and $A_1$ accepts $c$. By Claim \ref{claim:deep.reduction.avoids.constant}, $\D_a$ is the intersection of finitely many dense open subsets of $[a,A_1]$, hence is dense open. By Lemma \ref{lem:R.is.semiselective}, there exists some $B \in [a_0,A_1]$ that diagonalises $\{\D_b\}_{b\in\cal{AR}\restrictedto[a_0,A_1]}$. Then $B$ decides every $a \in \cal{AR}\restrictedto[a_0,B]$ in the same way $A_1$ does, so (1) holds. It remains to show that $B$ satisfies (2).

    Let $a,a'$ be $\sqsubseteq$-minimal among elements in $\cal{AR}\restrictedto[a_0,B]$ accepted by $B$, and suppose that $\Phi(a) = \Phi(a')$. Suppose for a contradiction that $\depth_B(a) < \depth_B(a')$. Let $b := r_{\depth_B(a)}(B)$. Since $\depth_B(a) \geq \depth_B(a_0) = \lh(a_0)$, we have that $a_0 \sqsubseteq b$, and $a \preceq_\fin b$. Let $c \sqsubseteq a'$ be $\sqsubseteq$-maximal such that $\depth_B(c) \leq \depth_B(a)$. Then $a_0 \sqsubseteq c$, $c \neq a'$ and $c \preceq_\fin b$, so $A_1$ rejects $c$ as $a'$ is $\sqsubseteq$-minimal among elements accepted by $B$. Let $C \in \D_b$ be such that $[b,B] \subseteq [b,C]$, and let $C' \in [a',B]$. Since $B$ accepts $a'$, $\Phi$ is constant on $[a',B]$ and $\Phi(C') = \Phi(a') = \Phi(a)$. However, $C' \in [c,C]$ and $\depth_C(r_{\lh(c)+1}(C')) > \lh(b)$, contradicting that $C \in \D_b \subseteq \D_{b,c,\Phi(a)}$.
\end{proof}

We define $(r_{k+})_n(X) := r_{k+n}(X)$ for all $X \in \cal{R}$ and $k,n \in \N$. Note that if $(\cal{R},\leq,r)$ is a closed triple satisfying \textbf{A1}--\textbf{A4}, then so is $([a,A],\leq,r_{\lh(a)+})$ for all $A \in \cal{R}$ and $a \in \cal{AR}\restrictedto A$.

\begin{lemma}
\label{lem:wp.Ellentuck.controlled}
    Let $(\cal{R},\leq,r)$ be a closed triple satisfying \textbf{A1}--\textbf{A4} that has finite amalgamable covers. Let $A \in \cal{R}$ and $a_0 \in \cal{AR}\restrictedto A$, and let $\mu : [a_0,A] \to [\N]^\infty$ be induced by a $\sqsubseteq$-preserving map $\alpha : \cal{AR}\restrictedto[a_0,A] \to [\N]^{<\infty}$ satisfying the following properties:
    \begin{enumerate}
        \item For all $B \in [a_0,A]$ and $a \in \cal{AR}\restrictedto[a_0,B]$ with $a \neq a_0$, $\alpha(a) \neq \emptyset$ and $\max(\alpha(a)) = \max(\alpha(r_{\depth_B(a)}(B)))$.
        
        \item For all $B \in [a_0,A]$ and $n \in \N$, there exists some $B' \in [a_0,B]$ such that $n \notin \mu(B')$.
    \end{enumerate}
    Let $h : [a_0,A] \to [\N]^\infty$ be a weak projection. Then there exist some $B \in [a_0,A]$ and an injective $\psi : \N \to \N$ such that $\psi``(\mu(C)) \cap h(C) \neq \emptyset$ for all $C \in [a_0,B]$.
\end{lemma}

\begin{proof}
    By Lemma \ref{lem:monotone.automatic.continuity} applied to $([a_0,A],\leq,r_{\lh(a_0)+})$, we may let $A_0 \in [a_0,A]$ be such that $h\restrictedto[a_0,A_0]$ is continuous. Let $\Phi := \min \circ (h\restrictedto[a_0,A_0])$, and note that $\Phi$ is also continuous. By Lemma \ref{lem:decision}, we may obtain some $A_1 \in [a_0,A_0]$ satisfying the conclusion of Lemma \ref{lem:decision}. For each $C \in [a_0,A_1]$, let $k_C$ be the least $k \geq \lh(a_0)$ such that $A_1$ accepts $r_k(C)$, and let $\Lambda(C) := \max(\alpha(r_{k_C}(C)))$.

    \begin{claim}
        $\Lambda$ is well-defined, and $A_1$ rejects $a_0$. Therefore, $\Lambda(C) \in \mu(C)$ for all $C \in [a_0,A_1]$.
    \end{claim}

    \begin{midproof}
        Let $C \in [a_0,A_1]$, and let $m := \Phi(C)$. Since $\Phi$ is continuous and $\N$ is equipped with the discrete topology, $\Phi^{-1}(m)$ is metrically open in $[a_0,A_0]$, so there exists some $k \geq \lh(a_0)$ such that $[r_k(C),C] \subseteq \Phi^{-1}(m)$. Then $\Phi$ is constant on $[r_k(C),C]$, so $C$ accepts $r_k(C)$, and since $A_1$ decides $r_k(C)$, we have that $A_1$ accepts $r_k(C)$ as well.

        If $A_1$ accepts $a_0$, then $\Phi$ is constant on $[a_0,A_1]$ with some value $m$, so $m \in h(C)$ for all $C \in [a_0,A_1]$. This is not possible as $h$ is Alaoglu--Birkhoff convergent. Thus, $A_1$ rejects $a_0$, so for all $C \in [a_0,A_1]$, we have that $r_{k_C}(C) \neq a_0$, hence $\alpha(r_{k_C}(C))$ is non-empty and $\Lambda(C) \in \alpha(r_{k_C}(C)) \subseteq \mu(C)$.
    \end{midproof}

    \begin{claim}
        There exists some $g : \N \to \N$ such that $\Lambda = g \circ \Phi$.
    \end{claim}

    \begin{midproof}
        It suffices to show that if $C,C' \in [a_0,A_1]$ are such that $\Phi(C) = \Phi(C')$, then $\Lambda(C) = \Lambda(C')$. By the minimality of $k_C$ and $k_{C'}$, $r_{k_C}(C)$ and $r_{k_{C'}}(C')$ are $\sqsubseteq$-minimal among the elements of $\cal{AR}\restrictedto[a_0,A_1]$ accepted by $A_1$, and $\Phi$ takes the same value on $[r_{k_C}(C),A_1]$ and $[r_{k_{C'}}(C'),A_1]$, so by Lemma \ref{lem:decision}(2), we have that $\depth_{A_1}(r_{k_C}(C)) = \depth_{A_1}(r_{k_{C'}}(C')) =: j$. Therefore, $\Lambda(C) = \max(\alpha(r_j(A_1))) = \Lambda(C')$.
    \end{midproof}

    Let:
    \begin{align*}
        \D := \{X \in [\N]^\infty : \text{$X$ is coinfinite and $g\restrictedto X$ is constant or injective}\}.
    \end{align*}
    It is clear that $\D$ is a dense open subset of $[\N]^\infty$. Since $h$ is Alaoglu--Birkhoff convergent, $h^{-1}(\D) \cap [a_0,A_1]$ is dense open in $[a_0,A_1]$, so there exists some $B \in [a_0,A_1]$ such that $h(B) \in \D$. Observe that $g\restrictedto h(B)$ is not constant --- otherwise, if $g\restrictedto h(B)$ takes the constant value $m$, then we may let $C \in [a_0,B]$ be such that $m \notin \mu(C)$, contradicting that $\Phi(C) \in h(C) \subseteq h(B)$ and $m = g(\Phi(C)) = \Lambda(C) \in \mu(C)$. Therefore, $g\restrictedto h(B)$ is injective, and since $h(B)$ is coinfinite, there exists some injective $\psi : \N \to \N$ such that $\psi(g(n)) = n$ for all $n \in h(B)$. Then for all $C \in [a_0,B]$, $\psi(\Lambda(C)) = \min(h(C)) \in h(C)$ and $\Lambda(C) \in \mu(C)$, so $\psi(\Lambda(C)) \in \psi``(\mu(C)) \cap h(C)$.
\end{proof}

\begin{lemma}
\label{lem:controlled.wps.give.isomorphic.generic.filters}
    Let $M$ be a ctm of $\ZFC$, and let $(\P,\leq_\P)$ be a poset in $M$. Let $G$ be $\P$-generic over $M$, and let $h_0,h_1 : \P \to [\N]^\infty$ be two weak projections in $M$. Suppose that there exist some $p \in G$ and some injective $\psi : \N \to \N$ in $M$ such that for all $q \leq p$, $\psi``(h_0(q)) \cap h_1(q) \neq \emptyset$. Then $M[h_0``(G)^\uparrow] = M[h_1``(G)^\uparrow]$, and in this model, $\psi$ witnesses that $h_0``(G)^\uparrow \cong h_1``(G)^\uparrow$ (i.e. they are \emph{Rudin--Keisler equivalent}).
\end{lemma}

We defer the definition of the Rudin--Keisler ordering and equivalence to \S\ref{sec:proof.of.CS.and.FIN}.

\begin{proof}
    Let $\U_0 := h_0``(G)^\uparrow$ and $\U_1 := h_1``(G)^\uparrow$. By Lemma \ref{lem:ab.equivalence}, $\U_0$ and $\U_1$ are $[\N]^\infty$-generic over $M$. Since $\psi \in M$ is injective, it suffices to show that:
    \begin{align*}
        \U_1 = \{W \in \Po(\N) \cap M : \psi^{-1}(W) \in \U_0\}.
    \end{align*}
    Suppose otherwise. Since $\U_0$ and $\U_1$ are $[\N]^\infty$-generic over $M$, there exists some $W \in \Po(\N) \cap M$ such that $\psi^{-1}(W) \in \U_0$, but $W \notin \U_1$, which implies that $W^c \in \U_1$. Let $q_0,q_1 \in G$ be such that $h_0(q_0) \subseteq \psi^{-1}(W)$ and $h_1(q_1) \subseteq W^c$. Choose $q \leq p$ in $G$ such that $q \leq q_0$ and $q \leq q_1$, and we have that $\psi``(h_0(q)) \subseteq W$ and $h_1(q) \subseteq W^c$, so $\psi``(h_0(q)) \cap h_1(q) = \emptyset$. This contradicts our choice of $p$ and $\psi$.
\end{proof}

A similar argument was given in the proof of \cite[Theorem 24]{RT12}. We are now ready to prove Theorem \ref{thm:no.wp.strong.subtrees.sigma-closed.projection}.

\begin{proof}[Proof of Theorem \ref{thm:no.wp.strong.subtrees.sigma-closed.projection}]
    We showed in Example \ref{ex:strong.subtrees.has.FACs} that $(\S_\infty,\subseteq,r)$ has finite amalgamable covers. Furthermore, the level set map $S : \S_\infty \to [\N]^\infty$ in Example \ref{ex:levels.strong.subtrees} satisfies the requirements of $\mu$ in Lemma \ref{lem:wp.Ellentuck.controlled}, where for all $a \in \S_{<\infty}$ we may let:
    \begin{align*}
        \alpha(a) := \{\dom(s) : \text{$s \in a$ splits or is terminal in $a$}\}.
    \end{align*}
    We leave the details of the verification of requirements to the reader.

    Suppose for a contradiction that there is a weak projection $f : \S_\infty \to \cal{R}$. Then $h := \pi \circ f : \S_\infty \to [\N]^\infty$ is a weak projection, so by Lemma \ref{lem:wp.Ellentuck.controlled}, there exist some $A \in \S_\infty$ and an injective $\psi : \N \to \N$ such that $\psi``(S(B)) \cap h(B) \neq \emptyset$ for all strong subtrees $B \subseteq A$. Let $M$ be the transitive collapse of a countable elementary submodel of some $H_\theta$ containing $f$, $\pi$, $A$, $\psi$ and everything else relevant, and let $G$ be $\S_\infty$-generic over $M$ such that $A \in G$. We consider the following filters in $M[G]$, which are generic over $M$ by Lemma \ref{lem:ab.equivalence}:
    \begin{align*}
        \U &:= S``(G)^\uparrow \tag*{($[\N]^\infty$-generic)} \\
        \cal{V} &:= f``(G)^\uparrow \tag*{($\cal{R}$-generic)} \\
        \cal{V}^* &:= \pi``(\cal{V})^\uparrow = h``(G)^\uparrow \tag*{($[\N]^\infty$-generic)}
    \end{align*}
    By Lemma \ref{lem:controlled.wps.give.isomorphic.generic.filters}, $M[\U] = M[\cal{V}^*]$. Since $(\cal{R},\leq,\leq^*,r)$ is $\sigma$-closed, $M[\cal{V}] \neq M[G]$, so $\cal{V} \in M[\U] = M[\cal{V}^*]$ by Lemma \ref{lem:strong.subtrees.forcing.minimal}. Therefore, $\cal{V} = i_{\cal{V}^*}(\tau)$ for some $[\N]^\infty$-name $\tau \in M$, so there exists some $X \in \cal{V}$ such that for every $\cal{R}$-generic filter $\cal{V}'$ over $M$ with $X \in \cal{V}'$:
    \begin{align*}
        \cal{V}' = i_{\pi``(\cal{V}')^\uparrow}(\tau).
    \end{align*}
    Thus, there exist some $X_0,X_1 \leq X$ in $M$ that are incompatible in $M$ and satisfy $\pi(X_0) = \pi(X_1)$. Let $\cal{V}_0$ be $\cal{R}$-generic over $M$ such that $X_0 \in \cal{V}_0$, and let $\U_0 := \pi``(\cal{V}_0)^\uparrow$. Note that $\pi(X_1) = \pi(X_0) \in \U_0$. Since $\pi$ is a projection in $M$, there exists some $\cal{R}$-generic filter $\cal{V}_1$ over $M$ such that $X_1 \in \cal{V}_1$ and $\pi``(\cal{V}_1)^\uparrow = \U_0$. Since $X \in \cal{V}_0 \cap \cal{V}_1$, we have that:
    \begin{align*}
        \cal{V}_0 = i_{\U_0}(\tau) = \cal{V}_1,
    \end{align*}
    contradicting that $X_0$ and $X_1$ are incompatible in $M$.
\end{proof}

We may now prove Corollary \ref{cor:no.wp.strong.subtrees.sigma-closed.examples} by verifying that the examples have projections that satisfy the hypothesis of Theorem \ref{thm:no.wp.strong.subtrees.sigma-closed.projection}:

\begin{proof}[Proof of Corollary \ref{cor:no.wp.strong.subtrees.sigma-closed.examples}]
    As discussed in the paragraph after Definition \ref{def:sigma-closed.quasi-order}, every example listed is $\sigma$-closed.
    \begin{enumerate}
        \item Recall from Example \ref{ex:min.max.FIN} the projection $\min : \FIN^{[\infty]} \to [\N]^\infty$. If $X = (x_n)_{n<\omega} \in \FIN^{[\infty]}$, then $X_0 := (x_{2n})_{n<\omega}$ and $X_1 := (x_{2n} \cup x_{2n+1})_{n<\omega}$ are incompatible, $X_0,X_1 \leq X$ and $\min(X_0) = \min(X_1)$. 

        \item[(2--4)] $\FIN_k^{[\infty]}$, $\FIN_{\pm k}^{[\infty]}$ and $W_{Lv}^{[\infty]}$ project onto $\FIN^{[\infty]}$ (Examples \ref{ex:tetris.operation}, \ref{ex:absolute.value.FIN_pmk} and \ref{ex:support.variable.words}), so a weak projection from $\S_\infty$ into any of these spaces would contradict (1).

        \item[(5)] Let $\pi_0 : \cal{R}_1 \to \cal{R}_0$ be the projection in Example \ref{ex:Dobrinen.spaces.projections} (where $\cal{R}_0$ is identified with $[\N]^\infty$, see Example \ref{ex:Dobrinen.space}). Let $X \in \cal{R}_1$. For each $n$, $X(2n+1)$ has exactly $2n + 2$ terminal nodes, so we may let $X_0(n)$ and $X_1(n)$ be the subtrees of $X(2n+1)$ whose terminal nodes are the $n + 1$ $\prec$-least and the $n + 1$ $\prec$-greatest terminal nodes of $X(2n+1)$ respectively. Let $X_0 := \bigcup_{n<\omega} X_0(n)$ and $X_1 := \bigcup_{n<\omega} X_1(n)$. Then $X_0,X_1 \in \cal{R}_1$, $X_0,X_1 \leq_1 X$, $\pi_0(X_0) = \pi_0(X_1)$ and $X_0$ and $X_1$ are incompatible. 

        \item[(6)] Let $p^* : \E_\infty \to [\N^+]^\infty$ be the projection in Example \ref{ex:minimal.representatives} (where $[\N^+]^\infty$ is identified with $[\N]^\infty$). Let $X \in \E_\infty$, and recall that $X(n)$ is the equivalence class of $X$ with the $n^\text{th}$-least representative. Let $X_0 \leq X$ be the equivalence relation whose equivalence classes are $X(2n) \cup X(2n+1)$ for $n < \omega$, and let $X_1 \leq X$ be the equivalence relation whose equivalence classes are $X(0) \cup X(1) \cup X(3)$, and $X(2n) \cup X(2n+3)$ for $n \geq 1$. Then $p^*(X_0) = p^*(X_1) = \{p_{2n}(X) : n \geq 1\}$, and the only equivalence relation on $\N$ coarser than both $X_0$ and $X_1$ is the trivial relation, so $X_0$ and $X_1$ are incompatible.
    \end{enumerate}
\end{proof}

We conclude this subsection with a brief discussion on the non-existence of finite amalgamable covers in $(\FIN^{[\infty]},\leq,r)$. This follows from Lemma \ref{lem:wp.Ellentuck.controlled} (with $\mu := \max$) and Proposition \ref{prop:no.wp.FIN.controlled} in the next subsection, as the identity map on $\FIN^{[\infty]}$ is a weak projection. It can also be shown directly: for any $A = (x_n)_{n<\omega} \in \FIN^{[\infty]}$, there is no finite amalgamable cover for $\emptyset$ in $[(x_0),A]$, which can be seen by considering, for each $1 \leq m < \omega$:
\begin{align*}
    B_m := (x_0 \cup x_m,x_{m+1},x_{m+2},\dots) \in [A].
\end{align*}
If $\{H_i\}_{i<k}$ is a finite amalgamable cover for $\emptyset$ in $[(x_0),A]$, then by Definition \ref{def:amalgamable.cover}(1), $B_m,B_{m'} \leq H_i$ for some $m < m'$ and $i < k$, which is impossible by Definition \ref{def:amalgamable.cover}(2). Alternatively, the conclusion of Lemma \ref{lem:decision} fails for $a_0 := \emptyset$ and the map $\Phi((x_n)_{n<\omega}) := \min(x_0)$: for any $B = (y_n)_{n<\omega} \in \FIN^{[\infty]}$ and $a \in \FIN^{[<\infty]}\restrictedto B$, $B$ accepts $a$ iff $\lh(a) \geq 1$. Let $a := (y_0)$ and $a' := (y_0 \cup y_1)$, and we have that $\Phi$ takes the same value $\min(y_0)$ on $[a,B]$ and $[a',B]$, but $\depth_B(a) = 1 \neq 2 = \depth_B(a')$.

\subsection{Carlson--Simpson space}
\label{sec:projections.Carlson-Simpson}
While the combinatorics of Ramseyan ultrafilters (see Example \ref{ex:Ramseyan.ultrafilter}) were studied in \cite{M86} and \cite{H01}, the Tukey theory of Ramsey ultrafilters of the Carlson--Simpson space $\E_\infty$ remains largely unexplored. However, it is possible to utilise other tools to prove the non-existence of weak projection maps into $\E_\infty$. The main theorems of this subsection are as follows:

\begin{theorem}
\label{thm:no.wp.Carlson-Simpson.and.FIN}
    $(\E_\infty,\leq,r)$ and $(\FIN^{[\infty]},\leq,r)$ do not weakly project into each other.
\end{theorem}

By considering composition with the projections in Examples \ref{ex:min.max.FIN}, \ref{ex:tetris.operation}, \ref{ex:absolute.value.FIN_pmk} and \ref{ex:support.variable.words}, we can conclude that:

\begin{corollary}
\label{cor:no.wp.Ellentuck.Carlson-Simpson}
    $([\N]^\infty,\subseteq,r)$ does not weakly project into $(\E_\infty,\leq,r)$.
\end{corollary}

\begin{corollary}
\label{cor:no.wp.Carlson-Simpson.examples}
    $(\E_\infty,\leq,r)$ does not weakly project into the following:
    \begin{enumerate}
        \item $(\FIN_k^{[\infty]},\leq,r)$ for $k \geq 2$.
        
        \item $(\FIN_{\pm k}^{[\infty]},\leq,r)$ for $k \geq 2$.
        
        \item $(W_{Lv}^{[\infty]},\leq,r)$, where $L$ is finite.
    \end{enumerate}
\end{corollary}

Recall that $(\cal{R},\leq_\cal{R},r)$ is \emph{unconditional} if $\cal{R}$ is either $[\N]^\infty$, $\FIN_k^{[\infty]}$, $\FIN_{\pm k}^{[\infty]}$, $W_{Lv}^{[\infty]}$ or $\FIN_{Lv}^{[\infty]}$ with $L$ finite, or $E^{[\infty]}$ --- see the end of \S\ref{sec:examples}. In each of these examples, there exists an ordering $<$ such that for a sequence $(x_n)_{n \leq N}$ of elements of $\cal{AR}_1$, $(x_n)_{n \leq N} \in \cal{AR}$ iff $(x_n)_{n<N} \in \cal{AR}$ and $(x_n)_{n<N} < x_N$. 

\begin{theorem}
\label{thm:no.wp.unconditional.Carlson-Simpson}
    Let $(\cal{R},\leq,r)$ be an unconditional \textbf{A2}-space.
    \begin{enumerate}
        \item There is a forcing notion $\Q * \dot{\P}$, where $\Q$ is $\sigma$-closed and $\forces_\Q$ $\dot{\P}$ is ccc, that forces ``$(\cal{R},\leq,r)$ does not weakly project into $(\E_\infty,\leq,r)$''.

        \item $(\cal{R},\leq,r)$ does not project into $(\E_\infty,\leq,r)$.

        \item Assuming the existence of large cardinals, $(\cal{R},\leq,r)$ does not weakly project into $(\E_\infty,\leq,r)$.
    \end{enumerate}
\end{theorem}

It is worth noting that the proof strategies of these results vary significantly. We outline them here:
\begin{itemize}
    \item For the $\E_\infty \nrightarrow \FIN^{[\infty]}$ direction of Theorem \ref{thm:no.wp.Carlson-Simpson.and.FIN}, we follow the first part of the proof of Theorem \ref{thm:no.wp.strong.subtrees.sigma-closed.projection} (more precisely, we use Lemmas \ref{lem:decision} and \ref{lem:wp.Ellentuck.controlled}). A contradiction is then obtained by composing any weak projection $f : \E_\infty \to \FIN^{[\infty]}$ with $\min$ and $\max$ respectively (Proposition \ref{prop:no.wp.FIN.controlled}).
        
    \item The converse direction uses a similar argument, but requires \"{O}zalp's classification, up to Rudin--Keisler equivalence, of the ultrafilters Tukey below a stable ordered-union ultrafilter (Theorem \ref{thm:Tukey.below.souu.RK.classification}).

    \item For Theorem \ref{thm:no.wp.unconditional.Carlson-Simpson}, we introduce an abstract notion of distributivity number $\frak{h}_\cal{R}$ for $\sigma$-closed \textbf{wA2}-spaces, and show that in a model by Brendle, $\frak{h}_{\E_\infty} < \frak{h}_\cal{R}$. Since weak projections preserve distributivity numbers (Lemma \ref{lem:wp.passes.distributivity}), Theorem \ref{thm:no.wp.unconditional.Carlson-Simpson}(1) follows, and Theorem \ref{thm:no.wp.unconditional.Carlson-Simpson}(2--3) follows from an absoluteness argument.
\end{itemize}

\subsubsection{Proof of Theorem \ref{thm:no.wp.Carlson-Simpson.and.FIN}}
\label{sec:proof.of.CS.and.FIN}
We first show that $(\E_\infty,\leq,r)$ has finite amalgamable covers (Definition \ref{def:finite.amalgamable.covers}), so that Lemmas \ref{lem:decision} and \ref{lem:wp.Ellentuck.controlled} apply.

\begin{example}
\label{ex:Carlson-Simpson.has.FACs}
    The Carlson--Simpson space $(\E_\infty,\leq,r)$ has finite amalgamable covers.
\end{example}

\begin{proof}
    The construction is similar to the strong binary subtrees case (Example \ref{ex:strong.subtrees.has.FACs}). Let $A \in \E_\infty$, $a \sqsubseteq A$ and $b \preceq_\fin a$. Let:
    \begin{align*}
        \F := \{c \in r_{\lh(b)+1}[b,A] : \depth_A(c) = \lh(a) + 1\}.
    \end{align*}
    Note that $\F$ is finite. For each $c \in \F$, let $H_c$ be the equivalence relation on $\N$ generated by $A \cup c$, so that $H_c$ is the $\leq$-maximal element of $[c,A]$. We shall show that $\{H_c\}_{c\in\F}$ is an amalgamable cover for $b$ in $[a,A]$. Recall that for $B \in \E_\infty$, $p_n(B)$ is the $n^\text{th}$-least representative of $B$ (see Example \ref{ex:Carlson-Simpson}).
    \begin{enumerate}
        \item Let $B \in [b,A]$ such that $\depth_A(r_{\lh(b)+1}(B)) > \lh(a)$, and let $c := B\restrictedto p_{\lh(a)+1}(A)$. The assumptions $b \preceq_\fin a$ and $\depth_A(r_{\lh(b)+1}(B)) > \lh(a)$ imply that:
        \begin{align*}
            \dom(b) \leq \dom(a) < p_{\lh(a)+1}(A) \leq p_{\lh(b)+1}(B),
        \end{align*}
        so $c = r_{\lh(b)+1}(B)\restrictedto p_{\lh(a)+1}(A)$ has exactly $\lh(b)+1$ equivalence classes. Therefore, $c \in \F$ and $B \leq H_c$.
        
        \item Let $A' \in [a,A]$ and $B \in [b,A']$ such that $B \leq H_c$ for some $c \in \F$. Since $b \sqsubseteq B$ and $B \leq H_c$, we have that $B\restrictedto\dom(c) = c$. Let $B'$ be the equivalence relation on $\N$ generated by $A'$ and all pairs $(\min(X),x)$ such that $X$ is an equivalence class of $B$ and $x \in X \setminus (A'(0) \cup \cdots \cup A'(\lh(a)))$. Then $A'(i) \subseteq B'(i)$ for all $i \leq \lh(a)$, which implies that $a \sqsubseteq B'$, so $B' \in [a,A']$. Since $B\restrictedto\dom(c) = c$, $B$ is the equivalence relation on $\N$ generated by $B' \cup c$, so $[b,B'] \cap [H_c] \subseteq [b,B]$.
    \end{enumerate}
\end{proof}

\begin{proposition}
\label{prop:no.wp.FIN.controlled}
    Let $(\cal{R},\leq,r)$ be a closed triple satisfying \textbf{A1}--\textbf{A4} that has finite amalgamable covers, and let $\mu : \cal{R} \to [\N]^\infty$ be a projection that satisfies the requirements of $\mu$ in Lemma \ref{lem:wp.Ellentuck.controlled} for $a_0 := \emptyset$. Then $(\cal{R},\leq,r)$ does not weakly project into $(\FIN^{[\infty]},\leq,r)$.
\end{proposition}

\begin{proof}
    Suppose for a contradiction that $g : \cal{R} \to \FIN^{[\infty]}$ is monotone Alaoglu--Birkhoff convergent. Let:
    \begin{align*}
        h_0 := \min \circ g, \quad h_1 := \max \circ g,
    \end{align*}
    where $\min,\max : \FIN^{[\infty]} \to [\N]^\infty$ are the maps defined in Example \ref{ex:min.max.FIN}. Then $h_0$ and $h_1$ are weak projections, so by applying Lemma \ref{lem:wp.Ellentuck.controlled} to $h_0$ and then to $h_1$, there exist some $B \in \cal{R}$ and injective $\psi_0,\psi_1 : \N \to \N$ such that $\psi_i``(\mu(B')) \cap h_i(B') \neq \emptyset$ for all $B' \leq B$ and $i < 2$. Since $\psi_0$ is injective, there exists some $\theta : \N \to \N$ such that $\theta \circ \psi_0 = \psi_1$. Let:
    \begin{align*}
        \K := \{X \in \FIN^{[\infty]} : \theta``(\min(X)) \cap \max(X) = \emptyset\}.
    \end{align*}

    \begin{claim}
        $\K$ is a dense open subset of $\FIN^{[\infty]}$.
    \end{claim}

    \begin{midproof}
        It is clear that $\K$ is open. Let $X = (x_n)_{n<\omega} \in \FIN^{[\infty]}$. By passing to a subsequence, we may assume that $n \mapsto \theta(\min(x_n))$ is constant or injective. 
        \begin{itemize}
            \item Suppose that it is constant with value $m$. Choose $N$ large enough so that $\max(x_N) > m$, and we have that $(x_n)_{n \geq N} \leq X$ is in $\K$. 
            
            \item Suppose that it is injective. We define $Y = (y_k)_{k<\omega} \leq X$ inductively as follows. Suppose that $y_i$ has been defined for all $i < k$, and that $\theta(\min(y_i)) \neq \max(y_j)$ for all $i,j < k$. Since both maps $n \mapsto \theta(\min(x_n))$ and $n \mapsto \max(x_n)$ are injective, we may choose $n_k < n_k'$, both greater than $n_i'$ for all $i < k$, such that $\theta(\min(x_{n_k})) \neq \max(x_{n_k'})$ and they are not $\theta(\min(y_i))$ or $\max(y_i)$ for any $i < k$. Let $y_k := x_{n_k} \cup x_{n_k'}$, and we have that $\theta(\min(y_k)) = \theta(\min(x_{n_k}))$ and $\max(y_k) = \max(x_{n_k'})$, completing the induction. Then it's clear from the construction that $Y \in \K$.
        \end{itemize}
    \end{midproof}

    Since $g$ is Alaoglu--Birkhoff convergent, there exists some $B' \leq B$ such that $g(B') \in \K$, i.e. $\theta``(h_0(B')) \cap h_1(B') = \emptyset$. For each $i < 2$, let:
    \begin{align*}
        P_i := \{n \in \mu(B') : \psi_i(n) \in h_i(B')\}.
    \end{align*}
    Observe that $P_0 \cap P_1 = \emptyset$ --- otherwise, we have that $\psi_0(n) \in h_0(B')$ and $\psi_1(n) \in h_1(B')$ for some $n \in \mu(B')$. But $\psi_1(n) = \theta(\psi_0(n)) \in \theta``(h_0(B'))$, so $\psi_1(n) \in \theta``(h_0(B')) \cap h_1(B')$, contradicting our choice of $B'$. Therefore, $\mu(B') \setminus P_i$ is infinite for some $i < 2$, and since $\mu$ is a projection, there exists some $B'' \leq B'$ such that $\mu(B'') \subseteq \mu(B') \setminus P_i$. But by the choice of $B$ at the beginning of the proof, we have that $\psi_i``(\mu(B'')) \cap h_i(B'') \neq \emptyset$, so there exists some $n \in \mu(B'')$ such that $\psi_i(n) \in h_i(B'') \subseteq h_i(B')$, contradicting that $\mu(B'') \cap P_i = \emptyset$.
\end{proof}

\begin{proof}[Proof of Theorem \ref{thm:no.wp.Carlson-Simpson.and.FIN}, $\E_\infty \nrightarrow \FIN^{[\infty]}$]
    It suffices to show that $(\E_\infty,\leq,r)$ satisfies the hypothesis of Lemma \ref{lem:wp.Ellentuck.controlled}, so that Proposition \ref{prop:no.wp.FIN.controlled} applies. We showed in Example \ref{ex:Carlson-Simpson.has.FACs} that $(\E_\infty,\leq,r)$ has finite amalgamable covers, and the map $p^* : \E_\infty \to [\N^+]^\infty$ in Example \ref{ex:minimal.representatives} satisfies the requirements of $\mu$ in Lemma \ref{lem:wp.Ellentuck.controlled} for $a_0 := \emptyset$, where for all $a \in \E_{<\infty}$ we may let:
    \begin{align*}
        \alpha(a) := \{\dom(r_k(a)) : 1 \leq k \leq \lh(a)\}.
    \end{align*}
    We leave the details of the verification of requirements to the reader.
\end{proof}

Note that, by Lemma \ref{lem:wp.Ellentuck.controlled} and Examples \ref{ex:Ellentuck.has.FACs} and \ref{ex:strong.subtrees.has.FACs}, Proposition \ref{prop:no.wp.FIN.controlled} also gives forcing-free proofs of Theorem \ref{thm:no.wp.Ellentuck.FIN} and Corollary \ref{cor:no.wp.strong.subtrees.sigma-closed.examples}(1).

To prove the converse direction of Theorem \ref{thm:no.wp.Carlson-Simpson.and.FIN}, we require \"{O}zalp's complete classification of ultrafilters Tukey below a stable ordered-union ultrafilter up to Rudin--Keisler equivalence. We first recall the notions needed to state \"{O}zalp's classification theorem. Given two countable sets $I,J$ and an ultrafilter $\cal{V}$ on $I$, if $h : I \to J$ is a function, then the \emph{RK-image} of $\cal{V}$ is the set:
\begin{align*}
    h(\cal{V}) := \{V \subseteq J : h^{-1}[V] \in \cal{V}\},
\end{align*}
which is an ultrafilter on $J$. Given another ultrafilter $\cal{W}$ on $J$, we say that $\cal{W}$ is \emph{Rudin--Keisler below} $\cal{V}$, denoted $\cal{W} \leq_{RK} \cal{V}$, if $\cal{W} = h(\cal{V})$ for some $h : I \to J$, and that $\cal{W}$ is \emph{Rudin--Keisler equivalent} to $\cal{V}$, denoted $\cal{W} \cong \cal{V}$, if $\cal{W} = h(\cal{V})$ for some $h : I \to J$ which is injective on an element of $\cal{V}$. Note that $\cal{W} \cong \cal{V}$ iff $\cal{W} \leq_{RK} \cal{V}$ and $\cal{V} \leq_{RK} \cal{W}$, so $\cong$ is an equivalence relation. 

If $\cal{W}$ is an ultrafilter on $I$ and $\cal{W}_x$ is an ultrafilter on $J$ for each $x \in I$, then the \emph{Fubini limit} is defined as the set:
\begin{align*}
    \lim_{x\to\cal{W}} \cal{W}_x := \{W \subseteq I \times J : \{x \in I : \{y \in J : (x,y) \in W\} \in \cal{W}_x\} \in \cal{W}\},
\end{align*}
which is an ultrafilter on $I \times J$. If $\cal{W}_x = \cal{W}'$ for some fixed ultrafilter $\cal{W}'$, then we write $\cal{W} \cdot \cal{W}' := \lim_{x\to\cal{W}} \cal{W}_x$. We also write $\cal{W}^2 := \cal{W} \cdot \cal{W}$, and observing that $\cdot$ is associative up to Rudin--Keisler equivalence, we inductively define $\cal{W}^{n+1} := \cal{W} \cdot \cal{W}^n$. We shall need the following basic facts about Fubini limits of ultrafilters:

\begin{lemma}
\label{lem:Fubini.limit.facts}
    Let $I,J$ be countable sets, let $\cal{W}$ be an ultrafilter on $I$, and let $\cal{W}_x$ and $\cal{V}_x$ be ultrafilters on $J$ for each $x \in I$.
    \begin{enumerate}
        \item $\cal{W} \leq_{RK} \lim_{x\to\cal{W}} \cal{W}_x$. Moreover, if $\cal{W}'$ is an ultrafilter and:
        \begin{align*}
            \{x \in I : \cal{W}' \leq_{RK} \cal{W}_x\} \in \cal{W},
        \end{align*}
        then $\cal{W}' \leq_{RK} \lim_{x\to\cal{W}} \cal{W}_x$.

        \item If:
        \begin{align*}
            \{x \in I : \cal{V}_x \leq_{RK} \cal{W}_x\} \in \cal{W}, 
        \end{align*}
        then $\lim_{x\to\cal{W}} \cal{V}_x \leq_{RK} \lim_{x\to\cal{W}} \cal{W}_x$.
    \end{enumerate}
\end{lemma}

Given a Ramsey ultrafilter $\U \subseteq \FIN^{[\infty]}$, the classes of ultrafilters of interest in \"{O}zalp's classification theorem are defined inductively as follows, where we identify every countable set with $\N$:
\begin{itemize}
    \item $\cal{C}_0 := \{\U_{\min},\U_{\max},\U_{\minmax},\U_{\FIN}\}$.

    \item $\cal{C}_{\alpha+1} := \{\lim_{x \to \cal{V}} \cal{W}_x : \cal{V} \in \cal{C}_0 \text{ and } \cal{W}_x \in \cal{C}_\alpha \text{ for all } x \in \N\}$.

    \item $\cal{C}_\alpha := \bigcup_{\beta<\alpha} \cal{C}_\beta$ if $\alpha$ is a limit ordinal.
\end{itemize}

\begin{theorem}[{\cite[Theorem 5.4]{O24}}]
\label{thm:Tukey.below.souu.RK.classification}
    Let $\U \subseteq \FIN^{[\infty]}$ be a Ramsey ultrafilter, and let $\cal{W}$ be a non-principal ultrafilter on a countable set such that $\cal{W} \leq_T \U$. Then $\cal{W} \cong \cal{W}'$ for some $\alpha < \omega_1$ and $\cal{W}' \in \cal{C}_\alpha$.
\end{theorem}

We may now prove the converse direction of Theorem \ref{thm:no.wp.Carlson-Simpson.and.FIN}.

\begin{proof}[Proof of Theorem \ref{thm:no.wp.Carlson-Simpson.and.FIN}, $\FIN^{[\infty]} \nrightarrow \E_\infty$]
    Suppose for a contradiction that there is a weak projection $f : \FIN^{[\infty]} \to \E_\infty$. Let $M$ be the transitive collapse of a countable elementary submodel of some $H_\theta$ containing $f$ and everything else relevant. Let $\U$ be $\FIN^{[\infty]}$-generic over $M$, and note that $M[\U]$ and $M$ share the same countable sequences of elements of $M$ (in particular, the same reals). We work in $M[\U]$ for the remainder of the proof. Recall the projections $\min,\max : \FIN^{[\infty]} \to [\N]^\infty$ (Example \ref{ex:min.max.FIN}) and $p^* : \E_\infty \to [\N]^\infty$ (Example \ref{ex:minimal.representatives}). We consider the following ultrafilters, which are generic over $M$ by Lemma \ref{lem:ab.equivalence} (hence Ramsey by Lemma \ref{lem:generic.filter.is.Ramsey}):
    \begin{align*}
        \cal{V} &:= f``(\U)^\uparrow \tag*{($\E_\infty$-generic)} \\
        \cal{V}^* &:= (p^*)``(\cal{V}) = (p^* \circ f)``(\U)^\uparrow \tag*{($[\N]^\infty$-generic)} \\
        \U_{\min} &:= \min``(\U)^\uparrow \tag*{($[\N]^\infty$-generic)} \\
        \U_{\max} &:= \max``(\U)^\uparrow \tag*{($[\N]^\infty$-generic)}
    \end{align*}
    We also define:
    \begin{align*}
        \tilde{\cal{V}} := \{r_2[X] : X \in \cal{V}\}^\uparrow,
    \end{align*}
    so $\tilde{\cal{V}} \leq_T \cal{V} \leq_T \U$, and $\tilde{\cal{V}}$ is an ultrafilter on the countable set $\E_2$ by Lemma \ref{lem:Ramsey.ultrafilter.finite.pigeonhole}. Note that $\tilde{\cal{V}}$ is non-principal, as $r_2[X]$ is infinite for all $X \in \E_\infty$.

    \begin{claim}
    \label{claim:Carlson--Simpson.ultrafilter.base.representation}
        $\cal{V} = \{X \in \E_\infty : r_2[X] \in \tilde{\cal{V}}\}$. Therefore, $\tilde{\cal{V}} \equiv_T \cal{V}$ and $M[\cal{V}] = M[\tilde{\cal{V}}]$. 
    \end{claim}

    \begin{midproof}
        If $r_2[X] \in \tilde{\cal{V}}$, then $r_2[Y] \subseteq r_2[X]$ for some $Y \in \cal{V}$, so it suffices to show that if $X,Y \in \E_\infty$ and $r_2[Y] \subseteq r_2[X]$, then $Y \leq X$. Suppose otherwise, so there exist some $i,j \in \N$ that are $X$-equivalent but not $Y$-equivalent. Without loss of generality, we may assume that $i \in Y(m)$ and $j \in Y(n)$ for some $m < n$. Let $k > n$ be large enough so that $p_k(Y) > \max\{i,j\}$, and let $Y' \leq Y$ be obtained from $Y$ by merging $Y(0),\dots,Y(n-1)$ into one equivalence class, and $Y(n),\dots,Y(k-1)$ into another. Then $i \in Y'(0)$, $j \in Y'(1)$ and $p_2(Y') = p_k(Y) > \max\{i,j\}$, so $i$ and $j$ are not $r_2(Y')$-equivalent. But $r_2(Y') \in r_2[Y] \subseteq r_2[X]$, which is a contradiction, as $i$ and $j$ are $X$-equivalent, and hence $a$-equivalent for every $a \in r_2[X]$ with $i,j \in \dom(a)$.
    \end{midproof}

    \begin{claim}
    \label{claim:Carlson-Simpson.generic.RK-images}
        Let $\cal{Y}$ be an ultrafilter on $\N$ that is $[\N]^\infty$-generic over $M$. If $\cal{Y} \leq_{RK} \tilde{\cal{V}}$, then $\cal{Y} \cong \cal{V}^*$.
    \end{claim}

    \begin{midproof}
        Recall from the proof of the $\E_\infty \nrightarrow \FIN^{[\infty]}$ direction that $p^*$ satisfies the requirements of $\mu$ in Lemma \ref{lem:wp.Ellentuck.controlled} for $a_0 := \emptyset$. Since $\cal{V}^* = (p^*)``(\cal{V})$ and $\cal{V}$ is $\E_\infty$-generic over $M$, by Lemmas \ref{lem:wp.Ellentuck.controlled} and \ref{lem:controlled.wps.give.isomorphic.generic.filters}, it suffices to show that there exist some $X_0 \in \cal{V}$ and some $h : \E_\infty \to \Po(\N)$ such that $h\restrictedto[X_0] : [X_0] \to [\N]^\infty$ is a weak projection, and $\cal{Y} = h``(\cal{V})^\uparrow$.
        
        Let $g : \E_2 \to \N$ be such that $\cal{Y} = g(\tilde{\cal{V}})$, and note that $g \in M$. Define $h : \E_\infty \to \Po(\N)$ by:
        \begin{align*}
            h(X) := g``(r_2[X]).
        \end{align*}
        It's clear that $h \in M$ is monotone, and $\cal{Y} = h``(\cal{V})^\uparrow$. Since $\cal{Y}$ is $[\N]^\infty$-generic over $M$, there exists some $X_0 \in \cal{V}$ such that $h``(\cal{V}')^\uparrow$ is $[\N]^\infty$-generic over $M$ for every $\E_\infty$-generic filter $\cal{V}'$ over $M$ with $X_0 \in \cal{V}'$. In particular, $h(X) \in [\N]^\infty$ for all $X \leq X_0$, and by Lemma \ref{lem:ab.equivalence}, $h\restrictedto[X_0] : [X_0] \to [\N]^\infty$ is a weak projection. 
    \end{midproof}

    If $\U_{\min} \leq_{RK} \tilde{\cal{V}}$ and $\U_{\max} \leq_{RK} \tilde{\cal{V}}$, then $\U_{\min} \cong \cal{V}^* \cong \U_{\max}$ by Claim \ref{claim:Carlson-Simpson.generic.RK-images}, contradicting that $\U_{\min} \not\cong \U_{\max}$ (cf. the paragraph before Theorem \ref{thm:souu.initial.Tukey.structure}). Thus, we may let $\cal{W},\cal{W}'$ be such that $\{\cal{W},\cal{W}'\} = \{\U_{\min},\U_{\max}\}$ and $\cal{W}' \not\leq_{RK} \tilde{\cal{V}}$.

    \begin{claim}
    \label{claim:W.V.isomorphic}
        $\cal{W} \cong \cal{V}^*$ and $\tilde{\cal{V}} \in M[\cal{W}]$. Therefore, $M[\tilde{\cal{V}}] = M[\cal{W}] = M[\cal{V}^*]$ (cf. Lemma \ref{lem:controlled.wps.give.isomorphic.generic.filters}).
    \end{claim}

    \begin{midproof}
        Since $\tilde{\cal{V}} \leq_T \U$, by Theorem \ref{thm:Tukey.below.souu.RK.classification}, there exist some $\alpha < \omega_1$ and some $\cal{Y} \in \cal{C}_\alpha$ such that $\cal{Y} \cong \tilde{\cal{V}}$, so $\cal{W}' \not\leq_{RK} \cal{Y}$. As $M[\U]$ and $M$ share the same reals, all functions witnessing Rudin--Keisler equivalences lie in $M$. Thus, by Claim \ref{claim:Carlson-Simpson.generic.RK-images}, it suffices to show, by induction on $\alpha$, that every $\cal{Y} \in \cal{C}_\alpha$ with $\cal{W}' \not\leq_{RK} \cal{Y}$ satisfies $\cal{W} \leq_{RK} \cal{Y}$, and is Rudin--Keisler equivalent to some ultrafilter in $M[\cal{W}]$.

        If $\alpha = 0$, then $\cal{Y} = \cal{W}$, as $\U_{\min},\U_{\max} \leq_{RK} \U_{\minmax} \leq_{RK} \U_\FIN$. The limit case is immediate. Now suppose that $\cal{Y} = \lim_{x\to\cal{Z}} \cal{Y}_x$, where $\cal{Z} \in \cal{C}_0$ and $\cal{Y}_x \in \cal{C}_\alpha$ for all $x \in \N$. By Lemma \ref{lem:Fubini.limit.facts}(1), we have that $\cal{Z} \leq_{RK} \cal{Y}$, so $\cal{W}' \not\leq_{RK} \cal{Z}$, and hence $\cal{Z} = \cal{W}$ by the case $\alpha = 0$. In particular, $\cal{W} \leq_{RK} \cal{Y}$. By Lemma \ref{lem:Fubini.limit.facts}(1) again, we have that:
        \begin{align*}
            I' := \{x \in \N : \cal{W}' \not\leq_{RK} \cal{Y}_x\} \in \cal{W}.
        \end{align*}
        By the induction hypothesis, for each $x \in I'$ there exists some ultrafilter $\cal{Y}_x' \in M[\cal{W}]$ such that $\cal{Y}_x \cong \cal{Y}_x'$, and we let $\cal{Y}_x' := \cal{W}$ for each $x \in \N \setminus I'$. Since $M[\U]$ and $M$ share the same countable sequences of elements of $M$, the sequence $(\cal{Y}_x')_{x \in \N}$ lies in $M[\cal{W}]$, and hence so does $\cal{Y}' := \lim_{x\to\cal{W}} \cal{Y}_x'$. Moreover, $\cal{Y} \cong \cal{Y}'$ by Lemma \ref{lem:Fubini.limit.facts}(2).
    \end{midproof}

    By Claims \ref{claim:Carlson--Simpson.ultrafilter.base.representation} and \ref{claim:W.V.isomorphic}, we have that $\cal{V} \in M[\tilde{\cal{V}}] = M[\cal{W}] = M[\cal{V}^*]$. By the proof of Corollary \ref{cor:no.wp.strong.subtrees.sigma-closed.examples}(6), for every $X \in \E_\infty$, there exist incompatible $X_0,X_1 \leq X$ such that $p^*(X_0) = p^*(X_1)$. Therefore, we obtain a contradiction as at the end of the proof of Theorem \ref{thm:no.wp.strong.subtrees.sigma-closed.projection}, with $\E_\infty$ and $p^*$ in place of $\cal{R}$ and $\pi$.
\end{proof}

\subsubsection{Proof of Theorem \ref{thm:no.wp.unconditional.Carlson-Simpson}}
\begin{definition}
\label{def:distributivity.number.qo}
    Let $(\P,\leq)$ be a quasi-order. The \emph{distributivity number of $\P$}, denoted by $\frak{h}_\P$, is the least cardinal $\kappa$ (if it exists) such that there exists a collection $\{\D_\alpha\}_{\alpha<\kappa}$ of dense open subsets in $(\P,\leq)$ for which $\bigcap_{\alpha<\kappa} \D_\alpha = \emptyset$.
\end{definition}

\begin{definition}
\label{def:distributivity.number}
    Let $(\cal{R},\leq,\leq^*,r)$ be a $\sigma$-closed \textbf{wA2}-space. The \emph{distributivity number of $\cal{R}$}, denoted by $\frak{h}_\cal{R}$, is the distributivity number of the quasi-order $(\cal{R},\leq^*)$. 
\end{definition}

It is clear that $\frak{h}_{[\N]^\infty} = \frak{h}$ --- see, for instance, \cite{B10} for details of the cardinal invariant $\frak{h}$. If $(\cal{R},\leq,\leq^*,r)$ is a $\sigma$-closed \textbf{wA2}-space, then $\omega_1 \leq \frak{h}_\cal{R}$. It is easy to check that for unconditional spaces, $\frak{h}_\cal{R}$ is well-defined and $\frak{h}_\cal{R} \leq \frak{c}$. The distributivity numbers for various topological Ramsey spaces have been previously studied.
\begin{itemize}
    \item The cardinal invariant $\frak{h}_\FIN := \frak{h}_{\FIN^{[\infty]}}$ was studied in \cite{BG17}, where it was shown, among other related results, that $\frak{h}_\FIN < \frak{h}$ is consistent. 
    
    \item The cardinal invariant $\frak{h}_{\E_\infty}$ was introduced and studied in \cite{CKMW00}. It was subsequently shown in \cite{B00} that $\ZFC + \MA + \frak{h}_{\E_\infty} = \omega_1$ is consistent.
\end{itemize}
A natural question therefore arises --- if $f : \cal{R} \to \S$ is a weak projection between two $\sigma$-closed topological Ramsey spaces, then must the inequality $\frak{h}_\cal{R} \leq \frak{h}_\S$ hold? It is clear that if there is an Alaoglu--Birkhoff convergent map from $\P$ to $\Q$, then $\frak{h}_\P \leq \frak{h}_\Q$. However, a weak projection $f : \cal{R} \to \S$ need not be Alaoglu--Birkhoff convergent with respect to $\leq_\cal{R}^*$ and $\leq_\S^*$, so it is unclear whether $\frak{h}_\cal{R} \leq \frak{h}_\S$ holds in general. It turns out that this is true if we assume that $(\S,\leq_\S^*)$ is separative.

\begin{lemma}
\label{lem:wp.passes.distributivity}
    Let $(\cal{R},\leq_\cal{R},\leq_\cal{R}^*,r)$ and $(\S,\leq_\S,\leq_\S^*,s)$ be two $\sigma$-closed \textbf{wA2}-spaces such that $(\S,\leq_\S^*)$ is separative. If there is a weak projection $f : \cal{R} \to \S$, then $\frak{h}_\cal{R} \leq \frak{h}_\S$.
\end{lemma}

\begin{proof}
    Let $\kappa < \frak{h}_\cal{R}$, and let $\{\D_\alpha\}_{\alpha<\kappa}$ be a collection of dense open subsets in $(\S,\leq_\S^*)$. 

    \begin{claim}
    \label{claim:do.equivalence}
        If $\D$ is a dense open subset in $(\S,\leq_\S^*)$, then $\D$ is a dense open subset in $(\S,\leq_\S)$.
    \end{claim}

    \begin{midproof}
        That $\D$ is open in $(\S,\leq_\S)$ follows from the fact that $\leq_\S^*$ coarsens $\leq_\S$. We shall show that $\D$ is dense in $(\S,\leq_\S)$. Let $X \in \S$, and let $Y \in \D$ such that $Y \leq^* X$. Since $(\S,\leq_\S^*)$ and $(\S,\leq_\S)$ share the same separative quotient, there exists some $Z \in \S$ such that $Z \leq Y$ and $Z \leq X$. Since $\leq_\S^*$ coarsens $\leq_\S$, we have that $Z \leq^* Y$, so $Z \in \D$ as $\D$ is open in $(\S,\leq_\S^*)$. Then $Z \in \D$ and $Z \leq X$, so $\D$ is dense in $(\S,\leq_\S)$.
    \end{midproof}

    Since $f$ is a weak projection, $f^{-1}(\D_\alpha)$ contains a dense open subset in $(\cal{R},\leq_\cal{R})$ for all $\alpha < \kappa$. Let:
    \begin{align*}
        \K_\alpha := \{B \in \cal{R} : B \leq_\cal{R}^* A \text{ for some } A \in f^{-1}(\D_\alpha)\}.
    \end{align*}
    It is clear that $\K_\alpha$ contains a dense open subset in $(\cal{R},\leq_\cal{R}^*)$. Since $\kappa < \frak{h}_\cal{R}$, we may let $B \in \bigcap_{\alpha<\kappa} \K_\alpha$. For each $\alpha < \kappa$, let $A_\alpha \in f^{-1}(\D_\alpha)$ be such that $B \leq_\cal{R}^* A_\alpha$. 

    \begin{claim}
    \label{claim:reduction.always.compatible}
        There exists some $Y \in \S$ such that for all $\alpha < \kappa$, every $Z \leq_\S Y$ is compatible with $f(A_\alpha)$ in $(\S,\leq_\S)$.
    \end{claim}

    \begin{midproof}
        Suppose otherwise. Then the set:
        \begin{align*}
            \cal{J} := \{Z \in \S : Z \perp f(A_\alpha) \text{ for some } \alpha < \kappa\}
        \end{align*}
        is dense open in $(\S,\leq_\S)$. Since $f$ is Alaoglu--Birkhoff convergent, there exist some $C \leq_\cal{R} B$ and $\alpha < \kappa$ such that $f(C) \perp f(A_\alpha)$. But $C \leq_\cal{R}^* A_\alpha$, so $C$ and $A_\alpha$ are compatible in $(\cal{R},\leq_\cal{R})$. Let $D \in \cal{R}$ such that $D \leq_\cal{R} C$ and $D \leq_\cal{R} A_\alpha$, and we have that $f(D) \leq_\S f(C)$ and $f(D) \leq_\S f(A_\alpha)$, contradicting that $f(C) \perp f(A_\alpha)$.
    \end{midproof}

    Let $Y \in \S$ be as in Claim \ref{claim:reduction.always.compatible}, and let $\alpha < \kappa$. Since $(\S,\leq_\S)$ and $(\S,\leq_\S^*)$ share the same separative quotient, every $Z \leq_\S^* Y$ is compatible with $f(A_\alpha)$ in $(\S,\leq_\S^*)$, so $Y \leq_\S^* f(A_\alpha)$ as $(\S,\leq_\S^*)$ is separative. Since $f(A_\alpha) \in \D_\alpha$, we have that $Y \in \D_\alpha$. Therefore, $Y \in \bigcap_{\alpha<\kappa} \D_\alpha$, so $\kappa < \frak{h}_\S$.
\end{proof}

We now compare the distributivity numbers of unconditional spaces and of the Carlson--Simpson space. For unconditional \textbf{wA2}-spaces, we have the following lower bound:

\begin{lemma}
\label{lem:unconditional.p.leq.h_R}
    If $(\cal{R},\leq,r)$ is an unconditional \textbf{wA2}-space, then $\frak{p} \leq \frak{h}_\cal{R}$.
\end{lemma}

\begin{proof}
    Let $A \in \cal{R}$ be arbitrary, and let $\kappa < \frak{p}$. Since a $\kappa$-closed quasi-order is $\kappa$-distributive, it suffices to show that $([A],\leq^*)$ is $\kappa$-closed --- that is, if $(A_\alpha)_{\alpha<\kappa}$ is a $\leq^*$-decreasing sequence in $[A]$, then there exists some $B \in \cal{R}$ such that $B \leq^* A_\alpha$ for all $\alpha$. We define a forcing notion $(\P,\preceq)$ as follows:
    \begin{itemize}
        \item $\P = \{(a,S) : a \in \cal{AR}\restrictedto A \wedge S \subseteq \kappa \text{ is finite}\}$.
        
        \item If $(a,S),(b,T) \in \P$, then $(b,T) \preceq (a,S)$ iff the following conditions hold:
        \begin{itemize}
            \item $a \sqsubseteq b$.
            
            \item $S \subseteq T$.
            
            \item If $b = (x_n)_{n<\lh(b)}$, then for all $\lh(a) \leq n < \lh(b)$ and $\alpha \in S$, $x_n \in \c{A_\alpha}$.
        \end{itemize}
    \end{itemize}
    Observe that $\P$ is $\sigma$-centred, as $\cal{AR}\restrictedto A$ is countable and for all $S,T \subseteq \kappa$ and $a \in \cal{AR}\restrictedto A$, $(a,S \cup T) \preceq (a,S)$ and $(a,S \cup T) \preceq (a,T)$. For each $\alpha < \kappa$ and $n < \omega$, let:
    \begin{align*}
        \D_\alpha &:= \{(a,S) \in \P : \alpha \in S\}, \\
        \K_n &:= \{(a,S) \in \P : \lh(a) \geq n\}.
    \end{align*}
    Observe that these sets are dense in $\P$ --- if $(a,S) \in \P$, then $(a,S \cup \{\alpha\}) \preceq (a,S)$ and $(a,S \cup \{\alpha\}) \in \D_\alpha$. For $\K_n$, if $n \leq \lh(a)$ then $(a,S) \in \K_n$. Otherwise, we may choose $a < x_{\lh(a)} < \cdots < x_{n-1} \in \bigcap_{\alpha \in S} \c{A_\alpha}$, and let $b := a^\frown (x_i)_{\lh(a) \leq i < n}$. Then $(b,S) \preceq (a,S)$ and $(b,S) \in \K_n$.

    Since $\kappa < \frak{p}$ and $(\P,\preceq)$ is $\sigma$-centred, there exists some filter $G \subseteq \P$ such that $G \cap \D_\alpha \neq \emptyset$ for all $\alpha$, and $G \cap \K_n \neq \emptyset$ for all $n$. Since $G \subseteq \P$ is a filter, we have that for all $(a,S),(b,T) \in G$, $a \sqsubseteq b$ or $b \sqsubseteq a$. Thus, we may let $B \in \cal{R}$ be the unique element such that for all $n$, $(r_n(B),S) \in G$ for some $S \subseteq \kappa$. Note that $B \in \cal{R}$ (and not $\cal{AR}$) as $G \cap \K_n \neq \emptyset$ for all $n$. We shall show that $B \leq^* A_\alpha$ for all $\alpha < \kappa$.

    Write $B = (x_i)_{i<\omega}$, and let $\alpha < \kappa$. Since $G \cap \D_\alpha \neq \emptyset$, there exist some $n < \omega$ and $S \subseteq \kappa$ with $\alpha \in S$ such that $(r_n(B),S) \in G$. Fix any $k \geq n$. By the choice of $B$, there is some $(r_k(B),T') \in G$, and since $G$ is a filter, we may fix $(c,U) \in G$ with $(c,U) \preceq (r_n(B),S)$ and $(c,U) \preceq (r_k(B),T')$. Then $r_k(B) \sqsubseteq c$, so $\lh(c) \geq k$, and $(c,U) \preceq (r_n(B),S)$ gives $x_i \in \c{A_\alpha}$ for all $n \leq i < \lh(c)$, in particular for all $n \leq i < k$. Since $k$ is arbitrary, we have that $x_i \in \c{A_\alpha}$ for all $i \geq n$, so $B/n \leq A_\alpha$ and hence $B \leq^* A_\alpha$.
\end{proof}

On the other hand, Brendle showed in \cite{B00} that $\frak{p} \leq \frak{h}_{\E_\infty}$ may fail.

\begin{theorem}[Brendle, \cite{B00}]
\label{thm:Carlson-Simpson.distributivity.omega_1}
    Let $M$ be a model of $\ZFC + \CH$, and let $\kappa \geq \aleph_2$ be a regular cardinal. Then there exists a ccc forcing $(\P,\leq)$ such that if $G$ is $\P$-generic over $M$, then $M[G]$ satisfies $\ZFC + \MA + \frak{c} = \kappa + \frak{h}_{\E_\infty} = \omega_1$. 
\end{theorem}

We turn our attention to the complexity of various properties of functions, which are crucial to our absoluteness arguments.

\begin{lemma}
\label{lem:projection.complexities}
    Let $(\cal{R},\leq_\cal{R},r)$ and $(\S,\leq_\S,s)$ be two \textbf{wA2}-spaces, with finitisations $\leq_{\cal{R},\fin}$ and $\leq_{\S,\fin}$ respectively. Assume that $\cal{AR}$ and $\cal{AS}$ are countable. Let $\alpha : \cal{AR} \to \cal{AS}$ be a $\sqsubseteq$-preserving map, and let $t$ be a real that codes $\alpha$, $\leq_{\cal{R},\fin}$ and $\leq_{\S,\fin}$.
    \begin{enumerate}
        \item ``$\alpha$ induces some $f : \cal{R} \to \S$'' is $\Pi_1^1(t)$.

        \item ``$f$ is monotone'' is $\Pi_1^1(t)$.

        \item ``$f$ is quasi-open'' is $\Pi_4^1(t)$.

        \item ``$f$ is a projection'' is $\Pi_2^1(t)$.
    \end{enumerate}
\end{lemma}

\begin{proof}
    We first note that:
    \begin{align*}
        A \leq_\cal{R} B &\iff \forall n \, \exists m \, \exists b \in \cal{AR}(b \sqsubseteq r_m(B) \wedge r_n(A) \leq_\fin b), \\
        X \leq_\S f(B) &\iff \forall n \, \exists m \, \exists y \in \cal{AS}(y \sqsubseteq \alpha(r_m(B)) \wedge s_n(X) \leq_\fin y), \\  
        f(A) \leq_\S Y &\iff \forall n \, \exists m \, \exists y \in \cal{AS}(y \sqsubseteq s_m(Y) \wedge \alpha(r_n(A)) \leq_\fin y), \\  
        X \leq_\S Y &\iff \forall n \, \exists m \, \exists y \in \cal{AS}(y \sqsubseteq s_m(Y) \wedge s_n(X) \leq_\fin y).
    \end{align*}
    All statements are arithmetic w.r.t. $t$.
    \begin{enumerate}
        \item Note that $\alpha$ induces some $f : \cal{R} \to \S$ iff $\lim_{n\to\infty} \lh(\alpha(r_n(A))) = \infty$ for all $A \in \cal{R}$. That is, it is equivalent to:
        \begin{align*}
            \forall A \in \cal{R} \, \forall x \sqsubseteq A \, \exists y \sqsubseteq A \,(x \sqsubseteq y \wedge \lh(\alpha(x)) < \lh(\alpha(y))).
        \end{align*}

        \item ``$f$ is monotone'' is equivalent to:
        \begin{align*}
            \forall A,B \in \cal{R}(A \leq_\cal{R} B \to f(A) \leq_\S f(B)).
        \end{align*}

        \item ``$f$ is quasi-open'' is equivalent to:
        \begin{align*}
            \forall A \in \cal{R} \, \exists X \in \S \, \forall Y \leq_\S X \, \exists B \leq_\cal{R} A(Y = f(B)).
        \end{align*}

        \item By (2), we know that ``$f$ is monotone'' is $\Pi_1^1(t)$. ``$f$ is cofinal'' is equivalent to:
        \begin{align*}
            \forall X \in \S \, \exists A \in \cal{R}(f(A) \leq_\S X),
        \end{align*}
        which is $\Pi_2^1(t)$. Similarly, the statement in Definition \ref{def:projection}(3) is equivalent to:
        \begin{align*}
            \forall A \in \cal{R} \, \forall X \in \S \, \exists B \in \cal{R}(X \leq_\S f(A) \to (B \leq_\cal{R} A \wedge f(B) \leq_\S X)).
        \end{align*}
    \end{enumerate}
\end{proof}

We are now ready to prove Theorem \ref{thm:no.wp.unconditional.Carlson-Simpson}.

\begin{proof}[Proof of Theorem \ref{thm:no.wp.unconditional.Carlson-Simpson}]
    Let $(\cal{R},\leq_\cal{R},r)$ be an unconditional \textbf{A2}-space.
    \begin{enumerate}
        \item Let $M$ be a ctm of $\ZFC$. Let $\Q = \Coll(\omega_1,\frak{c})$ be the L\'{e}vy collapse, whose conditions are countable partial functions from $\omega_1$ to $\frak{c}$. This forcing is $\sigma$-closed and forces $\CH$ (if $M$ already satisfies $\CH$, then let $\Q$ be any trivial $\sigma$-closed forcing). Let $\dot{\P}$ be the $\Q$-name of the ccc poset given by Theorem \ref{thm:Carlson-Simpson.distributivity.omega_1} for some regular $\kappa \geq \aleph_2$, and let $G$ be $\Q * \dot{\P}$-generic over $M$. Since $\MA$ implies that $\frak{p} = \frak{c}$, by Lemma \ref{lem:unconditional.p.leq.h_R}, $M[G] \models \omega_1 = \frak{h}_{\E_\infty} < \frak{h}_\cal{R} = \frak{c}$. Therefore, $(\cal{R},\leq,r)$ does not weakly project into $(\E_\infty,\leq,r)$ in $M[G]$, for otherwise $\frak{h}_\cal{R} \leq \frak{h}_{\E_\infty}$ holds by Lemma \ref{lem:wp.passes.distributivity} (as $(\E_\infty,\leq^*)$ is separative, see \cite[Lemma 2.2]{CKMW00}), a contradiction.
        
        \item Suppose for a contradiction that there is a projection $f : \cal{R} \to \E_\infty$. By Proposition \ref{prop:ab.induced.by.sqsubseteq-preserving}, there exists some $A \in \cal{R}$ such that $f\restrictedto[A]$ is continuous and induced by some $\sqsubseteq$-preserving $\alpha : \cal{AR}\restrictedto A \to \E_{<\infty}\restrictedto f(A)$. Since $f$ is a projection, so is $f\restrictedto[A] : [A] \to [f(A)]$. Note that $[A]$ is again unconditional, that $[f(A)]$ is again $\sigma$-closed, and that $\frak{h}_{[f(A)]} = \frak{h}_{\E_\infty}$. 
        \begin{itemize}
            \item Let $\frak{h}_{[A]}$ denote the distributivity number of $([A],\leq^*)$. The proof of Lemma \ref{lem:unconditional.p.leq.h_R} shows that $\frak{p} \leq \frak{h}_{[A]}$.
            
            \item There is a biprojection $H : [f(A)] \to \E_\infty$ defined as follows: Given $X \leq f(A)$ in $\E_\infty$, we let $H(X)$ be the unique equivalence relation on $\N$ such that for all $m,n \in \N$, $(m,n) \in H(X)$ iff $(p_m(f(A)),p_n(f(A))) \in X$. Therefore, $\frak{h}_{\E_\infty} = \frak{h}_{[f(A)]}$.
        \end{itemize}

        Let $M$ be the transitive collapse of a countable elementary submodel of some $H_\theta$ containing $\alpha$ and everything else relevant, and let $G$ be a $\Q * \dot{\P}$-generic filter over $M$, where $\Q * \dot{\P}$ is as in (1). By Lemma \ref{lem:projection.complexities}, ``$\alpha$ induces a projection $[A] \to [f(A)]$'' is $\mathbf{\Pi}_2^1$. It holds in $M$ by elementarity, so by Shoenfield's absoluteness it holds in $M[G]$. Thus, $M[G]$ contains a projection $[A] \to [f(A)]$. But the argument in (1) gives $\frak{h}_{[f(A)]} = \frak{h}_{\E_\infty} = \omega_1 < \frak{c} = \frak{h}_{[A]}$ in $M[G]$, a contradiction.

        \item Assume the existence of large cardinals so that $\mathbf{\Sigma}_4^1$ generic absoluteness holds, and suppose for a contradiction that there is a weak projection $f : \cal{R} \to \E_\infty$. We repeat the argument in (2) (with $\theta$ large enough for $M$ to satisfy the large cardinal hypothesis), replacing projections with weak projections. By Lemmas \ref{lem:ab.iff.qo.TRS} and \ref{lem:projection.complexities}, ``$\alpha$ induces a weak projection $[A] \to [f(A)]$'' is $\mathbf{\Pi}_4^1$. The remainder of the argument in (2) then applies verbatim.
    \end{enumerate}
\end{proof}

Removing the large cardinal hypothesis appears difficult, as the methods used in the proofs of the earlier theorems do not apply. Under $\CH$, the Carlson--Simpson space and all unconditional spaces are forcing equivalent, so models such as Brendle's, which distinguish the Carlson--Simpson space from the others, are needed, but such models are not well-studied. While canonisation theorems for equivalence relations exist (see \cite{ER52,PR82,PSV86,L07,D13,DT14,V14,DT15,DMT17}), canonisation of functions on unconditional spaces other than $\FIN^{[\infty]}$ is less studied.

\begin{qn}
\label{qn:no.wp.unconditional.Carlson-Simpson.ZFC}
    Does any of the following spaces weakly project into $(\E_\infty,\leq,r)$?
    \begin{enumerate}
        \item $(\FIN_k^{[\infty]},\leq,r)$, where $k \geq 2$.
        
        \item $(\FIN_{\pm k}^{[\infty]},\leq,r)$.
        
        \item $(W_{Lv}^{[\infty]},\leq,r)$, where $L$ is finite.
        
        \item $(E^{[\infty]},\leq,r)$, where $\mathbb{F} \neq \mathbb{F}_2$.
    \end{enumerate}
\end{qn}

\"{O}zalp's recent canonisation results in \cite{O26} may answer Question \ref{qn:no.wp.unconditional.Carlson-Simpson.ZFC}(1) for $k = 2$ in the negative.

The lack of existing results pertaining to the initial Tukey structure of Ramsey ultrafilters of $\E_\infty$ also makes Question \ref{qn:no.wp.unconditional.Carlson-Simpson.ZFC} difficult to answer. We saw in \S\ref{sec:non-existence.by.Tukey.type.theory} and \S\ref{sec:projections.Gowers.spaces} that the initial Tukey structure of Ramsey ultrafilters of topological Ramsey spaces is extremely helpful in determining the existence of weak projections. This leads to the following question:

\begin{qn}
    What is the initial Tukey structure of a Ramsey ultrafilter of $\E_\infty$?
\end{qn}

\subsubsection{Absoluteness of weak projections}
We discuss an alternative proof of Corollary \ref{cor:no.wp.Ellentuck.Carlson-Simpson}, which adapts the proof of Theorem \ref{thm:no.wp.unconditional.Carlson-Simpson}(2) by showing that being a weak projection from the Ellentuck space is a $\mathbf{\Pi}_2^1$ property, and hence is absolute. This suggests a possible approach to Question \ref{qn:no.wp.unconditional.Carlson-Simpson.ZFC}, namely to reduce the complexity of the property of being a weak projection between the spaces in question.

\begin{lemma}
\label{lem:projection.from.Ellentuck.complexity}
    Let $(\S,\leq,s)$ be a \textbf{wA2}-space such that $\cal{AS}$ is countable, and suppose that $(\S,\leq,s)$ satisfies the $(\varepsilon)$-properties. Let $\alpha : [\N]^{<\infty} \to \cal{AS}$ be a $\sqsubseteq$-preserving map, and let $t$ be a real that codes $\alpha$ and $\leq_\fin$. Then ``$\alpha$ induces a (continuous) monotone quasi-open $f : [\N]^\infty \to \S$'' is $\Pi_2^1(t)$.
\end{lemma}

\begin{proof}
    We saw from Lemma \ref{lem:projection.complexities} that ``$\alpha$ induces a (continuous) monotone $f : [\N]^\infty \to \S$'' is $\Pi_1^1(t)$. Thus, it suffices to show that ``$f$ is quasi-open'' is $\Pi_2^1(t)$. Note that $f$ is quasi-open iff the following statement holds:
    \begin{align*}
        \forall A \in [\N]^\infty \, \exists X \in \S \, [X] \subseteq f``[A].
    \end{align*}

    \begin{claim}
        $[X] \subseteq f``[A]$ iff for all $x \in \cal{AS}\restrictedto X$, there exists some $a \in [A]^{<\infty}$ such that $x \sqsubseteq \alpha(a)$.
    \end{claim}

    \begin{midproof}
        \underline{$\implies$:} Suppose that $[X] \subseteq f``[A]$. Let $x \in \cal{AS}\restrictedto X$, and let $Y \in [x,X]$. Then $x \sqsubseteq Y = f(B) = \lim_{n\to\infty} \alpha(r_n(B))$ for some $B \leq A$, so $x \sqsubseteq \alpha(r_n(B))$ for some $n$, with $r_n(B) \in [A]^{<\infty}$.

        \underline{$\impliedby$:} Suppose that for all $x \in \cal{AS}\restrictedto X$, there exists some $a \in [A]^{<\infty}$ such that $x \sqsubseteq \alpha(a)$. Fix some $Y \leq X$. By the $(\varepsilon_1)$-property, for each $n$, we may let $x_n \in \cal{AS}\restrictedto[s_n(Y),Y]$ be such that $x_n \neq s_{\lh(x_n)}(Y)$. By our assumption, there exists some $a_n \in [A]^{<\infty}$ such that $x_n \sqsubseteq \alpha(a_n)$. Consider the sequence $(a_n)_{n<\omega}$ in $[A]^{\leq\infty}$ under the product topology (by identifying $[A]^{\leq\infty}$ with $2^\N$, cf. the paragraph before Definition \ref{def:sqsubseteq-preserving}). $[A]^{\leq\infty}$ is compact under this topology, so there is a subsequence $(a_{n_k})_{k<\omega}$ that converges to some $B \in [A]^{\leq\infty}$. Observe that it suffices to show that $B \in [A]$ --- if this holds, then for all $m$, $r_m(B) \sqsubseteq a_{n_k}$ for all sufficiently large $k$, so $\alpha(r_m(B)) \sqsubseteq \alpha(a_{n_k})$ and $s_{n_k}(Y) \sqsubseteq x_{n_k} \sqsubseteq \alpha(a_{n_k})$. Choosing $k$ so that $n_k \geq \lh(\alpha(r_m(B)))$, we have that $\alpha(r_m(B)) \sqsubseteq s_{n_k}(Y) \sqsubseteq Y$. Since $\alpha$ induces $f$, we have that $f(B) = Y$. Therefore, it remains to show that $B$ cannot be finite.

        Suppose for a contradiction that $B$ is finite. Since $l_k := \lh(x_{n_k}) \to \infty$ and $x_{n_k} \sqsubseteq \alpha(a_{n_k})$, we have that $a_{n_k} \neq B$ for $k$ large enough. Since $a_{n_k} \to B$ and $B$ is finite,  we have that, for $k$ large enough, $a_{n_k} = B \cup b_k$ for some $\emptyset \neq b_k \in [A]^{<\infty}$ such that $\max(B) < \min(b_k)$, and $\min(b_k) \to \infty$. Pick $k < k'$ and $B' \in [A]$ such that $n_{k'} \geq \max\{l_k,\depth_Y(x_{n_k})\}$ and:
        \begin{align*}
            \max(b_k) < \min(b_{k'}) \leq \max(b_{k'}) < \min(B'),
        \end{align*}
        and let:
        \begin{align*}
            C &:= B \cup b_k \cup b_{k'} \cup B', \\
            C' &:= B \cup b_{k'} \cup B'.
        \end{align*}
        Since $f$ is monotone, $f(C') \leq f(C)$. Observe that $a_{n_k} \sqsubseteq C$, so $x_{n_k} \sqsubseteq \alpha(a_{n_k}) \sqsubseteq f(C)$. Similarly, $a_{n_{k'}} \sqsubseteq C'$, so $s_{n_{k'}}(Y) \sqsubseteq x_{n_{k'}} \sqsubseteq \alpha(a_{n_{k'}}) \sqsubseteq f(C')$. Since $l_k \leq n_{k'}$, we have that:
        \begin{align*}
            s_{l_k}(f(C')) = s_{l_k}(Y) \neq x_{n_k}.
        \end{align*}
        In other words, $x_{n_k} \sqsubseteq f(C)$, $f(C') \leq f(C)$ but $x_{n_k} \not\sqsubseteq f(C')$. By the $(\varepsilon_0)$-property, this implies that $x_{n_k} \notin \cal{AS}\restrictedto f(C')$. However, if $m := \depth_Y(x_{n_k})$, then $m \leq n_{k'}$, so $x_{n_k} \leq_\fin s_m(Y) = s_m(f(C'))$, which is a contradiction.
    \end{midproof}

    Hence, the statement ``$[X] \subseteq f``[A]$'' is arithmetic w.r.t. $t$, so ``$f$ is quasi-open'' is $\Pi_2^1(t)$.
\end{proof}

\begin{proof}[Proof of Corollary \ref{cor:no.wp.Ellentuck.Carlson-Simpson}]
    The proof of Theorem \ref{thm:no.wp.unconditional.Carlson-Simpson}(2) may be repeated, with projections replaced by weak projections. Indeed, if $f : [\N]^\infty \to \E_\infty$ is a weak projection, then by Proposition \ref{prop:ab.induced.by.sqsubseteq-preserving}, there exists some $A \in [\N]^\infty$ such that $f\restrictedto[A]$ is continuous and induced by some $\sqsubseteq$-preserving $\alpha : [A]^{<\infty} \to \E_{<\infty}\restrictedto f(A)$. Thus, by Lemmas \ref{lem:ab.iff.qo.TRS} and \ref{lem:projection.from.Ellentuck.complexity}, ``$\alpha$ induces a weak projection $[A] \to [f(A)]$'' is $\Pi_2^1$ in a real coding $\alpha$. Shoenfield's absoluteness therefore applies directly, and no large cardinal hypotheses are required.
\end{proof}

It is far from clear whether Lemma \ref{lem:projection.from.Ellentuck.complexity} holds for arbitrary \textbf{wA2}-spaces in place of the Ellentuck space. In the proof above, we relied on the fact that every sequence $(a_n)_{n=0}^\infty$ in $[\N]^{<\infty}$ has a subsequence $(a_{n_k})_{k=0}^\infty$ that either converges to an element of $[\N]^\infty$, or satisfies that for all $k < k'$, there exist $C' \subseteq C$ in $[\N]^\infty$ such that $a_{n_k} \sqsubseteq C$ and $a_{n_{k'}} \sqsubseteq C'$. In the proof, the latter case corresponds to $B$ being finite, and $C'$ is obtained by removing the finite set $b_k$ from $C$. The analogous property (with $\subseteq$ replaced by $\leq$) fails for arbitrary \textbf{A2}-spaces --- for instance, if $a_n := (\{0,n\})$ for each $n \geq 1$, then $(a_n)_{n=1}^\infty$ is a sequence in $\FIN^{[<\infty]}$ with no such subsequence, as if $C' \leq C$, $a_{n_k} \sqsubseteq C$ and $a_{n_{k'}} \sqsubseteq C'$, then the block of $C'$ containing $0$ would have to contain $n_k$. This leads to the following question, a positive answer to which for unconditional \textbf{wA2}-spaces would answer Question \ref{qn:no.wp.unconditional.Carlson-Simpson.ZFC} in the negative:

\begin{qn}
    Is the quasi-open property for continuous monotone functions $f : \cal{R} \to \S$ absolute?
\end{qn}

Given two \textbf{A2}-spaces $(\cal{R},\leq_\cal{R},r)$ and $(\S,\leq_\S,s)$, the absoluteness problem of the quasi-open property may be rephrased in the following manner: We may assign to each continuous monotone function $f : \cal{R} \to \S$ a unique code $t \in 2^\omega$, and write $f_t := f$. Let:
\begin{align*}
    \cal{C}_{\func{op}}(\cal{R},\S) := \{t \in 2^\omega : f_t \text{ is quasi-open}\}.
\end{align*}
Observe that the $\ZFC$-correctness of the quasi-open property for continuous monotone functions $f : \cal{R} \to \S$ is equivalent to the fact that $\cal{C}_{\func{op}}(\cal{R},\S)$ is $\mathbf{\Delta}_2^1$. On the one hand, if $\cal{C}_{\func{op}}(\cal{R},\S)$ is $\mathbf{\Delta}_2^1$, then the quasi-open property is absolute by Shoenfield absoluteness. On the other hand, if the quasi-open property is $\ZFC$-correct, then for any $t \in 2^\omega$:
\begin{align*}
    &\iffbreak \text{$f_t$ is quasi-open} \\
    &\iff \exists M \text{ ctm of $\ZFC$}(t \in M \wedge M \models \text{$f_t$ is quasi-open}) \\
    &\iff \forall M \text{ ctm of $\ZFC$}(t \in M \to M \models \text{$f_t$ is quasi-open}),
\end{align*}
Sabok showed in \cite{Sab12} that the set of codes for Ramsey positive analytic subsets of $[\N]^\infty$ is $\mathbf{\Sigma}_2^1$-complete, so it is not $\mathbf{\Delta}_2^1$. Wang proved the analogous result for $G_\delta$ subsets of \emph{well-indexed} topological Ramsey spaces in \cite{W25} (see also \S\ref{sec:strong.projections.and.well-indexed}). Should the analogue of Lemma \ref{lem:projection.from.Ellentuck.complexity} fail when the Ellentuck space is replaced by another topological Ramsey space, one might show this by finding a pair of topological Ramsey spaces for which $\cal{C}_{\func{op}}(\cal{R},\S)$ is not $\mathbf{\Delta}_2^1$.

\subsection{Carlson--Simpson space with variables}
\label{sec:projections.CS-variables}
We now turn our attention to the Carlson--Simpson space with variables. We begin with several simple observations. Since forcing with a $\sigma$-closed \textbf{A2}-space does not add a real, but forcing with $(\E_\infty^L,\leq)$ does for $|L| \geq 2$, we get the following:

\begin{theorem}
\label{thm:no.wp.sigma-closed.CS-variables}
    Let $(\cal{R},\leq,\leq^*,r)$ be a $\sigma$-closed \textbf{A2}-space, and let $L$ be a finite alphabet such that $|L| \geq 2$. Then $(\cal{R},\leq,r)$ does not weakly project into $(\E_\infty^L,\leq,r)$.
\end{theorem}

Furthermore, the Carlson--Simpson space with variables projects onto the Carlson--Simpson space without variables (see Example \ref{ex:CS-variables.merging.alphabets}), so we can also conclude the following from Corollary \ref{cor:no.wp.strong.subtrees.sigma-closed.examples}(6):

\begin{corollary}
\label{cor:no.wp.strong.subtrees.CS-variables}
    Let $L$ be a finite alphabet. Then $(\S_\infty,\subseteq,r)$ does not weakly project into $(\E_\infty^L,\leq,r)$.
\end{corollary}

In this subsection, we present two analogues of results in \S\ref{sec:projections.Carlson-Simpson} for the Carlson--Simpson space with variables. We first show that $(\E_\infty^L,\leq,r)$ has finite amalgamable covers (Definitions \ref{def:amalgamable.cover} and \ref{def:finite.amalgamable.covers}, Example \ref{ex:CS-variables.has.FACs}), which implies the following theorem:

\begin{theorem}
\label{thm:no.wp.CS-variables.FIN}
    Let $L$ be a finite alphabet. Then $(\E_\infty^L,\leq,r)$ does not weakly project into $(\FIN^{[\infty]},\leq,r)$. 
\end{theorem}

By considering composition with the projections in Examples \ref{ex:tetris.operation}, \ref{ex:absolute.value.FIN_pmk} and \ref{ex:support.variable.words}, we can conclude that:

\begin{corollary}
\label{cor:no.wp.CS-variables.FIN_k.Hales-Jewett}
    Let $L$ be a finite alphabet. Then $(\E_\infty^L,\leq,r)$ does not weakly project into the following:
    \begin{enumerate}
        \item $(\FIN_k^{[\infty]},\leq,r)$ for $k \geq 2$.
        
        \item $(\FIN_{\pm k}^{[\infty]},\leq,r)$ for $k \geq 2$.
        
        \item $(W_{Kv}^{[\infty]},\leq,r)$, where $K$ is finite.
    \end{enumerate}
\end{corollary}

Since $(\E_\infty,\leq,\leq^*,r)$ is $\sigma$-closed, by Theorem \ref{thm:no.wp.sigma-closed.trees}, it does not weakly project into $(\S_\infty,\subseteq,r)$. One may ask if the same holds for $(\E_\infty^L,\leq,r)$, for which forcing adds a real whenever $|L| \geq 2$. The following theorem answers this in the affirmative.

\begin{theorem}
\label{thm:no.wp.CS-variables.trees}
    Let $L$ be a finite alphabet. Then $(\E_\infty^L,\leq,r)$ does not weakly project into $(\S_\infty,\subseteq,r)$ or $(\mathbb{S},\subseteq,r)$.
\end{theorem}

We conclude this subsection by showing that the spaces $(\E_\infty^L,\leq,r)$ form a chain of topological Ramsey spaces ordered by the size of $L$. Recall, from Example \ref{ex:CS-variables.merging.alphabets}, that $(\E_\infty^L,\leq,r)$ projects onto $(\E_\infty^K,\leq,r)$ for all finite alphabets $K$ and $L$ such that $2 \leq |K| \leq |L|$. We show that this cannot be reversed.

\begin{theorem}
\label{thm:CS-variables.alphabet.hierarchy}
    Let $K$ and $L$ be finite alphabets such that $|K|,|L| \geq 2$. Then $(\E_\infty^L,\leq,r)$ weakly projects into $(\E_\infty^K,\leq,r)$ iff $|K| \leq |L|$.
\end{theorem}

As mentioned in Example \ref{ex:CS-variables}, if $|L| = 1$ then $(\E_\infty^L,\leq,r)$ is a reformulation of $(\E_\infty,\leq,r)$. Therefore, throughout this subsection, we fix a finite alphabet $L$ with $|L| \geq 2$.

\subsubsection{Existence of finite amalgamable covers}
\begin{example}
\label{ex:CS-variables.has.FACs}
    The Carlson--Simpson space with variables $(\E_\infty^L,\leq,r)$ has finite amalgamable covers.
\end{example}

\begin{proof}
    The construction is similar to the case of the Carlson--Simpson space (Example \ref{ex:Carlson-Simpson.has.FACs}). Let $A \in \E_\infty^L$, $a \sqsubseteq A$ and $b \preceq_\fin a$. Let:
    \begin{align*}
        \F := \{c \in r_{\lh(b)+1}[b,A] : \depth_A(c) = \lh(a) + 1\}.
    \end{align*}
    Note that $\F$ is finite. For each $c \in \F$, let $H_c$ be the equivalence relation on $\N \cup L$ generated by $A \cup c$, so that $H_c$ is the $\leq$-maximal element of $[c,A]$. We shall show that $\{H_c\}_{c\in\F}$ is an amalgamable cover for $b$ in $[a,A]$. Recall that for $B \in \E_\infty^L$, $p_n(B) = \min(B(n))$ (see Example \ref{ex:CS-variables}).
    \begin{enumerate}
        \item Let $B \in [b,A]$ such that $\depth_A(r_{\lh(b)+1}(B)) > \lh(a)$, and let $c := B\restrictedto(L \cup p_{\lh(a)}(A))$. Since $b \preceq_\fin a$ and $\depth_A(r_{\lh(b)+1}(B)) > \lh(a)$, we have that:
        \begin{align*}
            \dom(b) \subseteq \dom(a) \subsetneq L \cup p_{\lh(a)}(A) \subseteq L \cup p_{\lh(b)}(B),
        \end{align*}
        so $c = r_{\lh(b)+1}(B)\restrictedto(L \cup p_{\lh(a)}(A))$ has exactly $\lh(b)$ equivalence classes disjoint from $L$. Therefore, $c \in \F$ and $B \leq H_c$.

        \item Let $A' \in [a,A]$ and $B \in [b,A']$ such that $B \leq H_c$ for some $c \in \F$. Since $b \sqsubseteq B$ and $B \leq H_c$, we have that $B\restrictedto\dom(c) = c$. Let $Y$ be the union of $A'(0),\dots,A'(\lh(a)-1)$ and the equivalence classes of $A'$ meeting $L$. For each equivalence class $X$ of $B$, let $x_X$ be the unique element of $X \cap L$ if $X$ meets $L$, and let $x_X := \min(X)$ otherwise. Let $B'$ be the equivalence relation on $\N \cup L$ generated by $A'$ and all pairs $(x_X,x)$ such that $X$ is an equivalence class of $B$ and $x \in X \setminus Y$. Then $B'\restrictedto L$ is the equality relation and $A'(i) \subseteq B'(i)$ for all $i < \lh(a)$, which implies that $a \sqsubseteq B'$, so $B' \in [a,A']$. Since $B\restrictedto\dom(c) = c$, $B$ is the equivalence relation on $\N \cup L$ generated by $B' \cup c$, so $[b,B'] \cap [H_c] \subseteq [b,B]$.
    \end{enumerate}
\end{proof}

\begin{proof}[Proof of Theorem \ref{thm:no.wp.CS-variables.FIN}]
    It suffices to show that $(\E_\infty^L,\leq,r)$ satisfies the hypothesis of Lemma \ref{lem:wp.Ellentuck.controlled}, so that Proposition \ref{prop:no.wp.FIN.controlled} applies. We showed in Example \ref{ex:CS-variables.has.FACs} that $(\E_\infty^L,\leq,r)$ has finite amalgamable covers, and the map $p : \E_\infty^L \to [\N]^\infty$ in Example \ref{ex:minimal.representatives.not.L} satisfies the requirements of $\mu$ in Lemma \ref{lem:wp.Ellentuck.controlled} for $a_0 := \emptyset$, where for all $a \in \E_{<\infty}^L$ we may let:
    \begin{align*}
        \alpha(a) := \{\dom(r_k(a)) \setminus L : 1 \leq k \leq \lh(a)\}.
    \end{align*}
    Once again, we leave the details of the verification of requirements to the reader.
\end{proof}

\subsubsection{Forcing properties of $(\E_\infty^L,\leq)$}
Given a partition $\cal{L}$ of $L$, we let $q_{\cal{L}} : L \to \cal{L}$ be the unique map sending each $w \in L$ to the element of $\cal{L}$ containing it. The partition $\cal{L}$ also induces a natural map $\sigma_{\cal{L}} : \E_\infty^L \to \E_\infty^{\cal{L}}$ (Example \ref{ex:CS-variables.merging.alphabets}), where for each $A \in \E_\infty^L$, $\sigma_{\cal{L}}(A)$ is the equivalence relation on $\N \cup \cal{L}$ obtained by merging, for each $u \in \cal{L}$, the equivalence classes of all $w \in u$ into one, which is then identified with $u$. 

While $(\E_\infty^L,\leq)$ is not generally real-degree minimal (cf. Lemma \ref{lem:strong.subtrees.forcing.minimal.real}), the models generated by its new reals form a hierarchy indexed by the partitions of $L$ (see also Lemma \ref{lem:CS-variables.distinct.degrees}):

\begin{lemma}
\label{lem:CS-variables.forcing.reals}
    Let $M$ be a ctm of $\ZFC$, and let $G$ be $\E_\infty^L$-generic over $M$. If $x \in \omega^\omega \cap (M[G] \setminus M)$, then there exists some partition $\cal{L}$ of $L$ such that $|\cal{L}| \geq 2$ and $M[x] = M[\sigma_{\cal{L}}``(G)^\uparrow]$.
\end{lemma}

In particular, if $|L| = 2$ then the partition $\cal{L}$ must be discrete, so $(\E_\infty^L,\leq)$ is real-degree minimal if $|L| = 2$.

We recall from Examples \ref{ex:CS-variables.reformulation.variable-words} and \ref{ex:CS-variables.reformulation.trees} that $(\E_\infty^L,\leq,r)$ is equivalent to its tree formulation $(\cal{SE}_\infty^L,\leq,r)$ in the sense that the map $\Theta \circ \Phi : \E_\infty^L \to \cal{SE}_\infty^L$ is a biprojection (see Definition \ref{def:biprojection}), so $B \leq A$ iff $\Theta(\Phi(B)) \subseteq \Theta(\Phi(A))$. The elements of $\cal{SE}_\infty^L$ are subtrees of $L^{<\omega}$, and we endow $L^\omega$ with the product topology. For the rest of this subsection, we identify each $A \in \E_\infty^L$ with the tree $\Theta(\Phi(A))$, and we use the terminology for trees introduced in Example \ref{ex:strong.subtrees}. 

Recall also from Notation \ref{not:tree.restriction} that for a tree $T$, $\func{br}(T)$ is the set of all branches of $T$, and $T\restrictedto s = \{t \in T : t \sqsubseteq s \vee s \sqsubseteq t\}$ for each $s \in T$. Unravelling the definitions of $\Phi$ and $\Theta$, we have that for all $A \in \E_\infty^L$:
\begin{enumerate}
    \item $\func{br}(A)$ is the set of all $y \in L^\omega$ such that for all $m,n \in \N$ and $w \in L$, if $(w,n) \in A$ then $y(n) = w$, and if $(m,n) \in A$ then $y(m) = y(n)$.

    \item $p(A)$ is the set of all splitting levels of $A$, and for each $n < \omega$ and $s \in L^n$, there is a unique $t \in \func{split}_{n+1}(A)$ such that $t(p_k(A)) = s(k)$ for all $k < n$. In this case, we write $A[s] := A\restrictedto t$.
\end{enumerate}
(1) implies that for all $A \in \E_\infty^L$, partitions $\cal{L}$ of $L$ such that $|\cal{L}| \geq 2$ and $x \in L^\omega$, we have that:
\begin{align*}
\label{eq:branch.equivalence}
    x \in \func{br}(A) \implies q_{\cal{L}} \circ x \in \func{br}(\sigma_{\cal{L}}(A)). \tag{$*$}
\end{align*}
By (1), $A[s]$ is the equivalence relation coarsening $A$ obtained by merging the equivalence class $A(k)$ into that of $s(k)$ for each $k < n$, so $A[s] \in \E_\infty^L$ and $A[s] \leq A$. Therefore, by Lemma \ref{lem:tree.generic.real}, we have:

\begin{lemma}
\label{lem:CS-variables.generic.real.general}
    Let $M$ be a ctm of $\ZFC$, and let $G$ be $\E_\infty^L$-generic over $M$. Let:
    \begin{align*}
        g_G := \bigcup\{s \in L^{<\omega} : \forall A \in G \, s \in A\}.
    \end{align*}
    Then $g_G \in L^\omega \cap M[G]$ and for all $A \in \E_\infty^L \cap M$, $A \in G$ iff $g_G \in \func{br}(A)$, so $M[G] = M[g_G]$.
\end{lemma}

\begin{lemma}
\label{lem:CS-variables.merging.letters}
    Let $\cal{L}$ be a partition of $L$ such that $|\cal{L}| \geq 2$, and let $M$ be a ctm of $\ZFC$. If $G$ is $\E_\infty^L$-generic over $M$, then $G' := \sigma_{\cal{L}}``(G)^\uparrow$ is $\E_\infty^{\cal{L}}$-generic over $M$, and $g_{G'} = q_{\cal{L}} \circ g_G$.
\end{lemma}

\begin{proof}
    By Example \ref{ex:CS-variables.merging.alphabets} and Lemma \ref{lem:ab.equivalence}, $G'$ is $\E_\infty^{\cal{L}}$-generic over $M$. Since $g_G \in \func{br}(A)$ for all $A \in G$, we have that by (\ref{eq:branch.equivalence}), $q_{\cal{L}} \circ g_G \in \func{br}(X)$ for all $X \in G'$. Since $g_{G'}$ is the unique branch that lies in $\func{br}(X)$ for all $X \in G'$, $g_{G'} = q_{\cal{L}} \circ g_G$.
\end{proof}

\begin{lemma}
\label{lem:CS-variables.fusion}
    Let $(\D_n)_{n<\omega}$ be a sequence of dense open subsets of $\E_\infty^L$, and let $A \in \E_\infty^L$. Then there exists some $B \leq A$ such that $B\restrictedto t \in \D_n$ for all $n$ and $t \in \func{split}_{n+1}(B)$.
\end{lemma}

\begin{proof}
    This follows from a routine fusion argument similar to that in the proof of Lemma \ref{lem:strong.subtrees.forcing.minimal}. We provide an alternative proof by appealing to Lemma \ref{lem:R.is.semiselective}, working in the equivalence-relation formulation of $\E_\infty^L$, where the conclusion states that there exists some $B \leq A$ such that $B[s] \in \D_n$ for all $n$ and $s \in L^n$. For each $b \in \E_{<\infty}^L\restrictedto A$ and $s \in L^{\lh(b)}$, we let:
    \begin{align*}
        \cal{K}_{b,s} := \{C \in [b,A] : C[s] \in \D_{\lh(b)}\}.
    \end{align*}

    \begin{claim}
        $\cal{K}_{b,s}$ is dense open in $[b,A]$.
    \end{claim}

    \begin{midproof}
        Let $n := \lh(b)$. If $C \in \cal{K}_{b,s}$ and $C' \in [b,C]$, then $C(k) \subseteq C'(k)$ for all $k < n$ as $C' \leq C$ and $r_n(C') = r_n(C)$, so $C'[s] \leq C[s] \in \D_n$ and hence $C' \in \cal{K}_{b,s}$. To see that $\cal{K}_{b,s}$ is dense: let $B \in [b,A]$, and let $D \leq B[s]$ be such that $D \in \D_n$. Let $B'$ be the equivalence relation on $\N \cup L$ whose equivalence classes are:
        \begin{itemize}
            \item $B(k)$ for each $k < n$;

            \item $[w]_D \setminus \bigcup_{k<n} B(k)$ for each $w \in L$, where $[w]_D$ is the equivalence class of $w$ in $D$;

            \item the equivalence classes of $D$ disjoint from $L$.
        \end{itemize}
        Then $B' \in [b,B]$ and $B'[s] = D$, so $B' \in \cal{K}_{b,s}$.
    \end{midproof}

    Now let $\cal{K}_b := \bigcap_{s \in L^{\lh(b)}} \cal{K}_{b,s}$, which is dense open in $[b,A]$ as it is a finite intersection of dense open subsets of $[b,A]$. By Lemma \ref{lem:R.is.semiselective}, there exists some $B \in [A]$ that diagonalises $\vec{\cal{K}}$. Then $B \in \cal{K}_{r_n(B)}$ for all $n$, so $B[s] \in \D_n$ for all $n$ and $s \in L^n$.
\end{proof}

We shall use the following canonisation theorem of Doucha.

\begin{theorem}[Doucha, {\cite[Theorem 2.1]{D13}}]
\label{thm:CS-variables.canonisation}
    Let $L$ be a finite alphabet such that $|L| \geq 2$, let $A \in \E_\infty^L$, and let $E$ be an analytic equivalence relation on $\func{br}(A)$. Then there exist some $B \leq A$ and some partition $\cal{L}$ of $L$ such that for all $y,y' \in \func{br}(B)$, $(y,y') \in E$ iff $q_{\cal{L}} \circ y = q_{\cal{L}} \circ y'$.
\end{theorem}

\begin{lemma}
\label{lem:CS-variables.reading.of.names}
    Let $M$ be a ctm of $\ZFC$. In $M$, let $\dot{x}$ be an $\E_\infty^L$-name for an element of $\omega^\omega$, and let $A \in \E_\infty^L$. Then there exist (in $M$) some $B \leq A$, some partition $\cal{L}$ of $L$ and some continuous injection $\phi : \{q_{\cal{L}} \circ y : y \in \func{br}(B)\} \to \omega^\omega$ such that for every $\E_\infty^L$-generic filter $G$ over $M$ with $B \in G$:
    \begin{align*}
        i_G(\dot{x}) = \phi(q_{\cal{L}} \circ g_G).
    \end{align*}
\end{lemma}

\begin{proof}
    For each $n$, the set of all $C \in \E_\infty^L$ that decide $\dot{x}\restrictedto n$ is dense open. By Lemma \ref{lem:CS-variables.fusion}, there exists some $B_0 \leq A$ such that for all $n$ and $t \in \func{split}_{n+1}(B_0)$, $B_0\restrictedto t$ decides $\dot{x}\restrictedto n$, and we let $u_t \in \omega^n$ be such that $B_0\restrictedto t \forces \dot{x}\restrictedto n = u_t$. Observe that if $t \sqsubseteq t'$ both split in $B_0$, then $B_0\restrictedto t' \leq B_0\restrictedto t$, so $u_t \sqsubseteq u_{t'}$. Thus, we may define a continuous function $\phi_0 : \func{br}(B_0) \to \omega^\omega$ by:
    \begin{align*}
        \phi_0(y) := \bigcup_{n<\omega} u_{y\restrictedto p_n(B_0)}.
    \end{align*}

    \begin{claim}
        If $G$ is $\E_\infty^L$-generic over $M$ and $B_0 \in G$, then $i_G(\dot{x}) = \phi_0(g_G)$.
    \end{claim}

    \begin{midproof}
        Note that for each $n < \omega$, the set of all $C \leq B_0$ such that $C \leq B_0\restrictedto t$ for some $t \in \func{split}_{n+1}(B_0)$ is dense below $B_0$, so $B_0\restrictedto t \in G$ for some $t \in \func{split}_{n+1}(B_0)$. Then $g_G \in \func{br}(B_0\restrictedto t)$, so $t = g_G\restrictedto p_n(B_0)$. Since $B_0\restrictedto t \in G$ and $B_0\restrictedto t \forces \dot{x}\restrictedto n = u_t$, we have that:
        \begin{align*}
            i_G(\dot{x})\restrictedto n = u_t = \phi_0(g_G)\restrictedto n,
        \end{align*} 
        so the claim follows as $n$ is arbitrary.
    \end{midproof}

    By Theorem \ref{thm:CS-variables.canonisation}, there exist some $B \leq B_0$ and some partition $\cal{L}$ of $L$ such that for all $y,y' \in \func{br}(B)$, $\phi_0(y) = \phi_0(y')$ iff $q_{\cal{L}} \circ y = q_{\cal{L}} \circ y'$. Thus, we may define an injection $\phi : \{q_{\cal{L}} \circ y : y \in \func{br}(B)\} \to \omega^\omega$ by stipulating that for all $y \in \func{br}(B)$:
    \begin{align*}
        \phi(q_{\cal{L}} \circ y) := \phi_0(y).
    \end{align*}
    Note that $\phi$ is continuous, as the map $y \mapsto q_{\cal{L}} \circ y$ is closed and $\func{br}(B)$ is compact. If $G$ is an $\E_\infty^L$-generic filter over $M$ such that $B \in G$, then $g_G \in \func{br}(B)$, so:
    \begin{align*}
        i_G(\dot{x}) = \phi_0(g_G) = \phi(q_{\cal{L}} \circ g_G).
    \end{align*}
\end{proof}

We are now ready to prove Lemma \ref{lem:CS-variables.forcing.reals}.

\begin{proof}[Proof of Lemma \ref{lem:CS-variables.forcing.reals}]
    Let $\dot{x} \in M$ be an $\E_\infty^L$-name for an element of $\omega^\omega$ such that $i_G(\dot{x}) = x$. By Lemma \ref{lem:CS-variables.reading.of.names}, there exist (in $M$) some $B \in G$, some partition $\cal{L}$ of $L$ and some continuous injection $\phi : \{q_{\cal{L}} \circ y : y \in \func{br}(B)\} \to \omega^\omega$ such that $x = \phi(q_{\cal{L}} \circ g_G)$. We have that $|\cal{L}| \geq 2$ --- otherwise, if $\cal{L} = \{L\}$, then $q_{\cal{L}}$ is constant, so the domain of $\phi$ is a singleton and $x \in M$, a contradiction. Let $G' := \sigma_{\cal{L}}``(G)^\uparrow$. By Lemma \ref{lem:CS-variables.merging.letters}, $x = \phi(g_{G'})$, so $x \in M[G']$. Since $\phi$ is a continuous injection, $g_{G'} \in M[x]$, so $M[G'] \subseteq M[x]$ by Lemma \ref{lem:CS-variables.generic.real.general}.
\end{proof}

\subsubsection{Proof of Theorem \ref{thm:no.wp.CS-variables.trees}}
We are now ready to prove Theorem \ref{thm:no.wp.CS-variables.trees}. We recall several projections of interest here:
\begin{itemize}
    \item (Example \ref{ex:levels.strong.subtrees}) $S : \S_\infty \to [\N]^\infty$, where:
    \begin{align*}
        S(A) = \{\dom(s) : \text{$s$ splits in $A$}\}
    \end{align*}
    is the level set of splitting nodes in $A$.
    
    \item (Example \ref{ex:minimal.representatives}) $p^* : \E_\infty \to [\N^+]^\infty$, where $p^*(A) = \{p_n(A) : n \geq 1\}$ is the set of non-zero minimal representatives of equivalence classes.

    \item (Example \ref{ex:CS-variables.merging.alphabets}) $\sigma : \E_\infty^L \to \E_\infty$, where $\sigma(A) := \sigma_{\{L\}}(A)\restrictedto\N$.

    \item (Example \ref{ex:minimal.representatives.not.L}) $p : \E_\infty^L \to [\N]^\infty$, where $p(A) = \{p_n(A) : n < \omega\}$ is the set of minimal representatives of equivalence classes disjoint from $L$.
\end{itemize}

\begin{proof}[Proof of Theorem \ref{thm:no.wp.CS-variables.trees}]
    We shall follow the outline of the proof of Theorems \ref{thm:no.wp.strong.subtrees.Sacks} and \ref{thm:no.wp.strong.subtrees.sigma-closed.projection}.

    \underline{Sacks forcing:} Suppose for a contradiction that there is a weak projection $f : \E_\infty^L \to \mathbb{S}$. Let $M$ be the transitive collapse of a countable elementary submodel of some $H_\theta$ containing $f$ and everything else relevant, and let $G$ be $\E_\infty^L$-generic over $M$. By Lemma \ref{lem:ab.equivalence}, there exists a Sacks real $x \in M[G]$ over $M$. By Lemma \ref{lem:CS-variables.forcing.reals}, there exists some partition $\cal{L}$ of $L$ such that $|\cal{L}| \geq 2$ and $M[x] = M[G']$, where $G' := \sigma_{\cal{L}}``(G)^\uparrow$. By Lemma \ref{lem:CS-variables.merging.letters} and Example \ref{ex:CS-variables.merging.alphabets}, $G'$ is $\E_\infty^{\cal{L}}$-generic over $M$ and $\cal{W} := \sigma``(G')^\uparrow$ is $\E_\infty$-generic over $M$. Since $(\E_\infty,\leq,\leq^*,r)$ is $\sigma$-closed, $M[\cal{W}]$ and $M$ share the same reals, so $M \subsetneq M[\cal{W}] \subsetneq M[x]$, contradicting the minimality of Sacks forcing.

    \underline{Strong subtrees:} We showed in Example \ref{ex:CS-variables.has.FACs} that $(\E_\infty^L,\leq,r)$ has finite amalgamable covers. Furthermore, the map $p : \E_\infty^L \to [\N]^\infty$ in Example \ref{ex:minimal.representatives.not.L} satisfies the requirements of $\mu$ in Lemma \ref{lem:wp.Ellentuck.controlled} (see the proof of Theorem \ref{thm:no.wp.CS-variables.FIN}).

    Now suppose for a contradiction that there is a weak projection $f : \E_\infty^L \to \S_\infty$. Then $h := S \circ f : \E_\infty^L \to [\N]^\infty$ is a weak projection, so by Lemma \ref{lem:wp.Ellentuck.controlled}, there exist some $A \in \E_\infty^L$ and an injective $\psi : \N \to \N$ such that $\psi``(p(B)) \cap h(B) \neq \emptyset$ for all $B \leq A$. Let $M$ be the transitive collapse of a countable elementary submodel of some $H_\theta$ containing $f$, $A$, $\psi$ and everything else relevant, and let $G$ be $\E_\infty^L$-generic over $M$ such that $A \in G$. We consider the following filters in $M[G]$, which are generic over $M$ by Lemma \ref{lem:ab.equivalence}:
    \begin{align*}
        \U &:= p``(G)^\uparrow \tag*{($[\N]^\infty$-generic)} \\
        H &:= f``(G)^\uparrow \tag*{($\S_\infty$-generic)} \\
        \U' &:= S``(H)^\uparrow = h``(G)^\uparrow \tag*{($[\N]^\infty$-generic)}
    \end{align*}
    Let $x_H$ be as in Lemma \ref{lem:tree.generic.real}, so that $x_H \in \omega^\omega \cap (M[G] \setminus M)$ and $M[x_H] = M[H]$. As before, by Lemma \ref{lem:CS-variables.forcing.reals}, there exists some partition $\cal{L}$ of $L$ such that $|\cal{L}| \geq 2$ and $M[x_H] = M[\sigma_{\cal{L}}``(G)^\uparrow]$. This gives us following filters in $M[G]$, which are also generic over $M$ by Lemma \ref{lem:ab.equivalence}:
    \begin{align*}
        G' &:= \sigma_{\cal{L}}``(G)^\uparrow \tag*{($\E_\infty^{\cal{L}}$-generic)} \\
        \cal{W} &:= \sigma``(G')^\uparrow \tag*{($\E_\infty$-generic)} \\
        \cal{W}^* &:= p^*``(\cal{W})^\uparrow \tag*{($[\N^+]^\infty$-generic)}
    \end{align*}

    \begin{claim}
        $M[\cal{W}^*] = M[\cal{W}] = M[\U'] = M[\U]$.
    \end{claim}

    \begin{midproof}
        \underline{$M[\cal{W}^*] \subseteq M[\cal{W}]$:} Immediate.

        \underline{$M[\cal{W}] \subseteq M[\U']$:} Since $\cal{W} \in M[G'] = M[H]$, $H$ is $\S_\infty$-generic over $M$ and $M[\cal{W}] \neq M[H]$ (as $M[\cal{W}]$ and $M$ share the same reals), we have that $\cal{W} \in M[\U'] = M[S``(H)^\uparrow]$ by Lemma \ref{lem:strong.subtrees.forcing.minimal}.

        \underline{$M[\U'] \subseteq M[\U]$:} Immediate from Lemma \ref{lem:controlled.wps.give.isomorphic.generic.filters}.

        \underline{$M[\U] \subseteq M[\cal{W}^*]$:} Note that for all $B \in \E_\infty^L$, the sets $p(B)$ and $p^*(\sigma(\sigma_{\cal{L}}(B)))$ differ by only finitely many values ($0$ and the least $n \in \N$ such that $(w,n) \in B$ for some $w \in L$). Since $\U = p``(G)^\uparrow$ and $\cal{W}^* = p^*``(\cal{W})^\uparrow = (p^* \circ \sigma \circ \sigma_{\cal{L}})``(G)^\uparrow$ is $[\N^+]^\infty$-generic over $M$, we have that:
        \begin{align*}
            \U = \{Z \in \Po(\N) \cap M : Z \setminus \{0\} \in \cal{W}^*\},
        \end{align*}
        so $\U \in M[\cal{W}^*]$.
    \end{midproof}

    Thus, $\cal{W} \in M[\U] = M[\cal{W}^*]$. By the proof of Corollary \ref{cor:no.wp.strong.subtrees.sigma-closed.examples}(6), for every $X \in \E_\infty$ there exist incompatible $X_0,X_1 \leq X$ such that $p^*(X_0) = p^*(X_1)$. Therefore, we obtain a contradiction as at the end of the proof of Theorem \ref{thm:no.wp.strong.subtrees.sigma-closed.projection}, with $\E_\infty$, $\cal{W}$, $\cal{W}^*$ and $p^*$ in place of $\cal{R}$, $\cal{V}$, $\cal{V}^*$ and $\pi$.
\end{proof}

\subsubsection{Proof of Theorem \ref{thm:CS-variables.alphabet.hierarchy}}
We fix a finite alphabet $K$ such that $|K| \geq 2$. We shall prove Theorem \ref{thm:CS-variables.alphabet.hierarchy} entirely in the language of equivalence relations.

\begin{lemma}
\label{lem:CS-variables.distinct.degrees}
    Let $M$ be a ctm of $\ZFC$, and let $H$ be $\E_\infty^K$-generic over $M$. Let $\cal{K}$ and $\cal{K}'$ be partitions of $K$ with at least two elements each, and suppose that there exist $w_0,w_1 \in K$ which lie in the same element of $\cal{K}$ but in different elements of $\cal{K}'$. Then $\sigma_{\cal{K}'}``(H)^\uparrow \notin M[\sigma_{\cal{K}}``(H)^\uparrow]$.
\end{lemma}

\begin{proof}
    Suppose otherwise. By Example \ref{ex:CS-variables.merging.alphabets} and Lemma \ref{lem:ab.equivalence}, $\sigma_{\cal{K}}``(H)^\uparrow$ is $\E_\infty^{\cal{K}}$-generic over $M$, so there exists some $\E_\infty^{\cal{K}}$-name $\dot{z} \in M$ such that $\sigma_{\cal{K}'}``(H)^\uparrow = i_{\sigma_{\cal{K}}``(H)^\uparrow}(\dot{z})$. Thus, there exists some $A \in H$ such that for every $\E_\infty^K$-generic filter $H'$ over $M$ with $A \in H'$:
    \begin{align*}
        \sigma_{\cal{K}'}``(H')^\uparrow = i_{\sigma_{\cal{K}}``(H')^\uparrow}(\dot{z}).
    \end{align*}
    Let $B_0,B_1 \leq A$ be equivalence relations on $\N \cup K$ obtained from $A$ by merging $A(0)$ into the equivalence class of $w_0$ and of $w_1$ respectively. Note that for each $i < 2$, $A(0) \subseteq [w_i]_C$ for all $C \leq B_i$, where $[w_i]_C$ is the equivalence class of $C$ containing $w_i$. We define $\iota : [B_0] \to [B_1]$ by letting, for all $C \leq B_0$, $\iota(C)$ be the equivalence relation on $\N \cup K$ whose equivalence classes are:
    \begin{itemize}
        \item $[w_0]_C \setminus A(0)$;
        
        \item $[w_1]_C \cup A(0)$;
        
        \item all equivalence classes of $C$ other than $[w_0]_C$ and $[w_1]_C$.
    \end{itemize}
    Note that $\iota$ is a biprojection in $M$, and $\iota^{-1}$ is defined in the same way with $w_0$ and $w_1$ interchanged. Since $w_0$ and $w_1$ lie in the same element of $\cal{K}$, we have that $\sigma_{\cal{K}}(\iota(C)) = \sigma_{\cal{K}}(C)$ for all $C \leq B_0$.

    Now let $H_0$ be $\E_\infty^K$-generic over $M$ such that $B_0 \in H_0$. Define:
    \begin{align*}
        H_1 := \iota``(H_0 \cap [B_0])^\uparrow,
    \end{align*}
    so $H_1$ is $\E_\infty^K$-generic over $M$ and $B_1 \in H_1$. Since $B_0 \in H_0$ and $\sigma_{\cal{K}}(\iota(C)) = \sigma_{\cal{K}}(C)$ for all $C \leq B_0$, we have that:
    \begin{align*}
        \sigma_{\cal{K}}``(H_0)^\uparrow  = \sigma_{\cal{K}}``(H_1)^\uparrow.
    \end{align*}
    Since $B_0,B_1 \leq A$, we have that $A \in H_0 \cap H_1$, so:
    \begin{align*}
        \sigma_{\cal{K}'}``(H_0)^\uparrow = i_{\sigma_{\cal{K}}``(H_0)^\uparrow}(\dot{z}) = i_{\sigma_{\cal{K}}``(H_1)^\uparrow}(\dot{z}) = \sigma_{\cal{K}'}``(H_1)^\uparrow.
    \end{align*}
    However, this filter contains both $\sigma_{\cal{K}'}(B_0)$ and $\sigma_{\cal{K}'}(B_1)$, which are incompatible: if $C \leq \sigma_{\cal{K}'}(B_0)$ and $C \leq \sigma_{\cal{K}'}(B_1)$, then $(q_{\cal{K}'}(w_0),p_0(A)) \in C$ and $(q_{\cal{K}'}(w_1),p_0(A)) \in C$, so $C\restrictedto\cal{K}'$ is not the equality relation as $q_{\cal{K}'}(w_0) \neq q_{\cal{K}'}(w_1)$.
\end{proof}

\begin{proof}[Proof of Theorem \ref{thm:CS-variables.alphabet.hierarchy}]
    If $|K| \leq |L|$, then $(\E_\infty^L,\leq,r)$ projects onto $(\E_\infty^K,\leq,r)$ by Example \ref{ex:CS-variables.merging.alphabets}. Conversely, suppose for a contradiction that $|K| > |L|$ and there is a weak projection $f : \E_\infty^L \to \E_\infty^K$. Let $M$ be the transitive collapse of a countable elementary submodel of some $H_\theta$ containing $f$ and everything else relevant, and let $G$ be $\E_\infty^L$-generic over $M$. By Lemma \ref{lem:ab.equivalence}, $H := f``(G)^\uparrow \in M[G]$ is $\E_\infty^K$-generic over $M$. For each partition $\cal{K}$ of $K$ with at least two elements, we let:
    \begin{align*}
        H_{\cal{K}} := \sigma_{\cal{K}}``(H)^\uparrow.
    \end{align*}
    By Lemmas \ref{lem:CS-variables.merging.letters} and \ref{lem:CS-variables.generic.real.general}, we have that $M \neq M[H_{\cal{K}}] = M[g_{H_{\cal{K}}}]$. If $\cal{K}$ and $\cal{K}'$ are distinct such partitions, then, interchanging $\cal{K}$ and $\cal{K}'$ if necessary, there exist $w_0,w_1 \in K$ which lie in the same element of $\cal{K}$ but in different elements of $\cal{K}'$, so $M[H_{\cal{K}}] \neq M[H_{\cal{K}'}]$ by Lemma \ref{lem:CS-variables.distinct.degrees}. By Lemma \ref{lem:CS-variables.forcing.reals}, each $M[H_{\cal{K}}]$ is of the form $M[\sigma_{\cal{L}}``(G)^\uparrow]$ for some partition $\cal{L}$ of $L$ with at least two elements. Thus, there are at most as many partitions of $K$ with at least two elements as there are of $L$, contradicting $|K| > |L|$.
\end{proof}

\section{Strong projections and a dichotomy}
\label{sec:strong.projections}
This section addresses the following two observations one might make throughout the paper. Firstly, almost every example of a \textbf{wA2}-space in \S\ref{sec:examples} has a natural choice of projection \emph{onto} the Ellentuck space, which is related to the depth function. Secondly, every projection in \S\ref{sec:projection.examples} is ``uniform'' in the following sense: $f$ is induced by some $\alpha : \cal{AR} \to \cal{AS}$, and the restriction $f\restrictedto[a,A] : [a,A] \to [\alpha(a),f(A)]$ is still a projection. We shall introduce the notion of \emph{strong projections}, which are projections possessing this ``uniform'' property. We show that every \textbf{A2}-space $(\cal{R},\leq,r)$ in a large class either strongly projects into the Ellentuck space, or is ``trivial'' in a specific way. More specifically, we say that $A \in \cal{R}$ is an \emph{atom} if for all $B \leq A$, $A \leq B$ holds. The main theorem of this section is as follows:

\thmdichotomy

Note that the depth-uniform property cannot be removed from the theorem, as Sacks forcing $(\mathbb{S},\subseteq,r)$ is depth-preserving but not depth-uniform (Example \ref{ex:Sacks.forcing.not.depth-uniform}), and satisfies neither of the two cases above: it does not \emph{weakly} project into the Ellentuck space (Theorem \ref{thm:no.wp.Sacks.sigma-closed}), while it has no atoms at all.

\subsection{Setup and basic properties}
In this subsection, we fix two \textbf{A2}-spaces $(\cal{R},\leq_\cal{R},r)$ and $(\S,\leq_\S,s)$. 

\begin{definition}
\label{def:homomorphism}
    A \emph{homomorphism} $\alpha : \cal{AR} \to \cal{AS}$ is a function such that $\alpha(\emptyset) = \emptyset$, and for all $a,b \in \cal{AR}$:
    \begin{itemize}
        \item If $a \sqsubseteq b$, then $\alpha(a) \sqsubseteq \alpha(b)$.

        \item If $a \leq_\fin b$, then $\alpha(a) \leq_\fin \alpha(b)$.
    \end{itemize}
\end{definition}

In other words, $\alpha$ is a homomorphism between the structures $(\cal{AR},\leq_{\fin,\cal{R}},\sqsubseteq_\cal{R})$ and $(\cal{AS},\leq_{\fin,\S},\sqsubseteq_\S)$ satisfying $\alpha(\emptyset) = \emptyset$. Note that a map $f : \cal{R} \to \S$ induced by a homomorphism $\alpha$ is automatically monotone --- indeed, for all $A,B \in \cal{R}$, we have that:
\begin{align*}
    A \leq B &\implies \forall n \, \exists m \, r_n(A) \leq_\fin r_m(B) \\
    &\implies \forall n \, \exists m \, \alpha(r_n(A)) \leq_\fin \alpha(r_m(B)) \\
    &\implies \lim_{n\to\infty} \alpha(r_n(A)) \leq \lim_{m\to\infty} \alpha(r_m(B)) \\
    &\implies \forall x \sqsubseteq f(A) \, \exists y \sqsubseteq f(B) \, x \leq_\fin y \\
    &\implies f(A) \leq f(B).
\end{align*}

\begin{definition}
\label{def:strong.projection}
    A \emph{strong projection} is a function $f : \cal{R} \to \S$ such that the following properties hold:
    \begin{enumerate}
        \item There exists some homomorphism $\alpha : \cal{AR} \to \cal{AS}$ that induces $f$ (see Definition \ref{def:induce}).

        \item For all $a \in \cal{AR}$, the restricted map $f\restrictedto[a] : [a] \to [\alpha(a)]$ is a projection.

        \item For all $A \in \cal{R}$ and $a \in \cal{AR}\restrictedto A$, the restricted map $f\restrictedto[a,A] : [a,A] \to [\alpha(a),f(A)]$ is a projection.
    \end{enumerate}
\end{definition}

It is easy to verify that a composition of strong projections is again a strong projection, as the composition of the corresponding homomorphisms is again a homomorphism, and each restriction of the composition is a composition of projections. Since $\alpha(\emptyset) = \emptyset$ for any homomorphism $\alpha$, property (2) of Definition \ref{def:strong.projection} applied to $a = \emptyset$ states precisely that $f : \cal{R} \to \S$ is a projection. Therefore, by Lemma \ref{lem:projection.implies.ab}, we get the following:
\begin{gather*}
    \text{Strong projection} \implies \text{Projection} \implies \text{Weak projection}
\end{gather*}
We remark that if $f : \cal{R} \to \S$ is a function satisfying (1) and (3) of Definition \ref{def:strong.projection}, $\cal{R}$ has a $\leq_\cal{R}$-maximal element $\mathbb{A}$ and $f(\mathbb{A})$ is $\leq_\S$-maximal in $\S$, then $f$ satisfies (2) automatically, as we may apply (3) to $[a] = [a,\mathbb{A}]$ and $[\alpha(a)] = [\alpha(a),f(\mathbb{A})]$.

The following lemma provides us with a useful way to check whether an induced map is a strong projection.

\begin{lemma}
\label{lem:sp.one-point.extension.check}
    Let $\alpha : \cal{AR} \to \cal{AS}$ be a homomorphism inducing a function $f : \cal{R} \to \S$. Suppose that:
    \begin{enumerate}
        \item $f\restrictedto[a]$ is cofinal for all $a \in \cal{AR}$; and

        \item For all $A \in \cal{R}$, $a \in \cal{AR}\restrictedto A$ and $X \in [\alpha(a),f(A)]$, there exists some $b \in r_{\lh(a)+1}[a,A]$ such that $\alpha(b) \in \cal{AS}\restrictedto[\alpha(a),X]$.
    \end{enumerate}
    Then $f$ is a strong projection.
\end{lemma}

\begin{proof}
    Note that $f$ is monotone, as it is induced by a homomorphism. We first show that $f\restrictedto[a,A] : [a,A] \to [\alpha(a),f(A)]$ is cofinal for all $A \in \cal{R}$ and $a \in \cal{AR}\restrictedto A$ --- that is, for all $X \in [\alpha(a),f(A)]$, there exists some $B \in [a,A]$ such that $f(B) \in [\alpha(a),X]$. We shall define sequences $(b_n)_{n<\omega}$ and $(X_n)_{n<\omega}$ such that the following hold for all $n$:
    \begin{enumerate}
        \item $b_0 = a$ and $b_{n+1} \in r_{\lh(b_n)+1}[b_n,A]$.

        \item $X_{n+1} \in [\alpha(b_{n+1}),X_n] \subseteq [\alpha(b_{n+1}),X]$.
    \end{enumerate}
    If this construction can be done, we may let $B := \lim_{n\to\infty} b_n$. Then $B \in [a,A]$ and $f(B) \leq X$, as $\alpha(b_{n+1}) \in \cal{AS}\restrictedto X_n \subseteq \cal{AS}\restrictedto X$ for all $n$. Since $a \sqsubseteq B$, we also have $\alpha(a) \sqsubseteq f(B)$, so $f(B) \in [\alpha(a),X]$.

    For the base case, we let $X_0 := X$. Suppose that $b_n$ and $X_n$ have been defined for some $n$. Applying hypothesis (2) of the lemma to $X_n \in [\alpha(b_n),f(A)]$, we may let $b_{n+1} \in r_{\lh(b_n)+1}[b_n,A]$ such that $\alpha(b_{n+1}) \in \cal{AS}\restrictedto[\alpha(b_n),X_n]$. We complete the induction by letting $X_{n+1} \in [\alpha(b_{n+1}),X_n]$.

    It remains to show that $f\restrictedto[a,A]$ is a projection for all $A \in \cal{R}$ and $a \in \cal{AR}\restrictedto A$, and that $f\restrictedto[a]$ is a projection for all $a \in \cal{AR}$. To see that property (3) of Definition \ref{def:projection} holds for $f\restrictedto[a,A]$, we let $B \in [a,A]$, and let $X \in [\alpha(a),f(B)]$. Since $f\restrictedto[a,B]$ is cofinal, there exists some $C \in [a,B]$ such that $f(C) \in [\alpha(a),X]$. For $f\restrictedto[a]$, we assumed that it is cofinal in the hypothesis. If $A \in [a]$ and $X \in [\alpha(a),f(A)]$, then we may apply the cofinality of $f\restrictedto[a,A]$ again to obtain some $B \in [a,A]$ such that $f(B) \in [\alpha(a),X]$.
\end{proof}

Once again, hypothesis (1) of Lemma \ref{lem:sp.one-point.extension.check} is unnecessary if $f$ sends the $\leq_\cal{R}$-maximal element of $\cal{R}$ to the $\leq_\S$-maximal element of $\S$.

All examples listed in \S\ref{sec:projection.examples} are strong projections. It is easy to check that the corresponding maps $\alpha$ are homomorphisms, and the proofs that they are strong projections are similar to the proofs provided in \S\ref{sec:projection.examples} (see also Lemma \ref{lem:monotone.ab.stronger.properties}). For instance, in Example \ref{ex:min.max.FIN}, we can extend the maps min and max to $\FIN^{[<\infty]}$ pointwise --- that is, for each $a = (x_n)_{n<N} \in \FIN^{[<\infty]}$, we let:
\begin{align*}
    \min(a) &:= \{\min(x_n) : n<N\}, \\
    \max(a) &:= \{\max(x_n) : n<N\}. 
\end{align*}
It is straightforward to check that these maps are homomorphisms, and the induced $\min,\max : \FIN^{[\infty]} \to [\N]^\infty$ satisfy the properties of a strong projection.

\subsection{The depth-preserving and depth-uniform properties}
In this subsection, we define two sufficient properties for an \textbf{A2}-space to satisfy the dichotomy in Theorem \ref{thm:dichotomy}. We fix an \textbf{A2}-space $(\cal{R},\leq,r)$. 

\subsubsection{The depth-preserving property}
\begin{definition}
\label{def:depth-preserving}
    We say that $(\cal{R},\leq,r)$ is \emph{depth-preserving} if for all $A,B \in \cal{R}$ such that $B \leq A$, and all $a \in \cal{AR}\restrictedto B$, we have that:
    \begin{align*}
        \depth_A(r_{\depth_B(a)}(B)) = \depth_A(a).
    \end{align*}
\end{definition}

It is easy to check that all examples in \S\ref{sec:examples}, except for $\FIN_k^{[\infty]}$ and $\FIN_{\pm k}^{[\infty]}$ when $k \geq 2$, are depth-preserving. We give an example of the argument using the Hindman space $\FIN^{[\infty]}$ --- the same argument applies to all other examples.

\begin{example}
\label{ex:FIN.depth-preserving}
    $(\FIN^{[\infty]},\leq,r)$ is depth-preserving.
\end{example}

\begin{proof}
    Let $A = (x_n)_{n<\omega},B = (y_n)_{n<\omega} \in \FIN^{[\infty]}$ be such that $B \leq A$, and let $a = (z_n)_{n<N} \in \FIN^{[<\infty]}\restrictedto B$. We may assume that $N \geq 1$, as the case $a = \emptyset$ is trivial. Observe that:
    \begin{align*}
        \depth_B(a) = \min\{M : \max(y_{M-1}) = \max(z_{N-1})\}.
    \end{align*}

    Let $M := \depth_B(a)$. Then similarly, we have that:
    \begin{align*}
        \depth_A(r_M(B)) = \min\{K : \max(x_{K-1}) = \max(y_{M-1})\}.
    \end{align*}

    Letting $K := \depth_A(r_M(B))$, we have that:
    \begin{align*}
        \max(x_{K-1}) = \max(y_{M-1}) = \max(z_{N-1}),
    \end{align*}
    so $\depth_A(a) = K$.
\end{proof}

On the other hand, we have that:

\begin{example}
\label{ex:FIN_k.not.depth-preserving}
    For $k \geq 2$, $(\FIN_k^{[\infty]},\leq,r)$ and $(\FIN_{\pm k}^{[\infty]},\leq,r)$ are \textbf{not} depth-preserving.
\end{example}

\begin{proof}
    We give a counterexample in $(\FIN_k^{[\infty]},\leq,r)$ --- an example in $(\FIN_{\pm k}^{[\infty]},\leq,r)$ can be constructed similarly. For each $n$, let $e_n \in \FIN_k$ be the function such that $\supp(e_n) = \{n\}$ and $e_n(n) = k$. We let:
    \begin{itemize}
        \item $A = (e_n)_{n<\omega}$.
        
        \item $B = (e_0,e_1 + T(e_2))^\frown(e_n)_{n \geq 3}$.
        
        \item $a = (e_0 + T^{k-1}(e_1))$.
    \end{itemize}
    Note that $A$ is the $\leq$-maximal element, so $B \leq A$. Furthermore:
    \begin{align*}
        a = e_0 + T^{k-1}(e_1) = e_0 + T^{k-1}(e_1 + T(e_2)),
    \end{align*}
    so $a \in \FIN_k^{[<\infty]}\restrictedto B$. It is easy to see that $\depth_B(a) = \depth_A(a) = 2$, but $\depth_A(r_2(B)) = 3$.
\end{proof}

\subsubsection{The depth-uniform property} Given $A \in \cal{R}$ and $a \in \cal{AR}\restrictedto A$, we define:
\begin{align*}
    D_{a,A} := \{\depth_A(b) - \depth_A(a) - 1 : b \in r_{\lh(a)+1}[a,A]\}.
\end{align*}
Note that by axiom \textbf{A2}(3'), $D_{a,A} \subseteq \N$.

\begin{definition}
\label{def:depth-uniform}
    We say that $(\cal{R},\leq,r)$ is \emph{depth-uniform} if for all $A,B \in \cal{R}$, $a \in \cal{AR}\restrictedto A$ and $b \in \cal{AR}\restrictedto B$, $D_{a,A} = D_{b,B}$. In this case, we write $D(\cal{R}) := D_{a,A}$ for any $A \in \cal{R}$ and $a \in \cal{AR}\restrictedto A$.
\end{definition}

It is simple to check that every example $(\cal{R},\leq,r)$ in \S\ref{sec:examples}, except for Sacks forcing $\mathbb{S}$ and the singleton space (Example \ref{ex:singleton.space}), satisfies the property that $D_{a,A} = \N$ for all $A \in \cal{R}$ and $a \in \cal{AR}\restrictedto A$. As an example, we show that:

\begin{example}
    The Carlson--Simpson space $(\E_\infty,\leq,r)$ is depth-uniform, with $D(\E_\infty) = \N$.
\end{example}

\begin{proof}
    We shall show that $D_{a,A} = \N$ for all $A \in \E_\infty$ and $a \in \E_{<\infty}\restrictedto A$. Fix such $A$ and $a$. Let $k := \depth_A(a)$, and let $N := p_k(A) = \dom(a)$. Given $n \in \N$, we let $M := p_{k+n+1}(A)$, and we let $b$ be the equivalence relation on $\dom(b) := M$ generated by:
    \begin{align*}
        a \cup \{(i,j) : i < N \leq j < M \wedge (i,j) \in A\} \cup (S \times S),
    \end{align*}
    where $S := \{j : N \leq j < M \wedge \forall i < N \, (i,j) \notin A\}$.

    Observe that $b\restrictedto N = a$, and $b$ is coarser than $A\restrictedto M$, so $b \in \E_{<\infty}\restrictedto A$. Furthermore, $S$ is an equivalence class of $b$ with $\min(S) = N$ (as $N = p_k(A)$ is the least element of its equivalence class in $A$), while every other equivalence class of $b$ meets $\dom(a)$. This implies that $b$ has one more equivalence class than $a$, and $p_{\lh(a)}(b) = N$, so $a = b\restrictedto N = r_{\lh(a)}(b) \sqsubseteq b$ and $b \in r_{\lh(a)+1}[a,A]$. Since $\dom(b) = M = p_{k+n+1}(A)$, we have that $\depth_A(b) = k + n + 1$, so $n \in D_{a,A}$.
\end{proof}

On the other hand:

\begin{example}
\label{ex:Sacks.forcing.not.depth-uniform}
    Sacks forcing $\mathbb{S}$ is \textbf{not} depth-uniform.
\end{example}

\begin{proof}
    Fix some skewed perfect binary tree $A \in \mathbb{S}$, and let $\func{split}(A) = \{s_0,s_1,\dots\}$ be the $\prec$-increasing enumeration of the splitting nodes of $A$. We first show that $D_{\emptyset,A} = \N$. It is obvious that $D_{\emptyset,A} \subseteq \N$. To see that $\N \subseteq D_{\emptyset,A}$, observe that for each $n$:
    \begin{align*}
        B := \{t \in A : t \sqsubseteq s_n \text{ or } s_n \sqsubseteq t\}
    \end{align*}
    is a skewed perfect binary subtree of $A$ such that $\depth_A(r_1(B)) = n + 1$, where $r_1(B) = \{t \in A : t \sqsubseteq s_n\}$ (as for all $m$, $s_n \in r_m(A)$ iff $m \geq n + 1$).

    On the other hand, observe that $s_1 \sqsubseteq s_3$ (as $s_1$ and $s_3$ correspond to $(0)$ and $(0,0)$ respectively in $2^{<\omega}$, and the map in property (2) of Example \ref{ex:Sacks.forcing} is $\sqsubseteq$-preserving in both directions). Thus, if $B \in [2,A]$ then the $\prec$-third splitting node of $B$ cannot be $s_3$, as the map in property (2) of Example \ref{ex:Sacks.forcing} guarantees that if $\func{split}(B) = \{t_0,t_1,\dots\}$, then $t_1$ is not an initial segment of $t_2$, while $t_0 = s_0$ and $t_1 = s_1$ as $r_2(B) = r_2(A)$. Therefore, there is no $b \in r_3[2,A]$ such that $\depth_A(b) - \depth_A(r_2(A)) - 1 = 1$, so $1 \notin D_{r_2(A),A}$. However, $1 \in D_{\emptyset,A}$, so $\mathbb{S}$ is not depth-uniform.
\end{proof}

\subsection{\texorpdfstring{$D$-Ellentuck spaces}{D-Ellentuck spaces}}
We introduce a class of topological Ramsey spaces that will be relevant to the proof of Theorem \ref{thm:dichotomy}.  

\subsubsection{Notation and definitions}
Fix some $D \subseteq \N$. Given $A \in [\N]^\infty$ and $\beta \leq \omega$, we let $[A]^{\beta,D}$ be the set of all $B \in [A]^\beta$ (i.e. $B \subseteq A$ and $|B| = \beta$) such that if $A = \{m_n : n < \omega\}$ and $B = \{m_{n_k} : k < \beta\}$ are the increasing enumerations of $A$ and $B$ respectively, then $n_0 \in D$ and $n_{k+1} - n_k - 1 \in D$ whenever $k + 1 < \beta$. We also write the following, for $n < \omega$ and $a \in [A]^{<\infty,D}$:
\begin{align*}
    [A]^{<\infty,D} &:= \bigcup_{\beta<\omega} [A]^{\beta,D}, \\
    [A]^{\leq\infty,D} &:= [A]^{\infty,D} \cup [A]^{<\infty,D}, \\
    r_n(A) &:= \{m_i : i < n\}, \\
    [a,A]^{<\infty} &:= \{b \in [A]^{<\infty,D} : a \sqsubseteq b\}, \\
    [a,A]^{\leq\infty} &:= \{B \in [A]^{\leq\infty,D} : a \sqsubseteq B\}.
\end{align*}
Note that $[A]^{\beta,\N} = [A]^\beta$. It is clear that $A \in [\N]^{\infty,D}$ iff $r_n(A) \in [\N]^{n,D}$ for all $n$.

\begin{notation}
\label{not:subseteq.D}
    Let $b \subseteq \N$, and let $b = \{m_i : i < |b|\}$ be its increasing enumeration. For $a \subseteq \N$, we write $a \subseteq^D b$ if $a \subseteq b$ and, writing $a = \{m_{i_j} : j < |a|\}$ with $i_0 < i_1 < \cdots$, we have that $i_0 \in D$ and $i_{j+1} - i_j - 1 \in D$ whenever $j + 1 < |a|$.
\end{notation}

Note that $[A]^{\beta,D} = \{B \in [A]^\beta : B \subseteq^D A\}$ for all $A \in [\N]^\infty$ and $\beta \leq \omega$. While $\subseteq^D$ is not a partial order in general, it is one precisely when $(D,+)$ is a \emph{submonoid} of $(\N,+)$ (which we shall abbreviate as just \emph{submonoid}) --- that is, $0 \in D$ and for all $\delta,\eta \in D$, $\delta + \eta \in D$. Readers may also refer to \cite{RG09,ADG20} for more information about submonoids.

\begin{lemma}
\label{lem:cusd.iff.quasi-order}
    For any $D \subseteq \N$, $([\N]^{\leq\infty},\subseteq^D)$ is a partial order iff $D$ is a submonoid.
\end{lemma}

\begin{proof}
    \underline{$\implies$:} Since $\subseteq^D$ is reflexive, $\{0\} \subseteq^D \{0\}$, so $0 \in D$. Now let $\delta,\eta \in D$. We define:
    \begin{align*}
        c &:= \N, \\
        b &:= \{\eta + i : i \leq \delta\}, \\
        a &:= \{\delta + \eta\}.
    \end{align*}
    Then $a \subseteq^D b$ and $b \subseteq^D c$, so $a \subseteq^D c$ by the transitivity of $\subseteq^D$. Since the index of $\delta + \eta$ in $c$ is $\delta + \eta$, this implies that $\delta + \eta \in D$, as required.

    \underline{$\impliedby$:} Reflexivity is immediate as $0 \in D$, and antisymmetry is immediate as $\subseteq^D$ refines $\subseteq$. For transitivity, suppose that $a \subseteq^D b$ and $b \subseteq^D c$, and let $c = \{m_i : i < |c|\}$ be its increasing enumeration. It is clear that $a \subseteq c$. We may write $b = \{m_{i_j} : j < |b|\}$ and $a = \{m_{i_{j_k}} : k < |a|\}$, so that for all $j,k$ with $j + 1 < |b|$ and $k + 1 < |a|$:
    \begin{align*}
        i_0 \in D, \quad i_{j+1} - i_j - 1 \in D, \quad j_0 \in D, \quad j_{k+1} - j_k - 1 \in D.
    \end{align*}
    We need to show that $i_{j_0} \in D$ and $i_{j_{k+1}} - i_{j_k} - 1 \in D$ whenever $k + 1 < |a|$. This follows from the identities below, as their right-hand sides are sums of elements of $D$:
    \begin{align*}
        i_{j_0} &= i_0 + j_0 + \sum_{p<j_0} (i_{p+1} - i_p - 1), \\
        i_{j_{k+1}} - i_{j_k} - 1 &= (j_{k+1} - j_k - 1) + \sum_{j_k \leq p < j_{k+1}} (i_{p+1} - i_p - 1).
    \end{align*}
\end{proof}

\subsubsection{$D$-Ellentuck space and axioms \textbf{A1}--\textbf{A4}}
\begin{definition}
    Given some submonoid $D \subseteq \N$, the \emph{$D$-Ellentuck space} is the closed triple $([\N]^{\infty,D},\subseteq^D,r)$.
\end{definition}

Note that the $\N$-Ellentuck space is precisely the Ellentuck space. The main intermediate theorem required to prove Theorem \ref{thm:dichotomy} is the following:

\begin{theorem}
\label{thm:D-Ellentuck.dpdu.TRS}
    For any submonoid $D \subseteq \N$, the $D$-Ellentuck space is a depth-preserving depth-uniform closed triple satisfying \textbf{A1}--\textbf{A4}, with a $\subseteq^D$-maximal element $\N$, and satisfies $D([\N]^{\infty,D}) = D$.
\end{theorem}

Given some $g \in \N^+ := \N \setminus \{0\}$, we let $g\N := \{gn : n \in \N\}$. We require the following property of submonoids to prove Theorem \ref{thm:D-Ellentuck.dpdu.TRS}:

\begin{lemma}[{\cite[Lemma 2.1 and Proposition 2.2]{RG09}}]
\label{lem:cusd.implies.gcd}
    If $D \subseteq \N$ is a submonoid and $D \neq \{0\}$, then $D$ is a cofinite subset of $g\N$, where $g = \gcd(D)$ is the largest positive integer that divides every element of $D$.
\end{lemma}

\begin{proof}[Proof of Theorem \ref{thm:D-Ellentuck.dpdu.TRS}]
    \underline{\textbf{A1}--\textbf{A2}:} It is obvious that the $D$-Ellentuck space satisfies \textbf{A1}. It satisfies \textbf{A2} due to the finite version of $\subseteq^D$ (see Notation \ref{not:subseteq.D}), and note that $\{b : b \subseteq^D a\}$ is finite as $\subseteq^D$ refines $\subseteq$. For \textbf{A2}(3), note that if $a \subseteq^D b$ and $c \sqsubseteq a$, then the indices of $c$ in $b$ form an initial segment of those of $a$ in $b$, so $c \subseteq^D b$. 
    
    \underline{\textbf{A3}:} Observe first that for all $A \in [\N]^{\infty,D}$, with increasing enumeration $A = \{m_n : n < \omega\}$, and all non-empty $a \in [A]^{<\infty,D}$:
    \begin{align*}
        \depth_A(a) = n \iff \max(a) = m_{n-1}.
    \end{align*}

    \begin{claim}
    \label{claim:D-Ellentuck.extension}
        Let $A = \{m_n : n < \omega\} \in [\N]^{\infty,D}$ and $a \in [A]^{<\infty,D}$, and let $k := \depth_A(a)$. Then:
        \begin{align*}
            [a,A]^{\leq\infty} = \{a \cup \{m_{k+c} : c \in C\} : C \in [\N]^{\leq\infty,D}\}.
        \end{align*}
    \end{claim}

    \begin{midproof}
        Every $B \subseteq A$ with $a \sqsubseteq B$ is of the form $B = a \cup \{m_{k+c} : c \in C\}$ for a unique $C \subseteq \N$. Let $C = \{c_j : j < |C|\}$ be the increasing enumeration of such a $C$. Note that $\max(a) = m_{k-1}$ if $a \neq \emptyset$, while $k = 0$ if $a = \emptyset$. Since $a \subseteq^D A$, $(k + c_0) - (k - 1) - 1 = c_0$ and $(k + c_{j+1}) - (k + c_j) - 1 = c_{j+1} - c_j - 1$, we have that $B \subseteq^D A$ iff $C \in [\N]^{\leq\infty,D}$. Note also that if $B \subseteq^D A$, then $B \subseteq^D \N$ by the transitivity of $\subseteq^D$, so $B \in [\N]^{\leq\infty,D}$.
    \end{midproof}

    Now fix some $A = \{m_n : n < \omega\} \in [\N]^{\infty,D}$.
    \begin{enumerate}
        \item Let $a \in [A]^{<\infty,D}$, let $n := \depth_A(a)$, and let $B \in [n,A]$. Since $a \subseteq^D r_n(A) \sqsubseteq B$, we have that $a \in [B]^{<\infty,D}$ and $\depth_B(a) = n$. Therefore, $a \cup (B \setminus r_n(A)) \in [a,B]$ by Claim \ref{claim:D-Ellentuck.extension}, taking $C := \N$.

        \item Let $B \leq A$ and $a \in [B]^{<\infty,D}$. Let $n := \depth_A(a)$ and $k := \depth_B(a)$, so that $\depth_A(r_k(B)) = n$. By Claim \ref{claim:D-Ellentuck.extension}, we have that $B = r_k(B) \cup \{m_{n+c} : c \in C\}$ for some $C \in [\N]^{\infty,D}$, and that $B' := r_n(A) \cup \{m_{n+c} : c \in C\} \in [n,A]$. By (1), $a \in [B']^{<\infty,D}$, $\depth_{B'}(a) = n$ and $[a,B'] \neq \emptyset$. Since $B \setminus r_k(B) = B' \setminus r_n(A)$, Claim \ref{claim:D-Ellentuck.extension} implies that $[a,B'] = [a,B]$.
    \end{enumerate}

    \underline{\textbf{A4}:} Let $A = \{m_n : n < \omega\} \in [\N]^{\infty,D}$, $a \in [A]^{<\infty,D}$ and $\O \subseteq [\N]^{\lh(a)+1,D}$, and let $k := \depth_A(a)$. We let:
    \begin{align*}
        E := \{\delta \in D : a \cup \{m_{k+\delta}\} \in \O\}.
    \end{align*}
    It suffices to find some $C \in [\N]^{\infty,D}$ such that $D \cap C \subseteq E$ or $E \cap C = \emptyset$. Indeed, if $C = \{c_j : j < \omega\}$ is the increasing enumeration of such a $C$, then $B := r_k(A) \cup \{m_{k+c} : c \in C\} \in [k,A]$ by Claim \ref{claim:D-Ellentuck.extension}. If $b \in r_{\lh(a)+1}[a,B]$, then since $\depth_B(a) = k$, Claim \ref{claim:D-Ellentuck.extension} implies that $b = a \cup \{m_{k+c_\delta}\}$ for some $\delta \in D$. Since $\{c_\delta\} \subseteq^D C \subseteq^D \N$, we have that $c_\delta \in D \cap C$. Thus, $b \in \O$ if $D \cap C \subseteq E$, and $b \notin \O$ if $E \cap C = \emptyset$.

    The statement is trivial if $D = \{0\}$ (by taking $C := \N$), so we may assume that $D \neq \{0\}$, and hence that $D$ is infinite. Let $g := \gcd(D)$, and by Lemma \ref{lem:cusd.implies.gcd}, let $N$ be large enough so that $gn \in D$ for all $n \geq N$. If $E$ is finite, then we may let $\delta \in D$ be such that $\delta > \eta$ for all $\eta \in E$, and $C := \{\delta + n : n \in \N\} \in [\N]^{\infty,D}$ satisfies $E \cap C = \emptyset$. Thus, we may assume that $E$ is infinite, so we may let $\eta_0 < \eta_1 < \cdots$ be elements of $E$ such that $\eta_{i+1} - \eta_i \geq (N + 1)g$ for all $i$. We let:
    \begin{align*}
        C := \bigcup_{i<\omega} \{\eta_i + j : j < g\}.
    \end{align*}
    Since $E \subseteq D \subseteq g\N$, we have that $\eta_i + j \notin D$ whenever $0 < j < g$, so $D \cap C = \{\eta_i : i < \omega\} \subseteq E$. It remains to show that $C \in [\N]^{\infty,D}$. Note that $\min(C) = \eta_0 \in D$, and if $x < y$ are consecutive elements of $C$, then either $y - x - 1 = 0 \in D$, or $x = \eta_i + g - 1$ and $y = \eta_{i+1}$ for some $i$. In the latter case, $y - x - 1 = \eta_{i+1} - \eta_i - g \geq Ng$ and is a multiple of $g$, so it is in $D$.

    \underline{Depth-preserving:} Let $B,A \in [\N]^{\infty,D}$ with $B \subseteq^D A$, and let $a \in [B]^{<\infty,D}$. Let $n := \depth_B(a)$, so that $\max(a) = \max(r_n(B))$, and let $m := \depth_A(r_n(B))$, so that $\max(r_n(B)) = \max(r_m(A))$. Then $\max(a) = \max(r_m(A))$, so $\depth_A(r_n(B)) = m = \depth_A(a)$.

    \underline{Depth-uniform:} Let $A = \{m_n : n < \omega\} \in [\N]^{\infty,D}$ and $a \in [A]^{<\infty,D}$, and let $k := \depth_A(a)$. We need to show that $D_{a,A} = D$. By Claim \ref{claim:D-Ellentuck.extension}, we have that:
    \begin{align*}
        r_{\lh(a)+1}[a,A] = \{a \cup \{m_{k+\delta}\} : \delta \in D\}.
    \end{align*}
    Since $\depth_A(a \cup \{m_{k+\delta}\}) = k + \delta + 1$ for all $\delta \in D$, we have that $D_{a,A} = D$.

    Finally, $\N \in [\N]^{\infty,D}$ as $0 \in D$, and $A \subseteq^D \N$ for all $A \in [\N]^{\infty,D}$ since $[\N]^{\infty,D} = \{A \in [\N]^\infty : A \subseteq^D \N\}$.
\end{proof}

$D$-Ellentuck spaces serve as the natural ``intermediate space'' of the first part of Theorem \ref{thm:dichotomy} --- see Lemma \ref{lem:sp.dpdu.D-Ellentuck} and Theorem \ref{thm:sp.D-Ellentuck.Ellentuck}.

\subsection{Proof of Theorem \ref{thm:dichotomy}}
\begin{lemma}
\label{lem:D(R).is.cusd}
    If $(\cal{R},\leq,r)$ is a depth-preserving depth-uniform \textbf{A2}-space with a $\leq$-maximal element $\mathbb{A}$, then $D(\cal{R}) \subseteq \N$ is a submonoid.
\end{lemma}

\begin{proof}
    We write $D := D(\cal{R})$, and $\depth(a) := \depth_\mathbb{A}(a)$ for all $a \in \cal{AR}$. Fix $\delta,\eta \in D$. We may choose $\emptyset = a_0 \sqsubseteq a_1 \sqsubseteq \cdots \sqsubseteq a_{\delta+1} =: a$ such that $a_{i+1} \in r_{i+1}[a_i,\mathbb{A}]$ for all $i \leq \delta$, $\depth(a_1) = \eta + 1$, and $\depth(a_{i+1}) = \depth(a_i) + 1$ for all $1 \leq i \leq \delta$. This implies that $\lh(a) = \delta + 1$ and:
    \begin{align*}
        \depth(a) = \sum_{i \leq \delta} (\depth(a_{i+1}) - \depth(a_i)) = \delta + \eta + 1.
    \end{align*}
    Let $B \in [a,\mathbb{A}]$ be any element. Since $\delta \in D_{\emptyset,B}$, there exists some $b \in r_1[B]$ such that $\depth_B(b) = \delta + 1$. Since $(\cal{R},\leq,r)$ is depth-preserving, we have that:
    \begin{align*}
        \depth(b) = \depth(r_{\depth_B(b)}(B)) = \depth(r_{\delta+1}(B)) = \depth(a) = \delta + \eta + 1,
    \end{align*}
    so $\delta + \eta = \depth(b) - \depth(\emptyset) - 1 \in D_{\emptyset,\mathbb{A}} = D$.
\end{proof}

We are now ready to prove Theorem \ref{thm:dichotomy}. We shall split into two cases: $D(\cal{R}) = \{0\}$ and $D(\cal{R}) \neq \{0\}$.

\subsubsection{When $D(\cal{R}) \neq \{0\}$}
\begin{lemma}
\label{lem:sp.dpdu.D-Ellentuck}
    Let $(\cal{R},\leq,r)$ be a depth-preserving depth-uniform \textbf{A2}-space with a $\leq$-maximal element $\mathbb{A}$, and write $D := D(\cal{R})$. Then there exists a homomorphism $\alpha : \cal{AR} \to [\N]^{<\infty,D}$, inducing a map $f : \cal{R} \to [\N]^{\infty,D}$, satisfying the following properties:
    \begin{enumerate}
        \item For all $a \in \cal{AR}$, $\lh(\alpha(a)) = \lh(a)$.
        
        \item $f(\mathbb{A}) = \N$.

        \item For all $A \in \cal{R}$ and $a \in \cal{AR}\restrictedto A$, the restricted map $\alpha : \cal{AR}\restrictedto[a,A] \to [\alpha(a),f(A)]^{<\infty}$ is surjective.
        
        \item For all $A \in \cal{R}$ and $a \in \cal{AR}\restrictedto A$, $f``[a,A] = [\alpha(a),f(A)]$.
    \end{enumerate}
    In particular, $(\cal{R},\leq,r)$ strongly projects into $([\N]^{\infty,D(\cal{R})},\subseteq^{D(\cal{R})},r)$. 
\end{lemma}

\begin{proof}
    Let $\depth(a) := \depth_\mathbb{A}(a)$ for all $a \in \cal{AR}$. Recall, from Lemma \ref{lem:D(R).is.cusd}, that $D$ is a submonoid, so that the $D$-Ellentuck space $([\N]^{\infty,D},\subseteq^D,r)$ is defined and Theorem \ref{thm:D-Ellentuck.dpdu.TRS} applies. We define $\alpha : \cal{AR} \to [\N]^{<\infty,D}$ by:
    \begin{align*}
        \alpha(a) := \{\depth(a') - 1 : \emptyset \neq a' \sqsubseteq a\}.
    \end{align*}
    We shall show that this $\alpha$ works.

    Since $\depth$ is strictly increasing, we have that $\lh(\alpha(a)) = \lh(a)$ for all $a \in \cal{AR}$, so $\lim_{n\to\infty} \lh(\alpha(r_n(A))) = \infty$ for all $A \in \cal{R}$, and $\alpha(\emptyset) = \emptyset$. If $a \sqsubseteq b$, then the non-empty initial segments of $a$ form an initial segment of those of $b$, so $\alpha(a) \sqsubseteq \alpha(b)$ as $\depth$ is strictly increasing. Now suppose that $a \leq_\fin b$. We need to show that $\alpha(a) \subseteq^D \alpha(b)$. Let $\alpha(b) = \{m_n : n < \lh(b)\}$ be its increasing enumeration, and let $B \in [b,\mathbb{A}]$ be any element. Denote $n_i := \depth_B(r_i(a))$ for $i \leq \lh(a)$, and note that $n_i \leq \lh(b)$ for all $i$, so that $\depth(r_i(a)) - 1 = m_{n_i - 1}$ for $i \geq 1$ as $(\cal{R},\leq,r)$ is depth-preserving. Since $n_0 = 0$, to show that $\alpha(a) \subseteq^D \alpha(b)$, we need to show that $n_{i+1} - n_i - 1 \in D$ for $i < \lh(a)$. This follows from the fact that $r_{i+1}(a) \in r_{\lh(r_i(a))+1}[r_i(a),B]$ (which also implies that $\alpha(a) \in [\N]^{<\infty,D}$). Therefore, $\alpha$ is a homomorphism.

    Let $f : \cal{R} \to [\N]^{\infty,D}$ be the map induced by $\alpha$. It is clear from the definition of $\alpha$ that $f(\mathbb{A}) = \N$. We check that $\alpha$ satisfies property (3). Let $A \in \cal{R}$, $a \in \cal{AR}\restrictedto A$ and $x \in [\alpha(a),f(A)]^{<\infty}$. We may assume that $\lh(x) = \lh(a) + 1$ (for the general case, we construct the required $b$ inductively). Let $f(A) = \{m_n : n < \omega\}$ be its increasing enumeration, and let $k := \depth_A(a) = \depth_{f(A)}(\alpha(a))$ (the two agree as $(\cal{R},\leq,r)$ is depth-preserving). Observe that $x = \alpha(a) \cup \{m_{l-1}\}$, where $l := \depth_{f(A)}(x) > k$, and that $l - k - 1 \in D$ by Theorem \ref{thm:D-Ellentuck.dpdu.TRS}. Then $l - k - 1 \in D_{a,A}$, so there exists some $b \in r_{\lh(a)+1}[a,A]$ such that $\depth_A(b) - \depth_A(a) - 1 = l - k - 1$. This implies that $\depth_A(b) = l$, so $\depth(b) = \depth(r_l(A)) = m_{l-1} + 1$. Therefore, $\alpha(b) = \alpha(a) \cup \{m_{l-1}\} = x$.

    It remains to show that $f$ satisfies property (4), which follows from property (3). It is clear that $f``[a,A] \subseteq [\alpha(a),f(A)]$. Conversely, if $X \in [\alpha(a),f(A)]$, then by property (3), we may inductively pick $a = b_0 \sqsubseteq b_1 \sqsubseteq \cdots$ in $\cal{AR}\restrictedto[a,A]$ such that $\alpha(b_n) = r_{\lh(a)+n}(X)$. Let $B := \lim_{n\to\infty} b_n$, and we have that $f(B) = X$, as needed.

    Finally, $f$ is a strong projection: by property (4), $f``[a,B] = [\alpha(a),f(B)]$ for all $B \in [a,A]$, so each restriction $f\restrictedto[a,A]$ is a projection (cf. Lemma \ref{lem:monotone.ab.stronger.properties}), and property (2) of Definition \ref{def:strong.projection} follows from the remark after it, as $f(\mathbb{A}) = \N$ is the $\subseteq^D$-maximal element of $[\N]^{\infty,D}$.
\end{proof}

\begin{theorem}
\label{thm:sp.D-Ellentuck.Ellentuck}
    If $D \subseteq \N$ is a submonoid and $D \neq \{0\}$, then $([\N]^{\infty,D},\subseteq^D,r)$ strongly projects into the Ellentuck space $([\N]^\infty,\subseteq,r)$.
\end{theorem}

\begin{proof}
    Let $g := \gcd(D)$. By Lemma \ref{lem:cusd.implies.gcd}, we may let $N \in \N^+$ be large enough so that $gn \in D$ for all $n \geq N$. We define $\alpha : [\N]^{<\infty,D} \to [\N]^{<\infty}$ by:
    \begin{align*}
        \alpha(a) := \{k < \omega : gk \in a\}.
    \end{align*}
    It is clear that $\alpha$ is a homomorphism. If $A = \{m_n : n < \omega\} \in [\N]^{\infty,D}$, then since $D \subseteq g\N$, we have that $m_n \equiv n \pmod{g}$ for all $n$. In particular, $A$ contains infinitely many multiples of $g$, so we may let $f : [\N]^{\infty,D} \to [\N]^\infty$ be the map induced by $\alpha$. Since $\N$ is the $\subseteq^D$-maximal element of $[\N]^{\infty,D}$ and $f(\N) = \N$, hypothesis (1) of Lemma \ref{lem:sp.one-point.extension.check} holds automatically.

    Let $A = \{m_n : n < \omega\} \in [\N]^{\infty,D}$, $a \in [A]^{<\infty,D}$ and $X \in [\alpha(a),f(A)]$, and let $p := \depth_A(a)$. It suffices to find some $b \in r_{\lh(a)+1}[a,A]$ such that $\alpha(b) \in [\alpha(a),X]^{<\infty}$.
    \begin{itemize}
        \item If $p \not\equiv 0 \pmod{g}$, then $m_p$ is not a multiple of $g$, so we may let $b := a \cup \{m_p\}$, and we have that $\alpha(b) = \alpha(a)$.

        \item Suppose that $p \equiv 0 \pmod{g}$. Since $X \subseteq f(A)$, we may let $k \in X$ be large enough so that $gk = m_q$ for some $q \geq p + Ng$. Since $q \equiv m_q \equiv 0 \pmod{g}$, $q - p$ is a multiple of $g$ that is at least $Ng$, so $q - p \in D$. Therefore, $b := a \cup \{m_q\} \in r_{\lh(a)+1}[a,A]$ by Claim \ref{claim:D-Ellentuck.extension}, and $\alpha(b) = \alpha(a) \cup \{k\} \in [\alpha(a),X]^{<\infty}$.
    \end{itemize}
\end{proof}

\begin{proof}[Proof of Theorem \ref{thm:dichotomy}(1) if $D(\cal{R}) \neq \{0\}$]
    By Lemma \ref{lem:D(R).is.cusd}, $D(\cal{R})$ is a submonoid, so both Lemma \ref{lem:sp.dpdu.D-Ellentuck} and Theorem \ref{thm:sp.D-Ellentuck.Ellentuck} apply. 
\end{proof}

\subsubsection{When $D(\cal{R}) = \{0\}$}
We may now turn our attention to the case where $D(\cal{R}) = \{0\}$. Note that the singleton space (Example \ref{ex:singleton.space}) is a depth-preserving depth-uniform closed triple satisfying \textbf{A1}--\textbf{A4} with $D(\cal{R}) = \{0\}$, so case (2) of Theorem \ref{thm:dichotomy} is not vacuous.

\begin{proof}[Proof of Theorem \ref{thm:dichotomy}(2) if $D(\cal{R}) = \{0\}$]
    Fix a finitisation $\leq_\fin$ witnessing that $(\cal{R},\leq,r)$ is depth-preserving, depth-uniform and $D(\cal{R}) = \{0\}$. 

    \begin{claim}
        For all $A \in \cal{R}$, $[A]$ is metrically compact.
    \end{claim}
    
    \begin{midproof}
        Since $D_{a,A} = \{0\}$ for all $a \in \cal{AR}\restrictedto A$, one can show by induction on $\lh(a)$ shows that $\depth_A(a) = \lh(a)$ for all $a \in \cal{AR}\restrictedto A$. This implies that for all $n$:
        \begin{align*}
            \cal{AR}_n\restrictedto A \subseteq \{b \in \cal{AR} : b \leq_\fin r_n(A)\}
        \end{align*}
        and the latter set is finite by \textbf{A2}(1). Therefore, $[A]$ is compact as it is metrically closed.
    \end{midproof}

    Therefore, we may consider the partial order on the set $\F := \{[B] : B \leq A\}$, ordered by reverse containment. If $\cal{C} \subseteq \F$ is a $\subseteq$-decreasing chain, then $\bigcap\cal{C}$ is the intersection of a decreasing chain of compact sets, hence non-empty. Let $C \in \bigcap\cal{C}$, and we have that $[C]$ is a lower bound of $\cal{C}$. By Zorn's lemma, $[A]$ has a $\leq$-minimal element $B$ (that is, $C \leq B$ implies $B \leq C$), and $B$ is an atom.

    It remains to show that $(\cal{R},\leq,r)$ does not weakly project into the Ellentuck space. Suppose for a contradiction that we have a weak projection $f : \cal{R} \to [\N]^\infty$. Let $A \in \cal{R}$ be an atom. By Lemma \ref{lem:ab.iff.qo.TRS}, there exists some $X \leq f(A)$ such that $[X] \subseteq f``[A]$. Let $Y \subsetneq X$, and let $B \leq A$ such that $f(B) = Y$. Since $A$ is an atom, we have that $A \leq B$, so $X \subseteq f(A) \subseteq f(B) = Y \subsetneq X$, a contradiction.
\end{proof}

\subsubsection{Open problems}
As noted in Example \ref{ex:FIN_k.not.depth-preserving}, $(\FIN_k^{[\infty]},\leq,r)$ is not depth-preserving for $k \geq 2$ (though it is depth-uniform, since $D_{a,A} = \N$ always), yet it strongly projects into $([\N]^\infty,\subseteq,r)$ via the map $\max \circ T^{k-1}$ or $\min \circ T^{k-1}$. On the other hand, Sacks forcing $(\mathbb{S},\subseteq,r)$, which serves as the counterexample to the unrestricted version of Theorem \ref{thm:dichotomy}, is depth-preserving but not depth-uniform, so it witnesses that the depth-uniform hypothesis cannot be dropped. This suggests that the remaining hypothesis may be redundant.

\begin{qn}
\label{qn:dichotomy.depth-preserving.necessary}
    Is the depth-preserving hypothesis of Theorem \ref{thm:dichotomy} necessary?
\end{qn}

A key obstacle to adjusting the hypotheses of Theorem \ref{thm:dichotomy} lies in the difficulty of constructing strong projections that are not variants of the depth map. 

\begin{qn}
    What are some natural strong projections that are not the depth map (or some variants of it)?
\end{qn}

One may observe similarly that most natural weak projections are strong projections.

\begin{qn}
    What is a natural example of a weak projection that is not a strong projection?
\end{qn}

\begin{qn}
    Do there exist \textbf{A2}-spaces $(\cal{R},\leq_\cal{R},r)$ and $(\S,\leq_\S,s)$ such that $(\cal{R},\leq_\cal{R},r)$ weakly projects, but does not strongly project into $(\S,\leq_\S,s)$?
\end{qn}

\subsection{Strong projections and the well-indexed property}
\label{sec:strong.projections.and.well-indexed}
A notion closely related to strong projections into the Ellentuck space was considered by Wang in \cite{W25}. The following definition is slightly different from, but equivalent to, the original presentation (see \cite[Remark 5.2]{W25} and the proof of \cite[Proposition 5.3]{W25}).

\begin{definition}[{\cite[Definition 5.1]{W25}}]
\label{def:well-indexed}
    Let $(\cal{R},\leq,r)$ be an \textbf{A2}-space. We say that $(\cal{R},\leq,r)$ is \emph{well-indexed} if there exist some $\mathbb{A} \in \cal{R}$ and some $\sqsubseteq$-preserving map $\alpha : \cal{AR} \to [\N]^{<\infty}$ (Definition \ref{def:sqsubseteq-preserving}) satisfying the following properties:
    \begin{enumerate}
        \item 
        \begin{enumerate}
            \item $\alpha(\emptyset) = \emptyset$.
            
            \item $\alpha(a) \neq \emptyset$ for all $\emptyset \neq a \in \cal{AR}$.
            
            \item For all $A \in \cal{R}$ and $n < \omega$, $\lh(\alpha(r_{n+1}(A))) \leq \lh(\alpha(r_n(A))) + 1$.
        \end{enumerate}

        \item $\alpha$ induces some function $f : \cal{R} \to [\N]^\infty$ (see Definition \ref{def:induce}) satisfying the following properties:
        \begin{enumerate}
            \item $f(\mathbb{A}) = \N$.

            \item If $A \leq B \leq \mathbb{A}$, then $f(A) \subseteq f(B)$.

            \item For all $A \in \cal{R}$ and $n < \omega$, $[\alpha(r_n(A)),f(A)] \subseteq f``[r_n(A),A]$.
        \end{enumerate}
    \end{enumerate}
\end{definition}

Wang showed in \cite{W25} that if a closed triple $(\cal{R},\leq,r)$ satisfying \textbf{A1}--\textbf{A4}, with $\cal{AR}$ countable, is well-indexed, then the set of codes for Ramsey positive $G_\delta$ subsets of $\cal{R}$ is $\mathbf{\Sigma}_2^1$-complete. This generalises a result by Sabok in \cite{Sab12}, who proved the original result for just the Ellentuck space.

Applying Lemma \ref{lem:sp.dpdu.D-Ellentuck} in the case $D(\cal{R}) = \N$, we can conclude that:

\begin{theorem}
    If $(\cal{R},\leq,r)$ is a depth-preserving depth-uniform \textbf{A2}-space with a $\leq$-maximal element $\mathbb{A}$ such that $D(\cal{R}) = \N$, then $(\cal{R},\leq,r)$ is well-indexed.
\end{theorem}

\begin{proof}
    Let $\alpha$ and $f$ be as provided by Lemma \ref{lem:sp.dpdu.D-Ellentuck}. Since $\alpha$ is a homomorphism and $\lh(\alpha(a)) = \lh(a)$ for all $a \in \cal{AR}$, property (1) holds. Property (2)(a) follows from $f(\mathbb{A}) = \N$, property (2)(b) holds as $f$ is monotone, and property (2)(c) follows from property (4) of Lemma \ref{lem:sp.dpdu.D-Ellentuck}.
\end{proof}

Combining this result with the main theorem of \cite{W25}, we get:

\begin{corollary}
\label{cor:dpdu.DR-nontrivial.implies.complete.codes}
    If $(\cal{R},\leq,r)$ is a depth-preserving depth-uniform closed triple satisfying \textbf{A1}--\textbf{A4}, with a $\leq$-maximal element $\mathbb{A}$ such that $D(\cal{R}) = \N$, then the set of codes for Ramsey positive $G_\delta$ subsets of $\cal{R}$ is $\mathbf{\Sigma}_2^1$-complete.
\end{corollary}

A list of examples of well-indexed topological Ramsey spaces was also provided in \cite{W25}. The reader may verify that, in most of these examples, the witnessing map is precisely the depth map defined in the proof of Lemma \ref{lem:sp.dpdu.D-Ellentuck}.

We remark that if $\mathbb{A}$ and $\alpha$ witness that $(\cal{R},\leq,r)$ is well-indexed, then the restriction $f\restrictedto[\mathbb{A}] : [\mathbb{A}] \to [\N]^\infty$ is a projection (as in Definition \ref{def:projection}). For every $A \in \mathbb{S}$, $([A],\subseteq)$ is a dense subposet of Sacks forcing below $A$, so the proof of Theorem \ref{thm:no.wp.Sacks.sigma-closed} shows that $([A],\subseteq)$ does not weakly project into the Ellentuck space. Therefore, Sacks forcing $(\mathbb{S},\subseteq,r)$ is not well-indexed. This answers \cite[Problem 1]{W25} in the negative.

\section{Ramsey categories}
\label{sec:Ramsey.categories}
In this section, we discuss further possible connections between topological Ramsey theory and the Ramsey property of categories. In \S\ref{sec:Ramsey.quasi-orders}, we saw two examples of Ramsey categories (Examples \ref{ex:Ellentuck.space.category} and \ref{ex:Hindman.FIN.category}) that realise the finite Ramsey theorem and finite Hindman theorem. Given an \textbf{A2}-space $(\cal{R},\leq,r)$, we may define an analogous finite Ramsey-theoretic statement as follows: Given $a,b \in \cal{AR}$ such that $a \sqsubseteq b$ and a non-negative integer $l$, we define:
\begin{align*}
    \O_a(l,b) &:= \{d \in \cal{AR}_{\lh(a)+l}\restrictedto[a] : \exists b' \sqsubseteq b[d \leq_\fin b']\}, \\
    \O_a(b) &:= \bigcup_{l<\omega} \O_a(l,b).
\end{align*}
We shall assume that for all $a,b \in \cal{AR}$, if $a \leq_\fin b$ then $\lh(a) \leq \lh(b)$. The reader may verify that every example in \S\ref{sec:examples} satisfies this property.

\begin{definition}
    The statement $\varphi(\cal{R};a)$ asserts the following: For all positive integers $k,l,m$ with $l \leq m$, there is a positive integer $n$ such that for every $b \in \cal{AR}_{\lh(a)+n}\restrictedto[a]$ and $k$-colouring $\chi$ of $\O_a(l,b)$, there exists some $c \in \O_a(m,b)$ such that $\chi$ is monochromatic on $\O_a(l,c) \neq \emptyset$. We also let $\varphi(\cal{R})$ denote $\forall a \in \cal{AR} \, \varphi(\cal{R};a)$.
\end{definition}

Suppose that we wish to construct a Ramsey category $\mathbf{C}(\cal{R})$ such that the weak Ramsey property of $\mathbf{C}(\cal{R})$ is precisely the statement $\varphi(\cal{R})$. Following the construction in Examples \ref{ex:Ellentuck.space.category} and \ref{ex:Hindman.FIN.category}, we let $\Ob(\mathbf{C}(\cal{R})) = \N$, and we would like $\hom_{\mathbf{C}(\cal{R})}(l,n) = \O_a(l,b)$ for some $a \sqsubseteq b$ such that $\lh(b) - \lh(a) = n$. Since the choice of $a,b$ here is unclear, we introduce our first structural assumption, which removes the need to make a choice entirely.

\begin{definition}
\label{def:TRS.self-similar}
    An \textbf{A2}-space $(\cal{R},\leq,r)$ is \emph{self-similar} if for all $a,a' \in \cal{AR}$, $b \in \cal{AR}\restrictedto[a]$ and $b' \in \cal{AR}\restrictedto[a']$ such that $\lh(b') - \lh(a') = \lh(b) - \lh(a)$, the structures $(\O_a(b),\leq_\fin,\sqsubseteq)$ and $(\O_{a'}(b'),\leq_\fin,\sqsubseteq)$ are isomorphic. That is, there exists a function $i : \O_a(b) \to \O_{a'}(b')$ satisfying the following properties:
    \begin{enumerate}
        \item $i$ is bijective.
        
        \item For all $c,d \in \O_a(b)$, $c \sqsubseteq d$ iff $i(c) \sqsubseteq i(d)$.
        
        \item For all $c,d \in \O_a(b)$, $c \leq_\fin d$ iff $i(c) \leq_\fin i(d)$.
    \end{enumerate}
\end{definition}

\begin{example}
    The Ellentuck space $([\N]^\infty,\subseteq,r)$ and the Gowers spaces $(\FIN_k^{[\infty]},\leq,r)$ are self-similar:
    \begin{itemize}
        \item For $([\N]^\infty,\subseteq,r)$, let $j : b \setminus a \to b' \setminus a'$ be the unique strictly increasing bijection, and define $i : \O_a(b) \to \O_{a'}(b')$ by $i(c) := a' \cup j``(c \setminus a)$.

        \item For $(\FIN_k^{[\infty]},\leq,r)$, write $b/a = (x_q)_{q<n}$ and $b'/a' = (y_q)_{q<n}$ for the tails of $b$ and $b'$ beyond $a$ and $a'$ respectively. Define $i : \O_a(b) \to \O_{a'}(b')$ by stipulating that if $c = a^\frown(z_j)_{j<l}$ and for each $j$, $z_j = \sum_{q<n} T^{\lambda_q^j}(x_q)$, then let $w_j := \sum_{q<n} T^{\lambda_q^j}(y_q)$ and $i(c) := {a'}^\frown(w_j)_{j<l}$.
    \end{itemize}
\end{example}

Observe that if $i : \O_a(b) \to \O_{a'}(b')$ is an isomorphism between structures, then since $i$ preserves $\sqsubseteq$, $i$ also preserves the relative length of each $c \in \O_a(b)$ --- that is, $\lh(i(c)) - \lh(a') = \lh(c) - \lh(a)$ for all $c \in \O_a(b)$. Furthermore, for all $c \in \O_a(b)$, $l \leq \lh(c) - \lh(a)$:
\begin{align*}
    d \in \O_a(l,c) &\iff a \sqsubseteq d \wedge (\lh(d) = \lh(a) + l) \wedge \exists e \sqsubseteq c[d \leq_\fin e] \\
    &\iff a' \sqsubseteq i(d) \wedge (\lh(i(d)) = \lh(a') + l) \wedge \exists e \sqsubseteq i(c)[i(d) \leq_\fin e] \\
    &\iff i(d) \in \O_{a'}(l,i(c)).
\end{align*}
Thus, the restriction $i\restrictedto\O_a(l,c) : \O_a(l,c) \to \O_{a'}(l,i(c))$ is also an isomorphism between structures. This allows us to define $\mathbf{C}(\cal{R})$ concretely as follows: Fix some $A \in \cal{R}$, and define $\hom_{\mathbf{C}(\cal{R})}(l,n) := \O(l,n) := \O_\emptyset(l,r_n(A))$ for $l \leq n$ (if $l > n$, let $\hom_{\mathbf{C}(\cal{R})}(l,n) = \emptyset$). By self-similarity, the structure $(\O(l,n),\leq_\fin,\sqsubseteq)$ does not depend on the choice of $A$ (up to isomorphism), and is isomorphic to $(\O_a(l,b),\leq_\fin,\sqsubseteq)$ whenever $\lh(b) - \lh(a) = n$. 

Similar to the categories $\mathbf{E}$ and $\mathbf{BS}$ defined in Examples \ref{ex:Ellentuck.space.category} and \ref{ex:Hindman.FIN.category} respectively, we shall assume that there is a partial binary operation $\circ$, where $b \circ c \in \cal{AR}_l$ is defined for $b \in \cal{AR}_n$ and $c \in \O(l,n)$, such that for all $l \leq n$, $a \in \cal{AR}_l$ and $b \in \cal{AR}_n$:
\begin{align*}
    \exists b' \sqsubseteq b[a \leq_\fin b'] \iff \exists c \in \O(l,n) \, a = b \circ c.
\end{align*}
Note that if $a \in \O(m,n)$ and $b \in \O(l,m)$, then $a \circ b \in \O(l,n)$ --- indeed, by the equivalence above there exists some $a' \sqsubseteq a$ such that $a \circ b \leq_\fin a'$. Since $a \leq_\fin r_i(A)$ for some $i \leq n$ and $a \circ b \leq_\fin a' \sqsubseteq a$, we have that $a \circ b \leq_\fin r_j(A)$ for some $j \leq i \leq n$ by \textbf{A2}(3), so $a \circ b \in \O(l,n)$.

To complete our definition of the category $\mathbf{C}(\cal{R})$, we assume that $\circ$ is associative and every $\O(n,n)$ contains an identity morphism.

\begin{lemma}
\label{lem:Ramsey.category.iff.varphi}
    $\mathbf{C}(\cal{R})$ satisfies the weak Ramsey property iff $\varphi(\cal{R})$ holds.
\end{lemma}

\begin{proof}
    We fix $A \in \cal{R}$ such that $\hom_{\mathbf{C}(\cal{R})}(l,n) = \O_\emptyset(l,r_n(A))$ for all $l \leq n$. Given positive integers $l,m,n$ such that $l \leq m \leq n$, we observe that for all $b \in \O(m,n)$ and $a \in \cal{AR}_l$:
    \begin{align*}
        a \in b \circ \hom_{\mathbf{C}(\cal{R})}(l,m) &\iff \exists c \in \O(l,m) \, a = b \circ c \\
        &\iff \exists b' \sqsubseteq b[a \leq_\fin b'] \\
        &\iff a \in \O_\emptyset(l,b),
    \end{align*}
    so $b \circ \hom_{\mathbf{C}(\cal{R})}(l,m) = \O_\emptyset(l,b)$. Since $\hom_{\mathbf{C}(\cal{R})}(l,m) \neq \emptyset$ iff $l \leq m$, the weak Ramsey property of $\mathbf{C}(\cal{R})$ is equivalent to the following statement:
    \begin{enumerate}[label=($\varphi^*(\cal{R})$),leftmargin=*]
        \item For all positive integers $k,l,m$ with $l \leq m$, there is a positive integer $n$ such that for every $k$-colouring $\chi$ of $\O_\emptyset(l,r_n(A))$, there exists some $c \in \O_\emptyset(m,r_n(A))$ such that $\chi$ is monochromatic on $\O_\emptyset(l,c) \neq \emptyset$.
    \end{enumerate}
    It is clear that $\varphi^*(\cal{R})$ follows from $\varphi(\cal{R};\emptyset)$, which follows from $\varphi(\cal{R})$. Thus, it remains to show that $\varphi^*(\cal{R})$ implies $\varphi(\cal{R})$. 

    Fix some $a \in \cal{AR}$. Let $k,l,m$ be positive integers with $l \leq m$, and let $n$ be the positive integer witnessing that $\varphi^*(\cal{R})$ holds for $k,l,m$. Let $b \in \cal{AR}_{\lh(a)+n}\restrictedto[a]$, and let $\chi$ be a $k$-colouring of $\O_a(l,b)$. By self-similarity, there is a structure isomorphism $i : \O_\emptyset(r_n(A)) \to \O_a(b)$. Then $\chi \circ i$ is a $k$-colouring of $\O_\emptyset(l,r_n(A))$, so by $\varphi^*(\cal{R})$, there exists some $c \in \O_\emptyset(m,r_n(A))$ such that $\chi \circ i$ is monochromatic on $\O_\emptyset(l,c) \neq \emptyset$. Then $\chi$ is monochromatic on $i``(\O_\emptyset(l,c)) = \O_a(l,i(c)) \neq \emptyset$ and $i(c) \in \O_a(m,b)$, as desired.
\end{proof}

Thus, assuming self-similarity and the existence of a binary operation satisfying the required properties, the construction of a Ramsey category from an \textbf{A2}-space follows from verifying the statement $\varphi(\cal{R})$. The advantage in studying the statement $\varphi(\cal{R})$, independent of any category theory, is the non-requirement of a binary operation $\circ$ witnessing that $\mathbf{C}(\cal{R})$ is a category. As one may expect, for self-similar topological Ramsey spaces, $\varphi(\cal{R})$ follows from the abstract Ellentuck theorem (Theorem \ref{thm:abstract.Ellentuck.thm}) via a standard compactness argument:

\begin{theorem}
\label{thm:A4.implies.varphi}
    If $(\cal{R},\leq,r)$ satisfies \textbf{A4} and is self-similar, then $\varphi(\cal{R})$ holds.
\end{theorem}

\begin{proof}
    By the proof of Lemma \ref{lem:Ramsey.category.iff.varphi}, it suffices to show that $\varphi^*(\cal{R})$ holds. Fix some $A \in \cal{R}$ such that $\O(l,n) = \O_\emptyset(l,r_n(A))$ for all $l \leq n$. Fix positive integers $k,l,m$ with $l \leq m$, and suppose that for every positive integer $n$, there is a $k$-colouring $\chi_n$ of $\O_\emptyset(l,r_n(A))$ such that for all $c \in \O_\emptyset(m,r_n(A))$, $\chi_n$ is not monochromatic on $\O_\emptyset(l,c)$. For each $a \in \cal{AR}_l\restrictedto A$ and $i < k$, we let $P_{a,i}$ be a propositional variable representing ``$a$ has colour $i$'', and let $T$ be the theory containing the following formulas:
    \begin{enumerate}[label=(\arabic*)]
        \item $\bigvee_{i<k} P_{a,i}$ for each $a \in \cal{AR}_l\restrictedto A$.

        \item $\neg(P_{a,i} \wedge P_{a,j})$ for each $a \in \cal{AR}_l\restrictedto A$ and $i < j < k$.

        \item $\neg\bigwedge_{b \in \O_\emptyset(l,a)} P_{b,i}$ for each $a \in \cal{AR}_m\restrictedto A$ and $i < k$.
    \end{enumerate}
    Observe that $T$ is finitely satisfiable --- given a finite $T' \subseteq T$, we may let $n$ be greater than $\depth_A(a)$ for every $a \in \cal{AR}$ that appears as an index in some formula of $T'$, and let $P_{a,i}$ be true iff $\chi_n(a) = i$ for each $a \in \O_\emptyset(l,r_n(A))$ and $i < k$ (with arbitrary truth values for the remaining variables); then this assignment satisfies $T'$.
    
    By the compactness theorem, $T$ is satisfiable, say by the truth assignment $v$. We may thus define a colouring $\chi : \cal{AR}_l\restrictedto A \to k$ by letting, for each $a \in \cal{AR}_l\restrictedto A$, $\chi(a)$ be the unique $i < k$ such that $v(P_{a,i})$ is true. We define, for each $i < k$:
    \begin{align*}
        \X_i := \{B \leq A : \chi(r_l(B)) = i\}.
    \end{align*}
    Each $\X_i$ is a closed subset of $\cal{R}$, and the sets $\X_0,\dots,\X_{k-1}$ partition $[A]$, so by the abstract Ellentuck theorem there exists some $B \leq A$ such that $[B] \subseteq \X_i$ for some $i < k$. Then for all $b \in \O_\emptyset(l,r_m(B))$, there exists some $C \in [b,B] \subseteq [B] \subseteq \X_i$, and $\chi(b) = \chi(r_l(C)) = i$. This implies that $v(P_{b,i})$ is true for all $b \in \O_\emptyset(l,r_m(B))$, contradicting that $v$ satisfies (3) for $a = r_m(B)$.
\end{proof}

It is unclear whether the converse of Theorem \ref{thm:A4.implies.varphi} holds. This is not surprising, as the compactness argument that derives the finite Ramsey theorem from the infinite Ramsey theorem cannot be reversed. 

\begin{qn}
\label{qn:varphi.implies.A4}
    Let $(\cal{R},\leq,r)$ be a self-similar \textbf{A2}-space such that for all $a,b \in \cal{AR}$, if $a \leq_\fin b$ then $\lh(a) \leq \lh(b)$. If $\varphi(\cal{R})$ holds, does $(\cal{R},\leq,r)$ necessarily satisfy \textbf{A4}?
\end{qn}

Let $(\cal{R},\leq_\cal{R},r)$ be a closed triple satisfying \textbf{A1}--\textbf{A4}, and let $(\S,\leq_\S,s)$ be an \textbf{A2}-space, such that both triples satisfy all assumptions stated in this section (so that $\mathbf{C}(\cal{R})$ and $\mathbf{C}(\S)$ are categories). There are now two possible ways to transfer the weak (analytic) Ramsey property from $\cal{R}$ to $\S$:
\begin{itemize}
    \item A continuous quasi-open map $f : \cal{R} \to \S$ (by Lemma \ref{lem:cts.qo.transfers.weak.analytic.Ramsey.property}).
    
    \item A pre-adjunction $F : \Ob(\mathbf{C}(\S)) \rightleftarrows \Ob(\mathbf{C}(\cal{R})) : G$ (by Theorem \ref{thm:category.pre-adjunction.transfer.weak.Ramsey.property}).
\end{itemize}

\begin{qn}
    If there is a continuous quasi-open map $f : \cal{R} \to \S$, must there exist a pre-adjunction $F : \Ob(\mathbf{C}(\S)) \rightleftarrows \Ob(\mathbf{C}(\cal{R})) : G$? What about the converse?
\end{qn}

Since continuous quasi-open maps also transfer the weak analytic Ramsey property between abstract Ramsey spaces (more generally, Ramsey abstract orders, cf. Lemma \ref{lem:cts.qo.transfers.weak.analytic.Ramsey.property.ARS}), a framework detailing the construction of a Ramsey category from an abstract Ramsey space would also be interesting.

\section{Further remarks and future directions}
As discussed in \S\ref{sec:quasi-open}, it is unlikely that the monotonicity assumption can be removed from our framework without significant trade-offs. However, continuous quasi-open maps transfer the weak analytic Ramsey property (Lemma \ref{lem:cts.qo.transfers.weak.analytic.Ramsey.property}), so one may instead develop a framework of maps that transfer this weaker property, and study the (non-)existence of such maps between topological Ramsey spaces.

\begin{qn}
    Establish the non-existence or existence of Alaoglu--Birkhoff convergent maps (or quasi-open maps, pre-adjunctions etc.) between examples in \S\ref{sec:examples}.
\end{qn}

In the spirit of \S\ref{sec:projection.abstract.orders} and \S\ref{sec:non-existence.of.projections}, we also pose the following question:

\begin{qn}
    Is there a (continuous) quasi-open map:
    \begin{align*}
        (W_L^{[\infty]},W_{Lv}^{[\infty]},\leq,\leq^o,r,s) \to (\FIN_L^{[\infty]},\FIN_{Lv}^{[\infty]},\leq,\leq^o,r,s)?
    \end{align*}
\end{qn}

We remark that there are several obvious continuous quasi-open maps from $(\FIN_L^{[\infty]},\FIN_{Lv}^{[\infty]},\leq,\leq^o,r,s)$ or $(W_L^{[\infty]},W_{Lv}^{[\infty]},\leq,\leq^o,r,s)$ to some of the topological Ramsey spaces presented in \S\ref{sec:examples}, where a topological Ramsey space $(\cal{R},\leq,r)$ is viewed as the abstract order $(\cal{R},\cal{R},\leq,\leq)$. We state two examples, leaving the verification to the readers:
\begin{itemize}
    \item The maps $\max \circ \dom : \FIN_L \cup \FIN_{Lv} \to \N$ and $\min \circ \dom : \FIN_L \cup \FIN_{Lv} \to \N$, extended to sequences pointwise, are continuous quasi-open (cf. Example \ref{ex:min.max.FIN}).

    \item Since sequences in $W_L^{[\infty]}$ and $W_{Lv}^{[\infty]}$ are rapidly increasing, for each $X \in W_{Lv}^{[\infty]}$ we may define $\supp_X : \c{X}_L \cup \c{X}_{Lv} \to \FIN$ in the same way as Example \ref{ex:support.variable.words}. This map, extended to sequences pointwise, is also continuous quasi-open.
\end{itemize}

\begin{qn}
    Which \textbf{wA2}-spaces $(\cal{R},\leq,r)$ in \S\ref{sec:examples} admit a (continuous) quasi-open map from $(W_L^{[\infty]},W_{Lv}^{[\infty]},\leq,\leq^o)$ or $(\FIN_L^{[\infty]},\FIN_{Lv}^{[\infty]},\leq,\leq^o)$ to $(\cal{R},\cal{R},\leq,\leq)$ (or vice versa)?
\end{qn}

\section*{Acknowledgements}
The authors would like to thank Anton Bernshteyn, Natasha Dobrinen Jan Hubi\v{c}ka, S\l{}awomir Solecki, Stevo Todor\v{c}evi\'{c}, Jing Zhang and Andy Zucker for helpful conversations.

\section*{AI use statement}
A nearly complete draft was written without using AI, including all of \S\ref{sec:intro}, \S\ref{sec:Ramsey.quasi-orders.and.spaces}, \S\ref{sec:projections} and \S\ref{sec:strong.projections}, and most of \S\ref{sec:non-existence.of.projections}. Claude Fable 5.1 and Opus 5.5 were later consulted to verify the validity of proofs, check for typos, compile the bibliography and simplify the proofs of various results. AI also proposed proof strategies for Theorems \ref{thm:no.wp.strong.subtrees.sigma-closed.projection}, \ref{thm:no.wp.Carlson-Simpson.and.FIN}, \ref{thm:no.wp.CS-variables.FIN}, and \ref{thm:CS-variables.alphabet.hierarchy}, but the authors wrote the proofs themselves. 

\printbibliography[title={References}]

\end{document}